**W. B. VASANTHA KANDASAMY**

**FLORENTIN SMARANDACHE**

# FUZZY RELATIONAL MAPS AND NEUTROSOPHIC RELATIONAL MAPS

*2004*

# FUZZY RELATIONAL MAPS AND NEUTROSOPHIC RELATIONAL MAPS

**W. B. Vasantha Kandasamy**

Department of Mathematics
Indian Institute of Technology, Madras
Chennai – 600036, India
e-mail: **vasantha@iitm.ac.in**
web: **http://mat.iitm.ac.in/~wbv**

**Florentin Smarandache**
Department of Mathematics
University of New Mexico
Gallup, NM 87301, USA
e-mail: **smarand@gallup.unm.edu**

2004

# CONTENTS

## Chapter Two
## SOME APPLICATIONS OF FRE

## Chapter Three
## SOME NEW AND BASIC DEFINITIONS ON
## NEUTROSOPHIC THEORY

**Chapter Four**
**NEUTROSOPHIC RELATIONAL EQUATIONS AND THEIR PROPERTIES**

**Chapter Five**
**SUGGESTED PROBLEMS**

⓪⑥

*Dedicated to those*
*few, young, not-so-influential,*
*revolutionary scientists*
*and mathematicians who support*
*the newer paradigm shift*

①②

# Preface

The aim of this book is two fold. At the outset the book gives most of the available literature about Fuzzy Relational Equations (FREs) and its properties for there is no book that solely caters to FREs and its applications. Though we have a comprehensive bibliography, we do not promise to give all the possible available literature about FRE and its applications. We have given only those papers which we could access and which interested us specially. We have taken those papers which in our opinion could be transformed for neutrosophic study.

The second importance of this book is that for the first time we introduce the notion of Neutrosophic Relational Equations (NRE) which are analogous structure of FREs. Neutrosophic Relational Equations have a role to play for we see that in most of the real-world problems, the concept of indeterminacy certainly has its say; but the FRE has no power to deal with indeterminacy, but this new tool NRE has the capacity to include the notion of indeterminacy. So we feel the NREs are better tools than FREs to use when the problem under investigation has indeterminates. Thus we have defined in this book NREs and just sketched its probable applications.

This book has five chapters. The first chapter is a bulky one with 28 sections. These sections deal solely with FREs and their properties. By no means do we venture to give any proof for the results for this would make our book unwieldy and enormous in size. For proofs, one can refer the papers that have been cited in the bibliography.

The second chapter deals with the applications of FRE. This has 10 sections: we elaborately give the applications of FRE in flow rates in chemical industry problems, preference and determination of peak hour in the transportation problems, the social problems faced by bonded laborers etc.

Chapter three for the first time defines several new neutrosophic concepts starting from the notion of neutrosophic fuzzy set, neutrosophic fuzzy matrix, neutrosophic lattices, neutrosophic norms etc. and just indicate some of its important analogous properties. This chapter has six sections which are solely devoted to the introduction of several neutrosophic concepts which are essential for the further study of NRE.

Chapter four has eleven sections. This chapter gives all basic notions and definitions about the NREs and introduces NREs. Section 4.11 is completely devoted to suggest how one can apply NREs in the study of real world problems. We suggest many problems in chapter five for the reader to solve.

This is the third book in the Neutrosophics Series. The earlier two books are *Fuzzy Cognitive Maps and Neutrosophic Cognitive Maps* (http://gallup.unm.edu/~smarandache/NCMs.pdf) and *Analysis of Social Aspects of Migrant Labourers Living With HIV/AIDS Using Fuzzy Theory and Neutrosophic Cognitive Maps: With Specific Reference to Rural Tamil Nadu in India* (http://gallup.unm.edu/~smarandache/NeutrosophyAIDS.pdf).

Finally, we thank Meena Kandasamy for the cover design. We thank Kama Kandasamy for the layout of the book and for drawing all the figures used in this book. She displayed an enormous patience that is worthy of praise. We owe deep thanks to Dr.K.Kandasamy for his patient proof-reading of the book. Without his help this book would not have been possible.

W.B.VASANTHA KANDASAMY
FLORENTIN SMARANDACHE

**Chapter One**

# FUZZY RELATIONAL EQUATIONS: BASIC CONCEPTS AND PROPERTIES

The notion of fuzzy relational equations based upon the max-min composition was first investigated by Sanchez [84]. He studied conditions and theoretical methods to resolve fuzzy relations on fuzzy sets defined as mappings from sets to [0,1]. Some theorems for existence and determination of solutions of certain basic fuzzy relation equations were given by him. However the solution obtained by him is only the greatest element (or the maximum solution) derived from the max-min (or min-max) composition of fuzzy relations. [84]'s work has shed some light on this important subject. Since then many researchers have been trying to explore the problem and develop solution procedures [1, 4, 10-12, 18, 34, 52, 75-80, 82, 108, 111].

The max-min composition is commonly used when a system requires conservative solutions in the sense that the goodness of one value cannot compensate the badness of another value [117]. In reality there are situations that allow compensatability among the values of a solution vector. In such cases the min operator is not the best choice for the intersection of fuzzy sets, but max-product composition, is preferred since it can yield better or at least equivalent result. Before we go into the discussion of these Fuzzy Relational Equations (FRE) and its properties it uses and applications we just describe them.

This chapter has 28 sections that deal with the properties of FRE, methods of solving FRE using algorithms given by several researchers and in some cases methods of neural networks and genetic algorithm is used in solving problems. A complete set of references is given in the end of the book citing the names of all researchers whose research papers have been used.

## 1.1 Binary Fuzzy Relation and their properties

It is well known fact that binary relations are generalized mathematical functions. Contrary to functions from X to Y, binary relations R(X, Y) may assign to each element of X two or more elements of Y. Some basic operations on functions such as the inverse and composition are applicable to binary relations as well.

Given a fuzzy relation R(X, Y), its domain is a fuzzy set on X, dom R, whose membership function is defined by

$$\text{dom } (R(x)) = \max_{y \in Y} R \text{ } (x, y)$$

for each $x \in X$. That is, each element of set X belongs to the domain of R to the degree equal to the strength of its strongest relation to any member of set Y. The range of R (X, Y) is a fuzzy relation on Y, ran R whose membership function is defined by

$$\text{ran } R(y) = \max_{x \in X} R \text{ } (x, y)$$

for each $y \in Y$. That is, the strength of the strongest relation that each element of Y has to an element of X is equal to the degree of that elements membership in the range of R. In addition, the height of a fuzzy relation R(X,Y) is a number, h(R), defined by

$$h \text{ } (R) = \max_{y \in Y} \text{ } \max_{x \in X} (R \text{ } (x, y).$$

That is h(R) is the largest membership grade attained by any pair (x, y) in R.

A convenient representation of binary relation R(X, Y) are membership matrices $R = [r_{xy}]$ where $r_{xy} = R(x, y)$. Another useful representation of binary relation is a sagittal diagram. Each of the sets X, Y is represented by a set of nodes in the diagram nodes corresponding to one set is distinguished from nodes representing the other set.

Elements of $X \times Y$ with non-zero membership grades in R(X, Y) are represented in the diagram by lines connecting the respective nodes.

We illustrate the sagittal diagram of a binary fuzzy relation R(X, Y) together with the corresponding membership matrix in Figure 1.1.1.

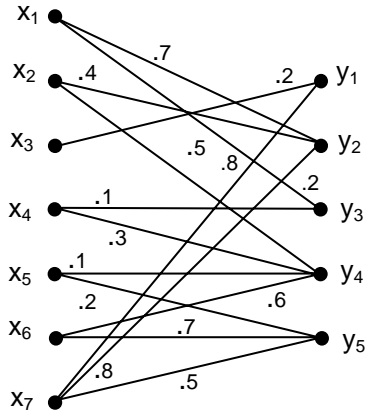

Figure: 1.1.1

The inverse of a fuzzy relation R(X, Y) denoted by $R^{-1}(Y, X)$ is a relation on $Y \times X$ defined by $R^{-1}(y, x) = R(x, y)$ for all $x \in X$ and for all $y \in X$. A membership matrix $R^{-1} = [r^{-1}_{yx}]$ representing $R^{-1}(Y, X)$ is the transpose of the matrix R for R (X, Y) which means that the rows of $R^{-1}$ equal the columns of R and the columns of $R^{-1}$ equal the rows of R.

Clearly $(R^{-1})^{-1} = R$ for any binary fuzzy relation. Thus a fuzzy binary relation can be represented by the sagittal diagram. The corresponding membership matrix:

$$
\begin{array}{c}
 \\
x_1 \\
x_2 \\
x_3 \\
x_4 \\
x_5 \\
x_6 \\
x_7
\end{array}
\begin{array}{ccccc}
y_1 & y_2 & y_3 & y_4 & y_5 \\
\left[\begin{array}{ccccc}
0 & .7 & .5 & 0 & 0 \\
0 & .4 & 0 & .1 & 0 \\
.2 & 0 & 0 & 0 & 0 \\
0 & 0 & .1 & 1 & 0 \\
0 & 0 & 0 & .3 & .7 \\
0 & 0 & 0 & .6 & .7 \\
.2 & 0 & .8 & 0 & .5
\end{array}\right]
\end{array}
$$

R is the membership matrix.

Consider now two binary fuzzy relations P(X, Y) and Q(Y, Z) with a common set Y. The standard composition of these relations, which is denoted by P(X, Y) ° Q(Y, Z), produces a binary relation R (X, Z) on $X \times Z$ defined by

$$R (x, z) = [P \circ Q] (x, z)$$
$$= \max_{y \in Y} \min [P(x, y), Q (y, z)]$$

for all $x \in X$ and all $z \in Z$. This composition, which is based on the standard t-norm and t-co-norm is often referred to as the max-min composition. It follows directly from the above equation that

$$[P (X, Y) \circ Q (Y, Z)]^{-1} = Q^{-1} (Z, Y) \circ P^{-1} (Y, X)$$
$$[P (X, Y) \circ Q (Y, Z)] \circ R (Z, W) = P (X, Y) \circ [Q (Y, Z) \circ R (Z, W)].$$

This is the standard (or max-min) composition, which is associative, and its inverse is equal to the reverse composition of the inverse relations.

However the standard composition is not commutative because Q(Y, Z) ° P(X, Y) is not even well defined when $X \neq Z$. Even if X = Z and Q (Y, Z) ° P(X, Y) are well defined, we may have P(X, Y) ° Q (Y, Z) $\neq$ Q (Y, Z ) ° P(X, Y). Compositions of binary fuzzy relations can be performed conveniently in terms of membership matrices of the relations. Let P = $[p_{ik}]$, Q = $[q_{kj}]$ and R = $[r_{ij}]$ be the membership matrices of binary relations such that R = P ° Q. We can then write using this matrix notations

$$[r_{ij}] = [p_{ik}] \circ [q_{kj}]$$

where $r_{ij} = \max_{k} \min (p_{ik} q_{kj})$. Observe that the same elements of P and Q are used in the calculation of R as would be used in the regular multiplication of matrices, but the product and sum operations are here replaced with max and min operations respectively.

A similar operation on two binary relations, which differs from the composition in that it yields triples instead of pairs, is known as the relational join. For fuzzy relations P(X, Y) and Q(Y, Z), the relational join P * Q, corresponding to the standard max-min composition is a ternary relation R(X, Y, Z) defined by

R (x, y, z) = [P * Q] (x, y, z) = min [P (x, y), Q (y, z)] for each x ∈ X, y ∈ Y and z ∈ Z.

The fact that the relational join produces a ternary relation from two binary relations is a major difference from the composition, which results in another binary relation.

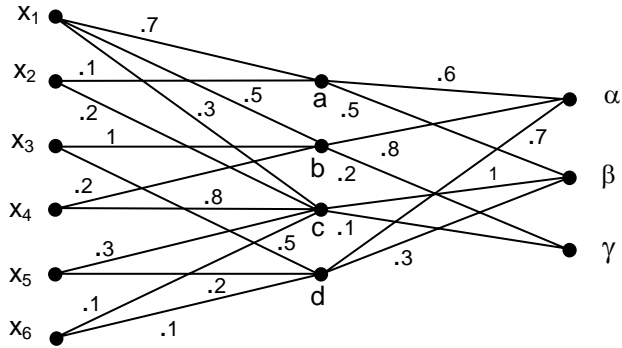

Figure: 1.1.2

Formally [P °Q] (x, z) = max [P * Q] (x, y, z) for each x ∈ X and z ∈ Z. Now we just see what happens if the binary relation on a single set. Binary relation R (X, X) can be expressed by the same forms as general binary relations.

The following properties are to be observed:

    i.    Each element of the set X is represented as a single node in the diagram.

    ii.    Directed connection between nodes indicates pairs of elements of X, with the grade of membership in R is non-zero.

    iii.    Each connection in the diagram is labeled by the actual membership grade of the corresponding pair in R.

An example of this diagram for a relation R (X, X) defined on X = {$x_1, x_2, x_3, x_4, x_5$} is shown in Figure 1.1.3. A crisp relation R (X, X) is reflexive if and only if (x, x) ∈ R. for each x ∈ R, that is if every element of X is related to itself, otherwise R(X, X) is called irreflexive. If (x, x) ∉ R for every x ∈ X the relation is called anti reflexive.

A crisp relation R(X, X) is symmetric if and only if for every (x, y) ∈ R, it is also the case that (y, x) ∈ R where x, y ∈ X. Thus whenever an element x is related to an element y through a symmetric relation, y is also related to x. If this is not the case for some x, y then the relation is called asymmetric. If both $(x, y) \in R$ and $(y, x) \in R$ implies $x = y$ then the relation is called anti symmetric. If either $(x, y) \in R$ or $(y, x) \in R$ whenever $x \neq y$, then the relation is called strictly anti symmetric.

$$
\begin{array}{c c}
 & \begin{array}{c c c c c} x_1 & x_2 & x_3 & x_4 & x_5 \end{array} \\
\begin{array}{c} x_1 \\ x_2 \\ x_3 \\ x_4 \\ x_5 \end{array} &
\left[ \begin{array}{c c c c c}
.2 & 0 & .1 & 0 & .7 \\
0 & .8 & .6 & 0 & 0 \\
.1 & 0 & .2 & 0 & 0 \\
.2 & 0 & 0 & .1 & .4 \\
.6 & .1 & 0 & 0 & .5
\end{array} \right]
\end{array}
$$

The corresponding sagittal diagram is given in Figure 1.1.3:

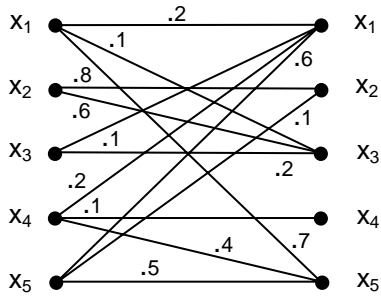

Figure: 1.1.3

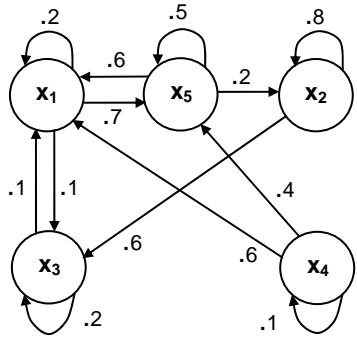

Figure: 1.1.4

Table

| x | y | R(x, y) |
|---|---|---------|
| $x_1$ | $x_1$ | .2 |
| $x_1$ | $x_3$ | .1 |
| $x_1$ | $x_5$ | .7 |
| $x_2$ | $x_2$ | .8 |
| $x_2$ | $x_3$ | .6 |
| $x_3$ | $x_1$ | .1 |
| $x_3$ | $x_3$ | .2 |
| $x_4$ | $x_1$ | .2 |
| $x_4$ | $x_4$ | .1 |
| $x_4$ | $x_5$ | .4 |
| $x_5$ | $x_1$ | .6 |
| $x_5$ | $x_2$ | .1 |
| $x_5$ | $x_5$ | .5 |

A crisp relation R (X, Y) is called transitive if and only if $(x, z) \in$ R, whenever both $(x, y) \in$ R and $(y, z) \in$ R for at least one $y \in$ X.

In other words the relation of x to y and of y to z imply the relation x to z is a transitive relation. A relation that does not satisfy this property is called non-transitive. If $(x, z) \notin$ R whenever both $(x, y) \in$ R and $(y, z) \in$ R, then the relation is called anti-transitive. The reflexivity, symmetry and transitivity is described by the following Figure 1.1.5:

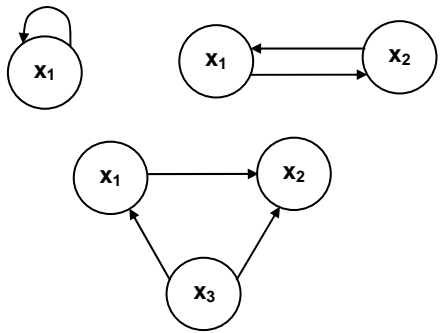

Figure: 1.1.5

A fuzzy relation R (X, X) is reflexive if and only if R(x,x) = 1 for all x ∈ X, if this is not the case for same x ∈ X, the relation is called irreflexive, if it is not satisfied for all x ∈ X, the relation is called anti reflexive. A weaker form of reflexivity referred to as ∈ - reflexivity denoted by R (x, x) ≥ ∈ where 0 < ∈ < 1. A fuzzy relation is symmetric if and only if

$$R(x, y) = R (y, x)$$

for all x, y ∈ X, if this relation is not true for some x, y ∈ X, the relation is called anti symmetric. Further more when R (x, y) > 0 and R (y, x) > 0 implies x = y for all x, y ∈ X the relation R is called anti symmetric.

A fuzzy relation R (X, X) is transitive if $R (x, z) \geq \max_{y \in Y} \min$ [R (x, y), R (y, z)] is satisfied for each pair (x, z) ∈ $X^2$. A relation failing to satisfy this inequality for some members of X is called non-transitive and if

$$R (x, z) < \max_{y \in Y} \min [R (x, y), R (y, z)]$$

for all (x,z) ∈ $X^2$, then the relation is called anti transitive.

## 1.2 Properties of Fuzzy Relations

In this section we just recollect the properties of fuzzy relations like, fuzzy equivalence relation, fuzzy compatibility relations, fuzzy ordering relations, fuzzy morphisms and sup-i-compositions of fuzzy relation. For more about these concepts please refer [43].

Now we proceed on to define fuzzy equivalence relation. A crisp binary relation R(X, X) that is reflexive, symmetric and transitive is called an equivalence relation. For each element x in X, we can define a crisp set $A_x$, which contains all the elements of X that are related to x, by the equivalence relation.

$$A_x = \{y \mid (x, y) \in R (X, X)\}$$

$A_x$ is clearly a subset of X. The element x is itself contained in $A_x$ due to the reflexivity of R, because R is transitive and symmetric each member of $A_x$, is related to all the other members of $A_x$. Further no member of $A_x$, is related to any element of X not included in $A_x$. This set $A_x$ is referred to an as equivalence class of R (X, X) with respect to x. The members of each equivalence class can be considered equivalent to each other and only to each other under the relation R. The family of all such equivalence classes defined by the relation which is usually denoted by X / R, forms a partition on X.

A fuzzy binary relation that is reflexive, symmetric and transitive is known as a fuzzy equivalence relation or similarity relation. In the rest of this section let us use the latter term. While the max-min form of transitivity is assumed, in the following discussion on concepts; can be generalized to the alternative definition of fuzzy transitivity.

While an equivalence relation clearly groups elements that are equivalent under the relation into disjoint classes, the interpretation of a similarity relation can be approached in two different ways. First it can be considered to effectively group elements into crisp sets whose members are similar to each other to some specified degree. Obviously when this degree is equal to 1, the grouping is an equivalence class. Alternatively however we may wish to consider the degree of similarity that the elements of X have to some specified element $x \in X$. Thus for each $x \in X$, a similarity class can be defined as a fuzzy set in which the membership grade of any particular element represents the similarity of that element to the element x. If all the elements in the class are similar to x to the degree of 1 and similar to all elements outside the set to the degree of 0 then the grouping again becomes an equivalence class. We know every fuzzy relation R can be uniquely represented in terms of its α-cuts by the formula

$$R = \bigcup_{\alpha \in (0,1]} \alpha. \,^{\alpha}R$$

It is easily verified that if R is a similarity relation then each α-cut, $^{\alpha}R$ is a crisp equivalence relation. Thus we may use any similarity relation R and by taking an α - cut $^{\alpha}R$ for any value $\alpha \in$ (0, 1], create a crisp equivalence relation that represents the presence of similarity between the elements to the degree α. Each of these equivalence relations form a partition of X. Let $\pi (^{\alpha}R)$ denote the partition corresponding to the equivalence relation $^{\alpha}R$. Clearly any two elements x and y belong to the same block of this partition if and only if R (x, y) ≥ α. Each similarity relation is associated with the set $\pi (R) = \{\pi (^{\alpha}R) \mid \alpha \in (0,1]\}$ of partition of

X. These partitions are nested in the sense that $\pi$ ($^{\alpha}$R) is a refinement of $\pi$ ($^{\beta}$R) if and only $\alpha \geq \beta$.

The equivalence classes formed by the levels of refinement of a similarity relation can be interpreted as grouping elements that are similar to each other and only to each other to a degree not less than $\alpha$.

Just as equivalences classes are defined by an equivalence relation, similarity classes are defined by a similarity relation. For a given similarity relation R(X, X) the similarity class for each x $\in$ X is a fuzzy set in which the membership grade of each element y $\in$ X is simply the strength of that elements relation to x or R(x, y). Thus the similarity class for an element x represents the degree to which all the other members of X are similar to x. Expect in the restricted case of equivalence classes themselves, similarity classes are fuzzy and therefore not generally disjoint.

Similarity relations are conveniently represented by membership matrices. Given a similarity relation R, the similarity class for each element is defined by the row of the membership matrix of R that corresponds to that element.

Fuzzy equivalence is a cutworthy property of binary relation R(X, X) since it is preserved in the classical sense in each $\alpha$-cut of R. This implies that the properties of fuzzy reflexivity, symmetry and max-min transitivity are also cutworthy. Binary relations are symmetric and transitive but not reflexive are usually referred to as quasi equivalence relations.

The notion of fuzzy equations is associated with the concept of compositions of binary relations. The composition of two fuzzy binary relations P (X, Y) and Q (Y, Z) can be defined, in general in terms of an operation on the membership matrices of P and Q that resembles matrix multiplication. This operation involves exactly the same combinations of matrix entries as in the regular matrix multiplication. However the multiplication and addition that are applied to these combinations in the matrix multiplication are replaced with other operations, these alternative operations represent in each given context the appropriate operations of fuzzy set intersections and union respectively. In the max-min composition for example, the multiplication and addition are replaced with the min and max operations respectively.

We shall give the notational conventions. Consider three fuzzy binary relations P (X, Y), Q (Y, Z) and R (X, Z) which are defined on the sets

$$X = \{x_i \mid i \in I\}$$

$$Y = \{y_j \mid j \in J\} \text{ and}$$
$$Z = \{z_k \mid k \in K\}$$

where we assume that $I = N_n$ $J = N_m$ and $K = N_5$. Let the membership matrices of P, Q and R be denoted by $P = [p_{ij}]$, $Q = [q_{ij}]$, $R = [r_{ik}]$ respectively, where $p_{ij} = P(x_i, y_j)$, $q_{jk} = Q(y_j, z_k)$ $r_{ij} = R(x_i, z_k)$ for all $i \in I$ ($=N_n$), $j \in J = (N_m)$ and $k \in K$ ($=N_5$). This clearly implies that all entries in the matrices P, Q, and R are real numbers from the unit interval [0, 1]. Assume now that the three relations constrain each other in such a way that P°Q = R where ° denotes max-min composition. This means that $\max_{j \in J} \min (p_{ij}, q_{jk})$

$= r_{ik}$ for all $i \in I$ and $k \in$ ˉ K. That is the matrix equation P° Q = R encompasses $n \times s$ simultaneous equations of the form $\max_{j \in J} \min (p_{ij}, q_{jk}) = r_{ik}$. When two of the components in each of the equations are given and one is unknown these equations are referred to as fuzzy relation equations.

When matrices P and Q are given the matrix R is to determined using P ° Q = R. The problem is trivial. It is solved simply by performing the max-min multiplication – like operation on P and Q as defined by $\max_{j \in J} \min (p_{ij}, q_{jk}) = r_{ik}$. Clearly the solution in this case exists and is unique. The problem becomes far from trivial when one of the two matrices on the left hand side of P ° Q = R is unknown. In this case the solution is guaranteed neither to exist nor to be unique.

Since R in P ° Q = R is obtained by composing P and Q it is suggestive to view the problem of determining P (or alternatively Q ) from R to Q (or alternatively R and P) as a decomposition of R with respect to Q (or alternatively with respect to P). Since many problems in various contexts can be formulated as problems of decomposition, the utility of any method for solving P ° Q = R is quite high. The use of fuzzy relation equations in some applications is illustrated. Assume that we have a method for solving P ° Q = R only for the first decomposition problem (given Q and R).

Then we can directly utilize this method for solving the second decomposition problem as well. We simply write P ° Q = R in the form $Q^{-1}$ o $P^{-1} = R^{-1}$ employing transposed matrices. We can solve $Q^{-1}$ o $P^{-1} = R^{-1}$ for $Q^{-1}$ by our method and then obtain the solution of P ° Q = R by $(Q^{-1})^{-1} = Q$.

We study the problem of partitioning the equations $P \circ Q = R$. We assume that a specific pair of matrices R and Q in the equations $P \circ Q = R$ is given. Let each particular matrix P that satisfies $P \circ Q = R$ is called its solution and let S (Q, R) = {P | P ∘ Q = R} denote the set of all solutions (the solution set).

It is easy to see this problem can be partitioned, without loss of generality into a set of simpler problems expressed by the matrix equations $p_i \circ Q = r_i$ for all $i \in I$ where

$$P_i = [p_{ij} \mid j \in J] \text{ and}$$
$$r_i = [r_{ik} \mid k \in K].$$

Indeed each of the equation in $\max_{j \in J} \min (p_{ij}q_{jk}) = r_{ik}$ contains unknown $p_{ij}$ identified only by one particular value of the index i, that is, the unknown $p_{ij}$ distinguished by different values of i do not appear together in any of the individual equations. Observe that $p_i$, Q, and $r_i$ in $p_i \circ Q = r_i$ represent respectively, a fuzzy set on Y, a fuzzy relation on $Y \times Z$ and a fuzzy set on Z. Let $S_i (Q, r_i)$ = $[p_i \mid p_i \circ Q = r_i]$ denote, for each $i \in I$, the solution set of one of the simpler problem expressed by $p_i \circ Q = r_i$.

Thus the matrices P in S (Q, R) = [P | P ∘ Q = R ] can be viewed as one column matrix

$$P = \begin{bmatrix} p_1 \\ p_2 \\ \vdots \\ p_n \end{bmatrix}$$

where $p_i \in S_i (Q, r_i)$ for all $i \in I (=N_n)$. It follows immediately from $\max_{j \in J} \min (p_{ij} q_{jk}) = r_{ik}$. That if $\max_{j \in J} q_{jk} < r_{ik}$ for some $i \in I$ and some $k \in K$, then no values $p_{ij} \in [0, 1]$ exists $(j \in J)$ that satisfy $P \circ Q = R$, therefore no matrix P exists that satisfies the matrix equation. This proposition can be stated more concisely as follows if

$$\max_{j \in J} q_{jk} < \max_{j \in J} r_{ik}$$

for some $k \in K$ then S (Q, R) = $\phi$. This proposition allows us in certain cases to determine quickly that $P \circ Q = R$ has no solutions its negation however is only a necessary not sufficient condition for the existence of a solution of P ° Q = R that is for S (Q, R) ≠ ϕ. Since P ° Q = R can be partitioned without loss of generality into a set of equations of the form $p_i$ ° Q = $r_i$ we need only methods for solving equations of the later form in order to arrive at a solution. We may therefore restrict our further discussion of matrix equations of the form P ° Q = R to matrix equation of the simpler form P ° Q = r, where p = [$p_j$ | j ∈ J], Q = [$q_{jk}$ | j ∈ J, k ∈ K] and r = {$r_k$ | k ∈ K}.

We just recall the solution method as discussed by [43]. For the sake of consistency with our previous discussion, let us again assume that p, Q and r represent respectively a fuzzy set on Y, a fuzzy relation on Y × Z and a fuzzy set on Z. Moreover let J = $N_m$ and K = $N_s$ and let S (Q, r) = {p | p ° Q = r} denote the solution set of

$$p ° Q = r.$$

In order to describe a method of solving p ° Q = r we need to introduce some additional concepts and convenient notation. First let $\wp$ denote the set of all possible vectors.

$$p = \{p_j \mid j \in J\}$$

such that $p_j$ ∈ [0, 1] for all j ∈ J, and let a partial ordering on $\wp$ be defined as follows for any pair $p^1$, $p^2$ ∈ $\wp$ $p^1$ ≤ $p^2$ if and only if $p_i^2 \le p_j^2$ for all j ∈J. Given an arbitrary pair $p^1$, $p^2$ ∈ $\wp$ such that $p^1$ ≤ $p^2$ let [$p^1$, $p^2$] = {p ∈ $\wp$ | $p^1$ ≤ p < $p^2$}. For any pair $p^1$, $p^2$ ∈ $\wp$ ({$p^1$, $p^2$} ≤ } is a lattice.

Now we recall some of the properties of the solution set S (Q, r). Employing the partial ordering on $\wp$, let an element $\hat{p}$ of S (Q, r) be called a maximal solution of p ° Q = r if for all p ∈ S (Q, r), p ≥ $\hat{p}$ implies p = $\hat{p}$ if for all p ∈ S (Q, r) p < $\tilde{p}$ then that is the maximum solution. Similar discussion can be made on the minimal solution of p ° Q = r. The minimal solution is unique if p ≥ $\hat{p}$ (i.e. $\hat{p}$ is unique).

It is well known when ever the solution set S (Q, r) is not empty it always contains a unique maximum solution $\hat{p}$ and it may contain several minimal solution. Let $\breve{S}$ (Q, r) denote the set of all minimal solutions. It is known that the solution set S (Q, r) is fully characterized by the maximum and minimal solution in the sense that it consists exactly of the maximum solution $\hat{p}$ all the minimal solutions and all elements of $\wp$ that are between $\hat{p}$ and the numeral solution.

Thus S (Q, r) = $\underset{p}{\cup}$ $\left[ \tilde{p}, \hat{p} \right]$ where the union is taken for all $\tilde{p} \in \check{S}$ (Q, r). When S (Q, r) $\neq \phi$, the maximum solution.

$\hat{p} = [ \hat{p}_j \mid j \in J]$ of p $\circ$ Q = r is determined as follows:

$$\hat{p}_j = \min_{k \in K} \sigma (q_{ik}, r_k) \text{ where } \sigma (q_{jk}, r_k) = \begin{cases} r_k & \text{if } q_{jk} > r_k \\ 1 & \text{otherwise} \end{cases}$$

when $\hat{p}$ determined in this way does not satisfy p $\circ$ Q = r then S(Q, r) = $\phi$. That is the existence of the maximum solution $\hat{p}$ as determined by $\hat{p}_j = \min_{k \in K} \sigma (q_{ik}, r_k)$ is a necessary and sufficient condition for S (Q, r) $\neq \phi$. Once $\hat{p}$ is determined by $\hat{p}_j = \min_{k \in K} \sigma$ $(q_{ik}, r_k)$, we must check to see if it satisfies the given matrix equations p $\circ$ Q = r. If it does not then the equation has no solution (S (Q, r) = $\phi$), otherwise $\hat{p}$ in the maximum solution of the equation and we next determine the set $\check{S}$ (Q, r) of its minimal solutions.

## 1.3 Fuzzy compatibility relations and composition of fuzzy relations

In this section we recall the definition of fuzzy compatibility relations, fuzzy ordering relations, fuzzy morphisms, and sup and inf compositions of fuzzy relations.

**DEFINITION [43]:** *A binary relation R(X, X) that is reflexive and symmetric is usually called a compatibility relation or tolerance relation. When R(X, X) is a reflexive and symmetric fuzzy relation it is sometimes called proximity relation.*

An important concept associated with compatibility relations is compatibility classes. Given a crisp compatibility relation R(X, X), a compatibility class is a subset A of X such that (x, y) $\in$ R for all x, y $\in$ A. A maximal compatibility class or maximal compatible is a compatibility class that is not properly contained with in any other compatibility class. The family consisting of all the maximal compatibles induced by R on X is called a complete cover of X with respect to R. When R is a fuzzy compatibility relation, compatibility classes are defined in terms of a specified membership degree $\alpha$. An $\alpha$-compatibility class is a subset A of X such that R $(x, y) \geq \alpha$ for all x, y $\in$ A. Maximal $\alpha$-compatibles and complete $\alpha$-cover are obvious generalizations of the corresponding concepts for crisp compatibility relations.

Compatibility relations are often conveniently viewed as reflexive undirected graphs contrary to fuzzy cognitive maps that are directed graphs. In this context, reflexivity implies that each node of the graph has a loop connecting the node to itself the loops are usually omitted from the visual representations of the graph although they are assumed to be present. Connections between nodes as defined by the relation, are not directed, since the property of symmetry guarantees that all existing connections appear in both directions. Each connection is labeled with the value corresponding to the membership grade R $(x, y) =$ R $(y,x)$.

We illustrate this by the following example.

**Example 1.3.1:** Consider a fuzzy relation R(X, X) defined on X $= \{x_1, x_2,\ldots, x_8\}$ by the following membership matrix:

$$
\begin{array}{c c}
 & \begin{array}{c c c c c c c c} x_1 & x_2 & x_3 & x_4 & x_5 & x_6 & x_7 & x_8 \end{array} \\
\begin{array}{c} x_1 \\ x_2 \\ x_3 \\ x_4 \\ x_5 \\ x_6 \\ x_7 \\ x_8 \end{array} &
\left[ \begin{array}{c c c c c c c c}
1 & .3 & 0 & 0 & .4 & 0 & 0 & .6 \\
.3 & 1 & .5 & .3 & 0 & 0 & 0 & 0 \\
0 & .5 & 1 & 0 & 0 & .7 & .6 & .8 \\
0 & .3 & 0 & 1 & .2 & 0 & .7 & .5 \\
.4 & 0 & 0 & .2 & 1 & 0 & 0 & 0 \\
0 & 0 & .7 & 0 & 0 & 1 & .2 & 0 \\
0 & 0 & .6 & .7 & 0 & .2 & 1 & .8 \\
.6 & 0 & .8 & .5 & 0 & 0 & .8 & 1
\end{array} \right]
\end{array}.
$$

Since the matrix is symmetric and all entries on the main diagonal are equal to 1, the relation represented is reflexive, and symmetric therefore it is a compatibility relation. The graph of the relation is shown by the following figure 1.3.1, its complete $\alpha$-covers for $\alpha >$ 0 and $\alpha \in \{0, .3, .1, .4, .6, .5, .2, .8, .7, 1\}$ is depicted. Figure 1.3.1 is the graph of compatibility relation given in example 1.3.1 while similarity and compatibility relations are characterized by symmetry, ordering relations require asymmetry (or anti-symmetry) and transitivity. There are several types of ordering relations.

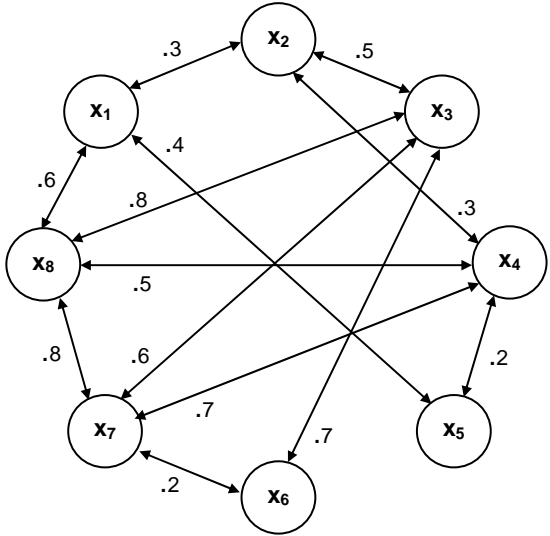

Figure: 1.2.1

A crisp binary relation R(X, X) that is reflexive, anti symmetric and transitive is called a partial ordering. The common symbol ≤ is suggestive of the properties of this class of relations. Thus x ≤ y denotes (x, y) ∈ R and signifies that x precedes y. The inverse partial ordering $R^{-1}$ (X, X) is suggested by the symbol ≥.

If y ≥ x including that (y, x) ∈ $R^{-1}$ then we say that y succeeds x. When x ≤ y; x is also referred to as a predecessor of y while y is called a successor of x. When x ≤ y and there is no z such that x ≤ y and z ≤ y, x is called an immediate predecessor of y and y is called an immediate successor of x. If we need to distinguish several partial orderings, such as P, Q and R we use the symbol $\overset{P}{\le}, \overset{Q}{\le}$ and $\overset{R}{\le}$ respectively.

Observe that a partial ordering '≤' on X does not guarantee that all pairs of elements x, y in X are comparable in the sense that either x ≤ y or y ≤ x. Thus, for some x, y ∈ X it is possible that x is neither a predecessor nor a successor of y. Such pairs are called non comparable with respect to ≤.

The following are definitions of some fundamental concepts associated with partial orderings:

1. If $x \in X$ and $x \leq y$ for every $y \in X$ then x is called the first member (or minimum) of X with respect to the relation denoted by $\leq$.

2. If $x \in X$ and $y \leq x$ for every $y \in X$, then x is called the last member (or maximum) of X with respect to the partial ordering relation.

3. If $x \in X$ and $y \leq x$ implies $x = y$, then x is called a minimal member of X with respect to the relation.

4. If $x \in X$ and $x \leq y$ implies $x = y$, then x is called a maximal member of X with respect to the relation [43]. Using these concepts every partial ordering satisfies the following properties:

    i. There exist at most one first member and at most one last member.
    ii. There may exist several maximal members and several minimal members.
    iii. If a first member exists then only one minimal member exists and it is identical with the first member.
    iv. If a last member exists, then only one maximal member exists and it is identical with the last member.
    v. The first and last members of a partial ordering relation correspond to the last and first members of the inverse partial ordering, respectively.

Let X again be a set on which a partial ordering is defined and let A be a subset of X $(A \subset X)$. If $x \in X$ and $x \leq y$ for every $y \in A$, then x is called a lower bound of A on X with respect to the partial ordering. If $x \in X$ and $y \leq x$ for every $y \in A$, then x is called an upper bound of A on X with respect to the relation. If a particular lower bound succeeds, every other lower bound of A, then it is called the greatest lower bound or infimum, of A. If a particular upper bound proceeds every other upper bound of A then it is called the least upper bound or superemum of A.

A partial ordering on a set X that contains a greatest lower bound and a least upper bound for every subset of two elements of X is called a lattice.

A partial ordering ≤ on X is said to be connected if and only if for all x, y ∈ X (x ≠ y) implies either x ≤ y or y ≤ x. When a partial ordering is connected all pairs of elements of X are comparable by the ordering such an ordering is usually called a linear ordering, some alternative names used in the literature are total ordering, simple ordering and complete ordering.

Partial ordering can be represented by diagrams and this sort of diagrams are called Hasse's diagrams.

A fuzzy binary relation R on a set X is a fuzzy partial ordering if and only if it is reflexive anti-symmetric and transitive under some form of fuzzy transitivity. Any fuzzy partial ordering based on max-min transitivity can be resolved into a series of crisp partial orderings in the same way in which this is done for similarity relations, that is by taking a series of α-cuts that produce increasing levels of refinement. When a fuzzy partial ordering is defined on a set X, two fuzzy sets are associated with each element x in X. The first is called the dominating class of x.

It is denoted by

$$R_{\geq [x]} \text{ and is defined by } R_{\geq [x]}(y) = R(x, y)$$

where y ∈ X. In other words the dominating class of x contains the members of X to the degree to which they dominate x. The second fuzzy set of concern is the class dominated by x, which is denoted by

$$R_{\leq [x]} \text{ and defined by } R_{\leq [x]}(y) = R(y, x)$$

where y ∈ X. The class dominated by x contains the elements of X to the degree to which they are dominated by x. An element x ∈ X is undominated if and only if R (x, y) = 0 for all y ∈ X and x ≠ y, an element x is undominating if and only if R (y,x) = 0 for all y ∈ X and y ≠ x. For a crisp subset A of a set X on which a fuzzy partial ordering R is defined, the fuzzy upper bound for A is the fuzzy set denoted by U (R, A) and defined by

$$U (R, A) = \bigcap_{x \in A} R_{\geq [x]}$$

where ∩ denotes an appropriate fuzzy intersection. This definition reduces to that of the conventional upper bound when the partial ordering is crisp. If a least upper bound of the set A exists it is the unique element x in U (R, A) such that

$$U (R, A) (x) > 0 \text{ and}$$

R (x, y) > 0 for all elements y in the support of U (R, A). Several other concepts of crisp orderings easily generalize to the fuzzy case. A fuzzy preordering is a fuzzy relation that is reflexive and transitive. Unlike a partial ordering, the preordering is not necessarily anti-symmetric.

A fuzzy weak ordering R is an ordering satisfying all the properties of a fuzzy linear ordering except anti-symmetry. Alternatively it can be thought of as a fuzzy preordering in which either R (x, y) > 0 or R (y, x) > 0 for all x ≠ y. A fuzzy strict ordering is anti-reflexive anti-symmetric and transitive, it can clearly be derived from any partial ordering R by replacing the values R(x, x) = 1 with zeros for all x.

Now we proceed on to recall the definition of fuzzy morphisms.

If two crisp binary relations R (X, X) and Q (Y, Y) are defined on sets X and Y, respectively then a function h : X → Y is said to be a homomorphism from (X, R) to (Y, Q) if $(x_1, x_2) \in R$ implies $(h(x_1), h(x_2)) \in Q$ for all $x_1, x_2 \in X$. In other words a homomorphism implies that for every two elements of the set X which are related under the relation R, their homomorphic images $h (x_1), h(x_2)$ in set Y are related under the relation Q.

When R (X, X) and Q (Y, Y) are fuzzy binary relations this implication can be generalized to R $(x_1, x_2) \le Q (h(x_1), h(x_2))$, for all $x_1, x_2 \in X$ and their images $h (x_1), h (x_2) \in Y$. Thus, the strength of relation between two elements under R is equated or exceeded by the strength of relation between their homomorphic images under Q.

Note that it is possible for a relation to exist under Q between the homomorphic images of two elements that are themselves unrelated under R. When this is never the case under a homomorphic function h, the function is called a strong homomorphism. It satisfies the two implications

$$\langle x_1, x_2 \rangle \in R \text{ implies } \langle h(x_1), h(x_2) \rangle \in Q$$

for all $x_1, x_2 \in X$ and $(y_1, y_2) \in Q$ implies $(x_1, x_2) \in R$, for all $y_1$, $y_2 \in Y$ where $x_1 \in h^{-1} (y_1)$ and $x_2 \in h^{-1} (y_2)$. Observe that when h is many to –one the inverse of h for each element of Y is a set of elements from X instead of a single element of X. If relations R (X, X) and Q (Y, Y) are fuzzy, then the criteria that a many-to-one function h must satisfy in order to be a strong homomorphism are somewhat modified. The function h imposes a partition. $\pi_h$ on the set X such that any two elements $x_1$, $x_2 \in X$ belong to the same block of the partition if and only if h maps them to the same element of Y. Let $A = \{a_1, a_2, \ldots, a_n\}$ and $B = \{b_1, b_2, \ldots, b_n\}$ be two blocks of this partition $\pi_h$ and let all elements of A be mapped to some element $y_1 \in Y$ and all elements of B to some element $y_2 \in Y$. Then the function h is said to be a strong homomorphism from $\langle X, R \rangle$ to $\langle Y, Q \rangle$ if and only if the degree of the strongest relation between any element of A and any element of B in the fuzzy relation R equals the strength of the relation between $y_1$ and $y_2$ in the fuzzy relation Q. Formally

$$\max_{ij} R(a_i\, b_j) = Q(y_1, y_2).$$

This equality must be satisfied for each pair of blocks of the partition $\pi_h$. If a homomorphism exists from (X, R) to (Y, Q) then the relation Q (X, Y) preserves some of the properties of the relation R (X, X) – namely, that all the pairs $(x_1, x_2) \in X \times X$ which are members of R have corresponding homomorphic images $\langle h(x_1), h(x_2) \rangle \in Y \times Y$ which are members of Q. Other members of Q may exist, however that are not the counterparts of any number of R. This is not the case when the homomorphism is strong. Here more properties of the relation R are preserved in relation Q. In fact Q represents a simplification of R in which elements belong to the same block of the partition $\pi_h$ created by the function h on the set X are no longer distinguished. These functions are useful for performing various kinds of simplifications of systems that preserve desirable properties in sets such as ordering or similarity.

If $h : X \rightarrow Y$ is a homomorphism from (X, R) to (Y, Q) and if h is completely specified, one to one and on to then it is called as isomorphism. This is effectively a translation or direct relabeling of elements of the set X into elements of the set Y that preserves all the properties of R in Q. If $Y \subseteq X$ then h is called an endomorphism. A function that is both an isomorphism and an endomorphism is called an automorphism. In this case the function maps the set X to itself and the relation R and Q are equal. Each of these terms applies without modification to fuzzy homomorphisms.

Now we just recall the notions of Sup-i-compositions of fuzzy relations. Sup-i compositions of binary fuzzy relations where i refers to a t-norm generalizes the standard max-min composition. The need to study these generalizations concepts from some applications such as approximate reasoning and fuzzy control. Given a particular t-norm i and two fuzzy relations P (X, Y) and Q (Y, Z), the Sup-i-composition of P and Q is a fuzzy relation $P^i \circ Q$ on X, Y, Z defined by

$$\{P^i \circ Q\} \, (x, z) = \underset{y \in Y}{Sup} \, i \, [P(x, y), Q(y,z)]$$

for all $x \in X$, $z \in Z$. When i is chosen to be the min operator $P \overset{i}{\underset{o}{}} Q$ becomes the standard composition $P \circ Q$.

Given fuzzy relations $P(X, Y)$, $P_j(X, Y)$, $Q(Y, Z)$, $Q_j(Y, Z)$ and $R(Z, V)$ where j takes values in an index set J, the following are basic properties of the Sup-i composition under the standard fuzzy union and intersections.

$$(P \overset{i}{\underset{o}{}} Q) \overset{i}{\underset{o}{}} R = P \overset{i}{\underset{o}{}} (Q \overset{i}{\underset{o}{}} R)$$

$$P \overset{i}{\underset{o}{}} \left( \bigcup_{j \in J} Q_j \right) = \bigcup_{j \in J} (P \underset{o}{} Q_j)$$

$$P \underset{o}{i} \left( \bigcap_{j \in J} Q_j \right) \subseteq \bigcap_{j \in J} (P \underset{o}{i} Q_j)$$

$$\bigcup_{j \in J} P_j \underset{o}{i} Q = \bigcup_{j \in J} (P_j \underset{o}{i} Q)$$

$$\bigcap_{j \in J} P_j \underset{o}{i} Q \subseteq \bigcap_{j \in J} (P_j \underset{o}{i} Q)$$

$$(P \underset{o}{i} Q)^{-1} = Q^{-1} \underset{o}{i} P^{-1}.$$

These properties follow directly from the corresponding properties of t-norms and their verification is left for the reader.

Sup i composition is also monotonic increasing that is for any fuzzy relation $P(X, Y)$, $Q$, $(Y, Z)$, $Q_2$ $(Y, Z)$ and $R(Z, V)$ if $Q_1 \subset Q$ then

$$P \underset{o}{i} Q_1 \subseteq P \underset{o}{i} Q_2$$

$$Q_1 \underset{o}{i} R \subseteq Q_2 \underset{o}{i} R.$$

The set of all binary fuzzy relations on $X^2$ forms a complete lattice ordered monoid $(\Im(X^2), \cap, \cup, \underset{o}{i})$ where $\cap$ and $\cup$ represent the meet and join of the lattice respectively and $\underset{o}{i}$ represents, the semi group operation of the monoid. The identity of $\underset{o}{i}$ is defined by the relations

$$E(x, y) = \begin{cases} 1 \ when \ x = y \\ 0 \ when \ x \neq y \end{cases}.$$

The concept of transitivity of a fuzzy relation, which is introduced in terms of the max-min composition can be generalized in terms of the Sup-i-compositions for the various t-norms i. We say the relation R on $X^2$ is i-transitive if and only if

$$R(x, z) \geq i[R(x, y), R(y, z)]$$

for all x, y, z $\in$ X. It is easy to show that a fuzzy relation R on $X^2$ is i-transitive if and only if $R \underset{o}{i} R \subseteq R$ which may be used as an alternative definition of i-transitivity.

When a relation R is not i-transitive, we define its i-transitive closure as a relation $R_{\tau(i)}$ that is the smallest i-transitive relation, containing R. To investigate properties of i-transitive closure let

$$R^{(n)} = R \underset{o}{i} R^{(n-1)}$$

n = 2, 3,… where R is a fuzzy relation on $X^2$ and $R^{(1)}$ = R. Using this notation the reader is expected to prove the following theorem:

**THEOREM [43]:** *For any fuzzy relation R on $X^2$, the fuzzy relation $R_{\tau(i)} = \bigcup_{n=1}^{\infty} R^{(n)}$ is the i-transitive closure of R.*

Prove if R be a reflexive fuzzy relation on $X^2$, where $|X| = n \geq 2$ then $R_{\tau(i)} = R^{(n-1)}$.

Now we proceed on to recall the notion of inf-$w_1$ compositions of fuzzy relations. Give a continuous t-norm i, let $w_i(a, b) = \sup \{x \in [0, 1] \ / \ i(a, x) \leq b\}$ for every a, b $\in$ [ 0, 1]. This operation referred to as operation $w_i$ plays an important role in fuzzy relation equations. While the t-norm i may be interpreted as logical conjunction, the corresponding operation $w_i$ may be interpreted as logical implication. The following basic properties of $w_i$ are left as an exercise for the reader to prove. For any a, $a_j$, b, d $\in$ [0, 1] where j takes values from an index set J prove $w_i$ has the following properties:

1.     $i\,(a, b) \leq d$   iff   $w_i(a, d) \geq b$
2.     $w_i\,(w_i\,(a, b) \geq a$
3.     $w_i\,(i\,(a, b), d) = w_i\,(a, w_i\,(b, d))$
4.     $a \leq d$ implies
    $$w_i\,(a, d) \geq w_i\,(b, d) \text{ and}$$
    $$w_i\,(d, a) \leq w_i\,(d, b)$$
5.     $i\,[w_i\,(a, b), w_i\,(b, d)] \leq w_i\,(a, d)$
6.     $w_i\,[\inf a_j, b] \geq \sup\limits_{j \in J} w_i\,(a_j, b)$
7.     $w_i\,(\sup a_j, b) = \inf\limits_{j \in J} w_i\,(a_j, b)$
8.     $w_i\,[b, \sup\limits_{j \in J} a_j] \geq \sup\limits_{j \in J} w_i\,(b, a_j)$
9.     $w_i\,(b, \inf\limits_{j \in J} a_j) = \inf\limits_{j \in J} w_i\,(b, a_j)$
10.    $i\,[a, w_i\,(a, b)] \leq b.$

Prove for the fuzzy relations P(X, Y), Q (Y, Z), R(X, Z) and S(Z, V) the following are equivalent

$$P \underset{o}{i} Q \subseteq R$$
$$Q \subseteq P^{-1} \underset{o}{w_i} R$$
$$P \subseteq (Q \underset{o}{w_i} R^{-1})^{-1}$$

Prove $P \underset{o}{w_i} (Q \underset{o}{w_i} S) = (P \underset{o}{w_i} Q) \underset{o}{w_i} S.$

Let P (X, Y), $P_j$ (X, Y), Q(Y, Z) and $Q_j$ (Y, Z) be fuzzy relations where j takes values in an index set J then prove.

$$\left( \bigcup_{j \in J} P_j \right) \underset{o}{w_i}\, Q = \bigcap_{j \in J} \left( P_j \underset{o}{w_i}\, Q \right)$$

$$\left( \bigcap_{j \in J} P_j \right) \underset{o}{w_i}\, Q = \supseteq \bigcup_{j \in J} P_j \underset{o}{w_i}\, Q$$

$$P \underset{o}{w_i} \left( \bigcap_{j \in J} Q_j \right) = \bigcap_{j \in J} P \underset{o}{w_i}\, Q_j$$

$$P \underset{o}{w_i} \bigcup_{j \in J} Q_j \supseteq \bigcup_{j \in J} P \underset{o}{w_i} Q_j .$$

Let P (X, Y), $Q_1$(Y, Z), $Q_2$ (Y, Z) and R (Z, V) be fuzzy relations. If $Q_1 \subseteq Q_2$ then prove.

$$P \underset{o}{w_i} Q_1 \subseteq P \underset{o}{w_i} Q_2 \text{ and}$$

$$Q_1 \underset{o}{w_i} R \supseteq Q_2 \underset{o}{w_i} R.$$

Now if P (X, Y), Q (Y, Z) and R (X, Z) be fuzzy relations prove

$$P^{-1} \underset{o}{i} (P \underset{o}{w_i} Q) \subseteq Q$$

$$R \subseteq P \underset{o}{w_i} (P^{-1} \underset{o}{i} R)$$

$$P \subseteq (P \underset{o}{w_i} Q) \underset{o}{w_i} Q^{-1}$$

$$R \subseteq (R \underset{o}{w_i} Q^{-1}) \underset{o}{w_i} Q.$$

## 1.4 Optimization Of FRE with Max-Product Composition

Jiranut Loelamonphing and Shu Cheng Fang [58] in their paper "optimization of fuzzy relation equation with max-product composition" (2001) has studied the solution set of fuzzy relation equations with max product composition and an optimization problem with a linear objective function subject to such FRE. By identifying the special properties of the feasible domain they determine an optimal solution without explicitly generating the whole set of minimal solutions.

The notion of FRE based upon the max-min composition was first investigated by [84]. He studied conditions and theoretical methods to resolve fuzzy relations on fuzzy sets defined as mappings from sets into complete Brouwerian lattices. Some theorems for existence and determination of solutions of certain basic fuzzy relation equations were presented in his work. However, the solution obtained in that work is only the greatest element (or the maximum solution) derived from the max-min (or min-max) composition of fuzzy relations.

Sanchez's work [84] has shed some light on this important subject. Since then, researchers have been trying to explore the problem and develop solution procedures [1, 11, 16, 24, 30, 82].

The "max-min" composition [117] is commonly used when a system requires conservative solutions in the sense that the goodness of one value cannot compensate the badness of another value. In reality, there are situations that allow compensatability among the values of a solution vector. In this case, the min operator is not the best choice for the intersection of fuzzy sets. Instead, the "max-product" composition is preferred since it can yield better, or at least equivalent, results. Note that when the intersection connector acts non-interactively, it can be uniquely defined by the min connector, but when the connector is interactive, it is application dependent and cannot be defined universally. Some outlines for selecting an appropriate connector has been provided by [112, 113].

Recently, researchers extended the study of an inverse solution of a system of FRE with max-product composition. They provided theoretical results for determining the complete solution sets as well as the conditions for the existence of resolutions. Their results showed that such complete solution sets can be characterized by one maximum solution and a number of minimal solutions. Since the total number of minimal solutions has a combinatorial nature in terms of the problem size, an efficient solution procedure is always in demand.

Motivated by the work of [24], we are interested in studying the optimization problem with a linear objective function subject to a system of FRE with the "max-product" composition.

Let $A = [a_{ij}]$, $0 \leq a_{ij} \leq 1$, be an $(m \times n)$ – dimensional fuzzy matrix, $b = (b_1, \ldots, b_n)$, $0 \leq b_j \leq 1$, be an n-dimensional vector, $I = \{1, 2, \ldots, m\}$ and $J = \{1, 2, \ldots, n\}$. A system of FRE defined by A and b is denoted by

$$X \text{ o } A = b, \tag{1}$$

where "o" represents the max-product composition. The resolution of (1) is a set of solution vectors $x = (x_1, \ldots, x_m)$, $0 \leq x_j \leq 1$, such that

$$\max_{i \in I} \{x_i . a_{ij}\} = b_j \text{ for } j \in J \tag{2}$$

Let $c = (c_1, \ldots, c_m) \in R^m$ be an m-dimensional vector where $c_i$ represents the weight (or cost) associated with variable $x_i$, for $i \in I$. The optimization problem we are interested in has the following form:

Minimize

$$Z = \sum_{i=1}^{m} c_i x_i \qquad (3)$$

Subject to x o A = b,

$$0 \le x_j \le 1.$$

Note that the characteristics of the solution sets obtained by using the max-min operator and the max-product operator are similar, i.e., when the solution set is not empty, it can be completely determined by a unique maximum solution and a finite number of minimal solutions [11, 34]. Since the solution set can be non-convex, traditional linear programming methods, such as the simplex and interior-point algorithms, cannot be applied to this problem.

[58] denote the solution set of problem (1) by X (A, b) = $\{(x_1, \ldots, x_m) \mid x_i \in [0, 1], i \in I$, and x o A = b$\}$. Define X = $[0, 1]^m$. For $x^1, x^2 \in X$, we say $x^1 \le x^2$ if and only if $x_i^1 \le x_i^2$, $\forall i \in I$. In this way, the operator "≤" forms a partial order relation on X and (X, ≤) becomes a lattice.

$\hat{x} \in$ X (A, b) is called a maximum solution if x ≤ $\hat{x}$ for all x ∈ X (A, b). Also $\bar{x} \in$ X (A, b) is called minimal solution if x ≤ $\bar{x}$, for any x ∈ X (A, b) implies x = $\bar{x}$. When X (A, b) is not empty, it can be completely determined by a unique maximum solution and a finite number of minimal solutions [11, 34].

The maximum solution can be obtained by applying the following operation [58]:

$$\hat{x} = A \, \Theta \, b = \left[ \bigwedge_{j=1}^{n} (a_{ij} \, \Theta \, b_j) \right]_{i \in I} \qquad (4)$$

where

$$a_{ij} \, \Theta \, b_j = \begin{cases} 1 & \text{if } a_{ij} \le b_j \\ b_j / a_{ij} & \text{if } a_{ij} > b_j \end{cases} \qquad (5)$$

$$a \wedge b = min \ (a, \ b).$$

Denote the set of all minimal solutions by $\breve{X}$ (A, b), the complete set of solution, X (A, b), is obtained by

$$X(A, b) = \bigcup_{\bar{x} \in X(A, b)} \{x \in X \mid \bar{x} \leq x \leq \hat{x}\} \qquad (6)$$

**DEFINITION 1.4.1:** *For a solution* $x \in X(A, b)$, *we call* $x_{i_0}$ *a binding variable if* $x_{i_0}$. $a_{i_0 j} = b_j$ *for* $i_0 \in I$ *and* $x_i a_{ij} \leq b_j$, *for all i* $\in I$.

*When the solution set of (1) is not empty, i.e.,* $X(A, b) \neq \emptyset$, *we define*

$$I_j = \{ i \in I./ \; \hat{x}_i , a_{ij} = b_j\}, \; \forall j \in J, \qquad (7)$$

$$\Lambda = I_1 \times I_2 \times ... \times I_n. \qquad (8)$$

*Here* $I_j$ *corresponds to a set of* $x_i$'s *that can satisfy constraint j of the fuzzy relation equations. And, the set* $\Lambda$ *represents all combinations of the binding variables such that every combination can satisfy every fuzzy relation constraint. Let each combination be represented by* $f = (f_1, f_2, ..., f_n) \in \Lambda$, *with* $f_j \in I_j$, $\forall j \in J$.

The optimization problem (3) can be decomposed into two problems, namely

$$\text{Minimize} \qquad Z' = \sum_{i=1}^{m} c'_i x_i \qquad (9)$$

Subject to     $x \circ A = b$,
$\qquad\qquad 0 \leq x_i \leq 1$

and

$$\text{Minimize} \qquad Z'' = \sum_{i=1}^{m} c''_i x_i \qquad (10)$$

Subject to     $x \circ A = b$,
$\qquad\qquad 0 \leq x_i \leq 1$

where $c' = (c'_1, c'_2, ..., c'_m)$ and $c'' = (c''_1, c''_2, ..., c''_m)$ are defined such that , $\forall i \in I$.

$$c'_i = \begin{cases} c_i & if \ c_i \geq 0, \\ 0 & if \ c_i < 0, \end{cases}$$

$$c''_i = \begin{cases} 0 & if \ c_i \geq 0, \\ c_i & if \ c_i < 0. \end{cases}$$

(11)

Apparently, the cost vector c = c' + c" and the objective value Z' = Z' + Z". Intuitively, when all the costs are non-positive, since $x_i$'s are non-negative and the problem is to minimize the objective value, we should make $x_i$ as large as possible.

Taking advantage of the special structure studied in the previous section, we now introduce some procedures to reduce the size of the original problem so that the effort to solve the problem is minimized. The key idea behind these reduction procedures is that some of the $x_i$'s can be determined immediately without solving the problem but just by identifying the special characteristic of the problem. Special cases which we can eliminate from consideration are as follows:

**Case I: $c_i \leq 0$.**

We know that $x^*_I = \hat{x}_i$, if $c_i \leq 0$. Hence, we can take any part that are associated with these $\hat{x}_i$'s out of consideration.

Here we define:

$$\hat{I} = \{i \in I \mid c_i \leq 0\}, \tag{12}$$

$$\hat{J} = \{j \in J \mid \hat{x}_i, a_{ij} = b_j, \forall \ i \in \hat{I} \}. \tag{13}$$

In other words, $\hat{J}$ is a set of indices of constraints which can be satisfied by a set of $\hat{x}_i$'s for $i \in \hat{I}$. Having defined $\hat{J}$ and $\hat{I}$, we now eliminate row i, $i \in \hat{I}$ and j, $j \in \hat{J}$ from matrix A as well as the j[th] element $j \in \hat{J}$, from vector b.

Let A' and b' be the updated fuzzy matrix and fuzzy vector, respectively.

Define J' = J \ $\hat{J}$, J' represents a reduced set of constraints.

**Case II : $I_j$ has only one element.**

Consider constraint $j \in J'$. If $I_j$ contains only one element, it means that only one $x_i$, $i \in I_j$, can satisfy the jth constraint. We have $x_i = b_j/a_{ij}$.

Define

$$\bar{I} = \{i \in I_j \mid \|I_j\| = 1; \, j \in J'\} \tag{14}$$

$$\bar{J} = \{j \in J' \mid x_i.a_{ij} = b_j; \, i \in \bar{I}\}. \tag{15}$$

Again, we can eliminate row i, $i \in \bar{I}$, and column j, $j \in \bar{J}$, from the updated fuzzy matrix A' as well as the $j^{th}$ element $j \in \bar{J}$, from the updated vector b'. Let A" and b" be the reduced fuzzy matrix and fuzzy vector corresponding to A' and b' respectively. We also need to update $\Lambda$. Define $J" = J' \setminus \bar{J}$, J" is an index set of constraints which need to be solved later by the branch-and-bound (B&B) method. The updated A" = $\Pi_{j \in J"} I_j$.

The branch-and-bounded algorithm will be performed on these A" and b". If b" is empty, then all constraints have been taken care of. Therefore, in order to minimize the objective value, since we are now left with positive $c_i$'s, we can assign the minimum value, i.e. zero, to all $x_i$'s whose values have not been assigned yet. When b" is not empty, we need to proceed further. Details will be discussed in the following:

In order to identity whether the problem is decomposable, consider a set of constraints, say B, which can be satisfied by a certain set of variables, say $X_B$. If the decision to choose which variable in the set $X_B$ to satisfy a constraint in B does not impact the decision on the rest of the problem, then we can extract this part from the whole problem.

Let k be the number of sub-problems, $1 \leq k \leq \|J"\|$.

Define

$$\Omega = \{I_j \mid j \in J"\}, \tag{16}$$

$$\Omega_l = \left\{ I_j \middle| \bigcap_{j \in J"} I_j \neq \phi \right\}, \, l = 1, \ldots, k, \tag{17}$$

$$\Omega_l \cap \Omega_{l'} = \phi, \, l \neq l', \tag{18}$$

$$\Omega = \Omega_1 \cup \Omega_2 \cup \ldots \cup \Omega_k, \tag{19}$$

$$\Lambda_l = \prod_{I_j \in \Omega_l} I_j, \tag{20}$$

$$I^{(1)} = \{ \, i \mid i \in I_j, \, I_j \in \Omega_1 \}, \tag{21}$$

$$J^{(1)} = \{ j \mid I_j \in \Omega_1 \}. \tag{22}$$

In this way, $\Omega_1$ contains sets of $I_j$'s which have some element(s) in common and we can decompose the original problem into k sub-problems. $I^{(1)}$ and $J^{(1)}$ correspond to sets of indices of variables and constraints, respectively, on which the B&B method is performed for sub-problem l [58].

Problem (9) can be transformed into the following 0-1 integer-programming problem

$$\text{minimize } Z' = \sum_{i=1}^{m} \left( c'_i \max_{j \in J} \left\{ \frac{b_i}{a_{ij}}, x_{ij} \right\} \right) \tag{23}$$

$$\text{subject to } \sum_{i=1}^{m} x_{ij} = 1 \qquad \forall j \in J$$

$x_{ij} = 0$ or $1 \; \forall \, i \in I, \, j \in J$.
$x_{ij} = 0 \; \forall \, i, j$ with $i \in j$.

### ALGORITHM 1 (ALGORITHM FOR FINDING AN OPTIMAL SOLUTION)

*Step 1*: Find the maximum solution of (1).
Compute $\hat{x} = A \otimes b = \left[ \Lambda_{j=1}^{n} (a_{ij} \otimes b_j) \right]_{i \in I}$ according to (4).

*Step 2*: Check feasibility.
If $\hat{x} \circ A = b$, continue. Otherwise, stop! $X(A, b) = \phi$ and problem (3) has no feasible solution.

*Step 3*: Compute index sets.
Compute $I_j = \left\{ i \in I \mid \hat{x}_i . a_{ij} = b_j \right\}$, $\forall \, j \in J$, which represents a set of $x_i$'s that can satisfy constraint j of the FRE.

*Step 4*: Arrange cost vector.
Define c' and c' according to (11).

*Step 5*: Perform problem reduction.
Compute $\hat{I} = \{i \in I \mid c_i \le 0\}$ and $\hat{J} = \{j \in J \mid \hat{x}_i . a_{ij} = b_j; i \in \hat{I}\}$.
Eliminate row $i \in \hat{I}$, and column j, $j \in \hat{J}$, from matrix A to obtain A'. Also eliminate the jth element, $j \in \hat{J}$, from vector b to obtain b'. Assign an optimal value $x^*_i = \hat{x}_i$, for $i \in \hat{I}$. If b' is empty, assign zero to unassigned $x^*_i = \hat{x}_i$, for $i \in \hat{I}$. If b' is empty assign zero to unassigned $x^*_i$ and go to Step 11. Otherwise, compute J' = J\ $\hat{J}$ and proceed to the next step.

*Step 6*: Find singleton $I_j$.
Compute $\bar{I} = \{i \in I_j \mid \|I_j\| = 1; j \in J'\}$ and $\bar{J}$ = {j ∈ J' | $x_i$. $a_{ij}$ = $b_j$; $i \in \bar{I}$ }. Eliminate row i, $i \in \bar{I}$ and column j, $j \in \bar{J}$ from matrix A' to obtain A''. Also eliminate the jth element, $j \in \bar{J}$, from vector b' to obtain b''. Assign $x^*_I = b_j / a_{ij}$, for $i \in \bar{I}$ and $i \in I_j$. If b'' is empty assign zero to unassigned $x^*_I$ and go to Step 11. Otherwise, compute J'' = J'\ $\bar{J}$ and proceed to the next step.

*Step 7*: Decompose the problem.
Decompose the problem by computing equations (16-22).

*Step 8*: Define sub-problems.
For each sub-problem l, define problem (11) and its corresponding 0-1 inter program using (23).

*Step 9*: Solve the integer program(s).
Solve each integer program by using the branch-and-bound method.

*Step 10*: Generate an optimal solution of sub-problem.
For each sub-problem l, define $f^l = (f_j)$, j ∈ $J^{(1)}$ with $f_j = i$ if $x_{ij} = 1$. Generate $F(f^*)$ via formula (10). Define $\bar{x}^* = \left(\bar{x}^*_1, ..., \bar{x}^*_m\right)$ with $\bar{x}^*_i = F_i(f^*)$.

*Step 11*: Generate an optimal solution.

Combine $\tilde{x}^*$ with the solution obtained from (5) and (6) to yield an optimal solution of problem (3).

For more refer [58].

## 1.5 Composite FRE-resolution based on Archimedean triangular norms

The resolution problem of FRE is one of the most important and widely studied problems in the field of fuzzy sets and fuzzy systems. The first step for the resolution of a FRE is to establish the existence of the solution. If the equation is solvable the solution set contains a maximum solution and possibly several minimum solutions. It has been proved that the finding of these solutions suffices for the finding of solution set.

[92] present sup t FREs. They prove in most practical cases the solution set of sup t FREs is non-empty and provide some new criteria for checking the existence of the solution.

Introducing the '*solution matrices*' formulation of the problem, we find an "if and only if" condition for the solution existence. Then, after a brief description of the most convenient algorithm for solving sup-t FREs proposed by [9], a fast algorithm is given for determining the solution set of sup-t FREs, where t is an Archimedean t-norm.

Let X, Y, Z be discrete crisp sets with cardinalities *n, m* and *k,* respectively, and A(X, Y), R(Y, Z), B(X, Z) be three binary fuzzy relations constraining its other with the relationship.

$$A \ o^t \ R = B, \tag{1}$$

where $o^t$ is the well-known sup-t composition (t is a triangular norm). Eq. (1) can be written in the matrix form

$$A \ o^t R = B, \tag{2}$$

where $A_{n \ x \ m}$, $R_{m \ x \ k}$, $B_{n \ x \ k}$ are the matrix representations of A, R, B, respectively. Eq. (2) is the typical form of a fuzzy relation equation (FRE) for which the following problems arise:

(i)     the resolution of (2) for R, when A and B are known, and

(ii)        the resolution of (2) for A, when R and B are known (inverse problem).

It is remarked that if we have a method for solving the first problem, using the same method for the equation $R^{-1} o^t A^{-1} = B^{-1}$ that employs transposed matrices, the second problem could be solved. Thus, without loss generality, we consider only the first problem. Moreover, (2) is actually a set of k simpler fuzzy relation equations that can be solved independently and so it suffices to consider only the equation

$$A \, o^t \, r = b, \tag{3}$$

Where $r_{m \times 1}$ and $b_{n \times 1}$ are column vectors of R and B, respectively. Clearly, (3) is a system of n equations of the form

$$a \, o^t \, r = b, \tag{4}$$

where a is a row vector of A. One can easily see that (3) has solution for r if and only if (iff) all the n equations of form (4) have at least one common solution for r.

Let S (A, b) be the solution set of (3), i.e. S(A, b) = {r : A $o^t$ r = b}. It is well known that if S(A, b) ≠ φ, then it contains a unique maximal solution $\hat{r}$ and may contain several minimal solutions $\check{r}$ [43]. The solution set is the union of all the lattices [$\check{r}$, $\hat{r}$] between each minimal and the maximum solution.

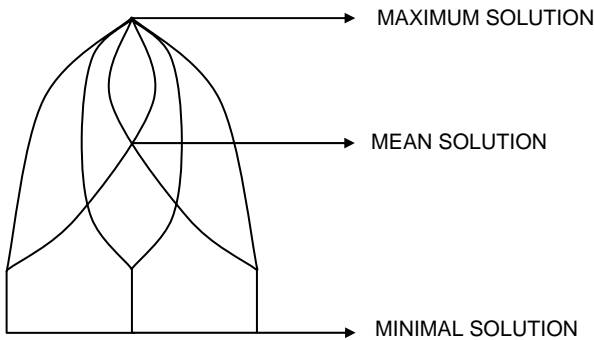

Figure: 1.5.1

In Figure 1.5.1, a view of the solution set is illustrated. We define the mean solution, as the minimal element of the intersection of the lattices $[\check{r},\hat{r}]$. Obviously, the mean solution always exists (since $\bar{r}$ always belongs to the intersection of the lattices $[\check{r},\hat{r}]$) and it is unique.

The maximum, the mean and the minimum solutions of (2) come from the respective solutions of (3) with the aid of the following:

(i) The maximum solution is the m × k matrix

$$[\hat{r}_1\,\hat{r}_2\ldots\hat{r}_k],\tag{5}$$

where $\hat{r}_i$ (i = 1,2,…,k) is the maximum solution of the ith equation of form (3).

(ii) The mean solution is the m × k matrix

$$\left[\bar{r}_1\ \bar{r}_2\ \cdots\ \bar{r}_k\right],\tag{6}$$

where $\bar{r}_i$ (i = 1,2,…k) is the mean solution of the i$^{th}$ equation of form (3).

(iii) The minimum solutions are the m × k matrices

$$[\check{r}_1\ \check{r}_2\ldots\check{r}_k],\tag{7}$$

where $\check{r}_i \in \breve{S}_i$ (i = 1,2,…, k) and $\breve{S}_i$ is the minimal solutions set of the i$^{th}$ equation of form (3).

The basis on which equations of form (3) can be solved is the simple equation of the form.

$$t\,(a,\,x) = b\tag{8}$$

where a and b are given Eq. (8) is actually a special case of (3), for n = m = 1. Our purpose in this section is the study of (8).

A function t = [0, 1] × [0, 1] $\to$ [0, 1] is a t-norm iff $\forall$a, b, d $\in$ [0, 1] it satisfies the following four axioms (axiomatic skeleton for t-norms):

*Axiom 1*: t (a, 1) = a and t (a, 0) = 0
(boundary condition).

*Axiom 2*: b ≤ d implies t (a, b) ≤ t (a, d) (monotonicity).

*Axiom 3*: t (a, b) = t (b, a) (commutativity).

*Axiom 4*: t (a, t (b, d)) = t (t (a, b), d) (associativity).

A  t-norm t is called Archimedean iff

*Axiom 5*: t is a continuous function.

*Axiom 6*: t (a, a) < a, $\forall a \in (0, 1)$ (subidempotency).

The class of t-norms has been widely studied by many researchers. It is proved in [43] that min (a, b) is the only idempotent t-norm. On the other hand, for any t-norm t it is true that [43]

$$t (a, b) \leq \min (a, b).$$

Thus, the only continuous t-norm that is not Archimedean is the min (a, b). This is a very interesting result, since it suggests a separate study of (3) when t is an Archimedean t-norm and when t is the min (a, b).

Equations of the form (8) do not always have a solution and when they do have one, it may not be unique. The following proposition can be easily proved.

**Proposition 1.5.1:** Let t be a continuous t-norm and a, b, x ∈ [0, 1]. The equation t(A, x) = b has a solution for x iff a ≥ b.

Let $a \,\widehat{\otimes}^{t}\, b$ and $a \,\widecheck{\otimes}^{t}\, b$ denote the maximal and the minimal, respectively, solution of (8) (if exist), i.e.

$$a \,\widehat{\otimes}^{t}\, b = \sup \{x \in [0, 1]: t (a, x) = b\},$$
$$a \,\widecheck{\otimes}^{t}\, b = \inf \{x \in [0, 1]: t (a, x) = b\}.$$

If the solution is unique, it is denoted by a $\otimes^t$ b. Based on the above notations, the maximal solution operator (max-SO) $\hat{\omega}_t$ and the minimal solution operator. (min-SO) $\breve{\omega}_t$ are defined as follows:

$$\hat{\omega}_t \ (a, b) = \begin{cases} 1, & a < b, \\ a \, \hat{\otimes}^t \, b, & a \geq b, \end{cases} \qquad (9)$$

$$\breve{\omega}_t (a, b) = \begin{cases} 0, & a < b, \\ a \, \breve{\otimes}^t \, b, & a \geq b, \end{cases} \qquad (10)$$

when (8) has a unique solution, min-SO and max-SO take the form

$$\hat{\omega}_t \ (a, b) = \begin{cases} 1, & a < b, \\ a \otimes^t b, & a \geq b, \end{cases}$$

$$\breve{\omega}_t (a, b) = \begin{cases} 0, & a < b, \\ a \otimes^t b, & a \geq b, \end{cases}$$

Max-SO and min-SO are extensions of the $\alpha$ and the $\sigma$ operators defined by Sanchez [84] for the min t-norm. The max-SO extension for any t-norm has been proposed by Di Nola et al [17]. [9] proposed the min-SO extension. The definition of max-SO given in [43] is

$$\hat{\omega}_t \ (a, b) = \sup \{x \in [0, 1]: t \, (a, x) \leq b\}. \qquad (11)$$

For any a, b $\in$ [0, 1], it is $\breve{\omega}_{\min} \ (a,b) \leq a$ and $\breve{\omega}_{\min} \ (a, b) \leq b$.

For any a, b, d $\in$ [0, 1] with a < b , it is $\breve{\omega}_{\min} \ (a,b) \leq \breve{\omega}_{\min} (b,d)$.

For any a, b $\in$ [0, 1] and any continuous t-norm t it is $t \, (a \, \breve{\omega}_t \, (a, b)) \leq b$.

For any a, b $\in$ [0, 1] and any continuous t-norm t it is $\hat{\omega}_t \ (a, b) \geq \breve{\omega}_t (a, b)$.

Here the resolution of (3) is studied. First, some conditions for the existence of the solution are given and then a resolution method is described. [9] has originally proposed the method. It is formulated in a different way, based on the solution matrices, in order to be more comprehensive. Moreover, some new results are given on the form of the solution matrices as well as on the solution of existence problem that help to clear up the underlying mechanism of the resolution process of FREs.

Let us first proceed with the solution existence problem. Eq. (3) has a solution iff all its equations of form (4) have a common solution. The following lemma can be established.

Let a $o^t$ r = b a FRE of form (4). We have S (a, b) $\neq \phi$ iff there exists j $\in N_m$ such that $a_j \geq b$.

Let A $o^t$ r = b be FRE of form (3) with m $\geq$ n. If for any i $\in N_n$ there exists j $\in N_m$ such that $A_{ij} \geq b_i$ and $A_{kj} \leq b_k$, $\forall k \in N_n - \{I\}$, then S (A, b) $\neq \phi$.

Let t be a continuous t-norm and A $o^t$ r = b be a FRE of form (3). The matrix $\overline{\Gamma}_{mxn}$ is the mean solution matrix (mean –SM) of the FRE, where

$$\widehat{\Gamma}_{ij} = \breve{\omega}_t \ (A_{ji,} \ b_j), \ \forall i \in N_m, j \in N_n,$$

The matrix $\widehat{\Gamma}$ is the maximal solution matrix (max-SM) of the FRM, where

$$\breve{\Gamma}_{ij} = \breve{\omega}_t \ (A_{ji,} \ b_j), \ \forall i \in N_m, j \in N_n,$$

The matrix $\widehat{\Gamma}$ is the minimal solution matrix (min-SM) of the FRM, where

$$\breve{\Gamma}_{ij} = \widehat{\omega}_t \ (A_{ji,} \ b_j), \ \forall i \in N_m, j \in N_n,$$

The matrix $\breve{\Gamma}$ is the minimal solution matrix (min-SM) of the FRE, where

$$\breve{\Gamma}_{ij} = \widehat{\omega}_{\min} \left( \inf_{k \in N_n} \widehat{\Gamma}_{ik}, \overline{\Gamma}_{ij} \right), \forall i \in N_m, j \in N_n,$$

Let t be a continuous t-norm and A $o^t$ r = b be a FRE of form (3). If S ($A_{j\bullet}$, $b_j$) $\neq \phi$ for any j $\in N_n$, then

$$\widehat{\Gamma}_{\bullet j} \in S \ (A_{j\bullet}, b_j), \forall j \in N_n$$
$$\overline{\Gamma}_{\bullet j} \in S \ (A_{j\bullet}, b_j), \forall j \in N_n$$

where $\hat{\Gamma}$ and $\overline{\Gamma}$ is the min-SM and the mean-SM of the equation, respectively.

Let A $o^t$ r = b be a FRE of form (3), where t is a continuous t-norm. The following propositions are equivalent:

(i) S (A, b) = $\phi$;
(ii) There exists j $\in$ $N_n$, such that $\breve{\Gamma}_{\bullet j} = 0$ and $b_j \neq 0$.

## Algorithm 1

*Step 1*: Write down the pseudo-polynomial form of min-SM:

$$P = \coprod_{\substack{j \in N_n \\ b_j \neq 0}} \sum_{i \in N_m} \frac{\breve{\Gamma}_{iij}}{i} \tag{12}$$

If $\breve{\Gamma}_{ij} = 0$, it is omitted from (12). All the operations involved in (12) (summation, multiplication, division) are symbolic.

*Step 2*: Calculate P according to the polynomial multiplication in symbolic form.

*Step 3*: Simplify P by multiplication and then summation.
Multiply the terms of the sum using the formula

$$\frac{a}{l}, \frac{b}{k} = \begin{cases} \dfrac{\max{(a,b)}}{l} & l = k, \\ unchanged, & l \neq k. \end{cases} \tag{13}$$

Sum among the terms using the formula

$$\frac{c_1}{l_1}, \frac{c_2}{l_2} \ldots \frac{c_p}{l_p} + \frac{d_1}{k_1} \cdot \frac{d_2}{k_2} \ldots \frac{d_p}{k_p}$$

$$= \begin{cases} \dfrac{c_1}{l_1}, \dfrac{c_2}{l_2} \ldots \dfrac{c_p}{l_p} & c_i \leq d_i, \forall i \in N_p, \\ unchanged & otherwise. \end{cases} \tag{14}$$

*Step 4*: Suppose that after the Step 3, P has s terms. Then (3) has s minimal solutions. They are computed using the following equation:

$$\breve{r}^{(i)} = \left( \breve{r}_1^{(i)}, \breve{r}_2^{(i)}, ..., \breve{r}_m^{(i)} \right) \qquad (15)$$

where $\breve{r}_j^{(i)} = c_j^{(s)}$, j = 1, 2, …, m, i = 1, 2,…, s.

Here some theoretical results are provided that lead to the simplification of the method described in the above whenever t is an Archimedean t-norm. Note that we have mentioned that the only continuous t-norm that is not Archimedean is the "minimum".

Let us first proceed with some issues on Archimedean t-norms. A very important way for generating Archimedean t-norms is based on the so-called decreasing generators. A decreasing generator is a continuous and strictly decreasing function $f : [0, 1] \to$ R such $f(1) = 0$.

The pseudo-inverse of $f$ is a function $f^{(-1)} :$ R $\to [0, 1]$ defined as

$$f^{(-1)}(a) = \begin{cases} 1, & a \in (-\infty, 0), \\ f^{(-1)}(a), & a \in [0, f(0)], \\ 0, & a \in (f(0), +\infty) \end{cases} \qquad (16)$$

where $f^{-1}$ is the classical inverse of $f$.

A binary operator t : $[0, 1] \times [0, 1] \to [0, 1]$ is an Archimedean t-norm iff there exists a decreasing generator $f$ such that $t(a, b) = f^{(-1)}(f(a) + f(b))$, $\forall$a, b $\in [0, 1]$.

Let t be an Archimedean t-norm and a, b, x $\in [0,1]$. The equation t(a, x) = b has a solution for x iff a $\geq$ b. If b $\in (0, 1)$ the solution $x_0$ is unique and $x_0 \in (0,1]$. Left as an exercise for the reader to prove.

**Algorithm 2**

*Step 1*: Write done the pseudo-polynomial form of Γ:

$$P = \coprod_{\substack{j \in N_n \\ b_j \neq 0}} \sum_{\substack{j \in N_m \\ \Gamma_{ij} = \breve{r}_j}} j \qquad (17)$$

*Step 2*: Calculate P according to the polynomial multiplication in symbolic form.

*Step 3*: Simplify P by multiplication and then summation.
    Multiply the sum using the formula

$$l.k = \begin{cases} l, & l = k, \\ l.k, & l \neq k \end{cases}. \qquad (18)$$

Sum among the terms using the formula

$$l_1. l_2 \ldots l_p + k_1.k_2 \ldots k_p$$

$$= \begin{cases} l_1.l_2 \ldots l_p & if \ \forall i \in N_p, \exists j N_q : l_i = k_j \\ unchanged & otherwise. \end{cases} \qquad (19)$$

Step 4: Suppose that after 3, P has s terms, then (3) has s minimal solutions, computed by

$$\breve{r}_j = \begin{cases} \bar{r}_j & j = l, \\ 0 & otherwise \end{cases} \qquad j \in N_{m.} \qquad (20)$$

We will now show the credibility of the above method.
    The column vector $\hat{r}$ computed is the maximum solution of (3), when t is an Archimedean t-norm.

*The column vector $\hat{r}$ so found is the mean solution of (3), when t is an Archimedean t-norm.*

*Algorithm 2 computes the minimal solution set of (3), if t is an Archimedean t-norm.*

## 1.6  Solving non-linear optimization problem with FRE constraints

J.Lu, S.C. Fang [61] have used fuzzy relation equation constraints to study the non-linear optimization problems. They have presented an optimization model with a nonlinear objective function subject to a system of fuzzy relation equations.
    The study of the fuzzy relation equations

$$x \circ A = b \qquad (1)$$

where $A = (a_{ij})_{m \times n}$, $0 \leq a_{ij} \leq 1$ is a fuzzy matrix 1b $= (b_1,\ldots, b_n)$, $0 \leq b_j \leq 1$ is an n-dimensional vector and '0' is the max-min composition [117].

The resolution of the equation $x \circ A = b$ is an interesting and on-going research topic. [61] in this paper instead of finding all solution of $x \circ A = b$, let $f(x)$ be the user's criterion function, they solve the following non linear programming model with fuzzy relation constrains

$$\min f(x) \text{ s.t } x \circ A = b \qquad (2)$$

A minimizer of Eq. (2) will provide a "best" solution to the user based on the objective function $f(x)$. some related applications of this model with different objective functions can be found in [107] for medical diagnosis, and in [60] for telecommunication equipment module test.

Contrary to the traditional optimization problems [62], problem (2) subjects to fuzzy relation constraints. From [34], we know that when the solution set of the fuzzy relation equations (1) is not empty, it is in general a non-convex set that can be completely determined by one maximum solution and a finite number of minimal solutions. Since the solution set is non-convex, conventional optimization methods [60] may not be directly employed to solve the problem (2). Recently, [24] studied problem (2) with a linear objective function subject to a system of FREs and presented a branch and bound procedure to find an optimal solution. Here, we focus on problem (2) with a nonlinear objective function and call it a nonlinear optimization problem with fuzzy relation constraints (NFRC). A genetic algorithm is proposed for solving MFRC. It is designed to be domain specific by taking advantage of the structure of the solution set of FREs. The individuals from the initial population are chosen from the feasible solution set. The genetic operations such as mutation and crossover are also kept within the feasible region. It is the beauty of this genetic algorithm to keep the search inside of the feasible solution set. The well-maintained feasibility of the population makes the search more efficient.

Genetic algorithms (GAs) are built upon the mechanism of the natural evolution of genetics. GAs emulate the biological evolutionary theory to solve optimization problems. In general,

GAs start with a randomly generated population and progress to improve solutions by using genetic operators such as crossover and mutation operators. In each interaction (generation), based on their performance (fitness) and some selection criteria, the relatively good solutions are retained and the relatively bad solutions are replaced by some newly generated off springs. An evaluation criterion (objective) usually guides the selection.

In the past few years, several methods were proposed to handle the constrained optimization problem using genetic algorithms. Although there was some variation in details among these algorithms, most of them used the penalty or barrier method [38, 41]. [64-66] introduced some special genetic operators to handle the constrained optimization problem, but those operators only work for the problems with a convex domain. A genetic algorithm for optimization problem with fuzzy relation constraints (GAOFRC) was proposed by them.

Unlike a general-purpose genetic algorithm, the proposed GAOFRC is designed specially for solving nonlinear optimization problem with fuzzy relation equations is non-convex in general, and the feasible domain is only a small portion of the convex hull of it, no existing method is readily available for solving NFRC, with the structure inside its feasible domain.

The proposed GAOFRC uses floating-point representation for individuals. Instead of randomly generating a population, the initialization process generates a feasible population utilizing the structure of fuzzy relation equation. The genetic operators are then designed to keep the feasibility of the individuals while they evolve. Those solutions with better objective function values will have higher opportunities to survive in the procedure. The algorithm terminates after it takes a pre-determined number of generations.

In GAOFRC, we use the floating point representation in which each gene or variable $x_i$ in an individual $x = (x_1, x_2,..., x_m)$ is real number from the interval [0 1] since the solution of fuzzy relation equations are nonnegative numbers that are less than one. More specifically individual $x \rightarrow (x_1, …, x_m)$ where $x_i \in [0, 1]$, $i = 1, 2, …, m$.

Compared to the GAs which have no feasibility requirement, GAOFRC's feasibility of individuals limits the search process to a much smaller space.

In general, a GA initializes the population randomly. It works well when dealing with unconstrained optimization problems. However, for a constrained optimization problem, randomly generated solutions may not be feasible. Since GAOFRC intends to keep the solutions (Chromosomes) feasible, we present an initialization module to initialize a population by randomly generating the individuals inside the feasible domain.

Since some elements will never play a role in determining the solutions to fuzzy relation equations. Therefore, we can modify the fuzzy relation matrix by identifying those elements and changing their values to us with the hope of easing the procedure of fining a new solution. To make it clear, we define some "equivalence operators".

**DEFINITION 1.6.1:** *If a value-changing in the element(s) of a given fuzzy relation matrix A has no effect on the solutions of fuzzy relation equations (1). This value changing is called an equivalence operation.*

***Lemma 1.6.1:*** For $j_1$, $j_2 \in \{1, 2,…,n\}$, if $b_{j_1} > b_{j_2}$, $a_{ij_1} \geq b_{j_2}$ and $a_{ij_2} > b_{j_2}$ for some i, then "resetting $a_{ij_1}$ to zero" is an equivalence operation.

Based on this idea, the initialization module originates a population consisting of a given number of randomly generated feasible solutions. The algorithm for initializing is described as follows:

- Get the matrix A, b, and size of population $P_{size}$.
- Compute the potential maximum solution $\hat{x}$ as follows:

$$\hat{x} = (A @ b) = \left[ \bigwedge_{j=1}^{n} (a_{ij} @ b_j) \right]_{1 \times m} \qquad (3)$$

where

$$a_{ij} @ b_j = \begin{cases} 1 & if\ a_{ij} \leq b_j \\ b_j & if\ a_{ij} > b_j. \end{cases}$$

- If $\hat{x}$ o A = b, continue. Otherwise, stop, the problem is infeasible.
- Simplify matrix A by the equivalence operations.
- For each element $a_{ij}$ of A,
  Initialize a lower bound parameter lb(i,j) = 0.
- For i = 1, …, m, j = 1, …, n

If $a_{ij} \geq b_j$, set lb $(i,j) = b_j$;

- For $i = 1, \ldots, m$,
  Set the maximal lower bound $LB_{max}(i) = \max^n_{j=1} lb(i, j)$.

- Set $k = 1$.

- WHILE ($k < P_{size}$)
  For $i = 1, \ldots, m$,
    Generate a random number pop $(k, i)$ in the interval
    $[LB_{max}(i), \hat{x}(i)]$
    Set $k \leftarrow k + 1$.

- Output matrix $[\text{pop }(k,i)]_{P_{size} \times m}$ as the initial population of size $P_{size}$.

Basically, we took up each individual equation and obtained the lower bound for the solutions for each variable $x_i$. Then, we compute the maximal lower bound for each $x_i$. $[LB_{max}(i), \hat{x}(i)]$ is usually an interval. The collection of these intervals we can generate a random number in the interval $[LB_{max}(i), \hat{x}(i)]$ for each $x_i$.

$$P(\text{Select the rth individual}) = q^1 (1-q)^{(r-1)},$$

where $q$ is the probability of selecting the best individual, $r$ is the rank of the individual, q' $= q/(1-(1-q)^{P_{size}})$, and $P_{size}$ is the population size.

Because GAOFRC would like to stay feasible, we cannot mutate the chromosomes randomly. Although various mutation operators handling the constrained optimization problems have been proposed in the literature [38, 41, 61, 65], they are all designed for the convex problems. There seldom is any mutation operator available.

**Table 1**
Elements for each columns such that $a_{ij} \geq b_j$

| Column 1 | Column 2 | Column 3 | Column 4 | Column 5 |
|---|---|---|---|---|
| $a_{21}$, $a_{31}$ | $a_{32}$ | $a_{43}$ | $a_{14}$, $a_{24}$, $a_{54}$, | $a_{45}$ |

for the non-convex problem. In what follows, we present a mutation operator for GAOFRC whose feasible domain is non-convex.

Note that a chromosome in GAOFRC is represented by a $1 \times m$ vector $x = (x_1, x_2,\ldots, x_m)$. For a given chromosome $x^1 = \left(x_1^1, x_2^1, \cdots, x_m^1\right)$, we define a feasible mutation operator that mutates the chromosome by randomly choosing $i_0$ from $1, 2,\ldots,m$ and decreasing $x_{i_0}^1$ to a random number between $[0, x_{i_0}^1]$, while this operation may make the chromosome $x^1$ infeasible we can adjust other $x_j^1$, $i \neq i_0$, to make $x^1$ feasible. In fact, the adjustment of making the infeasible solution become feasible in nothing but a process of finding a new solution. When the changing of $x_{i_0}^1$ pulls the $x^1$ outside of the convex hull of the feasible domain, decreasing $x_{i_0}^1$ results in an infeasible solution no matter how other $x_i^1$ s are adjusted. In this case, GAOFRC will neglect this decreasing operation and find another $x_i^1$ to decrease. Since both the choosing of $x_i^1$ and the extent of decreasing are randomly done, it is guaranteed that a feasible mutation is eventually attainable. We present a feasible mutation operation as follows:

1. Get the simplified matrix A, b and $x = (x_1, x_2,\ldots, x_m)$.

2. Find the decrease set D, a subset of $\{1, 2,\ldots,n\}$, such that there are more than one $a_{ij}$ at column j of A satisfying that $a_{ij} \geq b_j$, for $i \in D$.

3. Randomly choose an element k from D; generate a random number $x'_k$ from the interval $[0, x_k]$, set $x \leftarrow (x_1, x_2,.., x'_k,\ldots,x_m)$.

4. If $\bigvee\limits_{i=1}^{m} (x_i \wedge a_{ij}) = b_j$, $\forall j$, go to step 7; otherwise, go to next step.

5. Generate the increase set $N = D - \{k\}$.

6. For an equation j in which x is not satisfied, randomly choose an element $x_l$ from the increase set N such that $x_l < b_j$ and $a_{lj} \geq b_j$. Set $x_l = b_j$, go to step 4.

7. Go to crossover operation.

However, since the feasible domain of FRE is non-convex, the linear combination of two feasible individuals will very likely result in an infeasible one. Notice that the feasible domain of FRE is comprised of several connected convex sets that have a common maximum point (solution). We can take advantage of this special structure and call the maximal solution a "super" point.

**DEFINITION 1.6.2:** *If a non-convex set is a union of a number of connected, convex subsets and the intersection $S_0$ of these subsets is not empty, then any point s of $S_0$ is called a superpoint.*

From this definition, the maximum point of the feasible domain of FRE is a superpoint. In a connected set, S, for any two points of S, a linear contraction and extraction can be defined.

**DEFINITION 1.6.3:** *Given a connected set S and any two points $x^1$, $x^2$ of S, $0 \leq \lambda \leq 1$, $\gamma \geq 1$,*
*(i) A linear contraction of $x^1$ supervised by $x^2$ is defined by*

$$x^1 \leftarrow \lambda x^1 + (1- \lambda)x^2.$$

*(ii) A linear extraction of $x^1$ supervised by $x^2$ is defined by*

$$x^1 \leftarrow \gamma x^1 - (\gamma - 1)x^2.$$

Once the linear contraction and linear extraction are defined, we can present a "three-point crossover operator. Unlike most of the crossover found in the literature [4, 10, 12, 18, 38, 41, 61, 65, 111], the three-point crossover performs several operations for a point (parent). The operations on a parent will be both supervised by a superpoint and supervised by another parent.

Since the existing theory of genetic algorithms cannot provide measurement for the performance empirical computational testing is necessary. Test problems for computational experiments usually come from three different sources:

1. Published examples.
2. Problems taken from real world applications.
3. Randomly generated test problems.

Since the NFRC problem is in its early research, no published example is available. In this section, we propose a method for constructing test problems for NFRC.

To characterize the system of fuzzy relation equations (1), we introduce a pseudo-characteristic matrix P.

**DEFINITION 1.6.4:** *Given a system of fuzzy relation equations (1), a pseudo-characteristic matrix $P = (p_{ij})_{mxn}$ is defined as*

$$P_{ij} = \begin{cases} 1 & if \ a_{ij} > b_j, \\ 0 & f \ a_{ij} = b_j, \\ -1 & f \ a_{ij} < b_j, \end{cases}$$

With the help of p-matrix, we have:

**THEOREM 1.6.5:** *(Sufficient conditions for existence of solutions). For each column j of matrix A, if*

   *(i)      there is at least one $p_{ij} \neq 1$, and*
   *(ii)     $p_{ij} = 1$ and $b_{j'} > b_j$ implies that $p_{ij'} \neq 1$,*

*then $X(A, b) \neq \phi$.*

For more about the proof please refer [61].

## 1.7 Method of Solution to FRE in a complete Brouwerian lattice

Wang, X. [109] has given a method of solution to FRE in a complete Brouwerian lattice. Di Nola et al [17, 19] point out the problem of solving a FRE in a complete Brouwerian lattice.

Unfortunately, how to solve a FRE in a complete Brouwerian lattice is still an open problem [17]. To this problem, although [84] has given the sufficient and necessary condition to distinguish whether a FRE has a solution, and got the greatest solution in its solution set when it has a solution in [84], whether there exists a minimal element in the solution set and the determination of the minimal elements (if they exist) remain open in the finite case as well as in the infinite case. Here, we first show that there exists a minimal element in the solution set of a fuzzy relation equation $A \odot X = b$ (where $A = (a_1, a_2, \ldots, a_n)$ and b are known, and $X = (x_1, x_2,…, x_n)^T$ are unknown) when its solution set is nonempty and b has an irredundant finite join-decomposition. By the way, we give the method to solve $A \odot X = b$ in a complete Brouwerian lattice under the same conditions. Finally, a method to solve a more general FRE in a complete Brouwerian lattice when its solution set is nonempty is also given under similar conditions. All the works are completed in the case of finite domains.

It is assumed that $L = \langle L, \le, \vee, \wedge \rangle$ is a complete Brouwerian lattice with universal bounds 0 and 1, where $a \vee b = \sup\{1, b\}$, $a \wedge b = \inf\{a, b\}$,"$\le$" stands for the partial ordering of L. The formulas $a \nleq b$ and $b \ngeq a$ both mean that $a \le b$ does not hold. We also assume that X and Y are two finite sets. A mapping $A : X \to L$ is called a fuzzy set of X. A mapping $R : X \times Y \to L$ is called a fuzzy relation between X and Y. Let $X = \{x_1, x_2,…, x_n\}$, $Y = \{y_1, y_2,…, y_m\}$, $\underline{n} = \{1, 2,…n\}$, $\underline{m} = \{1, 2,…m\}$, $\underline{k} = \{1, 2,…k\}$, then a fuzzy set A of X can be denoted by a row vector $A = (a_1, a_2,…, a_n)$ or a column vector $A = (a_1, a_2,…, a_n)^T$ (the sign "T" denotes the "transpose"), where $a_i \in L$, $i \in \underline{n}$.

**DEFINITION 1.7.1 [17]:** *Let $R = (r_{ij})_{n \times m}$ and $A = (a_1, a_2,…, a_m)^T$, we define the max-min composition of R and A to be the fuzzy set $B = (b_1, b_2,…b_n)^T$, in symbols $B = R \odot A$, given by $b_i = \bigvee\limits_{j=1}^{m} (r_{ij} \wedge a_j)$ for any $i \in \underline{n}$.*

We propose three problems:

(q$_1$) Given R and B, determine $X = (x_1, x_2,…, x_m)^T$ such that

$$B = R \odot X \qquad (1)$$

holds. The solution set of (1) is denoted by $\aleph_1$.

(q$_2$) Given A and B determine $X = (x_{ij})_{nxm}$ such that

$$B = X \odot A \qquad (2)$$

holds. The solution set of (2) is denoted by $\aleph_2$.

($q_3$) Given $b \in L$ and $A = (a_1, a_2, \ldots, a_n)$, determine $X = (x_1, x_2, \ldots, x_n)^T$ such that

$$b = A \odot X$$

holds. The solution set of (3) is denoted by $\aleph_3$. Such fuzzy relation equations are called fuzzy elementary equations.

***Proposition 1.7.1:*** $\aleph_3 \neq \phi$ iff (if and only if) $(A \, \alpha \, b)^T \in \aleph_3$. Further, $(A \, \alpha \, b)^T \geq X$ for any $X \in \aleph_3$.

**THEOREM 1.7.2:** *If $\aleph_3 \neq \phi$, then $\aleph_3$ has minimal elements.*

## 1.8 Multi objective optimization problems with FRE constraints

[59] have studied a new class of optimization problems which have multiple objective functions subject to a set of FRE since the feasible domain of such a problem is in general non convex and the objective functions are not necessarily linear, traditional optimization methods become ineffective and inefficient.

Taking advantage of the special structure of the solution set, a reduction procedure is developed to simplify a given problem. Moreover, a genetic-based algorithm is proposed to find the "Pareto optimal solutions".

Let $X = [0,1]^m$, $I = \{1, 2, \ldots, m\}$ and $J = \{1, 2, \ldots, n\}$. Also, let $A$ be an $m \times n$ matrix, $[a_{ij}]_{m \times n}$, and $b$ be an n-dimensional vector $[b_j]_{1 \times n}$, such that $a_{ij} \in [0, 1]$, for all $i \in I$ and $j \in J$. Given $A$ and $b$, a system FRE is defined by

$$x \circ A = b \qquad (1)$$

where "o" represents the max-min composition [117]. A solution to (1) is a vector $x = (x_1, \ldots, x_m)$, $0 \leq x_i \leq 1$, such that

$$\max_{i \in I} \, [\min(x_i, a_{ij})] = b_j, \, \forall \, j \in J. \qquad (2)$$

In other words, the optimization problem we are interested in has the following form:

Minimize $\{f_1(x), f_2(x), \ldots, f_k(x)\}$. $\qquad (3)$

Such that

$$x \circ A = b$$
$$0 \leq x_i \leq 1, i \in I$$

where $f_k(x)$ is an objective function, $k \in K = \{1, \dots, p\}$.

This problem was first studied in [108] for medical applications, with $f_k$ being linear, for $k \in K$. The properties of efficient points were investigated and some necessary and sufficient conditions for identifying efficient points were provided. To facilitate decision making, a procedure was presented to transform the efficient point in an "interval-valued decision space" into a "constant-valued decision space" with a given level of confidence. This transformed problem becomes a multi-attribute decision making problem that can be evaluated by Yagar's method [112] to find an optimal alternative. Unfortunately, the work requires the objective functions to be linear and it also requires the knowledge of all minimal solutions of system (1), which is not trivial at all.

Here, $f_k$ is no longer required to be linear and the information of minimal solutions may be absent. A genetic algorithm (GA) is proposed to solve multi-objective optimization problems with FRE constraints. It is a stochastic searching method which explores the solution space by evaluating the population at hand and evolving the current population to a new one. Since each objective may not be commensurable, it is desirable to achieve a set of non-dominated criterion vectors.

For problem (3), let X be the feasible domain, i.e. $X = \{x \in R^m \mid x \circ A = b, 0 \leq x_i \leq 1, \forall i\}$. For each $x \in X$, we say x is a solution vector and define $z = (f_1(x), f_2(x), \dots, f_p(x))$ to be its criterion vector. Moreover, we define $Z = \{z \in R^p \mid Z = (f_1(x), f_2(x), \dots, f_p(x)), \text{ for some } x \in X\}$.

**DEFINITION 1.8.1:** *A point $\bar{x} \in X$ is an efficient or a Pareto optimal solution to problem (3) if and only if there does not exist any $x \in X$ such $f_k(x) \leq f_k(\bar{x})$, $\forall k \in K$, and $f_k(x) < f_k(\bar{x})$ for at least one k, otherwise, $\bar{x}$ is an inefficient solution.*

**DEFINITION 1.8.2:** *Let $z^1, z^2 \in Z$ be two criterion vectors. Then, $z^1$ dominates $z^2$ if and only if $z^1 \leq z^2$ and $z^1 \neq z^2$. That is, $z_k^1 \leq z_k^2 \quad \forall k \in K$, and $z_k^1 < z_k^2$ for at least one k.*

**DEFINITION 1.8.3:** *Let $\bar{z} \in Z$. Then, $\bar{z}$ is non-dominated if and only if there does not exist any $z \in Z$ that dominates, $\bar{z}$ is a dominated criterion vector.*

The idea of dominance is applied to the criterion vectors whereas the idea of efficiency is applied to the solution vectors. A point $\bar{x} \in X$ is efficient if its criterion vector is non-dominated in Z. That is, from an efficient point, it is not possible to move feasibly so as to decrease one of the objectives without necessarily increasing any other objective. The set of all efficient points is called the efficient set or Pareto optimal set. Also, the set of all non-dominated criterion vectors is called the non-dominated set. In the absence of a mathematical specification of the decision maker's utility function, we can only provide the decision maker with the Pareto optimal set for further analysis. For a problem with multiple linear objective functions, the concept of cones and related properties were used by [93] to characterize the Parato optimal solutions. [108] also used that concept to identify the efficient set.

A system of FRE may be manipulated in a way such that the required computational effort of the proposed genetic algorithm is reduced. Due to the requirement of x o A = b, some components of every solution vector may have to assume a specific value. These components can therefore be set aside from the problem. The genetic operators are then applied to this reduced problem.

Without loss of generality, we assume that the components of vector b are ordered in a decreasing fashion, i.e., $b_1 \geq b_2 \geq \ldots \geq b_n$. And, matrix A is rearranged correspondingly. Notice that the maximum solution $\hat{x}$ can be obtained by the following formula [34]:

$$\hat{x} = A \,@\, b = \left[ \bigwedge_{j=1}^{n} (a_{ij} \,@\, b_j) \right]_{i \in I} \tag{4}$$

where "$\wedge$" is the min operator and

$$a_{ij} \,@\, b_j = \begin{cases} 1 & if \ a_{ij} \leq b_j \\ b_j & if \ a_{ij} > b_{j.} \end{cases} \tag{5}$$

This $\hat{x}$ can then be used to check whether the feasible domain is empty. If

$$\max_{i \in I} [\min(\hat{x}_i, a_{ij})] = b_j, \, \forall \, j \in J \qquad (6)$$

then $\hat{x}$ is the maximum solution of (1). Otherwise, problem (3) is infeasible. For detection of zero procedure please refer [59].

Some of the elements in A play no role in the determination of solutions. These elements if used as a part of constraint satisfaction for some constraint may lead to the violation of another constraint. We shall detect such elements and modify them accordingly.

**DEFINITION 1.8.4:** *If a value-change in some element(s) of a given fuzzy relation matrix A has no effect on the solutions of the corresponding fuzzy relation equations, this value-change is called an equivalence operation.*

The reader is expected to prove the following lemma:

**Lemma 1.8.1:** For $j_1$, $j_2 \in J$, if $b_{j_1} > b_{j_2}$, and for some $i \in I$, $a_{ij_1} \geq b_{j_1}$ and $a_{ij_2} > b_{j_2}$, then changing $a_{ij_1}$ to be zero is an equivalence operation [61].

We now determine which components of every solution vector can assume only a specific value in order to satisfy system (1) to further simplify the system of FRE. The concept of "pseudo-characteristic matrix" [51] will be employed to detect such components. Notice that for constraint j, $j \in J$, if there exists only one $a_{ij}$ that is greater than or equal to $b_j$, then only the ith component of solution vectors can satisfy this constraint. In this case, the value of the ith component depends upon the value of $a_{ij}$. If $a_{ij} > b_j$, then $x_i$ can assume only a single value, $b_j$. If this is the case, the ith component of any solution vector has to be fixed and can be eliminated from the problem.

**DEFINITION 1.8.5:** *Given a system of FRE (1), a "pseudo-characteristic matrix" $P = [P_{ij}]_{mxn}$ is defined as*

$$P_{ij} = \begin{cases} 1 & if \ a_{ij} > b_j, \\ 0 & f \ a_{ij} = b_j, \\ -1 & f \ a_{ij} < b_j. \end{cases}$$

*This pseudo-characteristic matrix will be referred to as a p-matrix.*

**DEFINITION 1.8.6:** *Given the maximum solution $\hat{x}$, if there exists some $i \in I$ and some $j \in J$ such that $\hat{x}_i \wedge a_{ij} = b_j$, then the corresponding $a_{ij}$ of matrix A is called a critical element for $\hat{x}$.*

**Lemma 1.8.2:** If $a_{ij}$ is a critical element, then $p_{ij} \geq 0$.

It is left for the reader to supply proofs for the following lemmas:

**Lemmas 1.8.3:** For column, j if there is only one $i \in I$ such that $p_{ij} = 1$ and $p_{i'j} = -1 \; \forall \; i' \neq i$, then $\hat{x}_i = \breve{x}_i = b_j$.

Let C be the desired number of regions for containing efficient solutions and let P be the size of the efficient set E. To divide data points into C regions or clusters, we apply the concept of fuzzy clustering [4]. The idea is to find the degree of belonging of each data point to each cluster. Data points that belong to the same cluster should be "close" to each other. This closeness is measured by the membership value which is calculated from the distance of the current point from the centre of the cluster compared with those of other points. Given the membership value with respect to that cluster is the maximum value as compared to its membership values with respect to other clusters.

**DEFINITION 1.8.7:** *The matrix $\tilde{U} = [\mu_{cp}]$ is called a fuzzy c-partition if it satisfies the following conditions [4].*

1. *$\mu_{cp} \in [0, 1]$. $1 \leq c \leq C$, $1 \leq p \leq P$,*
2. *$\sum_{c=1}^{C} \mu_{cp} = 1, 1 \leq p \leq P$,*
3. *$0 < \sum_{p=1}^{p} \mu_{cp} < P, 1 \leq c \leq C$.*

*The fuzzy c-partition matrix, $\tilde{U}$, is a matrix of degrees of belonging to clusters of solutions, $\mu_{cp}$ represents the membership value to cluster c of solution p.*

*The location of a cluster is represented by its cluster center* $\upsilon^c = \left(\upsilon_1^c, \dots, \upsilon_q^c\right) \in R^q$, $c = 1, \dots, C$ *around which the data points are concentrated. Q is the number of dimensions of solutions.*

To determine the fuzzy c-partition matrix, U, we need to find the centers of clusters which can be obtained by using different methods. One of the frequently used criterion to identify the clusters in the so-called variance criterion [117]. The variance criterion measures the dissimilarity between points in a cluster and its cluster center by the Euclidean distance. This distance, $d_{cp}$ is calculated by [4]

$$
\begin{aligned}
d_{cp} \quad &= \quad d(x^p, v^c) \\
&= \quad \left\| x^p - v^c \right\| \\
&= \quad \left[ \sum_{j=1}^{q} \left( x_j^p - v_j^c \right)^2 \right]^{\frac{1}{2}}.
\end{aligned}
\tag{9}
$$

The variance criterion for fuzzy partitions corresponds to solving the following problem.

$$
(*)\ \text{Min z}\ (\tilde{U}, v) = \sum_{c=1}^{C} \sum_{p=1}^{P} (\mu_{cp})^m \left\| x^p - v^c \right\|^2
$$

$$
\text{s.t.}\ v^c = \frac{1}{\sum_{p=1}^{P} (\mu_{cp})^m} \sum_{p=1}^{P} (\mu_{cp})^m x^p,
$$

$c = 1, \dots, C,$ \hfill (8)

where $\mu_{cp}$ is determined by

$$
\mu_{cp} = \frac{\left( 1 / \left\| x^p - v^c \right\|^2 \right)^{1/(m-1)}}{\sum_{j=1}^{C} (1 / \left\| x^p - v^j \right\|^2)^{1/(m-1)}},
\tag{9}
$$

$c = 1, \dots, C;\ p = 1, \dots, P$ and $m > 1$ is given.

The system of (*) cannot be solved analytically. However, there exist some iterative algorithms, which approximate the optimal solution by starting from a given solution. One of the best-known algorithms for the fuzzy clustering problem is the fuzzy c-means algorithm [4]. For each $m \in (1, \infty)$, a fuzzy c-means algorithm iteratively solves the necessary conditions (8) as well as (9), and converges to a local optimum. It can be described as follows:

A goal of multi-objective optimization is to obtain the Pareto optimal set. We tested some optimization problems with both linear and nonlinear objective functions. For the linear case, we can theoretically obtain the Pareto optimal set given that the feasible domain is known. The theoretical Pareto optimal set is then used to compare with the results from the proposed genetic algorithm. For the non-linear case, we consider problems whose Pareto optimal sets can be identified. This allows us to precisely analyze the results. We show that our genetic algorithm is capable of finding the Pareto optimal set quite efficiently.

The Pareto optimal set may contain only one solution, a finite number of solutions, or an infinite number of solutions. This depends on the objective functions under consideration. Given different locations of the Pareto optimal set, we wish to investigate whether the proposed genetic algorithm can locate this set. For a problem with multiple linear objective functions, the concept of cones and their properties given by Steuer [93] and the results from Wang [108] will be used to identify the efficiency set.

## 1.9 Neural fuzzy relational system with a new learning algorithm

[92] has given to neural fuzzy relational systems a new learning algorithm. Fuzzy relational systems can represent symbolic knowledge in a formal numerical frame work with the aid of FRE.

It is actually a single layer of generalized neurons (compositional neurons) that perform the sup-t-norm composition. An on-line learning algorithm adapting the weights of its interconnections is incorporated into the neural network. These weights are actually the elements of the fuzzy relation representing the fuzzy relational system.

Fuzzy inference systems are extensions of crisp point-to-point into set-to-set mappings, i.e. mappings from the set of all the fuzzy subsets of the input space to the set of all the fuzzy subsets of the output space [43]. One of the widely used ways of constructing fuzzy inference systems is the method of approximate reasoning which can be implemented on the basis of compositional rule of inference. Different criteria have been proposed for the approximate reasoning to satisfy relevant conditions. The most useful is that of the perfect recall. Fuzzy inference systems that satisfy the perfect recall criterion can be implemented with the aid of max-min fuzzy relation equations (FREs) [17].

The need for more general research [3, 31-33, 37] lead to the representation of fuzzy inference systems on the basis of generalized sup-t-norm FREs [17, 91, 92]. A t-norm (triangular norm) is a function t: $[0,1]$ x $[0,1] \rightarrow [0, 1]$ satisfying for any a, b, c, d $\in [0, 1]$ the next four conditions:

1. $t(a, 1) = a$ and $t(a, 0) = 0$.
2. $b \leq d$ implies $t(a, b) \leq t(a, d)$.
3. $t(a, b) = t(b, a)$.
4. $t(a, t(b, d)) = t(t(a, b), d)$.

Moreover, it is called Archimedean iff

1. t is a continuous function.
2. $t(a, t(a, a) < a, \forall a \in (0, 1)$.

The class of t-norms has been studied by many researchers [32, 40, 92]. Their results are useful in the theory of FREs.

As previously explained, the union-intersection composition of fuzzy relations is one of the key issues of fuzzy set theory. In [91-92] a type of neuron that implements this operation is proposed. This type of neuron is referred to as compositional neuron. Pedrycz [79] and Wang [106] have proposed similar types of neurons, for the sup-min and the sup-product composition, respectively.

The general structure of a conventional neuron can be shown in Figure 1.9.1. The equation that describes this kind of neuron is as follows:

$$y = a \left( \sum_{k=1}^{n} w_i x_i + \partial \right)$$

$(1)$

where a is a non-linearity, $\partial$ is a threshold and $w_i$ are weights that can change on-line with the aid of a learning process.

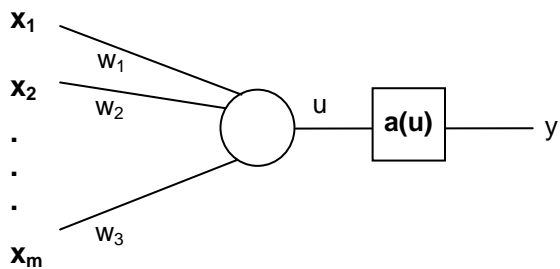

Figure: 1.9.1

The compositional neuron has the same structure with the neuron of Eq. (1) (Figure 1.9.1), but it can be described by the equation:

$$y = a \left( \underset{j \in N_n}{S} \, (x_i, w_i) \right) \qquad (2)$$

where S is a fuzzy union operator (an s-norm), t is a fuzzy intersection operator (a t-norm) and a is the activation function:

$$a\,(x) \;=\; \begin{cases} 0, & x \in (-\infty, 0), \\ x, & x \in [0, 1], \\ 1, & x \in (0, +\infty) \end{cases}$$

which is widely used in neural networks, From Eqs. (1) and (2), a similarity between the two neurons can be shown, since multiplication is a special case of an intersection operator and addition is a special case of a union operator. Compositional neurons can be used in more than one ways in order to construct a neural fuzzy system.

There are two different ways to use compositional neurons. Firstly, a general inference system that implements the relational equation of the approximate reasoning and secondly a general neural fuzzy system of arbitrary level of fuzziness implementing the interpolation method, are proposed. Here a single-layer neural network of compositional neurons is provided for the representation (identification) of a generalized fuzzy inference system. We first provide the reader with the formal problem statement.

Let $X = \{x_1, x_2, \ldots, x_m\}$ and $Y = \{y_1, y_2, \ldots, y_k\}$ be two finite crisp sets and let $D = \{A_i, B_i\}$, $i \in N_n$ be a set of input-output data with $A_i \in F(X)$ and $B_i \in F(Y)$, given sequentially and randomly to the system (some of them are allowed to reiterate before the first appearance of some others). The data sequence is described as $(A^{(\upsilon)}, B^{(\upsilon)})$, $\nu \in N$, where $(A^{(\upsilon)}, B^{(\upsilon)})$, $\in D$. The main problem that arise is the finding of the fuzzy relation R (the fuzzy system) for which the following equation holds:

$$A_i \ o^t \ R = B_i \text{ for each } i \in N_n \qquad (3)$$

where t is a continuous t-norm.

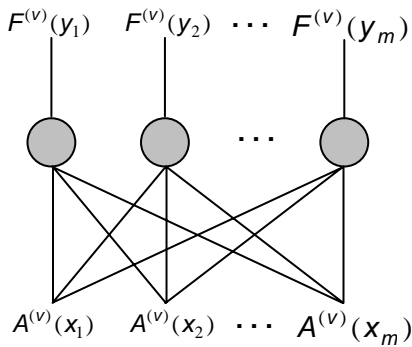

Figure: 1.9.2

In order to construct an efficient learning algorithm, the total error has to be determined. For this reason Pedrycz proposed [75-80] the maximization of an equality index $F^{(\upsilon)}(y)$ of the network output compared with the desired output $B^{(\upsilon)}$. The problem solved by Pedrycz has some differences with the problem solved here, since all the data is assumed to be known before the beginning of the training process. Our case is more difficult, since the data are unknown before the beginning of the training. On the other hand, we suppose that the respective FRE is solvable. The equality index is defined by:

$$F^{(\upsilon)}(y_i) \equiv B^{(\upsilon)}(y_i)] =$$
$$\frac{1}{2}\{[(F^{(\upsilon)}(y_i) \to B^{(\upsilon)}(y_i)) \wedge (B^{(\upsilon)}(y_i) \to (F^{(\upsilon)}(y_i))]$$
$$-[(\overline{F}^{(\upsilon)}(y_i) \to \overline{B}^{(\upsilon)}(y_i) \wedge (\overline{B}^{(\upsilon)}(y_i) \to (\overline{F}^{(\upsilon)}(y_i))] \qquad (4)$$

$$F^{(\upsilon)}(y_i) \equiv B^{(\upsilon)}(y_i)] =$$
$$\begin{cases} 1 + F^{(\upsilon)}(y_i) - B^{(\upsilon)}(y_i) & F^{(\upsilon)}(y_i) < B^{(\upsilon)}(y_i) \\ 1 - F^{(\upsilon)}(y_i) + B^{(\upsilon)}(y_i) & F^{(\upsilon)}(y_i) > B^{(\upsilon)}(y_i) \\ 1 & F^{(\upsilon)}(y_i) = B^{(\upsilon)}(y_i). \end{cases} \qquad (5)$$

The instant error is defined by:

$$E_i^{(\upsilon)} = 1 - [F^{(\upsilon)}(y_i) \equiv B^{(\upsilon)}(y_i)] \qquad (6)$$

which is actually the Hamming distance. We try to minimize the total error

$$E^{(\upsilon)} = \sum_i E_i^{(\upsilon)} \qquad (7)$$

at any time.

## 1.10 Unattainable Solution of FRE

[35] have obtained a necessary and sufficient condition for the existence of a partially attainable and an unattainable solution.

Let U and V be nonempty sets, and let L(U), L(V), and L(U × V) be the collections of fuzzy sets of U, V and U × V, respectively. Then an equation

$$X \text{ o } A = B \qquad (1)$$

is called a FRE, where $A \in L\,(U \times V)$ and $B \in L\,(V)$ are given and $X \in L\,(U)$ is unknown, and o denotes the ∧-∨ composition. A fuzzy set X satisfying the equation above is called a solution of the equation. If $\mu_x: U \to I$, $\mu_A : U \times V \to I$, and $\mu_B: V \to I$ are their membership functions where $I$ denotes the closed interval [0, 1] Eq. (1) is as follows:

$$(\forall \upsilon \in V) \left( \bigvee_{\upsilon \in U} (\mu_X(u) \wedge \mu_A(u, \upsilon)) = \mu_B(\upsilon) \right).$$

The solution set of FRE has been investigated by many researchers [34, 36, 39, 44, 46, 62, 68, 71, 84, 105], and several important properties are shown. Especially, in the case that U and V are both finite sets, it is shown that the solution set is completely determined by the greatest solution and the set of minimal solutions. However, when the cardinality of either U of V is infinite, a few properties about the solution set are investigated [39, 62, 105]. Here we use the concept of attainability to clarify some properties of the solution set of Eq. (1).

**DEFINITION 1.10.1:** *Let $\mu_X$ and $\mu_Y$ be membership functions of fuzzy set X, Y $\in$ L(U), respectively. Then, the partial order $\leq$, the join $\vee$, and the meet $\wedge$, are defined as follows:*

$$\mu_X \leq \mu_Y \Leftrightarrow (\forall u \in U)(\mu_X(u) \leq \mu_Y(u)),$$
$$\mu_X \vee \mu_Y: U \ni u \mapsto \mu_X(u) \vee \mu_Y(u) \in I,$$
$$\mu_X \wedge \mu_Y: U \ni u \mapsto \mu_X(u) \wedge \mu_Y(u) \in I,$$

*Note that $\mu_X \leq \mu_Y$ is equivalent to $X \subset Y$ for X, Y $\in$ L(U).*

**DEFINITION 1.10.2:** *Let $\aleph \subset L(U)$ be the solution set of Eq. (1). The greatest solution of Eq. (1) is an element $G \in \aleph$ such that $\mu_X \leq \mu_G$ (that is X, $\subset$ G) for all $X \in \aleph$. A minimal solution of Eq. (1) is an element $M \in \aleph$ such that $\mu_X < \mu_M$ (that is, X $\aleph$M) for no $X \in \aleph$. Moreover, $\aleph^o$ denotes the set of minimal solutions.*

**DEFINITION 1.10.3:** *For a, b $\in$ [0, 1]*

$$a \; \alpha \; b \; \underline{\Delta} \begin{cases} 1 & \text{if } a \leq b, \\ b & \text{otherwise} \end{cases}.$$

***Attainability of a solution***

When $X \in L(U)$ is a solution of Eq. (1)

$$(\forall \upsilon \in V)(\forall u \in U)(\mu x(u) \wedge \mu_A(u, \upsilon) \leq \mu_B(\upsilon))$$

holds. Moreover, when U and V are both finite sets,

$$(\forall \upsilon \in V) \, (\exists \, u_\upsilon \in U) \, (\mu_X (u_\upsilon) \wedge \mu_A (u_\upsilon, \upsilon) = \mu_B (\upsilon))$$

holds. Thus, we introduce the concepts of attainability and unattainability of a solution.

**DEFINITION 1.10.4:** *Let $X \in L(U)$ be a solution of Eq. (1), and let $V_1$ be a nonvoid subset of V, then*

*X is attainable for $V_1$*

$$\Leftrightarrow \left(\forall \upsilon_1 \in V_1\right)\left(\exists u_{\upsilon_1} \in U\right)\left(\mu_X\left(u_{u_1}\right) \wedge \mu_A\left(u_{\upsilon_1}, \upsilon_1\right) = \mu_B(\upsilon_1)\right)$$

*X is unattainable for $V_1$*

$$\Leftrightarrow \left(\forall \upsilon_1 \in V_1\right)\left(\forall \, u \in U\right)\left(\mu_X\left(u\right) \wedge \mu_A\left(u, \upsilon_1\right) < \mu_B(\upsilon_1)\right).$$

*Moreover, the set of solutions which is attainable for $V_1 \subset V$ is denoted by $\aleph_{\upsilon_1}^{(+)}$ and the set of solutions which is unattainable for $V_1 \subset V$ is denoted by $\aleph_{\upsilon_1}^{(-)}$.*

Note that when the set U and V are both finite, all solutions are attainable for V, that is,

$$\aleph = \aleph_{\upsilon_1}^{(+)}.$$

**DEFINITION 1.10.5:** *Let $X \in L(U)$ be a solution of Eq. (1), and let $V_1$ and $V_2$ be a nonvoid subsets of V satisfying $V_1 \cap V_2 = \phi$ and $V_1 \cup V_2 = V$, then*

*X is an attainable solution $\Leftrightarrow X \in \aleph_{\upsilon_1}^{(+)}$*

*X is a partially attainable solution*

$$\Leftrightarrow X \in \aleph_{\upsilon_1}^{(+)} \cap \aleph_{\upsilon_2}^{(-)}$$

*X is an unattainable solution $\Leftrightarrow X \in \aleph_{\upsilon_2}^{(-)}$.*

In [73], attainability is used for consideration about the extension principle for fuzzy sets. In a FRE the following properties about the set of attainable solutions are known.

**THEOREM 1.10.6 [62]:** *Let* $\hat{X}$ *be the greatest solution given in Theorem 1, then,*

$$\aleph_{\upsilon_1}^{(+)} \neq \phi \Leftrightarrow \hat{X} \in \aleph_{\upsilon_1}^{(+)}.$$

**THEOREM 1.10.7 [62]:**

$$X \in \aleph_{\upsilon_1}^{(+)} \Leftrightarrow (\exists X_g \in \aleph_{\upsilon_1}^{(+)})\mu_X \leq \mu_X \leq \mu_{\hat{X}}).$$

*A fuzzy set $X_g$ in Theorem 5 is called the reachable quasi-minimum solution of Eq. (1).*

**THEOREM 1.10.8 [39]:** *If V is a finite set, then,* $\hat{X} \in \aleph_{\upsilon_1}^{(+)} \Leftrightarrow \aleph^o \neq \phi.$

Some properties of the set of partially attainable and unattainable solutions, here, we show some properties about a partially attainable solution and an unattainable one. The following definition is useful for characterizing such kinds of solutions.

**DEFINITION 1.10.9:** *For a, b $\in$ I, we define $\beta$-operator as follows:*

$$\alpha \ \beta \ b \ \underset{=}{\triangle} \begin{cases} 1 & if \ a < b, \\ b & otherwise. \end{cases}$$

*Note that for a, b $\in$ I,*

$$a \ \alpha \ b = sup \ \{x \in [0, 1] \ a \wedge x \leq b\}$$

*and*

$$a \ \beta \ b = sup \ \{x \in [0, 1] \ a \wedge x < b\}.$$

## 1.11  Specificity shift in solving FRE

The specificity shift method can be classified as an approach situated in between analytical and numerical methods of solving

FRE. It relies on the original structure of the solution originating from the theory and simultaneously takes advantage of some optimization mechanisms available in the format of the parametric specificity shift affecting the relational constraints forming the FRE to be solved.

In this sense the optimal threshold values of the transformation functions provide with a better insight into the character of the data to be handled especially when it comes to their overall consistency level. In this study we are concerned with an important category of FRE with the max-t composition.

$$X \square R = y, \qquad (1)$$

where "t" is assumed to be a continuous t-norm while X, y and R are viewed as fuzzy sets and a fuzzy relation defined in finite universes of discourse. The problem of analytical solutions to these equations has been pursued in the depth; refer e.g. of the monograph by Di Nola et al. [17] as helpful source for the most extensive coverage of the area; see also [13]. On the applied side, these equations call for approximate solutions as quite often no analytical solutions can be generated. This pursuit has been handled with the aid of various techniques.

The approach introduced here falls under the category of data preprocessing by proposing the use to the well-known theoretical results to carefully preprocessed data (relational constraints). The emerging essence is what can be called a specificity shift of relational constraints being aimed at the higher solvability of the resulting FRE.

The problem is stated accordingly: Given is a collection of fuzzy data (treated as vectors in two finite unit hypercubes ) (x(1), y(1)), (x(2), y(2)),…,, (x(N), y(N)) where x (k) $\in [0, 1]^n$ and y (k) $\in [0, 1]^m$. Determine a fuzzy relation R satisfying the collection of the relational constraints (fuzzy relational equations)

$$x (k) \square R = y(k). \qquad (2)$$

Expressing (2) in terms of the corresponding membership functions of x(k), y(k) and R we derive

$$y_j (k) = \bigvee_{i=1}^{n} [x_i(k) tr_{ij}], \qquad (3)$$

where k = 1, 2, … , N, j = 1, 2, …, m. The emerging problem can be essentially classified as an interpolation task where the fuzzy relation R needs to go through all the already specified interpolation points (fuzzy sets).

Assuming that there exists a solution to (2), the theory [17] provides us with the solution to the fuzzy interpolation problem given in the form of the maximal fuzzy relation with the membership function equal to

$$R = \bigcap_{k=1}^{N}(x(k) \rightarrow y(k)).\qquad(4)$$

Note that the computations involve an intersection of the individual fuzzy relations determined via a psuedocomplement (residuation) associated with the t-norm standing in the original system of equations (2), namely

$$a \rightarrow b = \sup \{c \in [0, 1] \mid at\ c \leq b\}.$$

The main advantage lies in the simplicity, theoretical soundness, well-articulated properties and compactness of this solution. The major drawback originates from the fact that the determined result holds under a rather strong preliminary assumption about the existence of any solution to (2).

Now, if this assumption is evidently violated, the quality of the obtained solution could be very low. This is additionally aggravated by the fact that the derived solution is extremal (maximal) so that even a single relational constraint may contribute to the deterioration of the final aggregate result. The use of the optimization methods leads to better approximate solutions yet the entire procedure could be quite often time consuming.

Furthermore, as there could be a multiplicity of solutions, such approaches usually identify only one of them leaving the rest of them unknown.

In the investigated setting we are interested in making some repairs to the original relational constraints thus converting the original interpolation nodes into more feasible ones, meaning that there is a higher likelihood of finding a fuzzy relation capable of doing the interpolation of the modified constraints.

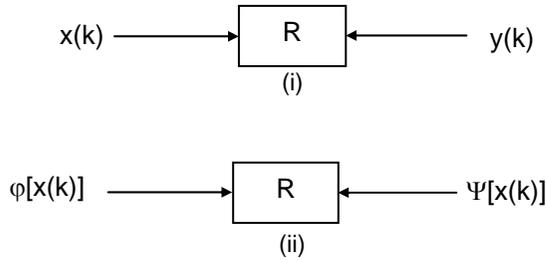

Figure: 1.11.1

Figure 1.11.1 depicts the determination of a fuzzy relation (i) original data (x(k), y(k)) and (ii) modified data.

Affecting the data and changing their membership values. To contrast the scheme of direct computation of R as implied by the theory refer to Figure 1.11.1 where both φ (x) and ψ(y) denote the membership functions resulting from these nonlinear transformations of the original fuzzy sets forming the interpolation nodes of the initial problem.

Thus, instead of (4), one proceeds with the following expression:

$$R = \bigcap_{k=1}^{N} (\varphi\,[x(k)] \rightarrow \psi\,[\,y(k)])$$ (5)

where both φ(x) and ψ(y) are defined point wise meaning that

$$\varphi(x) = [\varphi(x_1)\ \varphi\,(x_2)\ \ldots \varphi(x_n)]$$

and

$$\psi(x) = [\psi\,(x_1)\ \psi\,(x_2)\ \ldots \psi\,(x_m)].$$

There are two types of the transformation functions applied to the input and output data. The first one concerning the input fuzzy sets (x) is defined as a continuous mapping

$$\varphi\colon [0,\,1] \rightarrow [0,\,1]$$

such that

- φ is an increasing function of its argument,
- φ (1) = 1,
- φ(u) ≤ u.

The second operation of interest in this method applies to the output fuzzy sets (y) and is introduced as a continuous mapping

$$\psi: [0, 1] \rightarrow [0, 1]$$

such that

- $\psi$ is an increasing function of its argument,
- $\psi (1) = 1$,
- $\psi (u) \geq u$.

The algorithm combines the theory of the FRE with the heuristics of specificity shift applied to the input-output data (relational constraints). Let us briefly summarize the procedure:

1. select cut-off parameters of $\varphi$ and $\psi$,
2. transform data into a series of pairs ($\varphi(x(k))$, $\psi (y(k))$,
3. compute the fuzzy relation with use of (3),
4. verify the quality of the solution, e.g., by calculating a sum of squared distances (MSE criterion) between y(k) and x(k) □ R.

The entire process can be iterated with respect to the values of the cutoff parameters and these could be optimized so that they imply a minimal value of the MSE criterion.

Additionally, the obtained fuzzy relation can be viewed as a sound starting point for any finer optimization techniques, especially those relying on gradient-based mechanisms. Instead of being initialized from random fuzzy relations, one can start off the method from the fuzzy relation already known.

The numerical studies were completed to explore the efficiency of the introduced mechanism of specificity shift of the relational constraints. Let us emphasize that the choice of the cutoff parameters need to be carried out experimentally as such values are definitely problem-dependent. In general, the higher the values of $\alpha$ and the lower the values of $\beta$, the lower specificity of the computed fuzzy relation. In a limit case where $\alpha = 0$ and $\beta = 1$ the fuzzy relation becomes meaningless (R = 1). In this sense it satisfies all modified relational constraints but fails totally on the original data. As clearly emerges from the two boundary conditions, there are some plausible values of the threshold levels situated somewhere in-between.

Two possible ways of their determination could be sought: one forms an optimization problem involving $\alpha$ and $\beta$ or enumerates the values of the performance index being viewed now as a function of these two unknown parameters. Owing to a low dimensionality of the problem, the latter method sound realistic enough and it will be pursued in all the experiments reported below. Furthermore, to study the simplest possible scenario, the t-norm is specified as the minimum operation. Then the implication operator is a Godelian one governed by the formula

$$a \rightarrow b = \begin{cases} 1 & if \ a \leq b, \\ b & otherwise \end{cases}$$

$a, b \in [0, 1]$.

## 1.12  FRE with defuzzification algorithm for the largest solution

Kagei [42] has provided an algorithm for solving a new fuzzy relational equation including defuzzification. An input fuzzy set is first transformed into an internal fuzzy set by a fuzzy relation. That is the internal fuzzy relation is obtained from fuzzy input and defuzzified output. He classifies these problems into two types called type I and type II. There exists nontrivial largest solution for type I problem. For type II the largest solution is trivial. In addition when unique outputs are required there does not exist the largest solution for the set of solution is not closed set i.e. the supremum of the solutions does not give unique output.

Kagei [42] writes the simultaneous fuzzy relational equations as

$$q_\lambda = p_\lambda \ o \ R \qquad (1)$$

where $p_\lambda$ and $q_\lambda$ are fuzzy sets on X and Y, respectively ($\lambda$ is an index of the equation), R is a fuzzy relation from X to Y to be solved and o is a fuzzy composition operator. When the fuzzy set is defined as a mapping from a nonempty set to a complete Brouwerian lattice, the largest solution of R in Eq. (1) was solved by Sanchez [84] in 1976. After that, many works have been done on the fuzzy relational equations both in theory and in applications (for example, see [17]). The reported works use the fuzzy sets $q_\lambda$ as output data. However, some systems output defuzzified data (for example, refer to [95, 96]). We discuss how to solve the fuzzy relation R when the system includes defuzzification processes.

In these systems, membership values have to be compared with each other for the defuzzification of internal fuzzy sets. Therefore, a complete totally ordered set (a complete chain) should be used, instead of a complete Brouwerian lattice, as the range of membership functions. Such a set which appears in practical problems, like a closed interval and a finite set of real numbers, may be embedded into the unit interval [0, 1] with the order of real numbers. Although we employ a subset of the unit interval as the range of membership functions, the argument here can be applied to any complete totally ordered set which is isomorphic to a complete subset of the unit interval with the usual order.

Let U be a complete subset of the unit interval [0, 1] as a complete totally ordered set with the least element 0 and the greatest element 1. A totally ordered set is a lattice, where the min- and max-operations (meet and join) are given as $a \wedge b = a$ and $a \vee b = b$ for $a \leq b$. The existence of least upper bound is assured by the completeness. Here a fuzzy set is defined as a mapping from a nonempty set into U. Let p and q be fuzzy sets on nonempty sets X and Y, respectively, and R be a fuzzy relation between X and Y, respectively, R be a fuzzy relation between X and Y (i.e., a fuzzy set on X × Y). An input fuzzy set p is transformed to an internal fuzzy set q by the fuzzy relation R in the following two steps:

(1) For each y in Y, the fuzzy set on X is obtained as

$$\mu_p(x) \wedge \mu_R(x, y). \tag{2}$$

(2) The least upper bound operation $\vee$ is taken over X.

$$\mu_q(y) = \underset{x \in X}{\vee} \left\{ \mu_p(x) \wedge \mu_R(x, y) \right\} \tag{3}$$

On the other hand, defuzzification is to select an element with the largest membership value from the support set. Since the defuzzification process is the last step in ordinary systems, we can consider two types of systems including the fuzzy relation and the defuzzification process: (i) defuzzification of p ∩ R in Eq. (2), and (ii) defuzzification of q in Eq. (3).

### 1.12.1 Defuzzification of Eq. (2)

Assume that, for each y in Y. $\mu_R$ (x) $\wedge$ $\mu_R$ (x, y) takes a largest value at x = x* and the element x* is output to the input fuzzy set p. Then, x* satisfies the following equations:

$$\bigvee_{x \in X} \{\mu_p (x) \wedge \mu_R (x, y)\} = \mu_p (x^*) \wedge \mu_R (x^*, y) \qquad (4)$$

for each y in Y.

Note that, since x* depends on y, mapping x* (y) from Y to X is output for an input p. The problem is to find the fuzzy relation $\mu_R(x, y)$ to satisfy the above equation from the input fuzzy set p and the defuzzified output x*(y). Generally, there exist many pairs of inputs p and outputs x*. Each instance is distinguished by the suffix λ. Like as $p_\lambda$ and $x_\lambda$* (a set of λ is finite in engineering application, but this restriction is not needed in theoretical treatment). Since Eq. (4) can be solved in a same way for each y, it is a sufficient to solve the following problem for a fuzzy set $R^{(y)}$ on X, called type I problem, where $R^{(y)}$ is a fuzzy set on X whose membership values are

$$\mu_R(x, y), \text{ i.e.,} \mu_{R^{(y)}} (x) \equiv \mu_R(x, y).$$

Type I problem: For various pairs $(p_\lambda, x_\lambda^*)$ of input fuzzy sets $p_\lambda$ on X and defuzzified outputs $x_\lambda$* in X, obtain $\mu_{R^{(y)}} (x)$ such that

$$\bigvee_{x \in X} \{ \mu_{p_\lambda} (x) \wedge \mu_{R^{(y)}} (x, y)\} = \mu_{p_\lambda} ( x_\lambda^*) \wedge \mu_{R^{(y)}} ( x_\lambda^*) \qquad (5)$$

for all λ.

The superscript of $R^{(y)}$ is omitted in this section. It should be noted that R is regarded as fuzzy set on X for type I problem.

### 1.12.2 Defuzzification of Eq. (3)

Defuzzification is performed for $\mu_q(y)$ in Eq.(3). The output y* is an element of Y which gives the largest value of $\mu_q(y)$ for all y in Y. For many pairs of inputs and outputs, we reach the following type II problem:

*Type II problem:* For various pairs of input fuzzy sets $p_\lambda$ on X and defuzzified outputs $y_\lambda^*$ in Y, obtain $\mu_R(x, y)$ such that

$$\bigvee_{y \in Y} \left( \bigvee_{x \in X} \left\{ \mu_{p^\lambda}(x) \wedge \mu_R(x, y) \right\} \right) = \bigvee_{y \in Y} \left\{ \mu_{p^\lambda}(x) \wedge \mu_R(x, y_\lambda^*) \right\}$$

$$\text{for all } \lambda. \quad (6)$$

Figure 1.12.1(a) shows block diagrams for these systems. In this figure, the defuzzification process receives a fuzzy set, selects the element with the largest membership value, and outputs it.

In Figure 1.12.1(b), type I problem is solved for each y, $p_\lambda \cap R^{(y)}$ is a fuzzy set on X whose membership values are given as $\mu_{p_\lambda \cap R^{(y)}}(x) = \mu_{p_\lambda}(x) \wedge \mu_R(x, y)$.

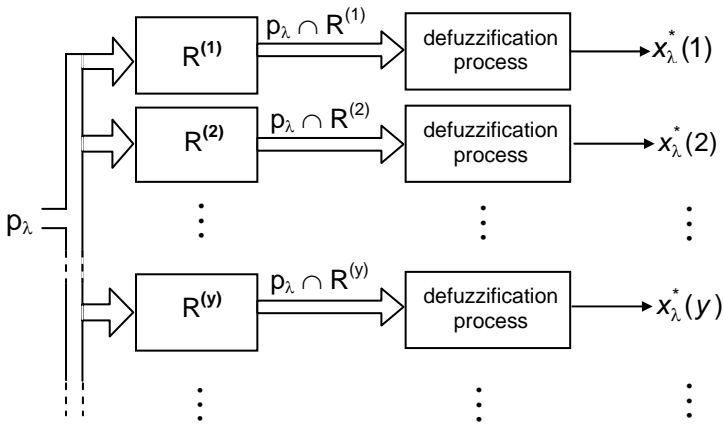

(a)

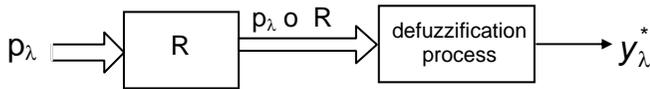

(b)

Figure: 1.12.1 (a) and (b)

### 1.12.3 Unique output

According to Eqs. (5) and (6), it is possible that the defuzzified output is not uniquely determined by R, i.e., the membership value may be the largest at more than one elements. When all defuzzified outputs are uniquely determined, Eqs. (5) and (6) should be replaced as follows:

For Type I problem, instead of Eq. (5)

$$\mu_{p_\lambda}(x) \wedge \mu_R(x) < \mu_{p_\lambda}(x_\lambda{}^*) \wedge \mu_R(x_\lambda{}^*)$$

for all $x \neq x_\lambda{}^*$ and for all $\lambda$.       (7)

For Type II problem, instead of Eq. (6)

$$\bigvee_{x \in X} \left\{ \mu_{p_\lambda}(x) \wedge \mu_R(x, y) \right\} < \bigvee_{x \in X} \left\{ \mu_{p_\lambda}(x) \wedge \mu_R(x, y_\lambda{}^*) \right\}$$

for all $y \neq y_\lambda{}^*$ for all $\lambda$.       (8)

If fuzzy relations $R_\xi$ satisfy Eq. (5) ($\xi$ is an index for distinguishing the solutions), the $U_\xi R_\xi$ also satisfies Eq. (5).

*There exists the largest solution for type I problem.*

For a fuzzy set p and an element x* of X, the largest solution $R^+$ of Eq. (5) is given as

$$\mu_{R^+}(x) = \mu_p(x) \, \alpha \, \mu_p(x^*) \qquad (9)$$

where $\alpha$ operation is defined as (refer to [84])

$$a \, \alpha \, b = \begin{cases} 1 & if \ a \leq b, \\ b & if \ a > b. \end{cases}$$

*Assume that $0 \leq \mu_R(x) \leq M_x$ for all x in X, where $M_x$'s are constants for each x such that $0 < M_x \leq 1$ ($M_x$ depends on x but does not on $\mu_R(x)$). The largest solution $R^+$ of type I problem for a single pair (p, x*) is given as*

$$\mu_{R^+}(x) = \left(\mu_p(x)\alpha\left\{\mu_p(x^*) \wedge M_{x^*}\right\}\right) \wedge M_{x^*}. \qquad (10)$$

In this subsection we assume that the support set has finite elements and the number of pairs (the number of λ) is also finite.

Let $R^+$ be the largest solution of Type I problem. For all pairs of $(p_\lambda, x_\lambda^*)$, we put

$$\mu_{R^+}(x) = \underset{\lambda}{\wedge}\ \left\{\mu_{p\lambda}(x)\ \alpha\ \mu_{p\lambda}(x_\lambda^*)\right\} \qquad (11)$$

**Algorithm 1.** *The largest solution of Type I problem:*

*Step 1:* Set $\mu_R(x) = 1$ for all x ∈ X.

*Step 2:* Repeat the Steps 2.1 and 2.2 for all pairs $(p_\lambda, x_\lambda^*)$.
*Step 2.1:* Set b = $\mu_{p\lambda}(x_\lambda^*) \wedge \mu_R(x_\lambda^*)$ .
*Step 2.2:* For all x not equal to $x_\lambda^*$,
 If $\mu_{p\lambda}(x) \wedge \mu_R(x) > b$ then set $\mu_R(x) = b$.
*Step 3:* When no alternation of $\mu_R(x)$ in Step 2.2 occurs in a single repetition of Step 2, then stop, otherwise, repeat Step 2.

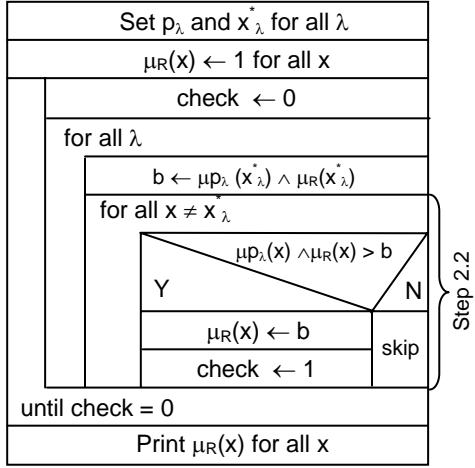

A quasi-largest solution R ‡ is defined as a solution expressed with the above symbols such that R ⊆ R ‡ (i.e., $\mu_R$ (x) ≤ $\mu_{R**}$(x) for all x) for any solution R.

**Algorithm 2.** The quasi-largest solution of type II problem with unique outputs:

*Step 1:* Set $\mu_R$(x, y) = 1 for all x ∈ X and all y ∈ Y.

*Step 2:* Repeat following Steps 2.1 and 2.2 for all pairs ($p_\lambda$, $y_\lambda$*).
*Step 2.1:* Set b = $V_{x' \in X}$ {$\mu_{p\lambda}$($x'$) ∧ $\mu_R$($x'$, $y_\lambda$*)}.
*Step 2.2:* For all x and all y not equal to $y_\lambda$*,
If $\mu_{p\lambda}$($x$) ∧ $\mu_R$($x$, $y$) ≥ b.
Then set $\mu_R$(x, y) = (b).

*Step 3:* When no alternation of $\mu_R$ (x, y) in Step 2.2 occurs in a single repetition of Step 2, then stop. Otherwise, repeat Step 2.

## 1.13 Solvability and Unique solvability of max-min fuzzy equations

Gavalec [27] has given a necessary and sufficient condition for the problem of solvability and for the problem of unique solvability of a fuzzy relation equation in an arbitrary max-min algebra.

By a max-min fuzzy algebra $\mathcal{B}$ we mean any linearly ordered set ($\mathcal{B}$, ≤) with the binary operations of maximum and minimum, denoted by ⊕ and ⊗. For any natural n > 0, $\mathcal{B}$ (n) denotes the set of all n-dimensional column vectors over $\mathcal{B}$, and $\mathcal{B}$ (m, n) denotes the set of all matrices of type m × n over $\mathcal{B}$. For x, y ∈ $\mathcal{B}$ (n), we write x ≤ y, if $x_i$ ≤ $y_i$ holds for all i ∈ N, and we write x< y, if x ≤ y and x ≠ y. The matrix operations over $\mathcal{B}$ are defined with respect to ⊕, ⊗ formally in the same manner as the matrix operations over any field.

$$A \otimes x = b, \qquad (1)$$

where the matrix A ∈ $\mathcal{B}$ (m, n) and the vector b ∈ $\mathcal{B}$ (m) are given, and the vector x ∈ $\mathcal{B}$ (n) is unknown. As matrices over $\mathcal{B}$ correspond to finite fuzzy relations, in the last section we apply our results to a FRE

$$A \otimes X = B, \qquad (2)$$

where $A \in \mathcal{B}(m, n)$, $B \in \mathcal{B}(m, p)$ are given fuzzy relations and the relation $X \in \mathcal{B}(n, p)$ is unknown.

An important special case of max-min algebra is Godel algebra, in which the underlying set is the closed unit interval with natural ordering of real numbers. This max-min algebra will be denoted by $\mathcal{B}_G$. An implication operator in $\mathcal{B}_G$ is defined by $\varphi_G(x, y): = 1$ for $x \le y$ and $\varphi_G(x, y): = y$, for $x > y$.

Here, $\mathcal{B}$ is supposed to be a general linearly ordered set which need not be dense nor bounded. An extension $\mathcal{B}^*$ is defined as the bounded algebra created from $\mathcal{B}$ by adding the least element, or the greatest element (or both), if necessary. If $\mathcal{B}$ itself is bounded, then $\mathcal{B} = \mathcal{B}^*$. The least element in $\mathcal{B}^*$ will be denoted by O, the greatest one by $I$. To avoid the trivial case, we assume $O < I$.

Let a matrix $A \in \mathcal{B}(m, n)$ and a vector $b \in \mathcal{B}(m)$ be fixed. We shall use the notation $M = \{1, 2, \ldots, m\}$, $N = \{1, 2, \ldots, n\}$. Further, we denote the solution sets.

$$S^* (A, b): = \{x \in \mathcal{B}^*(n); A \otimes x = b\},$$
$$S (A, b): = \{x \in \mathcal{B}(n); A \otimes x = b\},$$
$$\text{i.e. } S(A, b): = S^*(A, b) \cap \mathcal{B}(n).$$

The proofs of the following theorems are left as an exercise for the reader and also one can get the proof from [27].

**THEOREM 1.13.1:** *Let $A \in \mathcal{B}(m, n)$, $b \in \mathcal{B}(m)$. Equation $A \otimes x = b$ has a solution $x \in \mathcal{B}(n)$ if and only if $\bar{x}(A, b)$ is a solution in the extension $\mathcal{B}^*$, i.e. if $\bar{x}(A, b) \in S^*(A, b)$.*

The following notation will be useful in this section.
For $i \in M$, $j \in N$ we denote

$$F_{ij}: = \{x \in \bar{S}(A, b); a_{ij} \otimes x_j = b_i\}.$$

If $x \in F_{ij}$, then we say that $x_j$ fulfills the ith equation in $A \otimes x = b$. Of course, it does not mean that x is a solution.

***Lemma 1.13.1:*** Let $i \in M$, $j \in N$. If $x \in \bar{S}(A, b)$-$F_{ij}$, then $a_{ij} \otimes x_j < b_i$.

**Lemma 1.13.2:** Let $x \in \bar{S}$ (A, b). Then the following statements are equivalent:

      (i)                      $x \in S$ (A, b)

      (ii)                    $(\exists \varphi : M \to N)$ $(\forall i \in M)$ $x \in F_{i\varphi(i)}$.

**Lemma 1.13.3:** Let $i \in M$, $j \in N$. Then

    (i)    $F_{ij} = \{ x \in \bar{S}$ (A, b); $x_j = \bar{x}_j \}$, for every $i \in I_j$,

    (ii)   $F_{ij} = \{ x \in \bar{S}$ (A, b); $b_i \le x_j \le \bar{x}_j \}$, for every $i \in K_j$,

    (iii)  $F_{ij} = \emptyset$, for every $i \in M - (I_j \cup K_j)$.

For proof refer [27].

Unique solvability can conveniently be characterized using the notion of minimal covering. If S is a set and C'$\subseteq$ P(S) is a set of subsets of S, we say that ** is a covering of S, if UC' = S, and we say that a covering C' of S is minimal, if U(C' - {C}) $\neq$ S holds for every C $\in$ C'. In [8], a necessary condition for unique solvability is presented under assumption that max-min algebra $\mathcal{B}$ is bounded.

**THEOREM 1.13.2:**[8] *Let $A \in \mathcal{B}$ (m, n), $b \in \mathcal{B}$ (m), let $\mathcal{B}$ be bounded. If equation $A \otimes x = b$ has a unique solution $x \in \mathcal{B}$ (n), then the system $\{I_j \cup K_j, \ j \in N\}$ is a minimal covering of M.*

**THEOREM 1.13.3:** *Let $A \in \mathcal{B}$ (m, n), $b \in \mathcal{B}$ (m). The equation $A \otimes x = b$ has a unique solution $x \in \mathcal{B}$ (n), if and only if the system I is a minimal covering of the set $M - \bigcup K$.*

**THEOREM 1.13.4:** *Let $A \in \mathcal{B}$ (m, n), $b \in \mathcal{B}$ (m). Then*

    *(i)    $|S$ (A, b)$| \ge 1 \Leftrightarrow M = \bigcup (I \cup K),$*

    *(ii)   $|S$ (A, b)$| \le 1 \Leftrightarrow (\forall j \in N) M \neq \bigcup ((I - \{I_j\}) \cup K),$*

**THEOREM 1.13.5:** *Let $A \in \mathcal{B}$(m, n), $b \in \mathcal{B}$(m). Both questions, whether equation $A \otimes x = b$ is solvable, or uniquely solvable, respectively, can be answered in O(mn) time.*

**THEOREM 1.13.6:** *Let $A \in \mathcal{B}(m, n)$, $b \in \mathcal{B}(m)$. Then the greatest solution in S(A, b) exists if and only if $S(A, x) \neq \emptyset$ and $\bar{x}(A, b) \in \mathcal{B}(n)$. If this is the case, then $\bar{x} = \bar{x}(A, b)$ is the greatest solution.*

*Lemma 1.13.4:* Let $x \in S(A, b)$, then $\underline{x} \leq x$.

**THEOREM 1.13.7:** *Let $A \in \mathcal{B}(m, n)$, $b \in \mathcal{B}(m)$. Then the least solution in S(A, b) exists if and only if $\underline{x}(A, b) \in S(A, b)$. If this is the case then $\underline{x} = \underline{x}(A, b)$ is the least solution.*

The results of the previous sections can be applied to fuzzy relation equations of the form

$$A \otimes X = B \qquad (3)$$

where $A \in \mathcal{B}(m, n)$, $B \in \mathcal{B}(m, p)$ and $X \in \mathcal{B}(n, p)$. The relation equation (3) is equivalent to a set of p linear systems of the form (1). The systems use the columns of the matrix B as right-hand side vectors and their solutions form the column of the unknown matrix X. Thus, the results of the previous sections can be easily transferred to the case of (3). We bring here only the basic notation and formulation of theorems, without repeating any proofs.

In the notation $M = \{1, 2, \ldots, m\}$, $N = \{1, 2, \ldots, n\}$, $P = \{1, 2, \ldots, p\}$, we define a matrix $\bar{X} \in \mathcal{B}^*(n, p)$ by putting, for every $j \in N$, $k \in P$,

$$M_{jk} := \{i \in M; a_{ij} > b_{ik}\},$$
$$\bar{x}_{jk} : \min_{\mathcal{B}^*} \{b_{ik}; i \in M_{jk}\}.$$

**THEOREM 1.13.8:** *Let $A \in \mathcal{B}(m, n)$, $B \in \mathcal{B}(m, p)$. The equation $A \otimes X = B$ has a solution $X \in \mathcal{B}(n, p)$ if and only if the corresponding matrix $\bar{X} \in \mathcal{B}^*(n, p)$ fulfills $A \otimes \bar{X} = B$. If, moreover, $\bar{X} \in \mathcal{B}(n, p)$, then $\bar{X}$ is the maximum solution, otherwise, there is no maximum solution in $\mathcal{B}(n, p)$.*

*Further we define, for every $j \in N$, $k \in P$,*

$$I_{jk} = \{i \in M; \, a_{ij} \geq b_{ik} = \overline{x}_{jk} \, \},$$
$$I_K : = \{I_{jk} \, ; \, j \in N\},$$
$$K_{jk} : \{ \, i \in M; \, a_{ij} = b_{ik} < \overline{x}_{jk} \, \},$$
$$K_k : \{K_{jk} : j \in N\}.$$

**THEOREM 1.13.9:** *Let $A \in \mathcal{B}\,(m,\,n)$, $B \in \mathcal{B}\,(m,\,p)$, the equation $A \otimes X = B$ has a unique solution $X \in \mathcal{B}\,(n,\,p)$, if and only if, for every $k \in P$, the system $I_k$ is a minimal covering of the set $M - \bigcup K_k$.*

**THEOREM 1.13.10:** *Let $A \in \mathcal{B}\,(m,\,n)$, $B \in \mathcal{B}\,(m,\,p)$. Both questions, whether the equation $A \otimes X = B$ is solvable, or uniquely solvable, respectively, can be answered in $O(mnp)$ time.*

## 1.14  New algorithms for solving FRE

[63] have given a new algorithm to solve the fuzzy relation equation
$$P \circ Q = R \qquad\qquad (1)$$

with max-min composition and max-product composition. This algorithm operates systematically and graphically on a matrix pattern to get all the solutions of P.

**DEFINITION 1.14.1:** *If $p(Q,\,r)$ denotes the set of all solutions of $p \circ Q = r$, we call $\overline{p} \in \boldsymbol{p}(Q,\,r)$ the maximum solution of $\boldsymbol{p}(Q,\,r)$ if $p \leq \overline{p}$ for all $p \in \boldsymbol{p}(Q,\,r)$. Meanwhile, $\underline{p} \in \boldsymbol{p}(Q,\,r)$ is called a minimal solution of $p(Q,r)$, if $p \leq \underline{p}$ implies $p = \underline{p}$ for all $p \in \boldsymbol{p}\,(Q,\,r)$. The set of all minimal solutions of $\boldsymbol{p}(Q,\,r)$ is denoted by $\underline{\boldsymbol{p}}\,(Q,\,r)$ [34].*

**THEOREM 1.14.2:** $\boldsymbol{p}(\,\underline{p}\,,\,\overline{p}\,)\; \underline{\Delta}\, \Big\{\hat{p} \in \boldsymbol{p} /\;\; \underline{p} \leq \hat{p} \leq \overline{p}\Big\}$ *for each* $\underline{p} \in$ $\boldsymbol{p}$, $\boldsymbol{p}\,(\,\overline{p}\,)\; \underline{\underline{\Delta}}\, \Big\{\underline{\boldsymbol{p}} \in \boldsymbol{p} /\, p \leq \overline{p}\Big\}$ *and $\boldsymbol{p}_0$ denote the set of all minimal elements of $\boldsymbol{p}$, then*

1. $\boldsymbol{p}(Q,\,r) = \cup_p \boldsymbol{p}\,(\,\underline{p}\,,\,\overline{p}\,)$, *where* $\underline{p} \in \underline{\boldsymbol{p}}\;\,(Q,\,r)$.

2. $\boldsymbol{p}(Q,\,r) \neq 0 \Leftrightarrow \underline{\boldsymbol{p}}\;\,(Q,\,r) \neq 0,$

$p\ (Q,\ r) \neq 0 \Leftrightarrow \overline{p} \in p(Q,\ r)$,

$p(Q,\ r) \neq 0 \Leftrightarrow \overline{p}$ *is the maximum solution of Eq. (1) i.e.*
$\overline{p} = max\ p(Q,\ r)$.

3. $p_0 \subset \hat{p} \subset p \Leftrightarrow p\ (Q,\ r) = (\ \hat{p} \cap p(\ \overline{p}\ ))_0$, *where* $(\ \hat{p} \cap p(\ \overline{p}\ ))_0$
   *denotes the set of all minimal elements of* $(\ \hat{p}\ \cap p(\ \overline{p}\ ))$[34].

**Main results**

Following are the main algorithms for solving (1) with max-min (or max-product) composition:

*Step 1:* Check the existence of the solution refer [34, 43].

*Step 2:* Rank the elements of r with decreasing order and find the maximum solution $\overline{p}$ [34, 43].

*Step 3:* Build the table $M = [m_{jk}]$, $j = 1, 2, \ldots, m$; $k = 1, 2, \ldots, n$, where $m_{jk} \underline{\Delta} (\ \overline{p}_j, q_{jk})$. This matrix M is called "matrix pattern".

*Step 4:* Mark $m_{jk}$, which satisfies min $(\ \overline{p}_j, q_{jk}) = r_k$ (or $p_j. q_{jk} = r_k$), and then let the marked $m_{jk}$ be denoted by $\overline{m}_{jk}$.

*Step 5:* If $k_1$ is the smallest k in all marked $\overline{m}_{jk}$, then set $\underline{p}_{j1}$ to be the smaller one of the two elements in $\overline{m}_{jk_1}$ (or set $\underline{p}_{j1}$ to be $\overline{p}_{j1}$).

*Step 6.* Delete the $j_1$th row and the $k_1$th the column of M, and then delete all the columns that contain marked $\overline{m}_{j_1k}$, where $k \neq k_1$.

*Step 7:* In all remained and marked $\overline{m}_{jk}$, find the smallest k and set it to be $k_2$, then let $\underline{p}_{j_2}$ be the smaller one of the two elements in $\overline{m}_{j_2k_2}$ (or let $\underline{p}j_2$ be $\overline{p}_{j_2}$).

*Step 8:* Delete the $j_2$th row and the $k_2$th column of M, and then delete all columns that contain marked $\overline{m}_{j_2k}$, where $k \neq k_2$.

*Step 9:* Repeat steps 7 and 8 until no marked $\overline{m}_{jk}$ is remained.

*Step 10:* The other $\underline{p}_j$, which are not set in steps 5-8, are set to be zero.

**Lemma 1.14.1:** If the FRE is the form as (1), for giving m × n matrix Q and 1 × n vector r, the minimum solutions set ***p*** can be obtained by the above algorithm. Please refer [63] for proof.

## 1.15 Novel neural algorithms based on fuzzy S-rules for FRE

X. Li and D. Ruan [53, 54, 55] have given 3 papers in the years 1997, 1999 and 2000 as three parts on the same title. Their work is a commendable piece of work in the study of FRE and providing a novel neural algorithm based on fuzzy S-rules. In the year 1997 they have given a series of learning algorithm for max-min operator networks and max-min operator networks. These algorithms can be used to solve FRE and their performance and property which are strictly analyzed and proved better by mathematicians. An insight into their work is provided. For more please refer [53].

Any fuzzy system can be represented by a FRE system as

$$A \text{ o } W = B, \qquad (1)$$

where A and B are input and output, respectively, and the compositional operator o is generally a combination t-co-norm/t-norm. In addition to conventional methods [14, 15, 29, 84, 86], a new methodology to solve FRE using fuzzy neural networks [15, 33, 54, 55] is emerging. Since neural architectures were incorporated into the fuzzy field, it has been reasonable for us to think of using neural network architectures to find a solution of FRE.

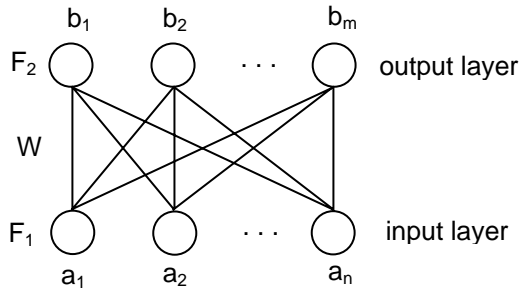

Figure: 1.15.1

A typical network architecture is shown in Figure 1.15.1, which contains an input layer, an output layer and some weighted connections. Its operation is to map an input vector (or pattern) A = $(a_1, a_2, \ldots, a_n)$ to an output vector (or pattern) B = $(b_1, b_2, \ldots, b_m)$, which can be expressed as the following formula:

$$b_j = f\left(\sum_{i=1}^{n} a_i w_{ij} - \theta\right) \qquad (2)$$

where W = $(w_{ij})_{n \times m}$ is the weight matrix of the network and $\theta$ is the threshold value of output neurons, and $f$ is the output transform functions [33]. If we omit the transform function $f$ and the threshold value $\theta$, which is reasonable (see the analysis by [6] in fact for our fuzzy neural network discussed later, formula (2) will become formula (3):

$$AW = B \qquad (3)$$

Still further if we change the operation $(+ , . )$ between A and W to a fuzzy operation e.g. $(\vee, \wedge)$ and confine all data to [0, 1] in Eq. (1). One main characteristic of neural networks is to be able to learn to generate the weighted connection matrix mentioned above from some known patterns which are often called sample data. This process corresponds naturally to that of solving FRE and implicates a new approach. In this case, a lot of researchers studied it [6, 53-55]. The most widely applied learning algorithm of neural networks is the error-based back propagation (BP). A BP algorithm to a two-layered network as shown in Fig.1.15.1 is the so-called $\delta$ rule which requires that the square of difference between the desired output $T_j$ and the actual output $O_j$ is derivable to each $w_{ij}$. So the key problem in realizing the new approach is how to efficiently fuzzify the $\delta$ rule. However a fuzzy system cannot guarantee the previous derivatives to exist. In order to obtain the derivatives, some ideas were proposed. A very representative work is in [6]. Blanco [6] limited the fuzzy operator to o = $(\vee, \wedge)$. His main idea was to define derivatives at those points where they do not exist originally, and he called his method "smooth derivative". Although he succeeded in getting a fuzzy $\delta$ rule, there are still some problems in this method:

- He made a table which contains some experiment result of the comparison between his method and Pedrycz's

[78]. In this comparison, he said he used the "best learning rate" in each case [6]. But why this is the best learning rate and how to decide it? No report.

- One conclusion of [6] is that "the weight matrix after the training is a possible solution for R". If a fuzzy relation equation has solutions, does this method guarantee to converge to a solution after training? No theoretical results and proofs exists.

- All training data in [6] were constructed by generating random inputs first and computing corresponding outputs then according to a known matrix R. It implicates that solutions do exist. But in practice we usually have training data pairs but do not know if a solution exist or not? In this case, could the method in [6] tell us whether a fuzzy relation equation has solution? No.

Due to the complexity in fact these problems exist not only in [6] but also in most other related fuzzy $\delta$ rules.

We describe the max-min operator networks and fuzzy $\delta$ rule

<u>The objective</u>

Assume the fuzzy relation equation is the following (4):

$$A \circ W = B, \qquad (4)$$

where $\circ = (\vee, \wedge)$ and

$$A = \begin{pmatrix} a_{11} & a_{12} & \cdots & a_{1n} \\ a_{21} & a_{22} & \cdots & a_{2n} \\ \vdots & \vdots & \ddots & \vdots \\ a_{p1} & a_{p2} & \cdots & a_{pn} \end{pmatrix},$$

$$B = \begin{pmatrix} b_{11} & b_{12} & \cdots & b_{1m} \\ b_{21} & b_{22} & \cdots & b_{2m} \\ \vdots & \vdots & \ddots & \vdots \\ b_{p1} & b_{p2} & \cdots & b_{pm} \end{pmatrix}.$$

It implies that we have a set of examples $(a_1, b_1)$, $(a_2, b_2)$,…, $(a_p, b_p)$ where $a_i = (a_{i1}, a_{i2},…, a_{in})$ and $b_i = (b_{i1}, b_{i2},…,b_{im}) = (i = 1, 2,…, p)$ are fuzzy vectors. Our objective is to use a neural network to solve the equation by training with these examples. Of course, our training algorithm should be different from any others and should offer strict theoretic results.

<u>Net topology</u>

The net topology in this section is same as that in [6].

A max-min operator neuron is founded by replacing the operator $(+, \bullet)$ of the traditional neuron with $(\vee, \wedge)$, and a network composed of such max-min operator neurons is called a max-min operator network. The inputs and weights of the max-min operator network are generally in [0,1], the output transform function is often $f(x) = x$ and the threshold $\theta = 0$ or we may consider no output transform function and no threshold value at all for all output neurons.

With a two-layered max-min operator network which has the same architecture as shown in Fig.1.15.1, we call it a fuzzy perceptron, whose every node in the input layer connects every node in the output layer. Here if we say a max-min operator network we always mean a two-layered max-min operator network. If input vector is $(a_1, a_2,…, a_n)$, output vector is $(b_1, b_2,…, b_m)$ and elements of the W matrix are $w_{ij}$, the outputs are obtained such that

$$b_j = \mathop{\vee}\limits_{i=1}^{n} (a_i \wedge w_{ij}) , j = 1, 2, …, m \qquad (5)$$

<u>Learning with fuzzy δ rule</u>

*Main idea*

The goal of training the network is to adjust the weights so that the application of set of inputs produces the desired set of outputs. This is driven by minimizing E, the square of the difference between the desired output $b_i$ and the actual outputs $b'_l$ for all the examples. Usually, a gradient descent method is used, which requires that $\delta E / \delta w_{ij}$ exists in (0, 1). However, this necessary condition is seldom satisfied in a fuzzy system. So we cannot directly move a conventional training method to our fuzzy neural network. But our fuzzy system possesses some following features:

- All elements of W are confined to [0, 1], so the greatest possible solution is a matrix with all elements of 1. This implies that if we expect to find the greatest solution of fuzzy relation equation we may initialize all weights to 1. This is very different from any other's methods which always initialize weights to random small positive real numbers.

- The solutions of Eq. (4) is the intersection of the solutions of the following sub equations:

$$a_i \, o \, W = b_i \quad (i = 1, 2, \ldots, p) \tag{6}$$

If we have a training algorithm to be able to find the greatest solution of any equation with the form of (6) but no bigger than the initial weight matrix, we may first obtain the greatest solution of $a_i \, o \, W = b_1$ by initializing all weights to 1 and training the network with first example pair $(a_1, b_1)$, and then obtain the greatest $a_2 \, o \, W = b_2$ but no more than the first solution by training the same network with the second example pair $(a_2, b_2)$, and so on, we will get the greatest solution of $a_p \, o \, W = b_p$ but no more than the $(p - 1)$ the solution by training the same network with the pth example pair $(a_p, b_p)$. We will prove later that the last solution will be the greatest solution of Eq. (4).

- From (5) we know that not all inputs and their weights play an important role to an output but only one or more that $a_i \wedge w_{ij}$ are the biggest which decide the last output. Therefore we need to adjust this weight or these weights if the actual output is bigger than the desired output. On the other hand, obviously those weights with $a_i \wedge w_{ij} > b_j$ should also be decreased.

Based on these features and ideas, now we present a training algorithm which we call fuzzy $\delta$ rule.

*Fuzzy $\delta$ rule*

*Step 1*: Initializing $_{wij} = 1$ for $\forall i, j$.

*Step 2*: Applying inputs and outputs $(a_{i1}, a_{i2}, \ldots, a_{in})$ is an input pattern, and $(b_{i1}, b_{i2}, \ldots, b_{im})$ is an output pattern.

*Step 3*: Calculating the actual outputs

$$(b_{ij})' = \bigvee_{k=1}^{n} (w_{kj} \wedge a_{ik}), j = 1, 2, \ldots, m, \qquad (7)$$

where $(b_{ij})'$ represents the actual output of the jth node when the ith data pair is being trained, and $w_{kj}$ is a connection weight from the kth input node to the jth output node, $a_k$ is the kth component of the input pattern.

*Step 4*: Adjusting weights. Let

$$\delta_{ij} = (b_{ij})' - b_{ij}, \qquad (8)$$

then

$$\begin{cases} w_{kj}(t+1) = w_{kj}(t) - \eta\delta_{ij} & if \ w_{kj}(t) \wedge a_{ik} > b_{ij} \\ w_{kj}(t+1) = w_{kj}(t) & else, \end{cases} \qquad (9)$$

where $\eta$ is a scale factor or a coefficient of step size, and $0 < \eta \leq 1$.

*Step 5*: Return to Step 3, until $w_{kj}(t+1) = w_{kj}(t)$ for $\forall k, j$.

*Step 6*: Repeat Step 2.

Obviously, this algorithm fully reflects the features mentioned above.

The proofs of the following theorem can be had from [53].

**THEOREM 1.15.1.** *If {W (t)} is a weight sequence of the fuzzy $\delta$ rule, then it is a monotone decreasing sequence.*

**THEOREM 1.15.2:** *The fuzzy $\delta$ rule is surely convergent.*

**Lemma 1.15.1.** *If $\exists$ a fuzzy vector $w = (w_1, w_2, \ldots, w_n)$ satisfying*

$$\boldsymbol{a} \ o \ \boldsymbol{w} = b, \qquad (10)$$

*where $\boldsymbol{a} = (a_1, a_2, \ldots, a_n)$ is a fuzzy vector, $0 \leq b \leq 1$. Then the fuzzy $\delta$ rule may converge to $\overline{w}$ which is the maximum solution of Eq.(10).*

**Lemma 1.15.2.** *Suppose the fuzzy matrix W is a solution of a o W = b, and $\overline{W}$ is another fuzzy matrix, $\overline{W} \supseteq W$. Let $\overline{\overline{W}}$ be the initial weight matrix of the fuzzy $\delta$ rule then the fuzzy $\delta$-rule may converge to $\overline{\overline{W}}$, where $\overline{\overline{W}}$ is the maximum solution which is smaller than or equal to $\overline{W}$.*

**THEOREM 1.15.3.** *If $\exists W$ makes the following equation tenable*

$$A \ o \ W = B, \qquad (11)$$

*The fuzzy $\delta$ rule may converge to the maximum solution W\* of this equation, where*

$$A = \begin{pmatrix} a_1 \\ a_2 \\ \vdots \\ a_p \end{pmatrix}, B = \begin{pmatrix} b_1 \\ b_2 \\ \vdots \\ b_p \end{pmatrix}.$$

**THEOREM 1.15.4**. *If Eq. (11) is solvable, the extended fuzzy $\delta$ rule B will converge to the maximum solution. If Eq. (11) has no solutions, the extended fuzzy $\delta$ rule B will converge to the maximum solution of A o W $\subset$ B.*

**THEOREM 1.15.5.** *The convergence matrix of the extended fuzzy $\delta$ rule B is*

$$w_{kj} = \bigwedge_i \left\{ b_{ij} \mid a_{ik} > b_{ij} \right\} \qquad (12)$$

*(appoint $\wedge_{\varphi = 1}$, $\varphi$ is the null set). The number of iteration steps is p which is the number of samples.*

**THEOREM 1.15.6.** *For a fuzzy relation equation*

$$X \ o \ R = S, \qquad (13)$$

*where $R \in \Im (V \times W)$ and $S \in \Im (U \times W)$ are known fuzzy relations, and $X \in \Im (U \times V)$ is an unknown fuzzy relation. Let*

$$\overline{X}(u, \upsilon) \underset{=}{\Delta} \bigwedge_W S(U \times W) \mid S(u, w) < R(u, w) \} \qquad (14)$$

*(appoint $\wedge_{\varphi} = 1$), then equation (13) is compatible if and only if $\overline{X} \; oR = S$, that is, $\underset{\upsilon \in V}{\vee} (\overline{X}(u, \upsilon) \wedge R(\upsilon, w)) = S(u, w) \quad (\forall (u, w))$*

*and $\overline{X}$ is the maximum solution of the equation.*

**THEOREM 1.15.7**. *For Eq. (11), the extended fuzzy $\delta$ rule B is equivalent to the solution method in Theorem 1.15.6.*

## 1.16 Novel neural network part I

In the year 1999 [53, 54, 55] have presented an extended fuzzy neuron and fuzzy neural network and a training algorithm which can be used to resolve some fuzzy relation equation. Their simulation results show if the equation has at least one solution the algorithm will converge to a solution; if an equation does not have a solution at all it can still converge to a matrix which most meets the equation.

**DEFINITION 1.16.1**: *(Fuzzy neuron operators). Suppose $\hat{+}$ and $\hat{\bullet}$ are a pair of binary operators defined on $R^2 \to R$ and for $\forall a, b, a', b', c \in R$, the following are satisfied:*

*1. Monotonicity:*
$$a \hat{+} b \leq a' \hat{+} b', \;\; a \hat{\bullet} b \leq a' \hat{\bullet} b'$$
*if $a \leq a'$ and $b \leq b'$.*
*2. Commutativity:*
$$a \hat{+} b = b \hat{+} a, \qquad a \hat{\bullet} b = b \hat{\bullet} a$$
*3. Associativity:*
$$(a \hat{+} b) \hat{+} c = a \hat{+} (b \hat{+} c), (a \hat{\bullet} b) \hat{\bullet} c = a \hat{\bullet} (b \hat{\bullet} c).$$
*4. Zero absorption of $\hat{\bullet}$ :*
$$a \hat{\bullet} \; 0 = 0.$$
*Then ( $\hat{+}$ , $\hat{\bullet}$ ) are called a pair of fuzzy neuron operators.*

*For example, $(+, \bullet)$, $(\Lambda, \bullet)$, $(\vee, \bullet)$, $(+, \Lambda)$, $(\Lambda, \Lambda)$ $(\vee, \Lambda)$ are all fuzzy neuron operators.*

**DEFINITION 1.16.2:** *(Extended fuzzy neutron). An extended fuzzy neuron is a memory system such as the following:*

$$Y = f \; \hat{+}_{i=1}^{n} \; (w_i \; \hat{\bullet} \; x_i) - \theta), \quad i = 1, 2, ..., n. \qquad (1)$$

*where ( $\hat{+}$ , $\hat{\bullet}$ ) are called a pair of fuzzy neuron operators. $x_i \in [0, 1]$, $i = 1, 2,..., n$, are inputs, $w_i \in [0, 1]$, $i = 1, 2,..., n$, are weights corresponding to the inputs. $\theta$ is a threshold, generally a positive real number or zero. $f$ is an output transform function and its range is [0, 1].*

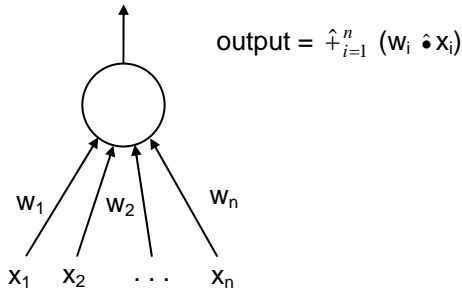

$$\text{output} = \hat{+}_{i=1}^{n} \ (w_i \ \hat{\bullet} \ x_i)$$

Figure 1.16.1

Figure 1.16.1 is a diagram of an extended fuzzy neuron without $\theta$ and $f$. The original reason to define the fuzzy neuron operators is based on some "axiomatic" assumptions of neurons. For example, for a neuron with n inputs we assume: (1) the more each input is, the more the total amount of activation is, and this corresponds to the monotonicity; (2) no input has priority in any order, that is, indifferent to the order in which the inputs to be combined are considered, and this corresponds to the commutativity and associativity of $\hat{+}$ ; (3) the contribution of an input without a weight is zero, this corresponds to the zero absorption of $\hat{\bullet}$ .

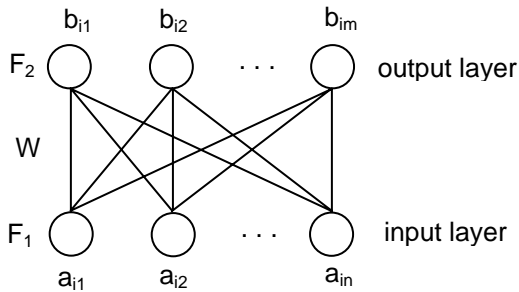

Figure: 1.16.2

We will however pay our attention on a two-layer network as shown in Figure 1.16.2. Input nodes have no computation and just for input, and all output neurons are the extended fuzzy neurons mentioned above without any transform function and threshold. We call this network an extended fuzzy neural network. It can be used to solve such fuzzy relation equations as the following:

$$A \ o \ W = B \qquad (2)$$

where

$$A = \begin{pmatrix} a_{11} & a_{12} & \cdots & a_{1n} \\ a_{21} & a_{22} & \cdots & a_{2n} \\ \vdots & \vdots & \ddots & \vdots \\ a_{p1} & a_{p2} & \cdots & a_{pn} \end{pmatrix},$$

$$B = \begin{pmatrix} b_{11} & b_{12} & \cdots & b_{1m} \\ b_{21} & b_{22} & \cdots & b_{2m} \\ \vdots & \vdots & \ddots & \vdots \\ b_{p1} & b_{p2} & \cdots & b_{pm} \end{pmatrix}.$$

As we know, Eq. (2) may be resolved by first solving m sub-equations as the following and then intersecting all solutions of a sub-equations.

$$\begin{pmatrix} a_{11} & a_{12} & \cdots & a_{1n} \\ a_{21} & a_{22} & \cdots & a_{2n} \\ \vdots & \vdots & \ddots & \vdots \\ a_{p1} & a_{p2} & \cdots & a_{pn} \end{pmatrix} o \ W = \begin{pmatrix} b_{1_i} \\ b_{2_i} \\ \vdots \\ b_{p_i} \end{pmatrix} \qquad (3)$$

where $i = 1, 2,\ldots, m$. So our task is to resolve equations of the form of Eq. (3). In order to simplify formulas, we let the right part of Eq. (3) be $B = (b_i, b_i,\ldots, b_p)^T$. Then we may use a simple network, which is a part of the extended fuzzy neural network in Figure 1.16.2 and called an extended fuzzy perceptron as shown in Figure 1.16.3 to resolve this fuzzy relation equation. In fact, if we regard $(a_{i1}, a_{i2},\ldots a_{in})$ and $b_i$ $(i = 1, 2,\ldots, p)$ as pair of sample patterns, we can find an appropriate W by training the perceptron.

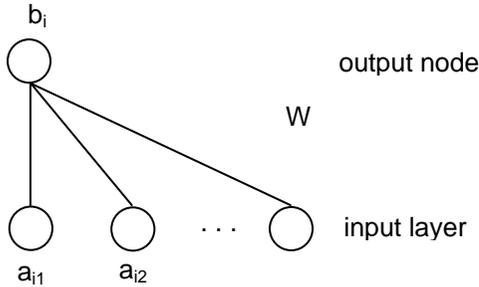

Figure: 1.16.3

The problem is what the training algorithm is.

Before we discuss the algorithm, we will explain the relationship between fuzzy neuron operators and t-norm and t-conorm.

According to [43], a t-norm /fuzzy intersection i is a binary operation on the unit interval that satisfies at least the following axioms for all a, b, d ∈ [0, 1]:

1.    $i(a, 1) = a$ (boundary condition),
2.    $b \leq d$ implies $i(a, b) \leq i(a, d)$ (monotonicity),
3.    $i(a, b) = i(b, a)$ (commutativity),
4.    $i(a, i(b, d) = i(i(a, b), d)$ (associativity).

If a t-norm i is a continuous function (continuity), then we call it a continuous t-norm.

The following are examples of some t-norms that are frequently used as fuzzy intersections (each defined for all a, b ∈ [0, 1]).

Standard intersection    :    $i(a, b) = \min (a, b) = a \wedge b$.
Algebraic product         :    $i(a, b) = ab$.
Bounded difference       :    $i(a, b) = \max(0, a + b - 1)$.
Drastic intersection       :

$$i(a, b) = \begin{cases} a & when\ b = 1, \\ b & when\ a = 1 \\ 0 & otherwise. \end{cases}$$

A t-co-norm / fuzzy union u is a binary operation on the unit interval that satisfies at least the following axioms for all a, b, d ∈ [0, 1],

1. u (a, 0) = a (boundary condition),
2. b ≤ d implies u (a, b) ≤ u (a, d) (monotonicity),
3. u (a, b) = u(b, a) commutativity),
4. u (a, u(b, d)) = u(u (a, b), d) (associativity).

If a t-co-norm u is a continuous function (continuity), then we call it a continuous t-co-norm.

The following are examples of some t-co-norms that are frequently used as fuzzy unions (each defined for all a, b ∈ [0, 1]).

Standard union     :     $u(a, b) = \max (a, b) = a \lor b$.
Algebraic sum     :     $u(a, b) = a + b - ab$.
Bounded sum     :     $u(a, b) = \min(1, a + b)$.
Drastic intersection     :

$$u(a, b) = \begin{cases} a & when\ b = 0, \\ b & when\ a = 0 \\ 1 & otherwise. \end{cases}$$

By comparison, we know most properties of t-norm or t-co-norm are same as those of fuzzy neuron operators, except that the boundary condition of t-norm or t-co-norm is different from the zero absorption of fuzzy neuron operators. In the following we will explain any neuron operators. In the following we will explain any combination of t-norm and t-norm or of t-co-norm and t-norm is a pair of fuzzy neuron operators.

For a t-norm i, we have i (0, 1) = 0 (boundary condition). Also from monotonicity and commutativity, we know i (a, 0) ≤ i (1, 0) = i (0, 1), therefore i (a, 0) = 0 (zero absorption). So, t-norm belongs to $\hat{\bullet}$ of fuzzy neuron operators. Of course, t-norm also belongs to $\hat{+}$. t-co-norm doesn't belong to $\hat{\bullet}$ of fuzzy neuron operators because its boundary condition is u(a, 1) = a and we cannot derive u (a, 0) = 0 (zero absorption). However, it definitely belongs to $\hat{+}$ of fuzzy neuron operators.

Why not just define the fuzzy neuron operators as a combination of t-co-norm and t-norm? It is possible, but the problem is that this definition is too narrow, for example, in this case (+, •) will be out of fuzzy neuron operators, and this is not what we have expected.

This section will give a training algorithm as the following:

*Step 1. Initialization.*
Let t = 0, and $w_i$ (t) = $w_i$(0) = 1.0, for i = 1, 2,…, n. $\hat{B} = (\hat{b}_1, \hat{b}_2, ... \hat{b}_p)^T$ = AW (0) and $\delta \vee_{i=1}^{p} (\hat{b}_i - b_i)$, where it should be noted that $\hat{b}_i > 1$, (i = 1, 2, … , p) often occurs when $\hat{+} = +$, but it does not affect solving Eq. (3). If $\delta \leq \varepsilon$, where $\varepsilon$ is a given very small positive real number, then W (0) is a solution of (1), and go to Step 5.

*Step 2. Calculation.*
For i = 1, 2,…, n, let
$W_i = (w_1(t), …, w_{i-1}(t), w_i(t) - \eta\delta.$
$\quad\quad w_{i+1}(t), …, w_n(t))^T,$
where $\eta$ is a step coefficient, $0 < \eta \leq 1$:
$B'_I = AW_i$, where $B'_I = (b'_{i1}, b'_{i2}, ..., b'_{ip})^T$

$\delta_{ji} = b'_{ij} - b_j,$
$E_{ji} = b''_j - b_{ij}$
for j = 1, 2, …, p, where $b''_j = b_j (t)$
and

$$\delta_i = \vee_{j=1}^{p} \delta_{ji} \quad \Delta_i = \sum_{j=1}^{p} \delta_{ji}$$

$$\delta_{min} = \wedge_{i=1}^{n} \delta_i \quad\quad \Delta_{min} = \wedge_{i=1}^{n} \Delta_i ,$$

$$E = \sum_{j=1}^{p} \sum_{i=1}^{n} E_{ji} .$$

*Step 3. Judgment and weight regulation.*

*Case 1:* ($\delta_{min} > \varepsilon$) is true. If E = 0, then let $w_i$ (t + 1) = $w_i$ (t) - $\eta\delta$, for i = 1, 2, …, n, and B" = (b"$_1$, …, b"$_p$) = AW (t +1), $\delta = \vee_{i=1}^{p} (b''_i - b_i)$
If E $\neq$ 0, then compute
Cos (B'$_I$ – B, $\hat{B}$ – B)

$$= \frac{\left(B_i^{'} - B\right)\left(\hat{B} - B\right)}{\left| B_i^{'} - B\right| \left|\hat{B} - B\right|} \ (i = 1, 2,\ldots, n), \tag{4}$$

and find an index k that makes

$$\cos\left(B_k^{'} - B, \hat{B} - B\right) = \max_{i=1}^{n} \cos\left(B_i^{'} - B, \hat{B} - B\right) \tag{5}$$

*(if the number of the index k satisfying the above equation is more than one, then the k which satisfies the following equation will be chosen,*

$$\delta_k = \delta_{min};$$

if the number of index k is still more than one, then the k which satisfies the following equation will be chosen,

$$\Delta_k = \Delta_{min}.$$

Normally, only one index k will remain. Otherwise, anyone is chosen randomly) and adjust weights with the following formula:

$$\begin{cases} w_i(t+1) = w_i(t) - \eta\delta & if \ i = k, \\ w_i(t+1) = w_i(t) & else. \end{cases}$$

*Case 2:* $(\delta_{min} > \varepsilon)$ is false. Find a index k that makes $\delta_{min} = \Lambda_{i=1}^{n} \delta_I = \delta_k$, and let $W = W_k$. W is a solution we ask. Go to Step 5.

*Step 4: Let t = t + 1 and return to Step 3.*

*Step 5: End.*

It seems very hard to understand this algorithm. In fact, if we analyze it in a fuzzy n-cube space, we will find that the above algorithm is very easy to understand. When we decrease W gradually from the top of the n-cube that is the greatest possible solution W = (1, 1,…,1), we observe the right part of (3), i.e., the output of the fuzzy perception in Figure 1.16.3, also decreases gradually from $\hat{B}$. When the output is near to B in time t, W(t) should be near to a solution of (3) too, because $(\hat{+}, \hat{\bullet})$ possesses monotonicity. After a while we will find simulation results support such a conclusion that if Eq. (3) has a solution at least the above algorithm many converge to a solution of Eq. (3).

## 1.17 Novel neural network part II

In 2000 [Li and Ruan [55]] have extended the fuzzy δ rule from (∨, ∧) to (∨, *) in which * is a general t-norm that is the fuzzy δ-rules J and K. A convergence theorem and an equivalence theorem point out respectively that the fuzzy δ-rule J can converge to the maximum solution and fuzzy δ-rule K is equivalent to the fuzzy method.

A general fuzzy relational equation systems may be expressed as

$$A \text{ o } W = B, \qquad (1)$$

where A and B are input and output, respectively, W is an unknown matrix, and o is a compositional operator which generally is a combination of t-co-norm/t-norm. In addition to conventional methods [67, 84, 86], a new methodology to solve fuzzy relation equations using fuzzy neural networks [7, 33] is emerging.

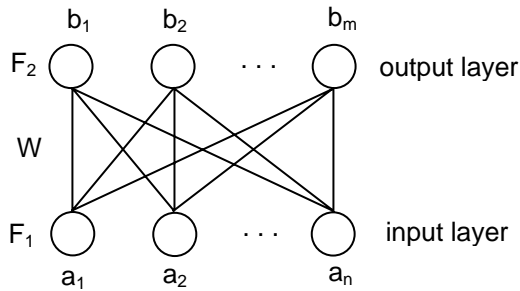

Figure: 1.17.1

A typical network architecture, as shown in Figure 1.17.1, contains an input layer, an output layer, and some weighted connections. Its operation is to map an input vector (or pattern) A = $(a_1, a_2, \ldots, a_n)$ to an output vector (or pattern) B = $(b_1, b_2, \ldots, b_m)$, which can be expressed as

$$b_j \;=\; f\left(\sum_{i=1}^{n} a_i w_{ij} - \theta\right), \qquad (2)$$

where $W = (w_{ij})_{nxm}$ is the weight matrix of the network and $\theta$ is the threshold value of output neurons, and $f$ is the output transform functions (see [33] in detail). If we omit both $f$ and $\theta$ which is reasonable (see the analysis by [6]) then formula (2) will become formula (3):

$$AW \;=\; B \qquad (3)$$

if we change the operation $(+,.)$ between A and W by a fuzzy operation e.g. $(\vee, \wedge)$ and confine all data to [0, 1], then Eq. (3) which will become a typical FRE like equation (1). Actually, we have discussed the cases when the compositional operator $o = (\vee, \wedge)$ and $o = (\vee, \bullet)$. Here we will extend o to a more general form, that is, $o = (\vee, *)$, where * is a t-norm. In the following sections, we will first give a training algorithm for the max-* operator networks and then give its convergence theorem.

Next, we will prove the neural method proposed here is equivalent to the fuzzy solving method in [17]. Afterwards, we will report some simulation and verify that its result is same as that of the old method in [17]. Finally, we will conclude this result and point out our future direction.

The objective

Assume the fuzzy relational equation

$$A \; o \; W = B, \qquad (4)$$

where $o = (\vee, *)$, in which * is a t-norm, and

$$A = \begin{pmatrix} a_{11} & a_{12} & \cdots & a_{1n} \\ a_{21} & a_{22} & \cdots & a_{2n} \\ \vdots & \vdots & \ddots & \vdots \\ a_{p1} & a_{p2} & \cdots & a_{pn} \end{pmatrix},$$

$$B = \begin{pmatrix} b_{11} & b_{12} & \cdots & b_{1m} \\ b_{21} & b_{22} & \cdots & b_{2m} \\ \vdots & \vdots & \ddots & \vdots \\ b_{p1} & b_{p2} & \cdots & b_{pm} \end{pmatrix}.$$

We have set of pairs $(a_1, b_1)$ $(a_2, b_2),\ldots,$ $(a_p, b_p)$, where $a_i = (a_{i1}, a_{i2},\ldots, a_{in})$ and $b_i = (b_{i1}, b_{i2},\ldots, b_{im})$ $(i = 1, 2,\ldots, p)$ are fuzzy vectors. Our objective is to use a neural network to solve the equation by training with the above set of pairs.

<u>Net topology</u>

A max-* operator neuron is founded by replacing the operator $(+, \bullet)$ of the traditional neuron with $(\vee, *)$, and the related network is called max-* operator network. The inputs and weights of the max-* operator network are generally in [0, 1], the output transform function is often $f(x) = x$ and the threshold $\theta = 0$ or we may consider no output transform function and no threshold value for all output neurons.

With a two-layered max-* operator network, which has the same architecture as shown in Figure 1.17.1, we call it a fuzzy perceptron, in which every node in the input layer connects every node in the output layer. Here for a max-* operator network we always mean a two-layered max-* operator network. If the input vector is $(a_1, a_2,\ldots, a_n)$, the output vector is $(b_1, b_2,\ldots, b_m)$, and the elements of the W matrix are $w_{ij}$, then

$$b_j = \bigvee_{i=1}^{n} (a_i * w_{ij}), \qquad j = 1, 2, \ldots, m. \qquad (5)$$

*Main idea*

The goal of training the network is to adjust the weight so that the application of a set of inputs produces the desired set of outputs. This is driven by minimizing E, the square of the difference between the desired output $\mathbf{b}_i$ and the actual output $\mathbf{b'}_i$ for all the pairs. Usually, a gradient descent method is used, which requires that $\partial E/\partial w_{ij}$ exists in (0, 1).

However this necessary condition is seldom satisfied in a fuzzy system. So we cannot directly move a conventional training method in our fuzzy neural network.

But our fuzzy system possesses the following features:

1.  All elements of W are confined to [0, 1], so the greatest possible solution is a matrix whose all elements are 1. This implicates that, if we expect to find the greatest solution of the fuzzy relation equation, we may initialize all weights to 1. This is very different from any other methods which always initialize weights to random small positive real numbers.

2.  The solutions of Eq. (4) is the intersection of the solutions of the following sub-equations:

    $$A_i \text{ o } W = b_i \quad (i = 1, 2, \ldots, p). \qquad (6)$$

3.  If we have a training algorithm to be able to find the greatest solution of any equation with the form of (6) but no bigger than the initial weight matrix, we may first obtain the greatest solution of $a_i \text{ o } W = b_1$ by initializing all weights to 1 and training the network with first example pair $(a_1, b_1)$, and then obtain the greatest solution of $a_2 \text{ o } W = b_2$ but no more than the first solution by training the same network with the second pair $(a_2, b_2)$, and so on, we will get the greatest solution of $a_p \text{ o } W = b_p$ but no more than the $(p - 1)$the solution by training the same network with the pth example pair $(a_p, b_p)$. We will prove later that the last solution will be the greatest solution of Eq. (4).

4.  From (5) we know that not all inputs and their weights play an important role to an output but only one or more that $a_i * w_{ij}$ are the biggest decide the last output. Therefore, we need to adjust this weight or these weights if the actual output is bigger than the desired output. On the other hand, obviously those weights with $a_i * w_{ij} > b_j$ should also be decreased.

*Fuzzy δ rule J*

Following fuzzy δ rules A, B,…, I in [53], we have now the following new training algorithm called fuzzy δ rule J for the (∨, *) operator network.

*Step 1. Initializing*
$w_{ij} = 1$ for $\forall I, j$.

*Step 2. Applying inputs and outputs: $(a_{i1}, a_{i2},…, a_{in})$ is an input pattern and $(b_{i1}, b_{i2},…, b_{im})$ is an output pattern.*

*Step 3. Calculating the actual outputs:*

$$(b_{ij})' = \bigvee_{k=1}^{n}(w_{kj} * a_{ik}), \qquad j = 1, 2, …, m \qquad (7)$$

where $(b_{ij})'$ represents the actual output of the jth node when the ith data pair is being trained, and $w_{kj}$ is a connection weight from the kth input node to the jth output node, and $a_{ik}$ is the kth component of the input pattern.

*Step 4. Adjusting weights: Let*

$$\delta_{ij} = (b_{ij})' - b_{ij}. \qquad (8)$$

then

$$\begin{cases} w_{kj}(t+1) = w_{kj}(t) - \eta\delta_{ij} \\ \quad if \ w_{kj}(t) * a_{ik} > b_{ij} \\ w_{kj}(t+1) = w_{kj}(t) \\ else, \end{cases} \qquad (9)$$

where η is a scale factor or coefficient of step size, and $0 < \eta \leq 1$.

*Step 5:* Return to Step 3, until $w_{kj}(t+1) = w_{kj}(t)$ for $\forall k, j$.

*Step 6:* Repeat Step 2.

The proof of the following theorem are left as an exercise for the reader if need be refer [55].

**THEOREM 1.17.1**: *If {W(t)} is a weight sequence of the fuzzy δ rule J, then it is a monotone decreasing sequence.*

**THEOREM 1.17.2.** *The fuzzy $\delta$ rule J is convergent.*

**Lemma 1.17.1.** *If $\exists$ a fuzzy vector $w = (w_1, w_2,\ldots, w_n)$ satisfying*

$$a \; o \; W = b, \qquad\qquad (10)$$

*where $a = (a_1, a_2,\ldots, a_n)$ is a fuzzy vector, $0 \leq b \leq 1$ and $o = (\lor, *)$, and a given $\eta$ is small enough, then the fuzzy $\delta$-rule J will converge to $\overline{W}$ which is the maximum solution of Eq.(8).*

**Corollary 1.17.1:** *If $\exists$ a fuzzy matrix $W$ satisfying*

$$a \; o \; W = b,$$

*where $a = (a_1, a_2,\ldots, a_n)$ and $b = (b_1, b_2,\ldots, b_m)$ are fuzzy vectors, then the fuzzy $\delta$ rule J will converge to the maximum solution $\overline{W}$ of (10) if $\eta$ is small enough.*

**Lemma 1.17.2:** *Suppose the fuzzy matrix $W$ is a solution of Eq. (9), and $\overline{W}$ is another fuzzy matrix, $\overline{W} \supseteq W$. Let $\overline{W}$ be the initial weight matrix of the fuzzy $\delta$ rule J, then the fuzzy $\delta$ rule J will converge to $\overline{\overline{W}}$ , if $\eta$ is small enough, where $\overline{\overline{W}}$ is the maximum solution which is smaller than or equal to $\overline{W}$ .*

**THEOREM 1.17.3:** *If $\exists W$ which makes the following equation tenable:*

$$A \; o \; W = B \qquad\qquad (11)$$

*The fuzzy $\delta$ rule J will converge to the maximum solution $W^*$ of this equation if $\eta$ is small enough,*
*where*

$$A = \begin{pmatrix} a_1 \\ a_2 \\ \vdots \\ a_p \end{pmatrix}, \quad B = \begin{pmatrix} b_1 \\ b_2 \\ \vdots \\ b_p \end{pmatrix}.$$

In this situation, the convergence value is $W = W(t)$ which is little smaller than the maximum solution. However, if $\eta$ is too small, the maximum solution will be pledged to get, but the training time will be increased too. How to get the maximum solution in the shortest time? Or how to decide "the best learning rate" $\eta$? This is a dilemma. Although we cannot decide what is "the best learning rate", may look for the lowest time cost. For example, we discover in the previous proof whatever $\eta$ is, the goal is to make $w_{kj}$ which satisfies $w_{kj} * a_k > b_j$ closer to $b_j$. So the method with the fastest speed is to seek for a $w_{kj}$ such that

$$w_{kj} * a_k = b_j. \tag{12}$$

Clearly, such $w_{kj}$ may not be unique, but according to the previous algorithm, we always choose the biggest one, that is,

$$\vee \{w_{kj} \mid w_{kj} * a_k = b_j\}. \tag{13}$$

Consequently, we have the following improved algorithm, which is called fuzzy $\delta$ rule K:

*Fuzzy $\delta$ rule K*

This algorithm is almost same to the fuzzy $\delta$ rule J, except that the Step 4 is changed to:

*Step 4*: Adjusting weights:

$$\begin{cases} w_{kj}(t+1) = V\{w_{kj} \mid w_{kj} * a_{ik} = b_{ij}\} \\ \quad if \ w_{kj}(t) * a_{ik} > b_{ij} \\ \quad w_{kj}(t+1) = w_{kj}(t) \\ \quad else. \end{cases} \tag{14}$$

If $\eta$ in the fuzzy $\delta$ rule J is not a constant but a variable parameter, by letting

$$\eta = \frac{w_{kj}(t) - (V\{w_{kj} \mid w_{kj} * a_k = b_j\})}{\delta_{ij}} \tag{15}$$

formula (10) will be equal to Eq. (14). In this case we can say $\eta$ is the best learning rate. The fuzzy $\delta$ rule K is the fastest method and each sample pattern is trained only once. With this algorithm we have

**THEOREM 1.17.4:** *If Eq. (11) is solvable, the fuzzy $\delta$ rule K will converge to the maximum solution.*

**THEOREM 1.17.5:** *If $R \neq \emptyset$, it is*

$$Q^{-1}\psi T = \overset{n}{\underset{i=1}{\wedge}}(Q_i \psi T_i),\qquad(16)$$

*where R is the set of solution, $\emptyset$ is the empty set, $Q^{-1} \psi T$ is the maximum solution, and $\psi$ is a defined composition operation of two fuzzy sets (see below).*

**THEOREM 1.17.6:** *If fuzzy relation equation (11) is solvable, then is maximum solution $\overline{W}$ can be obtained as follows:*

$$\overline{W} = \overset{n}{\underset{i=1}{\wedge}}(a_i \psi b_i).\qquad(17)$$

**THEOREM 1.17.7***: The fuzzy $\delta$ rule K is equivalent to the method of Theorem 1.17.5.*

## 1.18  Simple Fuzzy control and fuzzy control based on FRE

Vladimir P and Dusan Petro have studied fuzzy controllers based on fuzzy relational equations in the case when fuzzy controllers inputs are exact. Fuzzy relational equation with sup-t-composition results in plausible control and adjoint equation results in simple fuzzy control.

Fuzzy logic control is one of the expanding application fields of fuzzy set theory. Recent applications of fuzzy logic control spread over various areas of automatic control, particularly in process control [49-50, 94, 117]. Fuzzy logic controller (FLC) is used whenever conventional control methods are difficult to apply or show limitations in the performances, for example, due to complex system structure [74]. FLC allows simple and more human approach to a control design due to its ability to determine outputs for a given set of inputs without using conventional, mathematical models. FLC follows the general strategy of control worked out by a human being. Using set of control rules and membership functions, FLC converts linguistic variables into numeric values required in most applications. A typical closed-loop system with FLC is shown in Figure 1.18.1.

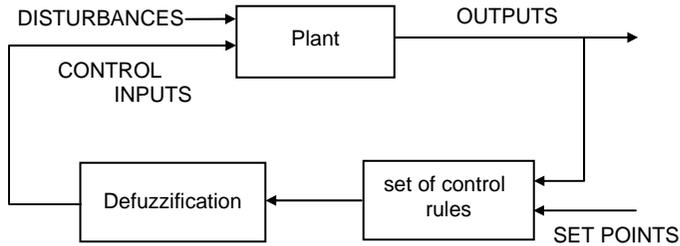

Figure: 1.18.1

A typical FLC is composed of three basic parts: an input signal fuzzification where continuous input signals are transformed into linguistic variables, a fuzzy engine that handles rule inference, and a defuzzification part that ensures exact and physically interpretable values for control variables. The design of FLC may include: the definition of input and output variables, the selection data manipulation method, the membership functions design and the rule (control) base design.

   The main source of knowledge to construct a set of control rules (control base) comes from the control protocol of the human operator. The protocol consists of a set of conditional "if-then" statements, where the first part (if) contains condition and the second part (then) deals with an action (control) that is to be taken. It conveys the human strategy, expressing which control is to be applied when a certain state of the process controlled is matched. For the reasons of simplicity, the following form of the control base is observed here.

$$\text{IF } X_i \text{ THEN } U_i, i = 1, \ldots, n. \qquad (1)$$

Conditions $X_i$ are expressed by membership functions $\mu(X_i(x))$, where x belongs to the space $[X]$. Control variables $U_i$ are expressed by membership functions $\mu(U_i(u))$, where u belongs to the space $[U]$. Statements like else or or are easy to incorporate in conditional part of control base (1).

   The most frequently used data manipulation method is "min-max gravity" method (simple fuzzy control). This Mamdani-type controller assumes min-max inference operators and centre of gravity defuzzification [117]. However, any t-norm and t-co-norm can be used as inference operators. Some properties of FLC using different inference operators can be found in [26, 70].

Simple control is a reasoning procedure based on modus ponens $(A \wedge (A \Rightarrow B)) \Rightarrow B$ tautology [117]. Modus ponens tautology reads:

| | | |
|---|---|---|
| *Implication* | *:* | if A then B |
| *Premise* | *:* | A is true |
| *Conclusion* | *:* | B is true |

where A and B are fuzzy statements or propositions.

Approximate reasoning based on another tautologies, such as modus tolens, syllogism or generalized modus ponens, which give $(A \Rightarrow (A \Rightarrow B)) \Rightarrow B$ [117], was also suggested [70, 115, 116]. Plausible control is observed here. Fuzzy control is plausible if it fulfills features given by F1-F4. Plausible control reads:

With implication if A then B

F1 Premise    A is true conclusion B is true
F2 Premise    A is not true conclusion B is unknown
F3 Premise    A is more fuzzy conclusion B is more fuzzy
F4 Premise    A is less fuzzy conclusion B is B

where feature F1 describes modus ponens tautology.

Inference methods can also be obtained by utilizing fuzzy relational equations with different implication functions. Fuzzy sets and fuzzy relations, calculated for the simple control, satisfy neither fuzzy relational equation with sup-t composition nor adjoint equation. Therefore, simple control is not "mathematically correct". However, solutions of fuzzy relational equations are not unique, because φ-operator and t-norm [26, 117] are not unique. It can be shown that combinations of different implication functions (φ-operator) and t-norms give plausible control. Here Godelian implication for φ-operator and minimum function (intersection) for t-norm are considered.

**Simple fuzzy control**

Using fuzzy relation notation and control base (1), simple fuzzy control is given by [76].

$$R_i = X_i * U_i, \ i = 1,\ldots, n$$

$$R = \bigcup_i R_i \,,$$

$$U = X \odot R, \qquad\qquad (2)$$

where * denotes Cartesian product operator, $\odot$ denotes sup-t (here sup-min) composition and $\cup$ denotes maximum (union) function [76, 117]. Controller input (condition-value obtained from the system), defined over space [X], is denoted by $X(x)$. Calculated fuzzy control, defined over [U], is denoted $U(u)$. Defuzzification (here centre to gravity [76, 117] is later applied to obtain exact control value.

Fuzzy relations $R_i$ are defined over product space [X, U] and are calculated as [76, 117]

$$R_i(x, u) = X_i(x) * U_i(u) = \min\{X_i(x), U_i(u)\} \qquad (3)$$

for all $x \in [X]$ and $u \in [U]$, where $X_i(x)$ and $U_i(u)$ are expressed by their membership functions. Relating to sup-min composition [76, 117] fuzzy control is

$$U(u) = X(x) \odot R(x, u) = \sup_{x \in [X]} \{\min\{X(x), R(x, u)\}\}. \qquad (4)$$

Subject to (3), this can be rewritten as

$$U(u) = \sup_{x \in [X]} \{\min\{X(x), \bigcup_i R_i(x, u)\}\}. \qquad (5)$$

Sup-min composition is distributed with respect to union [76, 117]:

$$U(u) = \sup_{x \in [X]} \{\bigcup_i \min\{X(x), R(x, u)\}. \qquad (6)$$

Subject to (3)

$$U(u) = \sup_{x \in [X]} \{\bigcup_i \min\{X(x), X_i(x) * U_i(u)\}. \qquad (7)$$

Subject to associativity this can be rewritten as

$$U(u) = \{\bigcup_i \sup_{x \in [X]} \min\{X(x), X_i(x) * U_i(u)\}$$

$$= \bigcup_i \wedge_I * U_i(u), \tag{8}$$

where $\wedge_I$ is a scalar value, called possibility of X with respect to $X_i$ [76], defined by

$$\wedge_i = \Pi(X / X_i) = \sup_{x \in [X]} \{\min \{X(x), X_i(x)\}\} \tag{9}$$

In simple fuzzy control, $\Lambda_I$ is a scalar value even in he case when controller input is fuzzy. However, a particular case, when X(x) is nonfuzzy, is observed here, which means

$$X(x) = \begin{cases} 1, & x = x_0, \\ 0, & otherwise \end{cases} \tag{10}$$

In that case, $\Lambda_i$ is calculated as

$$\wedge_i = \sup\{\min\{1, X_i(x_0)\}, \min\{0, X_i(x)\}\} = X_i(x_0). \tag{11}$$

*Fuzzy control based on fuzzy relational equation with sup-t composition*

Using fuzzy relation notation and control base (1) this type of fuzzy control is given by [76]

$$R_i = X_i \varphi U_i, \tag{12}$$

where
$$i = 1, \ldots, n,$$
$$R = \bigcap_i R_i,$$
$$U = X \odot R,$$

where $\cap$ denotes minimum (intersection) function [117]. Operator $\varphi$ represents implication function (here Godelian implication) [76].

$$(X \varphi Y)(x, y) = X(x) \varphi Y(y)$$

$$= \begin{cases} 1, & \mu(X(x)) \leq \mu(Y(y)), \\ \mu(X(x)), & \mu(X(x)) > \mu(Y(y)). \end{cases} \tag{13}$$

From (4) and (12) fuzzy control is

$$U(u) = \sup_{x \in [X]} \{\min \{X(x), \bigcap_i R_i (x, u)\}\} \qquad (14)$$

Sup-min composition is not distributed with respect to intersection [1]:

$$X \odot (Y \cap Z) \Leftarrow (X \odot Y) \cap (X \odot Z) \qquad (15)$$

Therefore, obtaining fuzzy control U(u) in (14), demands calculating fuzzy relation R (x, u)

**THEOREM 1.18.1:** *In the case when FLC input is exact (10), fuzzy control based on fuzzy relational equation with sup-t composition (14) can be obtained without calculation of fuzzy relation (12).*

*In view of (1), this fuzzy control is given by [76].*

$$R_i = X_i * U_i, \ \ i = 1,..., n,$$

*where*

$$R = \bigcup_i R_i,$$
$$U = X \, \varphi \, R, \qquad (16)$$

*where operators * and $\cup$ are already defined. Operator $\varphi$ in (18) denotes Godelian implication between fuzzy set and fuzzy relation [76].*

$$U(u) = \inf_{x \in [X]} \{X(x) \, \varphi \, R \, (x, u)\}. \qquad (17)$$

**THEOREM 1.18.2:** *In the case when FLC input is nonfuzzy (10), fuzzy control based on adjoint fuzzy relational equation (17) gives the same control algorithm as a simple fuzzy control (8).*

The proofs of the above theorem is left as an exercise for the reader.

## 1.19 A FRE in dynamic fuzzy systems

For a dynamic fuzzy system the fundamental method is to analyze its recursive relation of the fuzzy states. M. Kurano et al [47] gave the existence and the uniqueness of solution of a fuzzy relational equation.

We use the notations in [46]. Let X be a compact metric space. We denote by $2^X$ the collection of all subsets of X, and denote by $\zeta$ (X) the collection of all closed subsets of X. Let $\rho$ be the Hausdorff metric on $2^X$. Then it is well-known [48] that ($\zeta(X)$, $\rho$) is a compact metric space. Let $\zeta(X)$ be the set of all fuzzy sets $\bar{s} : X \rightarrow [0, 1]$ which are upper semi-continuous and satisfy $\sup_{x \in X} \bar{s}$ (x) = 1. Let $\bar{q} : X \times X \rightarrow [0, 1]$ be a continuous fuzzy relation on X such that $\bar{q}$ (x,.) $\in \zeta(X)$ for $x \in X$.

Here, we consider the existence and uniqueness of solution $\bar{p} \in \zeta$ (X) in the following fuzzy relation equation (1) for a given continuous fuzzy relation $\bar{q}$ on X (see [46]):

$$\bar{p} \ (y) = \sup_{x \in X} \{\bar{p}(x) \wedge \bar{q}(x, y)\}, \ y \in X, \tag{1}$$

where $a \wedge b := \min \{a, b\}$ for real numbers a and b. We define a map $\tilde{q}_\alpha : 2^X \rightarrow 2^X$ ($\alpha \in [0, 1]$ ) by

$$\tilde{q}_\alpha \ (D) := \begin{cases} \{y \mid \tilde{q}(x, y) \geq \alpha \ for \ some \ x \in D\} \\ \qquad for \ x \neq 0, D \in 2^X, D \neq \phi, \\ cl\{y \mid \tilde{q}(x, y) > 0 \ for \ some \ x \in D\} \\ \qquad for \ \alpha = 0, D \in 2^X, D \neq \phi, \\ X \ for \ 0 \leq \alpha \leq 1, D = \phi, \end{cases} \tag{2}$$

where cl denotes the closure of a set. Especially, we put $\tilde{q}_\alpha$ (x) := $\tilde{q}_\alpha$ ({x}) for $x \in X$. We note that $\tilde{q}_\alpha : \zeta$ (X) $\rightarrow \zeta(X)$.

For $\alpha \in [0, 1]$ and $x \in X$, a sequence $\left\{\tilde{q}_\alpha^k (x)\right\}_{k=0,1,...}$ is defined iteratively by

$$\tilde{q}_\alpha^0(x) := \{x\}, \ \tilde{q}_\alpha^1(x) := \tilde{q}_\alpha(x) \ and$$
$$\tilde{q}_\alpha^{k+1}(x) := \tilde{q}_\alpha \ \tilde{q}_\alpha^k(x) \ ) \ for \ k = 1, 2, \dots .$$

Then, let $G_\alpha$ (x) : = $\bigcup_{k=1}^\infty \tilde{q}_\alpha^k (x)$ and

$$F_\alpha \ (x) := \bigcup_{k=0}^\infty \tilde{q}_\alpha^k (x) = \{x\} \cup G_\alpha \ (x) \tag{3}$$

We now consider a class of invariant points for this iteration procedure, that is, x ∈ $G_\alpha$ (x). So put

$$R_\alpha := \{x \in X \mid x \in G_\alpha(x)\} \text{ for } \alpha \in [0, 1]. \quad (4)$$

Each state of $R_\alpha$ is called as an "α-recurrent" state. The following properties (i) and (ii) hold clearly:

(i) $\tilde{q}_\alpha$ ($F_\alpha(z)$) = $G_\alpha(z)$ for α ∈ [0, 1] and x ∈ X;

(ii) $R_\alpha \subset R_{\alpha'}$ for 0 ≤ α' < α ≤ 1.

**Lemma 1.19.1:** If z ∈ $R_1$, the following (i)-(ii) hold:

i. $\tilde{q}_\alpha$ ($F_\alpha$ (z)) = $F_\alpha$ (z) for α ∈ [0, 1].:

ii. $F_\alpha$ (z) ⊂ $F_{\alpha'}$ (z) for 0 ≤ α' < α ≤ 1;

**Lemma 1.19.2:** If z ∈ $R_1$, the following (i)-(iii) hold:

i. $\tilde{q}_\alpha$ ( $\hat{F}_\alpha$ (z)) = $\hat{F}_\alpha$ (z) for α ∈ [0, 1].:

ii. $\hat{F}_\alpha$ (z) ⊂ $\hat{F}_{\alpha'}$ (z) for 0 ≤ α' < α ≤ 1;

iii. $\hat{F}_\alpha$ (z) = $\lim_{\alpha' \uparrow \alpha}$ $\hat{F}_{\alpha'}$ (z) for α ≠ 0.

**THEOREM 1.19.1:**

*i. If $R_1 \neq \emptyset$, then there exists a solution of (1).*

*ii. Let z ∈ $R_1$. Define a fuzzy state*
$$\bar{s}^z(x) := \sup_{x \in [0,1]} \left\{ \alpha \wedge 1_{\hat{F}_z(z)}(x) \right\}, \, x \in X,$$

*Then $\bar{s}^z \in \mathfrak{I}(X)$ satisfies (1).*

**Lemma 1.19.3:** For $z_1$, $z_2$ ∈ $R_1$,

$$z_1 \sim z_2 \text{ if and only if } F_\alpha(z_1) = F_\alpha(z_2)$$

for all x ∈ [0, 1].

**THEOREM 1.19.2:** *Let $\tilde{p}^k \in P$ (k = 1, 2, …, 1). Then:*
*(i) Put*

$$\tilde{p}\,(x) := \max_{k=1,2,...,l} \tilde{p}^k\,(x) \text{ for } x \in X.$$

Then $\tilde{p} \in P$.

*(ii) Let $\{\alpha^k \in [0, 1] \mid k = 1, 2, ..., l\}$ satisfy $\max_{k=1, 2, ..., l} \alpha^k = 1$. Put*

$$\tilde{p}\,(x): \max_{k=1,2,...,l} \left\{\alpha^k \wedge \tilde{p}^k(x)\right\} \text{ for } x \in X.$$

Then $\tilde{p} \in P$.

Let B be a convex subset of an n-dimensional Euclidean space $\mathfrak{R}^n$ and $C_c(B)$ the class of all closed and convex subsets of B. Throughout this section, we assume that the state space X is a convex and compact subset of $\mathfrak{R}^n$. The fuzzy set $\bar{s} \in \mathfrak{I}(X)$ is called convex if its $\alpha$-cut $\tilde{s}_\alpha$ is convex for each $\alpha \in [0, 1]$. Let $\mathfrak{I}_c(X) := \{ \tilde{s} \in \mathfrak{I}(X) \mid \tilde{s} \text{ is convex}\}$.

By applying Kakutani's fixed point theorem [23], we have the following:

**Lemma 1.19.4.** Let $\alpha \in [0, 1]$ and $\tilde{q}_\alpha(x)$ is convex for each $x \in X$. Then, for any $A \in C_c(X)$ with $A = \tilde{q}_\alpha(A)$, there exists an $x \in X$ such that $\bar{q}\,(x, x) \geq \alpha$.

**Proposition 1.19.1.** Let $p \in \mathfrak{I}_c(X)$ be a solution of (1). Then, for each $\alpha \in [0, 1]$, there exists an $x \in p_\alpha$ with $\tilde{q}\,(x, x) \geq \alpha$.

Assumption A. The following A1-A3 hold.

A1.    The set $U_1$ is a one-point set, say u. That is, $U_1 = \{u\}$.

A2.    $U_\alpha \subset F_\alpha(u)$ for each $\alpha \in [0, 1]$, where u is given by A1 and $F_\alpha(u)$ is defined by (3).

A3.    Let $\alpha \in [0, 1]$ and $A \in C_c(X)$. If $A = \tilde{q}_\alpha(A)$, then
$$A = \bigcup_{x \in \cup \cap A} F_\alpha(x)$$

**THEOREM 1.19.3:** Under Assumption A, Eq. (1) has a unique solution in $\Im_c(X)$.

## 1.20  Solving FRE with a linear objective function

Fang and Li [24] has given an optimization model with a linear objective function subject to a system of fuzzy relational equations is present. Due to the non convexity of its feasible domain defined by fuzzy relation equations designing an efficient solution procedure for solving such problems is not a trivial job. Here they present a solution procedure.

Let $A = [a_{ij}]$, $0 \leq a_{ij} \leq 1$, be an $(m \times n)$-dimensional fuzzy matrix and $b = (b_1,\ldots, b_n)^T$, $0 \leq b_{ij} \leq 1$, be an n-dimensional vector, then the following system of fuzzy relation equations is defined by A and b:

$$x \circ A = b, \qquad (1)$$

where "o" denotes the commonly used max-min composition [117]. In other words, we try to find a solution vector $x = (x_1,\ldots, x_m)^T$, with $0 \leq x_i \leq 1$, such that

$$\max_{i=1,2,\ldots,m} \min (x_i, a_{ij}) = b_j \text{ for } j = 1,\ldots, n. \qquad (2)$$

The resolution of fuzzy relation equations (1) is an interesting and on-going research topic [1, 2, 11, 17, 30, 34, 52, 82, 84, 106]. Here, we study a variant of such problem.

Let $c = (c_1,\ldots, c_m)^T \in R^m$ be an m-dimensional vector where $c_i$ represents the weight (or cost) associated with variable $x_i$ for $i = 1,\ldots, m$. We consider the following optimization problem:

Minimize $\qquad Z = \sum_{i=1}^{m} c_i x_i \qquad (3)$

such that $\qquad x \circ A = b,$

$\qquad\qquad\qquad 0 \leq x_i \leq 1.$

Compared to the regular linear programming problems [25], this linear optimization problem subject to fuzzy relation equations has very different nature.

Note that the feasible domain of problem (3) is the solution set of system (1). We denote it by X (A, b) = {x = $(x_1,\ldots,x_m)^T \in R^m$| x o A = b, $x_i \in [0, 1]$}.

To characterize X (A, b), we define I = {1, 2, …, m}, J = {1, 2, …, n}, and X = {x ∈ $R^m$ | 0 ≤ $x_i$ ≤ 1, $\forall i \in I$}. For $x^1$, $x^2 \in X$, we say $x^1 \le x^2$ if and only if $x_i^1 \le x_i^2$, $\forall i \in I$. In this way, "≤" forms a partial order relation on X and (X, ≤) becomes a lattice. Moreover, we call $\hat{x} \in X$ (A, b) a maximum solution, if x ≤ $\hat{x}$, $\forall$x∈X (A, b). Similarly, $\hat{x} \in$ X (A, b) is called a minimum solution, if x ≤ $\hat{x}$ implies x = $\hat{x}$, $\forall$ x ∈ X(A, b). According to [11, 34], when X (A, b) ≠ ∅, it can be completely determined by one maximum solution and a finite number of minimum solutions. The maximum solution can be obtained by assigning.

$$\hat{x} = A @ b = \left[ \bigwedge_{j=1}^{n} (a_{ij} @ b_j) \right]_{i \in I} \qquad (4)$$

where

$$a_{ij} @ b_j = \begin{cases} 1 & if a_{ij} \le b_j, \\ b_j & if \ a_{ij} > b_j. \end{cases} \qquad (5)$$

Moreover, if we denote the set of all minimum solutions by $\check{X}$ (A, b), then

$$X (A, b) = \bigcup_{\check{x} \in \check{X}(A,b)} \{x \in X \mid \check{x} \le x \le \bar{x}\}. \qquad (6)$$

Now, we take a close look at X (A, b).

For proof of the following result please refer [24].

**Lemma 1.20.1:** If x ∈ X (A, b), then for each j ∈ J there exists $i_0 \in I$ such that $x_{io} \wedge a_{ioj} = b_j$ and $x_i \wedge a_{ij} \le b_j$, $\forall$i∈I.

**Lemma 1.20.2.** If X (A, b) ≠ ∅, then $I_j \ne \emptyset$ , $\forall$ j ∈ J.

**Lemma 1.20.3:** If X (A, b) ≠ ∅ then $\Lambda \ne \emptyset$ .

**THEOREM 1.20.1.** *Given that X (A, b) ≠ ∅* .

1. If $f \in \Lambda$, then $F(f) \in X(A, b)$.
2. For any $x \in X(A, b)$, there exists $f \in \Lambda$, such that $F(f) \le x$.

**Lemma 1.20.4:** If $c_i \le 0$, $\forall \, i \in I$, then $\hat{x}$ an optimal solution of problem (3).

**Lemma 1.20.5:** If $c_i \ge 0$, $\forall \, i \in I$, then one of the minimum solutions is an optimal solution of problem (3).

**THEOREM 1.20.2:** If $X(A, b)$ is non-empty and $x^*$ is defined according to

$$x^* = \begin{cases} \hat{x}_i^* & if \; c_i \ge 0, \\ \breve{x}_i & if \; c_i < 0, \end{cases} \qquad \forall \; i \in I.$$

then $x^*$ is an optimal solution of problem (3) with an optimal value

$$Z^* = c^T x^* = \sum_{i=1}^{m} \left( c_i^{''} \hat{x} + c_i^{'} \breve{x}_i^* \right).$$

Since $\breve{X}(A, b) \subset F(\Lambda)$ is implied by

$$\breve{X}(A, b) \subset F(\Lambda) \subset X(A, b)$$

when $c_i^{'} \ge 0$, $i = 1, \ldots, m$, solving

$$\begin{aligned} & \text{minimize } Z = \sum_{i=1}^{m} c_i^{'} x_i \\ & \text{subject to} \quad x \circ A = b \qquad\qquad (7) \\ & \qquad\qquad\quad x_i \in [0, 1] \end{aligned}$$

becomes equivalent to finding an $f^* \in \Lambda$ such that

$$\sum_{i=1}^{m} c_i^{'} F_i(f^*) = \min_{f \in \Lambda} \left\{ \sum_{i=1}^{m} c_i^{'} F_i(f) \right\}. \qquad (8)$$

Remembering the definition of $I_j = \{ i \in I \, / \, \hat{x}_i \, \wedge \, a_{ij} = b_j \}$ for all $j \in J$, we define variables

$$x_{ij} = \begin{cases} 1 & \text{if } i \text{ is chosen from } I_j \\ 0 & \text{otherwise}, \end{cases} \quad \forall\, i \in I, j \in J, \qquad (9)$$

and consider the following 0-1 integer programming problem:

$$\text{minimize } Z = \sum_{i=1}^{m}\left( c_i^{'}\, \max_{j \in J}\{b_j x_{ij}\} \right)$$

$$\text{subject to} \quad \sum_{i=1}^{m} x_{ij} = 1, \quad \forall\, j \in J, \qquad (10)$$

$$x_{ij} = 0 \text{ or } 1, \; \forall\, i \in I, j \in J,$$

$$x_{ij} = 0, \; \forall\, I, j \text{ with } i \notin I_i.$$

Note that the constraints of the above problem require that $\forall\, j \in J$, there exists exactly one $i \in I_j$, such that $x_{ij} = 1$. In this case, if we define $f = (f_1, ..., f_n)$ with $f_j = i$ whenever $x_{ij} = 1$, then $f \in \Lambda$. On the other hand, for any $f \in \Lambda$, by definition (9), we know that it corresponds to a feasible solution of problem (10). Moreover, from the definition of F, for any given $f \in \Lambda$, we have one feasible solution $x_{ij}$ of problem (10). Obviously, $\sum_{i=1}^{m} c_i^{'} F_i(f^*) = \sum_{i=1}^{m} c_i^{'}$ $\max_{j \in J}\{b_j x_{ij}\}$. Therefore, solving problem (10) is equivalent to finding an $f^* \in \Lambda$ for problem (8) via the relation defined by (9). In other words, solving problem (7) is equivalent to solving the 0-1 integer programming problem (10).

   While there are many different methods for solving integer programming problems, here we apply the commonly used branch-and-bound concept to solve problem (10). A branch-and-bound method implicitly enumerates all possible solutions to an integer programming problem. For our application, in the beginning, we choose one constraint to branch the original problem into several subproblems. Each subproblem is represented by one node. Then branching at each node is done by adding one additional constraint. New subproblems are created and represented by new nodes. Note that the more constraints added to a subproblem, the smaller feasible domain it has and, consequently, the larger optimal objective value Z it achieves.

   Based on the theory we built in previous sections, we propose an algorithm for finding an optimal solution of problem (3).

*Step1:*
Find the maximum solution of system (1). Compute $\hat{x} = A \text{ @ } b = \left[ \Lambda_{j=1}^{n} (a_{ij} \text{ @ } b_j) \right]_{i \in I}$ according to (4).

*Step 2:*
Check feasibility.
If $\hat{x} \text{ o } A = b$, continue. Otherwise, stop! $X(A, b) \neq \emptyset$ and problem (3) has no feasible solution.

*Step 3:*
Compute index sets.
Compute $I_j = \left\{ i \in I \mid \hat{x}_i \wedge a_{ij} = b_j \right\}$, $\forall \text{ j} \in \text{J}$.

*Step 4:*
Arrange cost vector.
Define c' and c'' according to

$$c'_i = \begin{cases} c_i & if \ c_i \geq 0, \\ 0 & if \ c_i < 0, \end{cases}$$

$$c''_i = \begin{cases} 0 & if \ c_i \geq 0, \\ c_i & if \ c_i < 0, \end{cases}$$

and define problem (7).

*Step 5:*
Define 0-1 integer program.
Define problem (10) via relation (9).

*Step 6:*
Solve integer program.
Use the branch-and-bound concept to solve (10).

*Step 7:*
Define $f = (f_1, ..., f_n)$ with $f_i = I$ if $x_{ij} = 1$.
Generate F(*f*) via formula

$$F_i(f) = \begin{cases} \max_{j \in J_f^i} b_j & if \ J_f^i \neq \phi \\ 0 & if \ J_f^i = \phi \ \forall \ i \in I. \end{cases}$$

*Define* $\breve{x}^* = \left( \breve{x}_1^*, ..., \breve{x}_m^* \right)^T$ *with* $\breve{x}_i^* = F_i$ *(f),* then $\breve{x}^*$ is an optimal solution of problem (7).

*Step 8*: Output.

1. Here, we have studied a linear optimization problem subject to a system of fuzzy relation equations and presented a procedure to find an optimal solution.

2. Due to the non-convexity nature of its feasible domain, we tend to believe that there is no polynomial-time algorithm for this problem. The best we can do here is that, after analyzing the properties of its feasible domain, we convert the original problem into a 0-1 integer programming problem, then apply the well-known branch-and-bound method to find one solution. The question of how to generate the whole optimal solution set is yet to be investigated.

3. From the analysis of Theorem 1.20.2, it is clearly seen that if all minimum solutions of a given system of fuzzy relation equations can be found, then an optimal solution of the optimization problem defined by (3) can be constructed. Therefore, solving this optimization problem is no harder than solving this optimization problem is no harder than solving a system of fuzzy relation equations for all minimum solution. Although it is not known whether these two problems are essentially equivalent or not, the basic concepts introduced have been further developed for solving fuzzy relation equations [52].

4. When a system of fuzzy relation equation (3) is derived for a particular application, such as the medical diagnosis [2], it is relatively easy to check if X(A, b) ≠ ∅ by taking the first two steps of the proposed algorithm. A liner function may not truly reflect the associate cost objective, but it can serve as an approximation. Extension to other types of objective functions is under investigation.

For more about these properties and proofs of the results please refer [24].

## 1.21 Some properties of minimal solution for a FRE

Fuzzy relation equation occurs in practical problems for example in fuzzy reasoning. Therefore it is necessary to investigate properties of the set of solutions. Here [39] have given a necessary and sufficient condition for existence of a minimal solution of a fuzzy relation equation defined on infinite index sets.

Let I and J be the index sets, and let $A = (a_{ij})$ be a coefficient matrix, $b = (b_j)$ be a constant vector where $i \in I$, $j \in J$. Then an equation

$$x \ o \ A = b, \qquad (1)$$

or,

$$\underset{i \in I}{V} (x_i \wedge a_{ij}) = b_j \ \text{ for all } j \in J. \qquad (1)$$

$x$ is called a fuzzy relation equation, where o denotes the sup-min composition, and all $x_i$ , $b_j$, $a_{ij}$'s are in the interval [0, 1]. An x which satisfies Eq. (1) is called a solution of Eq. (1). Fuzzy relation equation occurs in practice. For example, in fuzzy reasoning [70], when the inference rule and the consequences are known, a problem to determine antecedents to be used reduces to one of solving a FRE.

**DEFINITION 1.21.1:** *Let $(P, \leq)$ be a partially ordered set (poset) and $X \subset P$. A minimal element of X is an element $p \in X$ such that $x < p$ for $x \in X$. The greatest element of X is an element $g \in X$ such $x \leq g$ for all $x \in X$.*

**DEFINITION 1.21.2:** *Let $a = (a_i)$ and $b = (b_i)$ be vectors. Then the partial order $\leq$, the join $\vee$, and the meet $\wedge$ are defined as follows:*

$$a \leq b \Leftrightarrow a_i \leq b_i \quad \text{for all } i \in I$$
$$a \vee b \underset{=}{\Delta} (a_i \vee b_i), \qquad a \wedge b \underset{=}{\Delta} (a_i \wedge b_i).$$

**DEFINITION 1.21.3:** *Let $([0, 1]^I, \leq)$ be a poset with the partial order defined in Definition 1.21.2., and let $\aleph \subset [0, 1]^I$ be the solution set of Eq. (1). The greatest element of $\aleph$, a minimal* *element of* $\aleph$, *and* $\aleph^0$ *are called the greatest solution, a minimal solution, and a set of minimal solutions of Eq. (1), respectively,*

**DEFINITION 1.21.4:** *For a, b $\in$ [0, 1]*

$$a \ \alpha \ b \ \underline{\underline{\Delta}} \ \begin{cases} 1, & \text{if } a \leq b, \\ b, & \text{otherwise} \end{cases}$$

*Moreover*

$$A \ @ \ \boldsymbol{b^{-1}} \ \underline{\underline{\Delta}} \ \left[ \underset{j \in J}{\Lambda} \ a_{ij} \alpha b_j \right],$$

*where $b^{-1}$ denotes the transposition of vector b.*

**THEOREM 1.21.5 [84] :**
$$\aleph \neq \emptyset \ \Leftrightarrow A \ @ \ b^{-1} \in \aleph$$
*and then, A @ $b^{-1}$ is the greatest solution of Eq. (1).*

**THEOREM 1.21.6 [34]:** *When the index sets I and J are both finite $\aleph \neq \emptyset$ implies $\aleph \neq \emptyset$, and then*
$$x \in \aleph \Leftrightarrow (\exists \ \breve{x} \in \aleph) \ ( \bar{x} \leq x \leq \hat{x} ).$$

**THEOREM 1.21.7 [104]:** *Let the index set I be a metric compact space, and*
$$\aleph \underline{\underline{\Delta}} \ \{ \ x \in \aleph \ / \ x \text{ is upper semicontinuous on I} \}.$$
*If $\aleph_{use} \neq \emptyset$ holds, then $\aleph_{use} \neq \emptyset$, and for all $x \in \aleph$, there exists $x_{usc} \in \aleph_{use}^0$ with $x_{usc} \leq x$.*

**DEFINITION 1.21.8:** *The solution $x = (x_i) \in \aleph$ is attainable for $j \in J$ if there exists $I_j \in$ such that $x_i \wedge a_{i,j} = b_j$, and the solution $x = (x_i) \in \aleph$ is unattainable for $j \in J$ if $x_i \wedge a_{i,j} < b_j$ for all $i \in I$ .*

**DEFINITION 1.21.9:** *The solution $x \in \aleph$ is called an attainable solution if x is attainable for $j \in J$, the solution $x \in \aleph$ is called an unattainable solution if x is $\aleph_j^{(-)}$ for all $j \in J$, and the solution x $\in \aleph$ is called a partially attainable solution if $x \in \aleph$ is neither an attainable solution nor an unattainable solution. In other words*

$x \in \aleph$ *is an attainable solution* $\Leftrightarrow x \in \aleph_j^{(+)}$,

$x \in \aleph$ *is an unattainable solution* $\Leftrightarrow x \in \aleph_j^{(-)}$,

$x \in \aleph$, is a partially attainable solution

$$\Leftrightarrow x \in \aleph - \aleph_j^{(+)} - \aleph_j^{(-)}.$$

The set of all partially attainable solutions is denoted by $\aleph_j^{(*)}$.

**Remark 1.21.1:** It should be noticed that

$$\aleph = \aleph_j^{(+)} \cup \aleph_j^{(-)} \cup \aleph_j^{(*)}$$
$$\aleph_j^{(+)} \cap \aleph_j^{(-)} = \aleph_j^{(-)} \cap \aleph_j^{(*)} = \aleph_j^{(*)} \cap \aleph_j^{(+)} = \emptyset$$
$$J_1 \subset J_2 \Rightarrow \aleph_j^{(+)} \supset \aleph_j^{(+)}.$$

Note that when the index sets I and J are both finite, all solutions are attainable solutions, that is $\aleph_j^{(-)} = \aleph_j^{(*)} = \emptyset$.

**Remark 1.21.2 [69]:** Let x be the element of $\aleph$, then
$$(\exists j \in J) \, (b_j = 0) \Rightarrow x \in \aleph_j^{(+)}.$$

**Lemma 1.21.1:** Let x and y be the element of $\aleph$

    i.     $x \leq y$ and $x \in \aleph_j^{(+)} \Rightarrow y \in \aleph_j^{(+)}$

    ii.    $x \geq y$ and $x \in \aleph_j^{(-)} \Rightarrow y \in \aleph_j^{(-)}$

    iii.   $x \leq y$ and $x \in \aleph_j^{(*)} \Rightarrow y \in \aleph_j^{(*)} \cup \aleph_j^{(+)}$

    iv.   $x \geq y$ and $x \in \aleph_j^{(*)} \Rightarrow y \in \aleph_j^{(*)} \cup \aleph_j^{(-)}.$

**THEOREM 1.21.10:** *Let* $\hat{x}$ *be the greatest solution of Eq. (1) defined in Theorem 1.21.5. If J is finite set,*
$$\hat{x} \in \aleph_j^{(+)} \Leftrightarrow \aleph^0 \neq \emptyset.$$

**DEFINITION 1.21.11:**

$$I_j^1 \underline{\underline{\Delta}} \{i \in I \mid a_{ij} > b_j\}, \qquad I_j^2 \underline{\underline{\Delta}} \{i \in I \mid a_{ij} = b_j\},$$
$$I_j^3 \underline{\underline{\Delta}} \{i \in I \mid a_{ij} < b_j\}.$$

**Lemma 1.21.2.**

$$\aleph_j^{(+)} = \bigcup_{k \in I_j^1 \cup I_j^2} [\breve{x}_j(k), \hat{x}_j$$

where $\hat{x}_j = (\hat{x}_{ji})$ and $\breve{x}_j(k) = (\breve{x}_{ji}(k))$ are defined as follows:

$$\hat{x}_{ji} = \begin{cases} b_j, & if \ i \in I_j^1, \\ 1, & if \ i \in I_j^2 \cup I_j^3, \end{cases}$$

$$\breve{x}_{ji}(k) = \begin{cases} b_j, & if \ i = k, \\ 0, & if \ i \neq k. \end{cases}$$

**Lemma 1.21.3 [5]:** A finite poset has at least a minimal element.

**Lemma 1.21.4:** If J is a finite set and $\aleph_j^{(+)}$ is a nonvoid set, there exists at least one minimal solution $x_0 = (x_{0i}) \in \aleph_j^{(+)}$. Moreover, the cardinality of $\{i \in | x_{oi} > 0\}$ is finite.

**Lemma 1.21.5:** If x is the element of $\aleph^{(-)}$, there exists $y \in \aleph^{(-)}$ with x > y.

For proof please refer [39]. Each of these proofs can be taken up as an exercise by studious students.

## 1.22 Fuzzy relation equations and causal reasoning

D. Dubois and H. Prade have analyzed the fuzzy set approach to diagnosis problems and proposed a new model, more expressive for representing the available information, and where the intended meaning of the membership degrees has been clarified here, they pertain to uncertainty.

They have applied their model to a fault diagnosis problem in satellites. However this model have several limitation. They mention two of them.

The relational model given by them associates directly disorders and manifestations. This model is not capable of capturing the most general kind of incomplete information.

Let $S$ be a system whose current state is described by means of an n-tuple of binary attributes $(a_1,\ldots, a_i,\ldots, a_n)$. When $a_i = 1$ the manifestation $m_i$ is said to be present; when $a_i = 0$, it means that $m_i$ is absent. When there is no manifestation present, $S$ is said to be in its normal state and this state is described by the n-tuple $(0,\ldots, 0,\ldots, 0)$. Let $M$ denote the set of the n possible manifestations $\{m_1,\ldots, m_i,\ldots, m_n\}$. Let $D$ be a set of possible disorders $\{d_1,\ldots, d_j,\ldots, d_k\}$. A disorder can be present or absent. To each $d_i$ is associated the set $M(d_i)$ of manifestations which are entailed, or preferably, caused, produced, by the presence of $d_j$ alone. The completely informed case is first considered, where all the present manifestations are observed and where the set of manifestations which appear when a disorder is present is perfectly known. Thus if $m_i \notin M(d_j)$ it means that $m_i$ is not caused by $d_j$. A relation R on $D \times M$ is thus defined by $(d_j, m_i) \in R \Leftrightarrow m_i \in M(d_j)$, which associates manifestations and disorders.

Given a set $M^+$ of present manifestations which are observed, the problem is to find what disorder (s) may have produced the manifestations in $M^+$. It is supposed that the set $M^- = M - M^+ = \overline{M^+}$ is the set of manifestations which are absent, i.e. all manifestations which are present are observed. While deductive reasoning enables us to predict the presence of manifestations(s) from the presence of disorder(s), adductive reasoning looks for possible cause(s) of observed effects. In other words, one looks for plausible explanations (in terms of disorders) of an observed situation. Clearly while it is at least theoretically possible to find out all the possible causes which may have led to a given state of system $S$, the ordering of the possible solutions according to some levels of plausibility is out of the scope of logical reasoning, strictly speaking.

However one may for instance prefer the solutions which involve a small number of disorders, and especially the ones, if any, which rely on only one disorder. This is called the principle of parsimony. In case several disorders may be jointly present, the set of manifestations produced by the presence of a pair of disorders $\{d_i, d_j\}$ alone has to be defined, and more generally by a tuple of disorders. In the hypothesis that effect can be added and do not interfere, it holds that

$$M(\{d_i, d_j\}) = M(d_i) \cup M(d_j) \qquad (1)$$

and consequently

$$\overline{M(\{d_i, d_j\})} = \overline{M(d_i)} \cup \overline{M(d_j)} \ ,$$

i.e. the manifestations which are absent are those which are not produced by $d_i$ or $d_j$ separately. If this hypothesis is not acceptable, a subset M(D) of entailed manifestations should be prescribed for each subset $D \subseteq \mathcal{D}$ of disorders which can be jointly present. Under this new hypothesis, situations where disorder $d_i$ followed by $d_j$ has not the same effects in terms of manifestations as $d_j$ followed by $d_i$ are excluded. Since the set $\{d_i, d_j\}$ is not ordered. In other words, a relation on $2^D \times M$ is used, rather than on $D \times M$. If some disorders can never be jointly present, $2^D$ should be replaced by the appropriate set $A$ of associations of disorders which indeed make sense.

In the completely informed case described above, the following properties hold:

(i) $M^+ = \overline{M^-}$ , i.e. all present manifestations are observed, and equivalently all the manifestations which are not observed are indeed absent;

(ii) $\forall d, M(d) = M(d)^+ = \overline{M(d)^-}$ , where $M(d)^+$ (resp. $M(d)^-$) is the set of manifestations which are certainly present (resp. certainly absent) when disorder d alone is a event. When $M(D) \neq \cup_{d \in D} M(d)$, this condition is supposed to hold $\forall\ D \in 2$ in the completely informed case (and not only for D = {d}).

The potential set $\hat{D}$ of all disorders which can individually be responsible for $M^+$ is given by

$$\hat{D} = \{d \in D, M(d) = M^+\}. \tag{2}$$

Note that $M(d) = M^+ \Leftrightarrow \overline{M(d)} = M^-$. Namely no disorder outside $\hat{D}$ can cause $(M^+, M^-)$. Clearly, if $\hat{D} \neq \emptyset$, one must check for the set $\hat{D}\hat{D}$ of groups of disorders D, which together may have caused $M^+$,

$$\hat{D}\hat{D} = \{D \in A \subseteq 2^D, M(D) = M^+\} \tag{3}$$

Using the principle of parsimony, one might consider that the smaller the cardinality of D the more plausible it is. If M(D) can be obtained as $\cup_{d \in D} M(d)$, then the set $D_0$ of disorders which alone partially explains $M^+$,

$$D_0 = \{d \in D, M(d) \subseteq M^+\}, \qquad (4)$$

may be of interest for building elements of $\hat{D}\hat{D}$. Clearly $D_0 \supseteq \hat{D}$.

Now consider the case where information about the relation disorder–manifestation is still complete ($\forall\, d \in D$, M(d)= $M(d)^+$ = $\overline{M(d)^-}$ ), but where the present manifestations are not necessarily completely observed. When not all the information is available, the set $M^+$ of manifestations which are certainly present and the set $M^-$ of manifestations which are certainly absent no longer form a partition of $M$: indeed $M^+ \cap M^- \neq \emptyset$ but $M^+ \cup M^- \neq M$, or , equivalently, there exists a non–empty set $M^0 = M - (M^+ \cup M^-)$ of manifestations, the presence or absence of which are completely unknown. In other words, one may be unaware of the presence of some manifestations, and, perhaps, only a subset of the absent manifestations are known to be absent. Then (2) is changed into

$$\hat{D} = \{d \in D, M^+ \subseteq M(d) \subseteq \overline{M^-}\} \qquad (5)$$

since it is only known that the set of manifestations which are indeed present is lower bounded by the set $M^+$ and upper bounded by the set $\overline{M^-}$. Eq. (5) also writes

$$\hat{D} = \{d \in D, M(d) \cap M^- \neq \phi \text{ and } \overline{M(d)^-} \cap M^+ \neq \phi\}. \qquad (6)$$

The generalization of (5) – (6) to subsets D of joint disorders is straightforward.

Another particular case of incomplete information is when observations are complete ($M^+ = \overline{M^-}$ ), but some manifestations may sometimes be present or absent with a given disorder, i.e. for some d, it is sometimes not known if a manifestation m follows or not form d; in that case $m \notin M(d)^+$ and $m \notin M(d)^-$.

In other words, the union of the set $M(d)^+$ of manifestations which are certainly produced by d alone and the set $M(d)^-$ of manifestations which certainly cannot be caused by d alone, no longer covers $M$, i.e. $\exists d$, $M(d)^+ \cup M(d)^- \neq M$, but $M(d)^+ \cap M(d)^- \neq \emptyset$ always holds. Denoting $M(d)^0 = M$ $(M^+(d) \cup M^-(d))$, $m \in M(d)^0$ means that m is only a possible manifestation of d. In particular m may be absent or present when d is present. An unknown disorder can always be introduced, i.e., $d_0$ such that $M(d_0)^0 = M$, i.e. whose manifestations are unknown. Hence it is not a closed–world model. Then, the set $\hat{D}$ of potential disorders which alone can explain $M^+ = \overline{M^-}$ is given by

$$\hat{D} = \{d \in D, M(d)^+ \subseteq M^+ \text{ and } M(d) = \subseteq M^-\}. \qquad (7)$$

Obviously, when $M(d)^+ = \overline{M(d)^-} = M(d)$, (2) is recovered. Besides, (7) can be easily generalized over to a subset of joint disorders.

In the general case, both the information pertaining to the manifestations and the information relative to the association between disorders and manifestations is incomplete. Then d belongs to the set $\hat{D}$ of potential disorders each of which can alone explain both M+ and M– if and only if d does not produce with certainty any manifestation which is certainly absent in the evidence, and no observed manifestations are ruled out by d. Formally,

$$\hat{D} = \{d \in D, M(d)^+ \subseteq \overline{M^-} \text{ and } M(d)^- \subseteq \overline{M^+}\}. \qquad (8)$$

This also writes

$$\hat{D} = \{d \in D, M(d)^+ \cap M^- = \phi \text{ and } M(d)^- \cap M^+ = \phi\}. \qquad (9)$$

Clearly, (8) reduces to (2) in the completely informed case, since then $\overline{M^-} = M^+$ and $M(d)^- = \overline{M(d)^+}$. More generally, it is worth noticing that $\hat{D}$ gathers all the manifestations in $D$ which cannot be ruled out by the observations. If $M^- \neq \phi$ and $M(d)^- \neq \phi$ i.e. no information is available on the manifestations certainly absent, it can be verified that $\hat{D} = D$ and the whole set of possible disorders is obtained.

The unknown disorder $d_0$ whose manifestations are unknown is such that $M(d_0)^+ = \phi = M(d_0)^-$. Hence $d_0 \in \hat{D}$, i.e. $d_0$ can be always considered as potentially responsible for observed manifestations. The membership test (8) is thus very permissive, i.e. $\hat{D}$ can be very large, and contains what may look like irrelevant causes, since they cannot be ruled out. Namely a disorder d may belong to $\hat{D}$ defined by (8) – (9) even if $M(d)^+ \cap M^+ = \phi$ and $M(d)^- \cap M^- = \phi$. Indeed, in this case, nothing forbids $M(d)^+ \cup M(d)^- \subseteq M^0$. This means that $\hat{D}$ includes disorders, no sure manifestations of which are observed, and no forbidden manifestations are for sure absent. Such disorders, no sure manifestations of which are observed, and no forbidden manifestations are for sure absent. Such disorders may still be present when $M^+$ and $M^-$ are observed, since all the sure information we have about them pertains to manifestations in $M^0$ about which no observation is available (this is true for the unknown disorder $d_0$).

In other words, (9) achieves a consistency-based diagnosis only, that is, only the disorders which are incompatible with manifestations observed as present or absent are rejected. In order to be more selective, one may turn to adductive reasoning. First, among the disorders in $\hat{D}$, consider those, if any, which have at least a weak relevance to the observations, namely the subset $\hat{D}*$ of $\hat{D}$, defined by

$$\hat{D}* = \{d \in \hat{D}, \, M(d)^+ \cap M^+ \neq \emptyset \text{ or } M(d)^- \cap M^- \neq \emptyset \}. \qquad (10)$$

This eliminates the disorders such that $M(d)^- \cup M(d)^- \subseteq M^0$, i.e. the disorders which are not suggested by the observations while not being ruled out by them either (and $d_0$ particularly). It is worth noticing that the refinement of $\hat{D}$ by $\hat{D}*$ makes sense only in the general case of incomplete information. Indeed, $\hat{D} = \hat{D}*$, when $\hat{D}$ is defined by (5) (if $M^+ \neq \emptyset$ or $M^- \neq \emptyset$) or (7) (if $M(d)^+ \neq \emptyset$ or $M(d)^- \neq \emptyset$). Obviously, one may still refine $\hat{D}*$ by strengthening the requirements on d, by asking for a conjunction rather than a disjunction in (10) or by replacing the conditions by more demanding ones like $M(d)^+ \subseteq M^+$ (all sure manifestations of d are observed) or $M(d)^+ \supseteq M^+$ (all observed manifestations are among the ones which certainly accompany d), and similar conditions for the subsets pertaining to the absence of manifestations. In particular, the subset $\hat{D}**$ of $\hat{D}$ defined by

$$\hat{D}^{**} = \{d \in \hat{D}, M(d)^+ \supseteq M^+ \text{ and } M(d)^- \supseteq M^-\} \qquad (11)$$

gathers the disorders, if any, which offer a complete coverage of the observations, but which may also have some sure effects or sure absence of effects which have remained unobserved (i.e. which are in $M^0$). A genuine abductive task is performed by (11).

When $\hat{D} \neq \emptyset$, one can look for explanations in terms of subsets of disorders which are not singletons. Eq. (3) is then extended by

$$\overline{DD} = \{D \in \overline{DD} \subseteq 2^D, M(D)^+ \subseteq \overline{M^-} \text{ and } \subseteq M(D)^- \subseteq \overline{M^+} \qquad (12)$$

for the subsets of disorders which jointly may explain $M^+$ and $M^-$. $\overline{DD}$ can obviously be refined by extending the counterparts of (10) or (11) for defining "relevant" subsets of disorders; for instance (10) is generalized by

$$\overline{DD} = \{D \in \overline{DD}\ M(d)^+ \cap M^+ \neq \emptyset \text{ or } M(d)^- \cap M^- \neq \emptyset\}. \qquad (13)$$

As expected, what is present and what is absent play symmetrical roles, exchanging $+$ and $-$ in (8) – (12).

Note that if $M^- \neq \emptyset$ i.e. if only manifestations which are certainly present are known, (9) (or(12)) may yield a result $\hat{D} \neq D$ (or $\overline{DD} \neq A$) provided that $\overline{M(d)^-} \neq M$, i.e. we have non–trivial information on the set of manifestations $\overline{M(d)^-}$ (or $\overline{M(D)^-}$ which may be produced by a disorder d (or a subset of disorders D) alone; indeed $\overline{M(d)^-} \supseteq M(d)^+$, (resp. $\overline{M(d)^-} \supseteq M(D)^+$), and thus $\overline{M(d)^-}$ (resp. $\overline{M(D)^-}$ gathers the manifestations in $M(d)^+$ (resp. $M(D)^+$} which are certainly produced by d (resp.D) and the manifestations for which we do not know if they can or cannot follow from d (resp.D). For $M^- \neq \emptyset$, (9) and (12) write

$$\hat{D} = \{d \in D, \overline{M(d)^-} \supseteq M^+\} \qquad (14)$$

$$\hat{D}\hat{D} = \{D \in {}^* \subseteq 2^D, \overline{M(d)^-} \supseteq M^+\}. \qquad (15)$$

In the non–completely informed case, the hypothesis (1) that effects can be added and do not interfere writes {for two disorders)

$$M(\{d_i, d_j\})^+ = M(d_i)^+ \cup M(d_j)^+ \qquad (16a)$$

and $$M(\{d_i, d_j\})^- = M(d_i)^- \cap M(d_j)^-. \qquad (16b)$$

Clearly, (16) reduces to (1) in he completely informed case. Note that the second equality of (16) still writes

$$\overline{M(\{d_i, d_j\})^-} = \overline{M(d_i)^-} \cup \overline{M(d_j)^-}$$

which says that the possible manifestations of two simultaneous disorders gather the manifestations possibly produced by each disorder, as for certain manifestations.

According to the interpretation we have in mind, it will lead to different models with different interpretations of the results.

In this framework, modeling a quantifier like "most" comes down to assigning a high degree of importance to k manifestations, arbitrarily chosen in $M^+$, with k "close to" the number of observed effects in $M^+$ (this can be defined as a fuzzy set defined on the set of integers). In order to estimate to what extent the disorder $d_j$ explains "most" of the effects in $M^+$ we compute the maximum of $\mu_{\hat{D}^-}(d_j)$ on all permutations $\sigma$ of [1, n], i.e. we compute

$$\max_{\sigma} \min_{i} [\min(w_{\sigma(i)}, \mu_{M^+}(m_i)) \to \mu_R(d_j, m_i)],$$

where $w_{\sigma(i)} = 1$ if $\sigma(i) \in [1, k]$ and $w_\sigma(i) = 0$ otherwise; we get the set of causes which alone explain at least k manifestations in $M^+$. This readily extends to subsets of causes of a given cardinality. In the general case, the weights where $w_{\sigma(i)}$ may lie in [0, 1] and capture the idea of taking k as a fuzzy number. The analogy between this quantification problem and the one addressed by Yager's [113] OWA is worth noticing.

The modeling of uncertainty remains qualitative in the above approach. Indeed, we could use a finite completely ordered chain of levels of certainly ranging between 0 and 1, i.e. $l_1 = 0 < l_2 < \ldots < l_n = 1$ instead of [0, 1], with ,min ($l_i$, $l_k$) = $l_i$ and max ($l_i$, $l_k$) = $l_k$ if i ≤ k, and $1 - l_i = l_{n+1-i}$, and $l_i \to l_k = 1$ if $l_i \le l_k$, $l_i \to l_k = l_k$ if $l_i > l_k$. Taking into account the incomplete nature of the information about the presence or absence of manifestations decreases the discrimination power when going from the completely informed case (Eq. (2)) to the incomplete information case (Eqs. (8) or (9)), since then the number of possible disorders in $\hat{D}$ increases. This is due to the fact that now there are manifestations which are neither certain nor impossible and consequences of the presence of a given disorder d which are only possible, as picture in Figure 1.22.2., while $M(d)^+ = M(d)$ and $M(d)^- = \overline{M(d)}$ in Figure 122.1. (Similar figures can be drawn for $M^+$ and $M^-$).

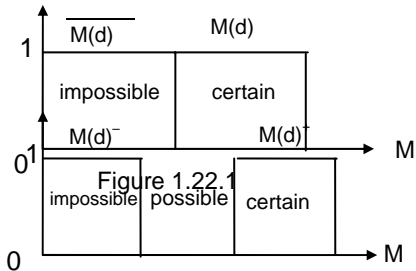

Figure 1.22.1

Figure 1.22.2

This suggests that in order to improve the discrimination power of the model, we have to refine the non–fuzzy model in such a way that consequences (resp. manifestations) previously expressed as certain (resp. certainly present) and impossible (resp. certainly absent) remain classified in the same way and where some possible consequences (resp. possibly present manifestations) are now allowed to be either somewhat certain (resp. somewhat certainly present) or somewhat impossible (resp. somewhat certainly absent); see Figure 1.22.3.

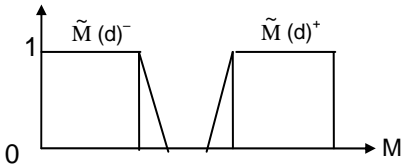

Figure 1.22.3

## 1.23 Identification of FRE of fuzzy neural networks

Blanco and et al [6] have established that any fuzzy system described by a max-min fuzzy relational equation may be identified by using a max-min fuzzy neural network.

The introduction of the smooth derivative by the authors have improved the effectiveness of the training process their method not only identifies the system but also solves the associated fuzzy relational equation. Any fuzzy system can be represented by a fuzzy relational equations system, and thus to identify, it forces us to solve equations like $X \oplus R = Y$, where X and Y are inputs and outputs, respectively, and where the composition operation $\oplus$ is generally a combination t-co-norm / t-norm.

*Identification of FRE by fuzzy networks without activation function*

### The problem
Our objective is to identify a fuzzy system through solving a fuzzy relational equation by a max-min fuzzy neural network. We will assume the fuzzy relational equation is $X \oplus R = Y$, $X \in [0, 1]'$, $Y \in [0, 1]'$, $R \in [0, 1]^{rxs}$. We will limit ourselves to the case $\oplus$ = max-min. We also suppose that we have a set of examples $[X^i, Y^i = i = 1,..., p]$ to solve R, and we will use a neural max-min for the identification (by using these examples, to train the neural network).The problem is to design the neural network (its topology) and the learning method.

### Net topology
We are going to consider a fuzzy network with the following topology: The input-out pairs are $(x_1,..., x_i,..., x_r)$ and $(Out_1,..., Out_j, ..., Out_s)$, where $Out_j$ is determined by $Out_j = \max[\min (x_i, w_{ij})]$, the $w_{ij}$ being the elements of the weight matrix W that assesses the synaptic connection (see Figure 1.23.1).

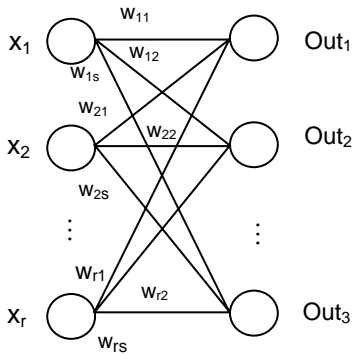

Figure: 1.23.1

Let us observe that no activation function is considered here. The objective of the activation function in an artificial neural is twice:

        i     To adapt the output to prefix range.

        ii    To fix a threshold.

In our case, we use the max-min operation, and obviously, the output range is fixed by the operation, it is [0, 1], so the first objective is reached.

On the other hand, the operation min is a threshold function, for each input $x_j$ and each weight $w_{ij}$, which represents the saturation level.

$$\min (x_j, w_{ij}) = \begin{cases} x_j & if \ x_i \leq w_{ij}, \\ w_{ij} & if \ x_i > w_{ij}. \end{cases}$$

Thus, use of any activation function is not needed afterwards to apply max-min to the input, since it is underlying in its own process. So, we are considering a neural network without a hidden layer, where the inputs are the values $X \in [0, 1]$' and the outputs $Y \in [0, 1]$' are obtained by $Y = \max(\min (W, X))$, W being the weight matrix. If $X = (x_1, x_2,\ldots, x_r)$, $Y = (Out_1,\ldots, Out_2, \ldots, Out_s)$ and the elements of the W matrix are $w_{ij}$, the outputs are obtained such that

$Out_1 = \max [\min (x_1, w_{11}), \min (x_2, w_{21}),\ldots, \min (x_r, w_{r1})],$
$\qquad\qquad \vdots$
$Out_s = \max [\min (x_1, w_{1s}), \min (x_2, w_{2s}),\ldots, \min (x_r, w_{rs})],$

The objective of training the network is to adjust the weights so that the application of a set of inputs produces the desired set of outputs. This is driven by minimizing the square of the difference between the desired output $T_j$ and the actual $O_j$, for all the patterns to be learnt, $E = \frac{1}{2} \Sigma(T_j - O_j)^2$, where $O_j = \max_i (\min (x_i, w_{ij}))$.

It is well known that

$$\frac{\partial E}{\partial w_{ij}} = \frac{\partial E}{\partial O_j} * \frac{\partial O_j}{\partial w_{ij}} \qquad\qquad (1)$$

in any interval where the derivatives are defined.

Let us expand the second factor of (1):

$$\frac{\partial O_j}{\partial w_{ij}} = \underbrace{\frac{\partial Max[Min(x_i, w_{ij}), \ i = 1,...,r]}{\partial Min \ (x_s, w_{sj})}}_{p1} * \underbrace{\frac{\partial Min \ (x_s, w_{sj})}{\partial w_{sj}}}_{p2}$$

$$P1 = \frac{\partial Max \left\{ Min(x_s, w_{sj}), \underset{i \neq s}{Max}[Min(x_i, w_{sj})] \right\}}{\partial Min \ (x_s, w_{sj})}.$$

Let us recall that the functions Min (y, p), Max (y, p) are derivable into the open intervals y < p and y > p but their derivative is not defined in y = p (see Figure 1.23.2).

$$\frac{\partial Min(y, \ p)}{\partial y} = \begin{cases} 1 & if \ y < p, \\ 0 & if \ y > p, \end{cases}$$

$$\frac{\partial Max(y, \ p)}{\partial y} = \begin{cases} 1 & if \ y > p, \\ 0 & if \ y < p, \end{cases}$$

As it is usual  we will assign in any case the value in y = p equal to the left or right derivative, respectively, that is, finally we will assume

$$\frac{\partial Min(y, \ p)}{\partial y} = \begin{cases} 1 & if \ y \leq p, \\ 0 & if \ y > p, \end{cases}$$

$$\frac{\partial Max(y, \ p)}{\partial y} = \begin{cases} 1 & if \ y \geq p, \\ 0 & if \ y < p, \end{cases}$$

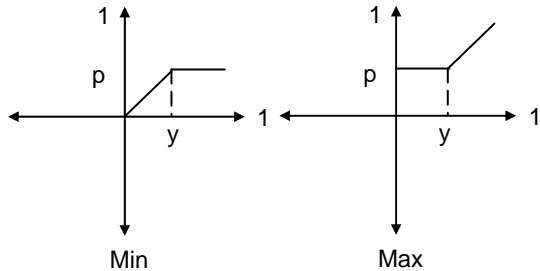

Figure: 1.23.2

According to these ideas

$$P1 = \begin{cases} 1 \; if \; Min \, (x_s, \, w_{ssj}) \geq \underset{i \neq s}{Max}(Min(x_i, w_{ij})), \\ 0 \; if \; Min \, (x_s, \, w_{sj}) < \underset{i \neq s}{Max} \, (Min \, (x_i, w_{ij})), \end{cases}$$

$$P2 = \frac{\partial Min \, (x_s, \, w_{sj})}{\partial w_{sj}} = \begin{cases} 1 & if \; x_s \geq w_{sj}, \\ 0 & if \; x_s < w_{sj}. \end{cases}$$

By combining P1 and P2 we will obtain the value of (2):

$$\frac{\partial O_j}{\partial w_{sj}} = \begin{cases} x_s < w_{sj} \begin{cases} \{ x_s \geq \underset{i \neq s}{Max}(Min(x_i, w_{ij})) \rightarrow C1 \\ \{ x_s < \underset{i \neq s}{Max}(Min(x_i, w_{ij})) \rightarrow C2 \end{cases} \\ x_s \geq w_{sj} \begin{cases} w_{sj} \geq \underset{i \neq s}{Max}(Min(x_i, w_{ij})) \rightarrow C3 \\ w_{sj} < \underset{i \neq s}{Max}(Min(x_i, w_{ij})) \rightarrow C4. \end{cases} \end{cases}$$

So, the value of (2), say C, will be 0 in the cases C1, C2 or C4 and 1 in C3.

Now we are going to expand the first factor of (1), for that we will denote - $\partial E / \partial O_j$ by $\delta_j$ $\partial E / \partial O_j = - (T_j - O_j)$.

Therefore, $\partial E / \partial w_{ij} = \delta_j C$. Finally, the changes for the weights will be obtained from a $\delta$-rule with expression: $\Delta w_{ij} = \mu \delta_j C$, where $\delta_j = (T_j - O_j)$.

Thus, the learning algorithm is similar to the classical back propagation but C is defined equal to $C_1 - C_3$ or $C_4$ according to the values of input-output pairs.

By applying this learning process, it is not guaranteed that the network will learn, obviously because the value of C is null in three of the four cases C1-C4.

To improve this behavior we will develop a new procedure, the so-called "smooth derivative" which is based on the following idea.

$$\frac{\partial Min(y, \, p)}{\partial y} = \begin{cases} 1 & if \; y \leq p, \\ 0 & if \; y > p, \end{cases}$$

is just the "crisp" truth value of the proposition "y is less than or equal to p". In the same way

$$\frac{\partial Max(y,\, p)}{\partial y} = \begin{cases} 1 & if \ \ y \geq p, \\ 0 & if \ \ y < p, \end{cases}$$

is the "crisp" truth value of the proposition "y is greater than or, equal to p".

Thus, to improve the performance of the learning process, we are interested in changing this "crisp" behavior by one "fuzzy" being able to capture the real meaning of $(y \leq p)$ or $(y \geq p)$ in a vague context. Taking into account that we are measuring the relative position of y with respect to p, we propose to measure for each y the inclusion degree p in y, which we will denote by $\|p \subset y\|$.

In turn, we can apply any implication function to assess the inclusion degree, taking into account that when $p \leq y$ then $p \subseteq y$ with degree equal to 1, whereas when $p > y$ it is reasonable to consider the inclusion degree of p into y to be equal to y, it intuitively results that Godel implication is the most suitable one. The experiences have confirmed this intuition. So, we will consider

$$\|p \subset y\| = p \xrightarrow{\ G\ } y = \begin{cases} 1 & if \ \ p \leq y, \\ y & if \ \ p > y. \end{cases}$$

On the other hand, similar to before when dealing with max $\{a, b,…, c\}$, if we denote "max 1 = max $\{a, b,…, c\}$" and "max 2 = max $\{\{a, b,…, c\} - \{max\ 1\}\}$". We are interested in knowing the inclusion degree of max 1 in max 2, $\|max\ 1 \subset max\ 2\|$.

$$\|max\ 1 \subset max\ 2 \| = max\ 1 \xrightarrow{\ G\ } max\ 2$$
$$= \begin{cases} 1 & if \ \ max\ 1 \leq max\ 2 \\ max\ 2 & if \ \ maxl > max\ 2. \end{cases}$$

So, when we use the inclusion degree, P1 and P2 not only have the values 0 or 1 but a value in [0, 1]. We observe that P2 in zero when $x_s, < w_{sj}$, but using the inclusion degree of $w_{sj}$ in $x_s$ $\|w \subset x\|$, a new value of P2 is obtained to be

$$P2 = \begin{cases} 1 & if \ \ w \leq x, \\ x & if \ \ w > x, \end{cases} \qquad\qquad (A)$$

We can treat the case for P1 similarly. Denoting max $2 = \max_{i \neq s}$ (Min $(x_i, w_{ij})$), note that P1 is obviously zero when min $(x_s, w_{sj}) <$ max 2. We take the inclusion degree of max 2 in min $(x_s, w_{sj})$:

$$\left\| \max 2 \subset \min(x_s, w_{sj}) \right\| = \begin{cases} 1 & \text{if } \max 2 \leq (\min(x_s, w_{sj})), \\ \min(x_x, w_{ij})) & \text{if } \max 2 > (\min(x_s, w_{sj})). \end{cases}$$

Finally, we obtain the P1 value:

$$P1 = \begin{cases} 1 & \text{if } \max 2 \leq (\min(x_s, w_{sj})), \\ \min(x_x, w_{ij})) & \text{if } \max 2 > (\min(x_s, w_{sj})). \end{cases} \qquad (B)$$

By combining (A) and (B), we will obtain

$$\frac{\partial O_j}{\partial w_{sj}} = \begin{cases} x_s < w_{sj} \begin{cases} \{ x_s \geq \underset{i \neq s}{Max}(Min(x_i, w_{ij})) \rightarrow C1 = x_s, \\ \{ x_s < \underset{i \neq s}{Max}(Min(x_i, w_{ij})) \rightarrow C2 = x_s * x_s \end{cases} \\ x_s \geq w_{sj} \begin{cases} w_{sj} \geq \underset{i \neq s}{Max}(Min(x_i, w_{ij})) \rightarrow C3 = 1, \\ w_{sj} < \underset{i \neq s}{Max}(Min(x_i, w_{ij})) \rightarrow C4 = w_{sj}. \end{cases} \end{cases}$$

Obviously, the values obtained from $\partial O_j / \partial w_{ij}$ depend on the implications chosen. We have made several trials with all the implications, except with those, which give us an inclusion degree of null.

## 1.24 Equations in classes of fuzzy relations

Drewniak [20] explains the problem of existence of concrete solutions and brings information about bounds of the family of all such solutions.

Fuzzy relation equations were introduced by Sanchez [84], [94] and have been investigated by many authors. We ask for the existence of solutions which belong to a concrete class of fuzzy relations [20, 21, 115]. The description of the family of all such solutions is also needed.

Let us consider the lattice $L = ([0, 1], \vee, \wedge, \rightarrow)$, where

$$a \vee b = \max (a, b), a \wedge b = \min (a, b), a \rightarrow b = \begin{cases} 1 & \text{if } a \leq b, \\ b & \text{if } a > b. \end{cases} \quad (1)$$

We deal with fuzzy relations R, S, T, U : $X \times X \rightarrow L$ on a finite set $X = \{x_1, \ldots, x_n\}$ with sup-inf composition:

$$(RS) (x, z) = \bigvee_{y \in X} (R(x, y) \wedge S(y, z)) \text{ for } x, z \in X. \quad (2)$$

Relation composition (2) is isotone, associative and

$$R(T \vee U) = RT \vee RU, \quad (3)$$
$$R(T \wedge U) \leq RT \wedge RU. \quad (4)$$

It has the identity element I,

$$IR = RI = R. \quad (5)$$

$$I(x, y) = \begin{cases} 1 & \text{if } x = y, \\ 0 & \text{if } x \neq y, \end{cases} \quad I'(x, y) = \begin{cases} 0 & \text{if } x = y, \\ 1 & \text{if } x \neq y. \end{cases} \quad (6)$$

And the null element $\phi$,

$$\phi R = R\phi = \phi, \phi (x, y) = 0. \quad (7)$$

The existence of solution of the relation equation

$$RU = T \quad (8)$$

(with unknown relation U) was characterized by Sanchez [84]:

**THEOREM 1.24.1:** *Eq. (8) has solution iff RU\* = T, where*

$$U^* (x, z) = \bigwedge_{y \in X} (R(y, x) \rightarrow T(y, z)) \quad \text{for } x, z \in X. \quad (9)$$

*If Eq. (8) has solutions, then the above formula gives the greatest one. In general, we always have*

$$RU^* \leq T. \quad (10)$$

*A description of the family of all solutions of (8) is more complicated* [11, 15, 20, 34]*:*

**THEOREM 1.24.2:** *If Eq. (8) has solutions, then there exist $m \in N$ and solutions $U_1, U_2,\ldots, U_m$ such that for any solution U there exist $k \le m$ fulfilling*

$$U_k \le U \le U^*.$$

*The family of all solutions of (8) has the form*

$$U (R, T) = \bigcup_{1 \le k \le m} [U_k, U^*], \tag{11}$$

*where $[U_k, U^*]$ is a lattice interval of relations.*

In examples and proofs the fuzzy relation R will appear as a square matrix $R = [r_{i,k}]$, where

$$r_{i,k} = R(x_i, x_k) \text{ for i, k} = 1,\ldots, n.$$

Matrix R will also be described as a sequence of columns:

$$R = (r_1,\ldots, r_n), \qquad r_k = \begin{bmatrix} r_{1,k} \\ \vdots \\ r_{n,k} \end{bmatrix}, k = 1, 2,\ldots, n.$$

**Lemma 1.24.1:** If Eq. (8) has reflexive solutions, then $U^*$ is reflexive.

**Lemma 1.24.2:** $U^*$ is reflexive iff $R \le T$.

**THEOREM 1.24.3:** *The solvable Eq. (8) has reflexive solutions iff $R \le T$.*

**Lemma 1.24.3:** If $R \le T$, then
$$R \le RU^* \le T. \tag{12}$$

**THEOREM 1.24.4:** *If $R = T$, then Eq. (8) is solvable and $[I, U^*]$ is the set of all the reflexive solutions of (8).*

**Lemma 1.24.4:** If $R \le T$ and U is a solution of (8), then $I \vee U$ is a reflexive solution of (8).

**THEOREM 1.24.5:** *If Eq. (8) is solvable and $R \leq T$, then the minimal reflexive solutions of (8) belong to the set $\{U_1 \vee I, \ldots, U_m \vee I\}$ and the family of all the reflexive solutions has the form*

$$\bigcup_{1 \leq k \leq m} [U_k \vee I, U^*].$$

**COROLLARY:** *If all the minimal solutions $U_i$ are reflexive, then any solution of (8) is reflexive.*

**Lemma 1.24.5:** If the relation U is irreflexive, then any solution of Eq. (8) is irreflexive.

**Lemma 1.24.6:** The relation $U^*$ is irreflexive iff

$$\forall_{x \in X} \; \exists_{y \in Y} \; (R\,(y,\,x) \neq 0, T\,(y,\,x) = 0). \qquad (13)$$

**THEOREM 1.24.6:** *If Eq. (8) is solvable and relations R, T fulfill condition (13), then all the solutions are irreflexive.*

**Lemma 1.24.7:** If the fuzzy relation R in Eq. (8) is reflexive, then $U^* \leq T$.

**THEOREM 1.24.7:** *If the fuzzy relation R in (8) is reflexive and T is irreflexive, then any solution of (8) is irreflexive.*

**Lemma 1.24.8:** For any fuzzy relation U, the relation $I' \wedge U$ is irreflexive.

**Lemma 1.24.9:** If U is an irreflexive solution of (8), then the relation $U^{ir} = I' \wedge U^*$ is also an irreflexive solution and $U \leq U^{ir}$.

**THEOREM 1.24.8:** *Eq. (8) has irreflexive solutions iff the relation $U^{ir}$ is a solution of (8). Moreover, if Eq.(8) has irreflexive solutions, then the solution $U^{ir}$ is the greatest one.*

*A fuzzy relation R is symmetric iff $R^{-1} = R$, where*

$$R^{-1}(x,\,y) = R\,(y,\,x) \text{ for } x,\,y \in X. \qquad (14)$$

Directly from (14) we get:

**Lemma 1.24.10:** For any fuzzy relation U, relations $U \vee U^{-1}$ and $U \wedge U^{-1}$ are symmetric.

**THEOREM 1.24.9:** *Eq. (8) has symmetric solutions iff the relation*

$$U^s = U^* \wedge (U^*)^{-1} \qquad (15)$$

*is a solution of (8). Moreover, if Eq. (8) has symmetric solutions, then formula (15) gives the greatest one.*

**THEOREM 1.24.10:** *If Eq. (8) has symmetric solutions, then for any solution $U \leq U^*$, the fuzzy relation $U \vee U^{-1}$ is a symmetric solution. In particular, any minimal symmetric solution of (8) belongs to the set*

$$\{U_1 \vee (U_1)^{-1},..., U_m \vee (U_m)^{-1}\}.$$

**Lemma 1.24.11:** If the relation U is anti-symmetric (asymmetric), then any relation $R \leq V$ is also anti-symmetric (asymmetric).

**THEOREM 1.24.11:** *Eq. (8) has anti-symmetric (asymmetric) solutions iff at least one of the minimal solutions is anti-symmetric (asymmetric). Moreover, if $U^*$ is antisymmetric (asymmetric), then any solution of (8) is anti-symmetric (asymmetric).*

Using Lemmas 4.1. and 2.3, we get:

**THEOREM 1.24.12:** *If the fuzzy relation R in (8) is reflexive and T is anti-symmetric (asymmetric), then any solution of (8) is anti-symmetric (asymmetric).*

**THEOREM 1.24.13:** *If Eq. (8) has anti-symmetric (asymmetric) solutions, then there exist $p \in N$ and maximal anti-symmetric (asymmetric) solutions $U_1^*, ..., U_p^*$ of (8) such that any anti-symmetric (asymmetric) solution U is bounded from above by some solution $U_k^*$, $k \in \{1, 2, ...p\}$.*

**Lemma 1.24.12:** If Eq. (8) has complete (strongly complete) solutions, then $U^*$ is complete (strongly complete).

**Lemma 1.24.13:** $U^*$ is complete iff

$$\mathop{\forall}_{1 \le i \ne k \le n} (r_i \le t_k \text{ or } r_k \le t_i) \qquad (16)$$

and strongly complete iff

$$\mathop{\forall}_{1 \le i, k \le n} (r_i \le t_k \text{ or } r_k \le t_i) \qquad (17)$$

**THEOREM 1.24.14:** *The solvable Eq. (8) has complete (strongly complete) solutions iff relations R and T fulfill condition (16) (condition (17)).*

**THEOREM 1.24.15:** *If Eq. (8) has complete (strongly complete) solutions, then there exist $q \in N$ and minimal complete (strongly complete) solutions $U_1^c, ..., U_q^c$ of (8). The family of all the complete solutions of (8) has the form $\cup_{1 \le k \le q} \left[ U_k^c, U* \right]$.*

**Lemma 1.24.14:** For any fuzzy relation R, the least transitive relation $R^\vee$ containing R is given by the formula

$$R^\vee = \mathop{\vee}_{1 \le k \le n} R^k.$$

Now we need a useful notation from formula (9) (cf. [6]: $R \to T = U*$.

**Lemma 1.24.15:** (cf. Wagenknecht [105]). For any fuzzy relation R, the relation $R \to R$ is transitive.

**Lemma 1.24.16:** If Eq. (8) has transitive (idempotent) solution U, then

$$U \le T \to T.$$

**Lemma 1.24.17:** If Eq. (8) has transitive solutions, then the relation

$$U^t = (T \to T) \wedge U* \qquad (18)$$

is also a transitive solution.

**THEOREM 1.24.16:** *Eq. (8) has transitive solutions iff the relation (18) is a solution of (8). If Eq. (8) has transitive solutions, then formula (18) gives the greatest transitive solution.*

**THEOREM 1.24.17:** *If R ≤ T then Eq. (8) has transitive solutions iff the relation T→ T is a solution of (8).*

**THEOREM 1.24.18:** *If Eq. (8) has transitive solutions, then the minimal transitive solutions belong to the set*

$$\{(U_1)^\vee, (U_2)^\vee, \dots, (U_m)^\vee\}.$$

The reader is expected to prove the above theorem. For more about these properties please refer [21].

## 1.25 Approximate solutions and FRE and a characterization of t-norms that define metrics for fuzzy sets.

Gottwald [29] gives a necessary and sufficient condition, which characterizes all those t-norms, which yield a particular metric for fuzzy sets.

One of the basic relations for fuzzy sets A, B ∈ F (ℵ) over a given universe of discourse ℵ is their inclusion relation.

$$A \subseteq B \Leftrightarrow \mu_A(x) \le \mu_B(x) \text{ for all } x \in \aleph. \qquad (1)$$

Using some ideas from many-valued logic it is quite natural to extend that relation to a fuzzified, i.e., graded inclusion relation ⊆ for fuzzy sets, based, e.g., on the Lukasiewicz implication. Having done this, it again is natural to extend this fuzzification of inclusion between fuzzy sets to the equality of fuzzy sets and to define a fuzzified, i.e. graded equality ≡ for fuzzy sets.

We assume that a fixed universe of discourse ℵ is given which contains at least two elements. Each fuzzy set A over ℵ, i.e. each A ∈ F (ℵ) is characterized by its membership function

$$\mu_A : \aleph \to [0, 1].$$

The union of two fuzzy sets A, B is denoted A ∪ B and as usual characterized by the maximum of the respective membership degrees, the usual intersection is A ∩ B and defined using the min-operator, and the usual complement is denoted $\overline{A}$. For each a ∈ ℵ and each u ∈ [0, 1] the fuzzy u-singleton of a , denoted by: ⟨⟨a⟩⟩ is the fuzzy set C with membership function:

$$C := \langle\langle a \rangle\rangle_u : \mu_C(x) = \begin{cases} u \ if \ x = a, \\ 0 \ otherwise. \end{cases} \quad (2)$$

The membership degrees $\mu_A(x)$ are considered as truth degrees of a generalized, i.e. many-valued membership predicate $\varepsilon$ . For membership of a point $a \in \aleph$ in a fuzzy set $A \in F(\aleph)$ we then write

$$a \ \varepsilon \ A, \quad (3)$$

but, as usual in formal logic, now have to distinguish between the well-formed formula (3) and its truth degree $[\![ \ a \ \varepsilon \ A \ ]\!]$ which is nothing else than the usual membership degree.

$$[\![ \ a\varepsilon A \ ]\!] \in A =_{def} \mu_A(a). \quad (4)$$

This notation $[\![ \ \dots \ ]\!]$ for the truth degree will also be used in case … is a more complex expression than simply (3).

As basic tools to build up more complex expressions than (3), connectives for conjunction, implication, negation and a quantifier for generalization are used. As usual, those logical operators are characterized by the way they operate on the truth degrees, i.e., by their truth functions. The simplest case is the negation operator $\neg$, which is determined by

$$[\![ \ \neg H \ ]\!] =_{def} 1 - [\![ \ H \ ]\!] \quad (5)$$

if H is any well-formed formula of the set theoretic languages we just are constituting. Quite standard, too, is the understanding of the generalization quantifier $\forall$ for which generalization $\forall x$ with respect to all $x \in \aleph$ means to consider the infimum of the corresponding truth degrees:

$$[\![ \ \forall x H(x) \ ]\!] =_{def} \inf_{x \in *} [\![ \ H(x) \ ]\!] \quad (6)$$

A wide variety of possibilities exists to interpret the conjunction connective $\wedge$. We will allow any t-norm to be used as truth function for $\wedge$. By a t-norm, we understand a binary operation t in the set [0, 1] of membership degrees, i.e. a function $t : [0, \ 1]^2 \rightarrow [0, \ 1]$ which is commutative, associative, monotonically nondecreasing, and has 1 as a neutral and 0 as zero element: that means each such t-norm r fulfills for all u, υ, w ∈ [0, 1]:

(T1)     utυ = utυ and ut(υtw) = (utυ)tw,

(T2)     $u \leq υ \Rightarrow utw \leq utw$,

(T3)     1tu = t(1, u) = u and   0tu = t(0, u) = 0;

cf. also [75].

We write $\wedge_t$ to indicate, that t is the truth function which characterizes $\wedge_t$. Hence one always has

$$[\![ H_1 \wedge_t H_2 ]\!] =_{def} [\![ H_1 ]\!] \, t \, [\![ H_2 ]\!].\tag{7}$$

There is a special case, the so-called Lukasiewicz conjunction and characterized via (7) using the t-norm $t_L$

$$t_L (u, υ) =_{def} \max \{0, u + υ - 1\}.\tag{8}$$

For t-norms a partial ordering $\leqq$ is considered which is pointwise defined by t

$$t_1 \leqq t_2 =_{def} t_1 (u,υ) \leqq t_2 (u,υ) \text{ for all } u,υ \in [0, 1].\tag{9}$$

Among the t-norms the left-continuous ones are of special interest. With them by the definition

$$uφ_t υ =_{def} \sup \{w \,|\, utw \leq υ\} \text{ for all } u,υ \in [0, 1].\tag{10}$$

a φ-operator $φ_t$ is connected which is the truth function of a suitable implication connective $\rightarrow_t$ to be considered together with $\wedge_t$; cf. [28, 29]. The left continuous t-norms t also have another important property.

$$uφ_t υ =1 \Leftrightarrow u \leq υ \text{ for all } u,υ \in [0, 1]\tag{11}$$

which becomes crucial for some results mentioned later on.

For the t-norm $t_L$ the corresponding φ-operator $φ_L$ is the well-known Likasiewicz implication:

$$uφ_L υ = \min \{1, 1 – u + υ\}.\tag{12}$$

And for the t-norm $t = t_G = \min$ the corresponding φ-operator is the Godel implication

$$u\varphi_L \upsilon = \begin{cases} 1 & \textit{if } u \le \upsilon, \\ \upsilon & \textit{if } u > \upsilon. \end{cases} \qquad (13)$$

With those preliminaries the fuzzified inclusion is defined in essentially the same way as in classical set theory by

$$A \subseteq_t B =_{\text{def}} \forall x\, (x\, \varepsilon A \rightarrow_t x\, \varepsilon\, B), \qquad (14)$$

which means in more traditional notation

$$\llbracket A \subseteq_t B \rrbracket = \inf_{x \in X} \sup \{w \llbracket x\, \varepsilon\, A \rrbracket \, tw \le \llbracket x\, \varepsilon\, B \rrbracket \qquad (15)$$

and thus especially for $t = t_L$

$$\llbracket A \subseteq_t B \rrbracket = \inf_{x \in X}\ (\llbracket x\, \varepsilon\, A \rrbracket\, \varphi_L \llbracket x\, \varepsilon\, B \rrbracket$$

$$= \inf_{x \in X}\ \min \{1, \llbracket x \not\varepsilon A \rrbracket + \llbracket x\, \varepsilon\, B \rrbracket \}.$$

And the fuzzified identity is defined as

$$A \equiv_t B =_{\text{def}} A \subseteq_t B \wedge_t B \subseteq_t A. \qquad (16)$$

Introducing also the notation $\models$ for (generalized) logical validity as in [28, 29] by

$$\models A \equiv_t A. \qquad (17)$$

Whose proof partly uses property (11), of symmetry

$$\models A \equiv_t B \rightarrow_t B \equiv_t A, \qquad (18)$$

and of transitivity

$$\models A \equiv_t B \wedge_t B \equiv_t C \rightarrow_t A \equiv_t C, \qquad (19)$$

cf. [28, 29]. Additionally, again using (11), here in an essential way, one has

$$A = B \Leftrightarrow \models A \equiv_t B \Leftrightarrow \models A \equiv_t B. \qquad (20)$$

To see how these definitions work let us look at the t-norms $t_G = \min$ and $t_L$. For a readable formulation of the following results we use besides the support:

$$\text{Supp}(A) = \{x \in \aleph \mid \mu_A(x) > 0\} = \{x \in \aleph \mid [\![ x \, \varepsilon \, A ]\!] \neq 0\}$$

of $A \in F(\aleph)$ for any fuzzy sets $A, B \in F(\aleph)$ also the crisp sets

$$\{A > B =_{\text{def}} \{x \in \aleph \mid [\![ x \, \varepsilon \, A ]\!] > [\![ x \, \varepsilon \, B ]\!]\}, \qquad (21)$$

$$\{A \neq B\} =_{\text{def}} \{x \in \aleph \mid [\![ x \, \varepsilon \, A ]\!] \neq [\![ x \, \varepsilon \, B ]\!]\}. \qquad (22)$$

These sets generalize the support in the sense that

$$\text{Supp}(A) = \{A > \emptyset\} = \{A \neq \emptyset\}.$$

Holds true for each $A \in F(\aleph)$.

Straightforward calculations give for $t_G = \min$ the results with inf over empty set is 1 and sup over empty set is 0.

$$[\![ A \subseteq_{t_G} B ]\!] = \inf_{x \in \{A > B\}} [\![ x \, \varepsilon \, B ]\!] \qquad (23)$$

with the corollary

$$\text{supp}(A) \backslash \text{supp}(B) \neq \emptyset \Rightarrow [\![ A \subseteq_{t_G} B ]\!] = 0. \qquad (24)$$

Then immediately one also has by definition (16)

$$[\![ A \equiv_{t_G} B ]\!] = \inf_{x \in \{A \neq B\}} \min [\![ x \, \varepsilon \, A ]\!], [\![ x \, \varepsilon \, B ]\!]$$

$$= \inf_{x \in \{A \neq B\}} [\![ x \, \varepsilon \, A \cap B ]\!] \qquad (25)$$

now with the corollary

$$\text{supp}(A) \neq \text{supp}(B) \Rightarrow [\![ A \equiv_{t_G} B ]\!] = 0, \qquad (26)$$

which indicates that $\equiv_{t_G}$ is quite a strong fuzzified equality.

The other case $t = t_L$ again by elementary calculations first gives

$$\llbracket \, A \subseteq_{t_1} B \, \rrbracket = 1 - \sup_{x \in \{A > B\}} (\llbracket \, x \, \varepsilon \, A \, \rrbracket - \llbracket \, x \, \varepsilon \, B \, \rrbracket) \qquad (27)$$

and therefore with the auxiliary notation

$$\Delta \, (A, B) =_{\text{def}} \sup_{x \in \{A > B\}} (\llbracket \, x \, \varepsilon \, A \, \rrbracket - \llbracket \, x \, \varepsilon \, B \, \rrbracket)$$

quite directly

$$\llbracket \, A \equiv_t B \, \rrbracket = \max \{0, \, 1 - \Delta \, (A, B) + \Delta \, (B, A))\}. \quad (28)$$

This time, contrary to the results (24) and (26), both of the claims

$$\llbracket \, A \subseteq_{t_1} B \, \rrbracket \neq 0 \text{ and supp } (A) \nsubseteq \text{ supp } (B)$$

as well as both of the claims

$$\llbracket \, A \equiv_{t_1} B \, \rrbracket \neq 0 \text{ and supp } (A) \neq \text{ supp } (B),$$

can be true at once.

To some extent therefore $\subseteq_{t_L}$ and $\equiv_{t_L}$ better meet the intuition behind the fuzzified inclusion and equality than $\subseteq_{t_G}$ and $\equiv_{t_G}$, namely the intuition that "small" deviations from the "true", i.e. complete inclusion or equality do not completely falsify the generalized inclusion or equality. In addition, deviations from supp $(A) \subseteq$ supp $(B)$ and supp $(B)$ should surely count as "small" as long as combined with small differences in the membership degrees over the "critical" regions supp $(A) \setminus$ supp $(B)$ and supp$(A) \Delta$ supp $(B)$[2].

Finally, we need the notion of a metric in $F(\aleph)$. A dyadic function Q from $F(\aleph)$ into the nonnegative reals $R^+$ is a metric iff for all A, B, C $\in F(\aleph)$ the following conditions hold true:

(M1)    Q (A, B) = 0 iff A = B            (identity property),
(M2)    Q (A, B) = Q (B, A)               (symmetry),
(M3)    Q (A, C) + Q (C, B) ≥ Q (A, B)    (triangle inequality).

Sometimes the identity condition (M1) is weakened to the condition

$$(M1^p) \quad Q(A, A) = 0.$$

By a pseudo-metric Q then a function is meant that fulfills conditions $(M1^p)$, M2), (M3) .

In a general setting, in referring to a fuzzy (relational) equation one quite often has in mind an equation describing a relationship between fuzzy sets in two (possibly different) "space", i.e. universes of discourse. Such a form of relationship is supposed to be represented by a fuzzy relation between the elements of those "spaces" i.e., over the (crisp) Cartesian product of those universe of discourse. More precisely, let be given two fuzzy sets $A \in F(\aleph)$ and $B \in F(\wp)$ as well as a fuzzy relation R $\in F(\aleph \times \wp)$. Then a fuzzy relational equation can be written down in a general form as

$$\Theta (R, A) = B, \qquad (29)$$

where $\Theta$ is a suitable operator producing a fuzzy set B out of a fuzzy set A and a fuzzy relation R. And the case of a system.

$$\Theta (R, A_i) = B_i, i = 1,\ldots, n \qquad (30)$$

of fuzzy relational equations fits into those considerations as well.

Even more general, of course, is to consider $\Theta$ as an operator (of some finite parity) which maps fuzzy sets and relations onto fuzzy sets or fuzzy relations – and for which some of the arguments has (have) to be determined. We will not discuss the problem of fuzzy equations in this generality here.

The standard examples of fuzzy equations are fuzzy relational equations like

$$R" A = B, \text{ i.e. } A o_t R = B$$

for given fuzzy sets A, B or "pure" relational equations like

$$R o_t S = T$$

for given fuzzy relations S, T, both with an unknown fuzzy relation. But also again fuzzy set equations

$$R" A = A, \text{ i.e. } A \, o_t \, R = A$$

for an unknown fuzzy set A, cf. [84], or even fuzzy arithmetical equations for given fuzzy numbers A, B like

$$A \oplus X = B,$$
$$A \otimes X = B.$$

With an unknown fuzzy number X; cf. [28, 29].

Our strategy toward discussing the solvability behavior of fuzzy (relational) equations essentially is that one of a "many-valued translation" – we change, with respect to some t-norm t, from traditional equations like (29): $B = \Theta(R, A)$ which are to be solved for R or A to their many-valued counterparts $B \equiv_t \Theta(R, A)$ which – for lower semicontinuous t-norms t- have the property that always.

$$\llbracket \Theta(R, A) \equiv_t B \rrbracket = 1 \Leftrightarrow \Theta(R, A) = B.$$

What we reach to different levels of satisfaction are characterizations of the truth degrees $\llbracket \exists X(\Theta(R, X) \equiv_t B \rrbracket$ and $\llbracket \exists X(\Theta(X, A) \equiv_t B \rrbracket$ which do not involve the variable X, i.e. which are built up using only the "given data" R, B or A, B, respectively.

For system of fuzzy equations the situation is almost the same: only the truth degree to be determined now is e.g., of the form

$$\llbracket \exists X \prod_{i=1}^{n} (\Theta(X, A_i) \equiv_t B_i) \rrbracket,$$

i.e. has to be taken as the truth degree of the sentence. The system of fuzzy equations

$$\Theta(X, A_i) \equiv_t B_i, \ i = 1, \ldots, n$$

has a solution.

Instead of discussing directly the problem of solvability of fuzzy relational Eqs. (29) or of systems (30) of such equations we consider the truth degrees which we just mentioned as solvability degrees indicating the solvability behavior of our (systems of) fuzzy equations.

**DEFINITION 1.25.1:** *For each one of the fuzzy (relational Eqs. (29) and of the systems (30) of such equations, which are*

*supposed to be solved with respect to the fuzzy relation R, their solvability degree is the truth degree*

$$\xi_0 = {}_{def} \left[\!\left[ \exists X(\Theta(X, A_i) \equiv_t B_i) \right]\!\right],$$

*in the case of one equation, and it is in the case of system of equations the truth degree*

$$\xi_0 = {}_{def} \left[\!\left[ \exists X \prod_{i=1}^{n} (\Theta(X, A_i) \equiv_t B_i) \right]\!\right].$$

**COROLLARY 1.25.1:** *For all A, $A_i \in F(\aleph)$ and B, $B_i \in F(\wp)$ the solvability degree $\xi_0$ of equation (29) and the solvability degree $\xi$ of the system (30) of equations are*

$$\xi_0 = sup \{ \left[\!\left[ \Theta(R, A) \equiv_t B \right]\!\right] / R \in F(\aleph \times \wp) \},$$

$$\xi = sup \{ \underset{i=1}{\overset{n}{T}} \left[\!\left[ \Theta(R, A_i) \equiv_t B_i \right]\!\right] / R \in F(\aleph \times \wp) \}.$$

**Proposition 1.25.1:** If a fuzzy equation (29) or a system (30) of such equations has a solution then its solvability degree is = 1.

Given a continuous t-norm t and finite t-clan L, for all A, $A_i \in F_L(\aleph)$ and B, $B_i \in F_L(\wp)$ we consider also the relative solvability degrees

$$\xi_0^{(L)} = {}_{def} sup \{ \left[\!\left[ \Theta(R, A) \equiv_t B \right]\!\right] | R \in F_L(\aleph \times \wp) \},$$

$$\xi^{(L)} = {}_{def} sup \{ \underset{i=1}{\overset{n}{T}} \left[\!\left[ \Theta(R, A_i) \equiv_t B_i \right]\!\right] | R \in F_L(\aleph \times \wp) \},$$

of course, using bounded quantification and writing $\Re = F_L(\aleph \times \wp)$ one has

$$\xi_0^{(L)} = \left[\!\left[ \exists_\Re X(\Theta(X, A) \equiv_t B) \right]\!\right],$$

$$\xi^{(L)} = \left[\!\left[ \exists_\Re X \prod_{i=1}^{n} (\Theta(X, A_i) \equiv_t B_i) \right]\!\right].$$

**Proposition 1.25.2:** Suppose that L is finite t-clan with respect to a continuous t-norm t. Then one has

$$\xi_0^{(L)} \leq \xi_0 \text{ and } \xi^{(L)} \leq \xi$$

for all such (systems of) fuzzy equations with A, $A_i \in F_L(\aleph)$ and .
B, $B_i \in F_L(\wp)$.

***Proposition 1.25.3:*** Suppose that L is finite t-clan with respect to a continuous t-norm t and that $\xi_0^{(L)} = 1$ or $\xi^{(L)} = 1$ for some (system of fuzzy equation(s). Then this fuzzy equation or this system of fuzzy equations has a solution.

Given a left continuous t-norm t, a binary "distinguishability" function $Q_t$ is defined on $F(\aleph)$ by

$$Q_t(A, B) =_{def} 1 - [\![ A \equiv_t B ]\!] \qquad (*)$$

i.e. by always putting $Q_t(A, B) = [\![ \neg(A \equiv_t B) ]\!]$.

For fuzzy sets A, $B \in F(\aleph)$ in case of the t-norm $t_G = \min$ one gets, using the Godel implication (13), by simple calculations the corresponding "distinguishability" function $Q_G$ as

$$Q_G = \sup_{\substack{x \in \aleph \\ \|x \varepsilon A\| \neq \|x \varepsilon B\|}} (1 - [\![ x \varepsilon A \cap B ]\!])$$

$$= \sup_{\substack{x \in \aleph \\ \|x \varepsilon A\| \neq \|x \varepsilon B\|}} ([\![ x \varepsilon \overline{A \cap B} ]\!])$$

and in the case of the t-norm $t = t_L$ one gets, using the Lukasieqicz implication (12), after some elementary transformations

$$Q_L(A,B) = \min \{1, \max \{0, \sup_{x \in X} ([\![ x \varepsilon A ]\!] - [\![ x \varepsilon B ]\!] + \max \{0, \sup_{x \in X} ([\![ x \varepsilon B ]\!] - [\![ x \varepsilon A ]\!]\}\}.$$

This function $Q_L(A, B)$ is loosely related to the Cebysev distance of the membership functions $\mu_A$, $\mu_B$ defined as

$$d_C(\mu_A, \mu_B) = \sup_{x \in X} |\mu_A(x) - \mu_B(x)|.$$

In the sense that one always has

$$d_C(\mu_A, \mu_B) \leq Q_L(A, B) \leq 2. \, d_C(\mu_A, \mu_B)$$

and especially

$$(\models A \subseteq_{t_L} B) \Rightarrow Q_L (A, B) = d_C (\mu_A, \mu_B),$$

which means, using the original crisp implication relation $\subseteq$ for fuzzy sets as mentioned in (1),

$$A \subseteq B \Rightarrow Q_L (A, B) = d_C (\mu_A, \mu_B).$$

As mentioned in the introduction, the fuzzified identity $\equiv_t$ is considered as a kind of graded indistinguishability or similarity. Indeed, each graded relation $\equiv_t$ is reflexive in the sense of (17), it is symmetric in the sense of (18) and it is sup-t-transitive in the sense of property (19).

The main problem of the present section is to find a necessary and sufficient condition for t to yield via (*) a metric with properties (M1),…, (M3) as distinguishability function $Q_t$.

As simple consequences of the definitions (2) and (14), (15) first notice that for all $a \in \chi$ and $u, \upsilon \in [0, 1]$ one has

$$\llbracket \langle\!\langle a \rangle\!\rangle_u \subseteq_t \langle\!\langle a \rangle\!\rangle_v \rrbracket = u \, \varphi_t \, \upsilon = \sup \{ w \mid utw \leq \upsilon \}$$

and hence especially for $u = 1$:

$$\llbracket \langle\!\langle a \rangle\!\rangle_1 \subseteq_t \langle\!\langle a \rangle\!\rangle_v \rrbracket = \nu.$$

Now we can formulate and prove our characterization result.

**THEOREM 1.25.2:** *Suppose t is a left continuous t-norm. Then the function $Q_t$ of (36) is a metric in $F(\chi)$ iff $t \geqq t_L$, i.e. iff for all $u, \upsilon \in [0, 1]$:*

$$Max \{0, u + \upsilon - 1\} \leq ut\upsilon.$$

For more please refer [29].

## 1.26 Solvability criteria for systems of fuzzy relation equations

Dorte Neundorf and Rolf Bohm [72] have given a solvability criteria for a system of relational equations with two different composition methods. They have proved under certain conditions the system of relational equations is always solvable.

The motivation here is to obtain fuzzy relation models of real processes [84]. A set of rules is created, which describes the behavior of the process. The rules have the format "IF $A_i$ THEN $B_i$ ($i = 1,…, n$)", with $A_i$ and $B_i$ as fuzzy numbers. $A_i$ and $B_i$ are called the premise of a rule and the conclusion, respectively. Each rule can be transformed into a relation equation "$B_i = A_i$ o $R_i$" . $R_i$ is also fuzzy number, the relation. The relation performs a mapping form $A_i$ to $B_i$. It describes the "transformation" that is the meaning of the ith rule. The operator "o" is a symbol for any kind of evaluation that may be useful. It will be specified later.

The idea is to get a particular $R_i$ as the description or model of the subprocess described by the rule with index i. $R_i$ is a solution of the relation equation $B_i = A_i$ o $R_i$. The solution of a system of relation equations (and hence the model of the entire process) is a fuzzy number R that fulfills every relation equation of the given set. In some cases more than one solution can be determined.

It will be determined which parts of the membership functions $A_i$ and $B_i$ may be changed by the solution process without changing the resulting relation $R_i$. A combination of these conditions for all existing rules will result in intersection conditions for the membership functions that can still be tested easily. Different composition possibilities are included. The first step will be to look at these problems using a general expression for the $\underline{t}$-norm. A detailed analysis for a specific $\underline{t}$-norm, the minimum (u $\underline{t}$ υ = min {u, υ}), is added.

To ease the reading of the text, the term "relation equation" will be abbreviated by RE and the term "system of REs" by SRE, when it is helpful. The expression "supp A" means the support of function A. The support is the subset of the domain, where A does not vanish. The possibilities to interpret and solve a system of fuzzy relation equations vary with the method to calculate the operator "o". More detailed analysis is given in [84]. The necessary definitions are given here.

**DEFINITION 1.26.1** ($\underline{t}$-NORM)**:** *A mapping $\underline{t}$ is $\underline{t}$-norm, if the following properties hold: $\underline{t}$ :[0, 1] $\rightarrow$ [0, 1] with*

| | |
|---|---|
| *commutativity:* | $u \underline{t} v = v \underline{t} u$ |
| *associativity:* | $u \underline{t} (v \underline{t} w) = (u \underline{t} v) \underline{t} w$ |
| *monotonicity:* | $w \leq u \wedge z \leq u \Rightarrow w \underline{t} z \leq u \underline{t} v$ |
| *neural element:* | $u \underline{t} 1_t = u,$ |
| *zero element:* | $u \underline{t} 0_t = 0.$ |

**DEFINITION 1.26.2** (DUAL $\underline{t}$ -CONORM) : $\underline{s}_t$ *is the dual $\underline{t}$ -co-norm of t if u $\underline{s}$ v := 1 – (1 – u) $\underline{t}$ (1 – v).*

**DEFINITION 1.26.3** ($\underline{\rho}$ -OPERATOR): *The operator $\rho$ is defined as u $\underline{\rho}$ w: = sup {z ∈ [0, 1]/ u $\underline{t}$ z ≤ w}.*

**DEFINITION 1.26.4** (INTERSECTION AND UNION)**:** *Intersection and union are defined by the $\underline{t}$ -norm:*

$D: = A \cap_t B \Leftrightarrow \mu_D(x) := \mu_A(x) \underline{t} \mu_B(x) \, \forall x$

$C: = A \cap_t B \Leftrightarrow \mu_C(x) : = \mu_A(x) \underline{s}_t \mu_B(x) \, \forall x.$

The first step is to discuss a general system of relation equations that is composed with sup $\underline{s}_t$ -algorithm. Proposing the relation equation B = A ° R to describe a rule, the membership function of R can be written as

$$\mu_R (x, y) = \mu_A (x) \underline{\rho} \mu_B (y)$$
$$= \sup \{z \in [0, 1] \mid \mu_A (x) \underline{t} z \le \mu_B(y)\}$$

Refer table given below

Types and solutions for RE and SRE [2]

| Type of equation | Solution of single RE In case of solvability | Solution of SRE In case of solvability |
|---|---|---|
| Sup $\underline{t}$ -composition $B_i = A_i \circ_t R_i$ $\mu_B(y) = \sup \{\mu_A(x) \underline{t} \mu_R(x, y)\}$ | $R_i = A_i \rho B_i$ Greatest solution in the sense of inclusion | $R = \cap(A_i \rho B_i)$ $= \cap \sup R_i.$ |
| inf $\underline{t}$ -composition $B_i = A_i \rho R_i$ $\mu_B(y) = \inf \{\mu_A(x) \rho \mu_R(x, y)\}$ | $R_i = A_i \underline{t} B_i$ smallest solution in the sense of inclusion | $R = \cup (A_i \underline{t} B_i)$ $= \cup \inf R_i.$ |

On the other hand, starting with A and R, B can be calculated as

$$\mu_B(y) = \sup_x \mu_A(x) \ \underline{t} \ \mu_R(x, y)$$

$$= \max\left\{ \sup_{x \notin \text{supp } A} \mu_A(x) \, \underline{t} \, \mu_R(x, y), \ \sup_{x \in \text{supp } A} \mu_A(x) \, \underline{t} \, \mu_R(x, y) \right\}$$

$$= \max\left\{ 0, \ \sup_{x \in \text{supp } A} \mu_A(x) \, \underline{t} \, \mu_R(x, y) \right\}$$

$$= \sup_{x \in \text{supp } A} \mu_A(x) \, \underline{t} \, \mu_R(x, y) . \qquad\qquad (+)$$

Because of the definition of $\mu_R(x, y)$ in above given table, the latter is always smaller than or equal to $\mu_B(y)$,. If $\mu_A(x) < \mu_B(y)$, $\mu_R(x, y)$ is equal to one, so that

$$\mu_A(x) \ \underline{t} \ \mu_R(x, y) = \mu_A(x) < \mu_B(y).$$

It is obvious the supremum is an element of the set $\{(x, y)|\ \mu_A(x) \geq \mu_B(y)\}$. That is why the membership function $\mu_R(x, y)$ may be decreased in the set $\{(x, y)|\ \mu_A(x) < \mu_B(y)\}$ without any influence on $\mu_B(y)$.

Hence it is possible to identify parts of the support of R, where any changes of the membership function of R affect the result – in the sequel titled as fixed parts, in all other regions R may be decreased without any effect on the result. Figure 1 gives an example for the fixed parts of one rule with trapezoidal membership functions on a one-dimensional domain.

The application of a relation on a given input membership function is a supremum process over the support of A. Because of the monotonicity of the $\underline{t}$-norm it will be possible to decrease the relation even in the area $\{(x, y)|\ \mu_A(x) \geq \mu_B(y)\}$ without modifying the resulting B, if the maximum values of this set are not influenced.

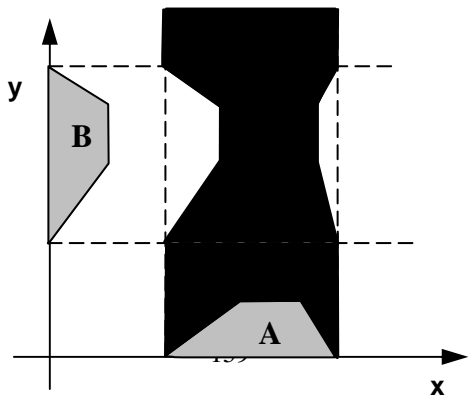

Figure 1.26.1

The supremum process in (+) is applied on a line, which is parallel to the x-axis. Its elongation meets the fixed y and its end are fixed by the border of the fixed parts of the relation equation (the grey area for the example givcen in Figure 1.26.1). On this line the values of the relation may only be changed in a way that the supremum of the $\underline{t}$-norms of these values is not touched.

In case of continuous premises with finite support, some additional information can be given:

$$\forall\, c \in [0, 1]: \exists\, x_1, x_2: [\mu_A(x_1) = \mu_A(x_2) = c]$$
$$\wedge\, [c = \max_x \mu_A(x) \vee x_1 \neq x_2].$$

Figure 1.26.2 illustrates these ideas for a triangular premise: the maximum $c_2$ is only reached once. On the contrary, $c_1$ is not a maximum; because of the continuity of the premise membership function it will be reached once "on the way up" and "once on the way down".

For $x_1$ and $x_2$ $A(x_1)$ is equal to $A(x_2)$.y is constant, so the relation R is identical for $x_1$ and $x_2$. To reproduce B correctly we will only need one of these two points if the other value is not allowed to increase. In other words, the relation must remain unchanged at least on an area that covers the codomain of A once.

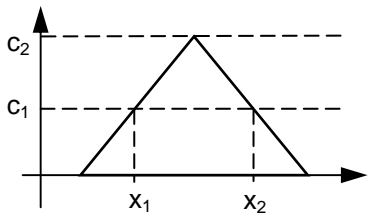

Figure 1.26.2

If the membership function is convex on a connected support (this holds e.g. for trapezoidal premises), the support can be divided in an increasing and decreasing part. Then it is sufficient to keep the relation on one of these parts. Then it is sufficient to keep the relation on one of these parts. All other parts of the membership function may be decreased at all points. Of course, there is more than one possibility to choose the fixed parts.

For a SRE the solutions of the REs are composed by a minimum process. An application of this composition can only result in decreasing, which has been shown not to change the result, if it does not occur on fixed parts.

**THEOREM 1.26.5** (SOLVABILITY CRITERIA FOR GENERAL SRE USING SUP-$\underline{t}$-COMPOSITION ): *A system of relation equations*

$$B_i = A_i \mathbin{\overset{\circ}{{}_t}} R_i,\ i = 1,\ ...,\ n$$

*evaluated by using the sup-$\underline{t}$- composition, is solvable if for every $A_i$ there is a subset $S \subset supp\ A_i$ with $A_i(S) = codomain\ (A_i)$ and $A_j(s) = 0$ for all $s \in S$ and $j = 1,\ ...,\ n$ and $j \neq i$.*

**THEOREM 1.26.6** (SOLVABILITY CRITERIA FOR GENERAL SRE USING INF-$\rho$-COMPOSITION): *A system of relation equations*

$$B_i = A_i \mathbin{\overset{\circ}{{}_t}} R_i,\ i = 1,\ ...,\ n$$

*evaluated by using the inf -$\rho$- composition, is solvable if for every $A_i$ there is a point $s \in supp\ A_i$ with $A_j(s) = 0$ for all $s \in S$ and $j = 1,\ ...,\ n$ and $j \neq i$.*

**THEOREM 1.26.7** (SOLVABILITY CRITERIA FOR SRE USING SUP-$\underline{t}$-COMPOSITION AND $\underline{t} =$ MIN). *A system of relation equations*

$$B_i = A_i \mathbin{\overset{\circ}{{}_t}} R_i,\ i = 1,\ ...,\ n$$

*evaluated by using the sup-$\underline{t}$-composition with $\underline{t} =$ min, is solvable if for every $A_i$ there is a subset $S \subset supp\ A_i$ with $A_j(s) = 0$ for all $s \in S$ and $j = 1,\ ...,\ n$ and $j \neq i$.*

**THEOREM 1.26.8** (SOLVABILITY CRITERIA FOR GENERAL SRE USING INF-$\rho$-COMPOSITION AND $\underline{t} =$ MIN): *A system of relation equations*

$$B_i = A_i \mathbin{\overset{\circ}{{}_t}} R_i,\ i = 1,\ ...,\ n$$

*evaluated by using the inf -$\rho$- composition, with $t =$ min is always solvable.*

### 1.27 Infinite FRE to a complete Brouwerian lattice

Wang [109] have proved and obtained a method to construct a minimal solution from any given solution with in finite steps.

Let $I$ and $J$ be index sets, and let $A = (a_{ij})_{I \times J}$ be a coefficient matrix, $B = (b_i)_{i \in I}$ be a constant column vector. Then an equation.

$$A \odot X = B \qquad (1)$$

or

$$\underset{j \in J}{V} (a_{ij} \wedge x_j) = b_i$$

for all $i \in I$

is called a fuzzy relational equation assigned on a complete Brouwerian lattice L, where $\odot$ denote the sup-inf composite operation, and all $x_j$, $b_i$, $a_{ij}$'s are in L. An X which satisfies (1) is called a solution of (1), the solution set of (1) is denoted by $\aleph = \{X : A \odot X = B\}$. A special case of (1) is as follows:

$$A \odot X = b \qquad (2)$$

or

$$\underset{j \in J}{V} (a_{ij} \wedge x_j) = b_i \, ,$$

where $b \in L$, $A = (a_i)_{i \in J}$ is a row vector. The solution set of (2) is denoted by $\aleph = \{X : A \odot X = b\}$.

**DEFINITION 1.27.1 [5]:** *In a distributive lattice L, if p is join-irreducible if $b \vee c = a$ implies $b = a$ or $c = a$.*

**PROPOSITION 1.27.1 [5]:** *In a distributive lattice L, if p is join-irreducible, then*

$$p \leq V_{i=1}^{k} x_i$$

*implies $p \leq x_i$ for some i.*

**DEFINITION 1.27.2:** *For an element a of a lattice L, if there are join-irreducible elements $p_1$, $p_2$,..., $p_n$ such that a $V_{i=1}^{n} p_i$, we say that a has a finite join-decomposition.*

Further, if for any $j \in \{1, 2,...,n\}$, we have moreover $a \neq V_{i=1, i \neq j}^n p_i$, then the decomposition is called *irredundant*, and we say that $a$ has an *irredundant finite join-decomposition*.

**DEFINITION 1.27.3 [5]:** *A Brouwerian lattice is a lattice L in which, for any given elements a and b, the set of all $x \in L$ such that $a \wedge x \leq b$ contains a greatest elements, denoted by $a \alpha b$, the relative pseudo-complement of a in b.*

***Remark 1.27.1 [84].*** If $L = [0, 1]$, then it is easy to see that for any given $a, b \in L$,

$$a \alpha b = \begin{cases} 1 \ a \leq b, \\ b \ a > b. \end{cases}$$

**PROPOSITION 1.27.2 [5]:** Any Brouwerian lattice L is distributive

**DEFINITION 1.27.4 [5]:** *Let $(P, \leq)$ be a partially ordered set and $X \subseteq P$. A minimal element of X is an element $p \in X$ such that there exists no $x \in X$ with $x < p$. The greatest element of X is an element $g \in X$ such that $x \leq g$ for all $x \in X$.*

**DEFINITION 1.27.5 (SANCHEZ [11]):** *Let $A = (A_{ij})_{I \times J}$ and $B = (b_{ij})_{I \times J}$ be two matrices. Then the partial order $\leq$, the join $\vee$, and the meet $\wedge$ are defined as follows:*

$$A \leq B \text{ if and only if } a_{ij} \leq b_{ij} \text{ for all } i \in I, j \in J,$$
$$A \vee B = (a_{ij} \vee b_{ij})_{I \times J}, A \wedge B = (a_{ij} \wedge b_{ij})_{I \times J}.$$

**PROPOSITION 1.27.3:** For each $i \in I$, let $A_i = (a_{ij})_{j \in J}$ be a row vector, and $\aleph_2$ be the solution set of $b_i = A_i \odot X$, then:

(a) $\aleph_1 \neq \phi$ if and only if $\bigcap_{i \in I} \aleph_{i2} \neq \phi$. Further $\aleph_1 \bigcap_{i \in I} \aleph_2$

(b) If $\aleph_1 \neq \phi$, then $X^* = A^T \alpha B$ is the greatest solution of (1), where $A^T$ is the transpose of A, $A^T \alpha B = \left( \wedge_{i \in I} (a_{ij} \alpha b_i) \right)_{j \in J}$ is a column vector.

**DEFINITION 1.27.6:** *An element a in a complete lattice L is called compact if whenever $a \leq \vee S$ there exists a finite subset $T \subseteq S$ with $a \leq \vee T$.*

**Lemma 1.27. 1:** Let J be a finite index set. If $\aleph_2 \neq \phi$, and b has an irredundant finite join-decomposition, then for each $X \in \aleph_2$, there exists a minimal element $X_*$ of $\aleph_2$ such that $X_* \leq X$.

**PROPOSITION 1.27.4:** If $\aleph_2 \neq \phi$, and $X_* = (x_{j*})_{j \in J}$ is a minimal element of $\aleph_2$, then $b = V_{j \in J} \, x_j^*$.

**THEOREM 1.27.7:** *If $\aleph_2 \neq \phi$, then for each $X \in \aleph_2$, there exists a minimal element $X_*$ of $\aleph_2$ such that $X_* \leq X$ if and only if there is a subset B of L with B satisfying:*

*(i)*      *$\vee B = b$;*
*(ii)*     *For each $p \in B$, if $p \neq 0$, then $b \neq V (B \setminus \{p\})$;*
*(iii)*    *For each $X = (x_j)_{j \in J} \in \aleph_2$ and each $p \in B$ there is $j \in J$ such that $p \leq a_j \wedge x_j$*

**THEOREM 1.27.8:** *If $\aleph_1 \neq \phi$, I is a finite index set, and every component $b_i$, $i \in I$, of B is a compact element and for each $b_i$, $i \in I$, there exists a subset $B_i$ of L such that*

*i.*      *$VB_i = b_i$;*
*ii.*     *For each $p_{it} \in B_i$, if $p_{it} \neq 0$, then $b_i \neq V (B_i \setminus \{p_{it}\})$;*
*iii.*    *For each $X = (x_j)_{j \in J} \in \aleph_{i2}$ and each $p_{it} \in B_i$ there is $j \in J$ such $p_{it} \leq a_{ij} \wedge x_j$*
*iv.*    *For each $p \in \bigcup_{i \in I} B_i$, if $p \neq 0$, then there is no subset Q of $\bigcup_{i \in I} B_i$ such that $p \leq V (Q \setminus \{p\})$.*

*Then for each $X \in \aleph_1$, there exists a minimal element $X_*$ of $\aleph_1$ such that $X_* \leq X$.*

For proof please refer [109].

## 1.28 Semantics of implication operators and fuzzy relational product

After a brief discussion of the need for fuzzy relation theory in practical systems work, here we explain the new triangle products of relations and the sort of results to be expected from them, starting from a crisp situation. The asymmetry of these products, in contrast to correlation, is noted as essential to the investigation of hierarchical dependencies. The panoply of multi-valued implication operators, with which the fuzzification of these products can be accomplished, is presented, and a few of their properties noted.

For checklist paradigm please refer [43]. Using a well-known psychological test in an actual situation, so that the finer structure is in fact available, a comparison is made between a checklist measure and several of the operator values, showing the interrelationship concretely. Finally, some products and their interpretations are presented, using further real-world data.

The difficulties of saying anything meaningful about a system increase enormously with its complexity. The vogue for, and success of, statistical methods are evidence of one way of doing this. Here we are concerned with quite another, the possibilistic [43], rather than the probabilistic way.

In any real-world situation our information about a system is too voluminous and intricate, and needs to be summarized; or it is approximate from the very beginning. A scientist, attempting to analyze such a system, implicitly asserts his belief that a number of significant things can be said about the system – could they only be found! In his attempt to analyze a real-world system, he is working with a model of its, simplified so as to be manageable and comprehensible. The danger of the assumption that this model can always be deterministic has been demonstrated in [43].

In general, it can be said that unwarranted structural assumptions imposed on the working model can lead to dangerous artifacts that do not reflect anything that is contained in the real-world data; this leads consequently to totally meaningless result of the analysis masquerading as "scientific truth".

On the other hand, rejecting such strong unwarranted assumptions, we may still be able to provide some meaningful answers to our questions such as: What structural relationships between the individual items of the analyzed data must exist? Which cannot exist? Which may exist perhaps if …? These modal terms in which we all think, but which we usually rule out in our "scientific discourse", are in fact the proper terms for possibilistic systems.

Possibility theory can be crisp; any given structure, say, may be utterly (1) or not at all (0) contained in another structure. More attractive and more consonant with summarized data from the real world, however, is fuzzy possibilistic theory: here the degree to which X can be contained in Y is (estimated as) some number from 0 to 1 inclusive. This may sound like a probability, but it is not. The quickest way to see this is from the fact that entirely different operations are performed on these fuzzy degrees than are performed on probabilities; this reflects, of course, a deeper semantic and epistemological difference, on which there is a large literature, of which Zadeh [115, 116] are particularly illuminating.

Chapter Two

# SOME APPLICATIONS OF FUZZY RELATIONAL EQUATIONS

In this chapter we give several of the applications of fuzzy relational equations in studies like chemical engineering, transportation, medicine etc. The fuzzy relational equations happen to be a seemingly simple method but in reality it can be used to solve many complicated problems, problems that even do not have solutions by using linear equations. This chapter is completely devoted to the applications of fuzzy relational equations. There are 11 sections in this chapter, which gives how FREs are applied to special problem. [81, 85, 100-101,111, 114]

## 2.1 Use of FRE in chemical engineering

The use of fuzzy relational equations (FRE) for the first time has been used in the study of flow rates in chemical plants. They have only used the concept of linear algebraic equations to study this problem and have shown that use of linear equations does not always guarantee them with solutions. Thus we are not only justified in using fuzzy relational equation but we are happy to state by adaptation of FRE we are guaranteed of solutions to the problem. We have adapted the fuzzy relational equations to the problem of estimation of flow rates in a chemical plant, flow rates in a pipe network and use of FRE in a 3 stage counter current exaction unit [99].

Experimental study of chemical plants is time consuming expensive and need intensive labor, researchers and engineers prefer only theoretical approach, which is inexpensive and effective. Only linear equations have been used to study: (1). A typical chemical plant having several inter-linked units (2). Flow distribution in a pipe network and (3). A three stage counter current extraction unit. Here, we tackle these problems in 2 stages.

At the first stage we use FRE to obtain a solution. This is done by the method of partitioning the matrix as rows. If no solution exists by this method we as the second stage adopt Fuzzy Neural Networks by giving weightages. We by varying the weights arrive at a solution which is very close to the predicted value or the difference between the estimated value and the predicted value is zero. Thus by using fuzzy approach we see that we are guaranteed of a solution which is close to the predicted value, unlike the linear algebraic equation in which we may get a solution and even granted we get a solution it may or may not concur with the predicted value.

To attain both solution and accuracy we tackle the problems using Fuzzy relational equations at the first stage and if no solution is possible by this method we adopt neural networks at the second stage and arrive at a solution.

Consider the binary relation $P(X, Y)$, $Q(Y, Z)$ and $R(X, Z)$ which are defined on the sets $X = \{x_i / i \in I\}$ $Y = \{y_i / j \in J\}$ and $Z\{z_k / k \in K\}$ where we assume that $I = N_n$, $J = N_r$ and $K = N_s$. Let the membership matrices of P, Q and R be denoted by $P = [p_{ij}]$, $Q = [q_{ik}]$ and $R = [r_{ik}]$ respectively, where $p_{ij} = P(x_i, y_j)$, $q_{ik} = Q(y_j, z_k)$ and $r_{ik} = R(x_i, z_k)$ for $i \in I (= N_n)$, $j \in J (= N_m)$ and $k \in K (= N_s)$. Entries in P, Q and R are taken from the interval [0, 1]. The three matrices constrain each other by the equation

$$P \circ Q = R \qquad (1)$$

(where $\circ$ denotes the max-min composition) known as the fuzzy relation equation (FRE) which represents the set of equation

$$\text{Max } p_{ij}q_{jk} = r_{ik} \qquad (2)$$

for all $i \in N_n$, $k \in N_s$. If after partitioning the matrix and solving the equation (1) yields maximum of $q_{jk} < r_{ik}$ for some $q_{jk}$, then this set of equation has no solution. So at this stage to solve the equation 2, we use feed-forward neural networks of one layer with n-neurons with m inputs shown in Figure 2.1.1.

Inputs of the neuron are associated with real numbers $W_{ij}$ referred as weights. The linear activation function f is defined by

$$f(a) = \begin{cases} 0 & \text{if } a < 0 \\ a & \text{if } a \in [0, 1] \\ 1 & \text{if } a > 1 \end{cases}$$

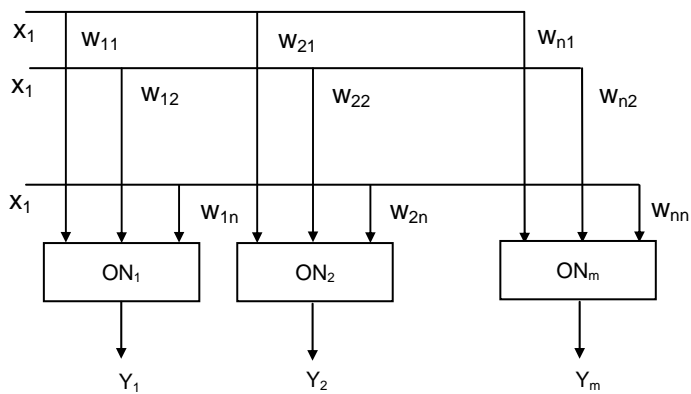

Figure: 2.1.1

The output $y_i = f(\max W_{ij}x_j)$, for all $i \in N_n$ and $j \in N_m$. Solution to (1) is obtained by varying the weights $W_{ij}$ so that the difference between the predicted value and the calculated value is zero.

### *FRE to estimate flow rates in a chemical plants*

A typical chemical plant consists of several interlinked units. These units act as nodes. The flowsheet is given in Figure 2.1.2.

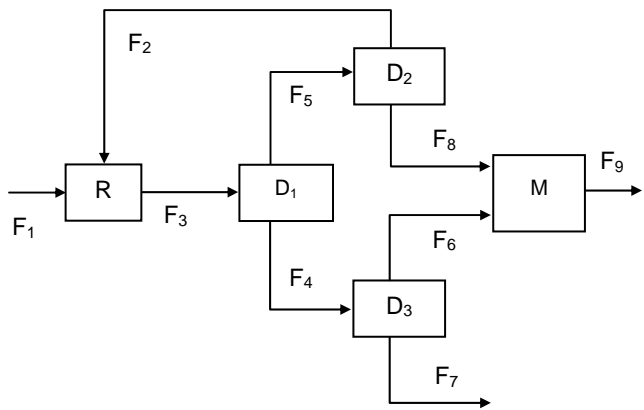

Figure: 2.1.2

An experimental approach would involve measuring the nine flow-rates to describe the state of the plant which would involve more money and labor. While studying this problem in practice researchers have has neglected density variations across each stream. The mass balance equations across each node at steady state can be written as

$$F_3 - F_2 = F_1,$$
$$F_2 - F_4 = F_5,$$
$$F_4 - F_7 = F_6,$$
$$F_2 + F_8 = F_5,$$
$$F_8 = F_9 - F_6. \qquad (3)$$

Here $F_i$ represents the volumetric flow rate of the $i^{th}$ stream. In equation (3) at least four variables have to be specified or determined experimentally. The remaining five can then be estimated from the equation (3), which is generated by applying the principle of conservation of mass to each unit. We assume $F_1$, $F_5$, $F_6$ and $F_9$ are experimentally measured, equation (3) reads with known values on the right-hand side as follows:

$$\begin{bmatrix} -1 & 1 & 0 & 0 & 0 \\ 0 & 1 & -1 & 0 & 0 \\ 0 & 0 & 1 & -1 & 0 \\ 1 & 0 & 0 & 0 & 1 \\ 0 & 0 & 0 & 0 & 1 \end{bmatrix} \begin{bmatrix} F_2 \\ F_3 \\ F_4 \\ F_7 \\ F_8 \end{bmatrix} = \begin{bmatrix} F_1 \\ F_5 \\ F_6 \\ F_5 \\ F_9 - F_6 \end{bmatrix} \qquad (4)$$

$$P \circ Q = R \qquad (5)$$

where, P, Q and R are explained. Using principle of conservation of mass balanced equation we estimate the flow rates of the five liquid stream. We in this problem aim to minimize the errors between the measured and the predicted value. We do this by giving suitable membership grades $p_{ij} \in [0, 1]$ and estimate the flow rates by using these $p_{ij}$'s in the equation 3. Now the equation 4 reads as follows:

$$\begin{bmatrix} p_{11} & p_{12} & 0 & 0 & 0 \\ 0 & p_{22} & p_{23} & 0 & 0 \\ 0 & 0 & p_{33} & p_{34} & 0 \\ p_{41} & 0 & 0 & 0 & p_{45} \\ 0 & 0 & 0 & 0 & p_{55} \end{bmatrix} \begin{bmatrix} F_2 \\ F_3 \\ F_4 \\ F_7 \\ F_8 \end{bmatrix} = \begin{bmatrix} F_1 \\ F_5 \\ F_6 \\ F_5 \\ F_9 - F_6 \end{bmatrix} \qquad (6)$$

where     $P = (p_{ij})$,
          $Q = (q_{ik}) = [F_2 \ F_3 \ F_4 \ F_7 \ F_8]^t$ and
          $R = (r_{ik}) = [F_1 \ F_5 \ F_6 \ F_5 \ F_9 - F_6]^t$.

We now apply the partitioning method of solution to equation (6). The partitioning of P correspondingly partitions R, which is give by a set of give subsets as follows:

$$[p_{11} \ p_{12} \ 0 \ 0 \ 0] \begin{bmatrix} F_2 \\ F_3 \\ F_4 \\ F_7 \\ F_8 \end{bmatrix} = \begin{bmatrix} F_1 \\ F_5 \\ F_9 \\ F_5 \\ F_9 - F_6 \end{bmatrix}, \dots$$

$$[0 \ 0 \ 0 \ 0 \ p_{55}] \begin{bmatrix} F_2 \\ F_3 \\ F_4 \\ F_5 \\ F_8 \end{bmatrix} = \begin{bmatrix} F_1 \\ F_5 \\ F_6 \\ F_5 \\ F_9 - F_6 \end{bmatrix}.$$

Suppose the subsets satisfies the condition max $q_{ik} < r_{ik}$ then it has no solution. If it does not satisfy, this condition, then it has a final solution. If we have no solution we proceed to the second stage of solving the problem using Fuzzy Neural Networks.

**NEURAL NETWORKS TO ESTIMATE THE FLOW RATES IN A CHEMICAL PLANT**

When the FRE has no solution by the partition method, we solve these FRE using neural networks. This is done by giving weightages of zero elements as 0 and the modified FRE now reads as

$$P_1 \circ \begin{bmatrix} F_2 \\ F_3 \\ F_4 \\ F_7 \\ F_8 \end{bmatrix} = \begin{bmatrix} F_1 \\ F_5 \\ F_6 \\ F_5 \\ F_9 - F_6 \end{bmatrix}.$$

The linear activation function f defined earlier gives the output $y_i$ = f (max $W_{ij} x_j$) (i $\in N_n$) we calculate max $W_{ij}x_j$ as follows:

1. $W_{11}x_1 = 0.02F_2$, $W_{12}x_2 = 0F_2$, $W_{13}x_3 = 0F_2$ $W_{14}x_4 = 0.045F_2$, $W_{15}x_5 = 0F_2$
        $y_1 = f (max_{j \in Nm} W_{ij}x_j) = f (0.02F_2, 0F_2, 0.045F_2, 0F_2)$

2. $W_{21}x_1 = 0.04F_3$, $W_{22}x_2 = 0.045F_3$, $W_{23}x_3 = 0F_3$, $W_{24}x_4 = 0.0F_3$, $W_{15}x_5 = 0F_3$
        $y_2 = f (max_{j \in Nm} W_{ij}x_j) = f (0.04F_3, 0.045F_3, 0F_3, 0.0F_3, 0F_3)$

3. $W_{31}x_1 = 0.0F_4$, $W_{32}x_2 = 0.085F_4$, $W_{33}x_3 = 0.15F_4$, $W_{34}x_4 = 0.0F_4$ $W_{35}x_5 = 0F_4$
        $y_3 = f (max_{j \in Nm} W_{ij}x_j) = f (0F_4, 0.085F_4, 0.15F_4, 0F_4, 0F_4)$

4. $W_{41}x_1 = 0.0F_7$, $W_{42}x_2 = 0F_7$, $W_{43}x_3 = 0.2F_7$, $W_{44}x_4 = 0.0F_7$, $W_{45}x_5 = 0F_7$
        $y_4 = f (max_{j \in Nm} W_{ij}x_j) = f (0F_7, 0F_7, 0.2F_7, 0.0F_7, 0F_7)$

5. $W_{51}x_1 = 0.0F_8$, $W_{52}x_2 = 0F_8$, $W_{53}x_3 = 0F_8$, $W_{54}x_4 = 0.45F_8$, $W_{55}x_5 = 0.5F_8$
        $y_5 = f (max_{j \in Nm} W_{ij}x_j) = f (0F_8, 0F_8, 0F_8, 0.45F_8, 0.5F_8)$

shown in Figure 2.1.2. Suppose the error does not reach 0 we change the weights till the error reaches 0. We continue the process again and again until the error reaches to zero.

Thus to reach the value zero we may have to go on giving different weightages (finite number of time) till say $s^{th}$ stage $P_s \circ Q_s$ whose linear activation function f, makes the predicted value to be equal to the calculated value. Thus by this method, we are guaranteed of a solution which coincides with the predicted value.

FUZZY NEURAL NETWORKS TO ESTIMATE VELOCITY OF FLOW DISTRIBUTION IN A PIPE NETWORK

In flow distribution in a pipe network of a chemical plant, we consider liquid entering into a pipe of length T and diameter D at a fixed pressure $P_i$, The flow distributes itself into two pipes each of length $T_1(T_2)$ and diameter $D_1(D_2)$ given in Figure 2.1.3.

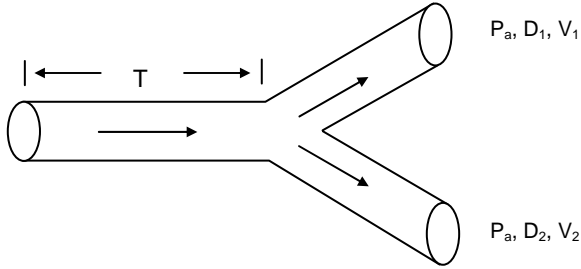

Figure: 2.1.3

The linear equation is based on Ohm's law, the drop in voltage V across a resistor R is given by the linear relation $V = iR$ (Ohm's law). The hydrodynamic analogue to the mean velocity v for laminar flow in a pipe is given by $\nabla_p = v\ (32\mu T/D^2)$. This is classical-Poiseulle equation. In flow distribution in a pipe network, neglecting pressure losses at the junction and assuming the flow is laminar in each pipe, the macroscopic momentum balance and the mass balance at the junction yields,

$$P_1 - P_a = (32\mu T/D^2)v + (32\mu T_1 D_1^2)v_1,$$
$$P_i - P_a = (32\mu T/D^2)v + 32\mu T/D_2^2)v_2,$$
$$D^2v = D_1^2v_1 + v_2D_2^2\ . \qquad (7)$$

Hence $P_a$ is the pressure at which the fluid leaves the system at the two outlets. The set of three equation in (7) can be solved and we estimate $v, v_1, v_2$ for a fixed $(P_i - P_a)$. The system reads as

$$\begin{bmatrix} 32\mu T/D^2 & 32\mu T_1/D_1^2 & 0 \\ 32\mu T/D^2 & 0 & 32\mu T_2/D_2^2 \\ -D^2 & D_2^1 & D_2^2 \end{bmatrix} \begin{bmatrix} v \\ v_1 \\ v_2 \end{bmatrix} = \begin{bmatrix} p_i - p_a \\ p_i - p_a \\ 0 \end{bmatrix}.$$

We transform this equation into a fuzzy relation equation. We use a similar procedure described earlier and obtain the result by fuzzy relation equation. We get max $(0.2v, 0.025v, 0.03v)$, max $(0.035v, 0v_1, 0.04v_1)$, max $(0v_2, 0.04v_2, 0.045v_2)$ by using neural networks for fuzzy relation equation described in section 3. Suppose the error does not reach to 0, we change the weights till the error reaches 0. We continue the process again and again till the error reaches zero.

## FUZZY NEURAL NETWORKS TO ESTIMATE THREE STAGE COUNTER CURRENT EXTRACTION UNIT

*Three-stage counter extraction unit is shown in Figure 2.1.4. The components A present in phase E (extract) along with a nondiffusing substance as being mixture.*

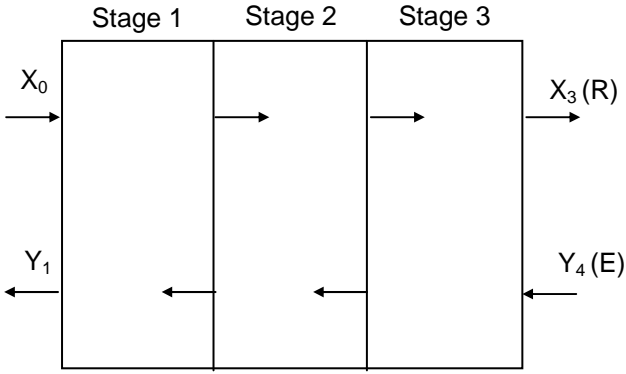

Figure: 2.1.4

*It is extracted into R by a nondiffusing solvent. The 3 extraction stage is given by the three equation.*

$$E_s Y_4 + R_s X_s = R_s X_3 + E_s Y_3,$$
$$E s Y_3 + R s X_1 = E_s + R_s X_2,$$
$$E s Y_2 + R_s X_0 = E_s Y_1 + R_s X_1 \qquad (8)$$

$Y_i(X_i)$ = moles of A, The flow of each stage is denoted by $E_s(R_s)$ and this constant does not vary between the different stages. The assumption of a linear equilibrium relationship for the compositions leaving the $i^{th}$ stage equations

$$Y_i = K X_i \qquad (9)$$

for i = 1, 2, 3 reads as

$$
\begin{bmatrix}
R_s & E_s & 0 & -E_s & 0 & 0 \\
K & -1 & 0 & 0 & 0 & 0 \\
-R_s & 0 & R_s & E_s & 0 & -E_s \\
0 & 0 & K & -1 & 0 & 0 \\
0 & 0 & -R_2 & 0 & R_s & E_s \\
0 & 0 & 0 & 0 & K & -1
\end{bmatrix}
\begin{bmatrix}
X_1 \\ Y_1 \\ X_2 \\ Y_2 \\ X_3 \\ Y_3
\end{bmatrix}
=
\begin{bmatrix}
R_s X_0 \\ 0 \\ 0 \\ 0 \\ E_s Y_4 \\ 0
\end{bmatrix}
$$

where $\{X_1, Y_1, X_2, Y_2, X_3, Y_3\}$ can be obtained for a given $E_s$, $R_s$ and $K$. Since use of linear algebraic equation does not result in the closeness of the measured and predicted value, we use neural networks for fuzzy relation equations to estimate the flow-rates of the stream, moles of the three-stage counter extraction unit and velocity of the flow distribution in a pipe network. As neural networks is a method to reduce the errors between the measured value and the predicted value. This allows varying degrees of set membership (weightages) based on a membership function defined over the range of value. The (weightages) membership function usually varies from 0 to 1. We use the similar modified procedure described earlier and get result by fuzzy relation equation. We get max $(0.2X_1, 0.25X_1, 0.3X_1, 0X_1, 0X_1, 0X_1)$, max $(0.35Y_1, 0.4Y_1, 0Y_1, 0Y_1, 0Y_1, 0Y_1)$ max $(0X_2, 0X_2, 0.45X_2, 0.5X_2, 0.55X_2, 0X_2)$, max $(0.6Y_2, 0Y_2, 0.65Y_5, 0.7Y_2, 0Y_2, 0Y_2)$ max $(0X_3, 0X_3, 0X_3, 0X_3, 0.75X_3, 0.8X_3)$, max $(0Y_3, 0Y_3, 0.85Y_3, 0Y_3, 0.9Y_3, 0.95Y_3)$ by neural networks for fuzzy relation equation. We continue this process until the error reaches zero or very near to zero.

*Thus we see that when we replace algebraic linear equations by fuzzy methods to the problems described we are not only guaranteed of a solution, but our solution is very close to the predicted value.*

## 2.2  New FRE to estimate the peak hours of the day for transport system

In this section we just recall the notion of new fuzzy relational equations and study the estimation of the peak hour problem for transport systems using it we have also compared our results with the paper of W.B.Vasantha Kandasamy and V. Indra where FCMs are used. We establish in our study which yields results, which are non over lapping and unique solutions.

MODIFIED FUZZY RELATION EQUATION HAS BEEN DERIVED FOR ANALYZING PASSENGER PREFERENCE FOR A PARTICULAR HOUR IN A DAY

Since any transport or any private concern which plys the buses may not in general have only one peak hour in day, for; the peak hours are ones where there is the maximum number of passengers traveling in that hour. The passengers can be broadly classified as college students, school going children, office going people, vendors etc. Each category will choose a different hour according to their own convenience. For example the vendor group may go for buying good in the early morning hours and the school going children may prefer to travel from 7.00 a.m. to 8 a.m., college students may prefer to travel from 8 a.m. to 9 a.m. and the office going people may prefer to travel from 9.00 a.m. to 10.00 a.m. and the returning hours to correspond to the peak hours as the school going children may return home at about 3.00 p.m. to 4.00 p.m., college students may return home at about 2.00 p.m. to 3.30 p.m. and the office going people may return home at about 5.00 p.m. to 6.00 p.m. Thus the peak hours of a day cannot be achieved by solving the one equation $P \circ Q = R$. So we reformulate this fuzzy relation equation in what follows by partitioning $Q_i$'s. This in turn partition the number of preferences depending on the set Q which correspondingly partitions R also. Thus the fuzzy relation equation say $P \circ Q = R$ reduces to a set of fuzzy relations equations $P_1 \circ Q_1 = R_1$, $P_2 \circ Q_2 = R_2$, …, $P_s \circ Q_s = R_s$ where $Q = Q_1 \cup Q_2 \cup … \cup Q_s$ such that $Q_i \cap Q_j = \phi$ for $i \neq j$. Hence by our method we get s preferences. This is important for we need at least 4 to 5 peak hours of a day. Here we give a new method by which we adopt the feed forward neural network to the transportation problem.

We briefly describe the modified or the new fuzzy relation equation used here. We know the fuzzy relation equation can be represented by neural network. We restrict our study to the form

$$P \circ Q = R \qquad\qquad (1)$$

where $\circ$ is the max-product composition; where $P = [p_{ij}]$, $Q = [q_{jk}]$ and $R = [r_{ik}]$, with $i \in N_n$, $j \in N_m$ and $k \in K_s$. We want to determine P. Equation (1) represents the set of equations.

$$\max_{j \in J_m} p_{ij}\, q_{jk} = \tau_{ik} \qquad\qquad (2)$$

for all $i \in N_n$ and $k \in N_s$.

To solve equation (2) for $P_{ij}$ ($i \in N_n$, $j \in N_m$), we use the feed forward neural network with m inputs and only one layer with n neurons.

First, the activation function employed by the neurons is not the sigmoid function, but the so-called linear activation function $f$ defined for all $a \in R$ by

$$f(a) = \begin{cases} 0 & \text{if } a < 0 \\ a & \text{if } a \in [0,\,1] \\ 1 & \text{if } a > 1. \end{cases}$$

Second, the output $y_i$ of neuron i is defined by

$$y_i = f\!\left( \max_{j \in N} W_{ij} X_j \right) \quad \left( i \in N_n \right).$$

Given equation 1, the training set of columns $q_k$ of matrix Q as input ($x_j = q_{ik}$ for each $j \in N_m$, $k \in N_s$) and columns $r_k$ of matrix R as expected output ($y_i = r_{jk}$ for each $i \in N_n$ and $k \in N^s$). Applying this training set to the feed forward neutral network, we obtain a solution to equation 1, when the error function reaches zero. The solution is then expressed by the weight $w_{ij}$ as $p_{ij} = w_{ij}$ for all $i \in N_n$ and $j \in N_m$. Thus $p = (w_{ij})$ is a n × n matrix.

It is already known that the fuzzy relation equation is in the dominant stage and there is lot of scope in doing research in this area, further it is to be also tested in real data.

Here we are transforming the single equation P ∘ Q = R into a collection of equations. When the word preference is said, there should be many preferences. If only one choice is given or the equation results in one solution, it cannot be called as the preference. Further, when we do some experiment in the real data, we may have too many solutions of preferences. For a unique solution sought out cannot or may not be available in reality, so to satisfy all the conditions described above, we are forced to reformulate the equation P ∘ Q = R. We partition the set Q into number of partition depending on the number preferences. When Q is partitioned, correspondingly R also gets partitioned, hence the one equation is transformed into the preferred number of equations.

Thus Q and R are given and P is to be determined. We partition Q into s sets, say $Q_1, Q_2, \ldots, Q_s$ such that $Q = Q_1 \cup Q_2 \cup \ldots \cup Q_s$, correspondingly R will be partitioned as $R = R_1 \cup R_2 \cup \ldots \cup R_s$. Now the resulting modified fuzzy equat6ions are $P_1 \circ Q_1 = R_1, P_2 \circ Q_2 = R_2, \ldots, P_s \circ Q_s = R_s$ respectively. Hence by our method, we obtain s preferences.

Since in reality it is difficult to make the error function $E_p$ to be exactly zero, we in our new fuzzy relation equation accept for the error function $E_p$ to be very close to zero. This is a deviation of the formula. Also we do not accept many stages in the arriving of the result. So once a proper guess is made even at the first stage we can get the desired solution by making $E_p$ very close to zero.

We are to find the passenger preference for a particular hour. The passenger preference problem for a particular hour reduces to finding the peak hours of the day (by peak hours of the day, we mean the number of passengers traveling in that hour is the maximum). Since the very term, preference by a passenger for a particular hour is an attribute, we felt it would be interesting if we adopt the modified fuzzy relation equation to this problem

So in our problem, we use the fuzzy relation equation $P \circ Q = R$, where P denotes the preference of a passenger to a particular hour, Q denotes the specific hour under study say $h_i$, $i = 1, 2, \ldots, 17$ where $h_i$ denotes the hour ending at 6 a.m., $h_2$ denotes the hour ending at 7 a.m., $\ldots$, $h_{17}$ denotes the hour ending at 10 p.m. and $R_i$ denotes the number of passengers traveling during that particular hour $h_i$, for $i = 1, 2, \ldots, 17$.

Here we use the fuzzy relation equation to determine P. We formulate the problem as follows:

If $h_i$, for $i = 1, 2, \ldots, n$ are the n-hour endings, $R_i$, for $i = 1, 2, \ldots, n$ denote the number of passengers traveling during hour $h_i$, for $i = 1, 2, \ldots, n$. We denote by R the set $\{R_1, R_2, \ldots, R_n\}$ and $Q = \{h_1, h_2, \ldots, h_n\}$. To calculate the preference of a passenger to a particular hour we associative with each $R_i$, a weight $w_i$. Since $R_i$ correspond to the number of passenger traveling in that hour $h_i$, is a positive value and hence comparison between any two $R_i$ and $R_j$'s always exist. Therefore, if $R_i < R_j$, then we associate a weight $w_i$ to $R_i$ and $w_j$ to $R_j$ such that $w_i < w_j$, where $w_i$ and $w_j$ take values in the interval [0, 1].

Now we solve the matrix relation equation $P \circ Q = R$ and obtain the preference of the passenger P for a particular time period, which is nothing but the maximum number of passengers traveling in that hour.

If we wish to obtain only one peak hour for the day, we take all the n elements and form a matrix equation,

$$P_{max} \bullet \begin{bmatrix} h_1 \\ h_2 \\ \vdots \\ h_n \end{bmatrix}_{n \times 1} = \begin{bmatrix} R_1 \\ R_2 \\ \vdots \\ R_n \end{bmatrix}_{n \times 1}$$

and find the n × n matrix $P = (w_{ij})$ using the method described in the beginning. We choose in several steps the weight function $w_1$, $w_2$, …, $w_n$ so that the error function $E_p$ reaches very near to zero. It is pertinent to mention here for our passengers preference problem we accept a value other than zero but which is very close to zero as limit which gives us the desired preference.

If we wish to have two peaks hours, we partition Q into $Q_1$ and $Q_2$ so that correspondingly R gets partitioned in $R_1$ and $R_2$ and obtain the two peak hours using the two equations $P_1 \circ Q_1 = R_1$ and $P_2 \circ Q_2 = R_2$ respectively. Clearly $P_1 = \begin{pmatrix} 1 \\ w_{ij} \end{pmatrix}$ the weights associated with the set $R_1$ and $P_2 = \begin{pmatrix} 2 \\ w_{ij} \end{pmatrix}$ the weights associated with the set $R_2$.

If we wish to have a peak hours, s < n, then we partition $h_i$ for i = 1, 2, …, n into s disjoint sets and find the s peak hours of the day. This method of partitioning the fuzzy relation equation can be used to any real world data problem, though we have described in the context of the transportation problem.

We have tested our hypothesis in the real data got from Pallavan Transport Corporation.

Hour ending Q; 6, 7, 8, 9, 10, 11, 12, 13, 14, 15, 16, 17, 18, 19, 20, 21, 22.

Passengers per hour R: 96, 71, 222, 269, 300, 220, 241, 265, 249, 114, 381, 288, 356, 189, 376, 182, 67.

We have partitioned the 17 hours of the day Q.

i) by partitioning Q into three elements each so as to get five preferences,

ii) by partitioning Q into five elements each so as to get three preferences and

iii) by arbitrarily partitioning Q into for classes so as to get four preferences.

In all cases from these real data, our predicated value coincides with the real preference value.

Since all the concepts are to be realized as fuzzy concepts, we at the first state make the entries of Q and R to lie between 0 and 1. This is done by multiplying Q by $10^{-2}$ and R by $10^{-4}$ respectively.

We partition Q into three elements each by taking only the first 15 elements from the table. That is $Q = Q_1 \cup Q_2 \cup Q_3 \cup Q_4 \cup Q_5$ and the corresponding $R = R_1 \cup R_2 \cup R_3 \cup R_4 \cup R_5$.

$$\text{For } Q_1 = x_i \qquad R_1 = r_{ik}$$

| | |
|---|---|
| 0.06 | 0.0096 |
| 0.07 | 0.0071 . |
| 0.08 | 0.0222 |

The fuzzy relation equation is

$$P1 = \begin{bmatrix} 0.06 \\ 0.07 \\ 0.08 \end{bmatrix} = \begin{bmatrix} 0.0096 \\ 0.0071 \\ 0.0222 \end{bmatrix}.$$

We employ the same method described earlier, where, the linear activation function $f$ is defined by

$$f(a) = \begin{cases} 0 & if \ a < 0 \\ a & if \ a \in [0,\ 1] \\ 1 & if \ a < 1 \end{cases}$$

for all $a \in \mathbb{R}$ and output $y_i$ of the neuron i is defined by

$$y_i = f\left( \max_{j \in N_m} W_{ij} \ X_j \right) \ \left( i \in N_n \right)$$

calculate $\max\limits_{j \in N_m} W_{ij}$ $X_j$ as follows:

(i) $w_{11}x_1$ = $0.03 \times 0.06$ = $0.0018$
    $w_{12}x_2$ = $0.0221875 \times 0.07$ = $0.001553125$
    $w_{13}x_3$ = $0.069375 \times 0.08$ = $0.00555$
    ∴ Max $(0.0018, 0.001553125, 0.00555) = 0.00555$

$$f\left(\max\limits_{j \in N_m} W_{ij} \ X_j\right) = f(0.00555) = 0.00555 \text{ (Since } 0.00555 \in [0, 1])$$
$$\therefore y_1 = 0.00555.$$

(ii) $w_{21}x_1$ = $0.06 \times 0.06$ = $0.0036$
     $w_{22}x_2$ = $0.044375 \times 0.07$ = $0.00310625$
     $w_{23}x_3$ = $0.13875 \times 0.08$ = $0.0111$
     ∴ Max $(0.0036, 0.00310625, 0.0111) = 0.0111$

$$f\left(\max\limits_{j \in N_m} W_{ij} \ X_j\right) = f(0.0111) = 0.0111 \text{ (Since } 0.0111 \in [0, 1])$$
$$\therefore y_2 = 0.0111$$

(iii) $w_{31}x_1$ = $0.12 \times 0.06$ = $0.0072$
      $w_{32}x_2$ = $0.08875 \times 0.07$ = $0.0062125$
      $w_{33}x_3$ = $0.2775 \times 0.08$ = $0.0222$
∴ Max $(0.0072, 0.0062125, 0.0222) = 0.0222$ (Since $0.0222 \in [0, 1]$
$\therefore y_3 = 0.0222$

$$f\left(\max\limits_{j \in N_m} W_{ij} \ X_j\right) = f(0.0222) = 0.0222 \text{ (Since } 0.0222 \in [0, 1])$$
$$\therefore y_3 = 0.0222.$$

Feed Forward Neural Network representing the solution is shown above.

$$\therefore P_1 = \begin{bmatrix} 0.03 & 0.06 & 0.12 \\ 0.0221875 & 0.044375 & 0.08875 \\ 0.069375 & 0.13875 & 0.2775 \end{bmatrix}.$$

Verification:

Consider, $P \circ Q = R$

that is $\max_{j \in N_m} p_{ij} q_{jk} = r_{ik}$

$\therefore$ Max (0.0018, 0.0042, 0.0096) $=$ 0.0096
Max (0.00133125, 0.00310625, 0.0071) $=$ 0.0071
Max (0.0041625, 0.0097125, 0.0222) $=$ 0.0222.

Similarly by adopting the above process, we have calculated the passenger preferences $P_2$, $P_3$, $P_4$ and $P_5$ for the pairs $(Q_2, R_2)$, $(Q_3, R_3)$, $(Q_4, R_4)$ and $(Q_5, R_5)$.

For

$$
\begin{array}{cc}
Q_2 & R_2 \\
0.09 & 0.0269 \\
0.10 & 0.0300 \\
0.11 & 0.0220
\end{array}
\quad \text{we have } P_2 =
\begin{bmatrix}
0.1345 & 0.269 & 0.06725 \\
0.15 & 0.3 & 0.075 \\
0.11 & 0.22 & 0.00605
\end{bmatrix}.
$$

For

$$
\begin{array}{cc}
Q_3 & R_3 \\
0.12 & 0.0241 \\
0.13 & 0.0265 \\
0.14 & 0.0249
\end{array}
\quad \text{we have } P_3 =
\begin{bmatrix}
0.2008 & 0.1004 & 0.0502 \\
0.2208 & 0.1104 & 0.0552 \\
0.2075 & 0.10375 & 0.051875
\end{bmatrix}.
$$

For

$$
\begin{array}{cc}
Q_4 & R_4 \\
0.15 & 0.0114 \\
0.16 & 0.0381 \\
0.17 & 0.0288
\end{array},
\quad P_4 =
\begin{bmatrix}
0.035625 & 0.07125 & 0.0178125 \\
0.1190625 & 0.23125 & 0.05953125 \\
0.09 & 0.18 & 0.045
\end{bmatrix},
$$

and for

$$
\begin{array}{cc}
Q_5 & R_5 \\
0.18 & 0.0356 \\
0.19 & 0.0189 \\
0.20 & 0.0376
\end{array}
\quad \text{we have } P_5 =
\begin{bmatrix}
0.0445 & 0.089 & 0.178 \\
0.023625 & 0.04725 & 0.0945 \\
0.047 & 0.094 & 0.188
\end{bmatrix}.
$$

On observing from the table, we see the preference $P_1$, $P_2$, $P_3$, $P_4$, and $P_5$ correspond to the peak hours of the day, $h_3$ that is 8 a.m. with 222 passengers, $h_5$ that is 10 a.m. with 300 passengers, $h_8$ that is 1 p.m. with 265 passengers, $h_{11}$ that is 4 p.m. with 381 passengers and $h_{15}$ that is 8 p.m. with 376 passengers. Thus this partition gives us five preferences with coincides with the real data as proved by the working.

max (0.003108, 0.00444, 0.00666, 0.0222, 0.01221) = 0.0222
max (0.003766, 0.00538, 0.0080694, 0.0269, 0.014795) = 0.0269
max (0.0042, 0.006, 0.009, 0.03, 0.0165) = 0.03
max (0.00308, 0.0044, 0.0065997, 0.0222, 0.0121) = 0.0222

Similarly we obtain the passenger preference P for the other entries using the above method.

For     $Q_2$          $R_2$
        0.12           0.0241
        0.13           0.0265
        0.14           0.0249
        0.15           0.0114
        0.16           0.0381

we have

$$P_2 = \begin{bmatrix} 0.030125 & 0.03765625 & 0.05020833 & 0.0753125 & 0.150625 \\ 0.033125 & 0.04140625 & 0.0552083 & 0.0828125 & 0.165625 \\ 0.031125 & 0.03890625 & 0.051875 & 0.0778125 & 0.155625 \\ 0.01425 & 0.0178125 & 0.02375 & 0.035625 & 0.07125 \\ 0.047625 & 0.05953 & 0.079375 & 0.1190625 & 0.238125 \end{bmatrix}$$

and for   $Q_3$              $R_3$
          0.17               0.0288
          0.18               0.0356
          0.19               0.0189
          0.20               0.0376
          0.21               0.0182 we have

$w_{54} x_4$ = $0.15 \times 0.10$ = 0.015
$w_{55} x_5$ = $0.11 \times 0.11$ = 0.0121.

$\therefore$ Max = (0.002485, 0.00888, 0.012105, 0.015, 0.0121) = 0.015

$$f\left(\max_{j \in N_m} W_{ij} \; X_j\right) = f(0.015) = 0.015 \text{ (Since } 0.015 \in [0.1])$$

$\therefore y_5$ = 0.015.

Feed forward neural network representing the solution is shown above

$$\therefore P_1 = \begin{bmatrix} 0.0142 & 0.01775 & 0.02366 & 0.071 & 0.03555 \\ 0.0444 & 0.0555 & 0.074 & 0.222 & 0.111 \\ 0.0538 & 0.06725 & 0.08966 & 0.269 & 0.1345 \\ 0.06 & 0.075 & 0.1 & 0.3 & 0.15 \\ 0.044 & 0.055 & 0.07333 & 0.220 & 0.11 \end{bmatrix}.$$

Verification :

Consider, $P \circ Q = R$

that is $\max_{j \in N_m} p_{ij} \, q_{jk} = r_{ik}$

max (0.000994, 0.00142, 0.0021294, 0.0071, 0.003905) = 0.0071

$w_{24} \, x_4$ = $0.075 \times 0.10$ = $0.0075$
$w_{25} \, x_5$ = $0.055 \times 0.11$ = $0.00605$
$\therefore$ Max (0.0012425, 0.00444, 0.0060525, 0.0075, 0.00605) = 0.0075

$f\left(\max_{j \in N_m} W_{ij} \, X_j\right) = f(0.0075) = 0.0075$ (Since $0.0075 \in [0,1]$)

$\therefore y_2 = 0.0075$

(iii)
$w_{31} \, x_1$ = $0.02366 \times 0.07$ = $0.016566$
$w_{32} \, x_2$ = $0.074 \times 0.08$ = $0.00592$
$w_{33} \, x_3$ = $0.08966 \times 0.09$ = $0.00807$
$w_{34} \, x_4$ = $0.1 \times 0.10$ = $0.010$
$w_{35} \, x_5$ = $0.07333 \times 0.11$ = $0.008066$
$\therefore$ Max (0.016566, 0.00592, 0.00807, 0.010, 0.008066) = 0.016566

$f\left(\max_{j \in N_m} W_{ij} \, X_j\right) = f(0.016566) = 0.016566$ (Since $0.016566 \in [0, 1]$)

$\therefore y_3 = 0.016566$

(iv)
$w_{41} \, x_1$ = $0.071 \times 0.07$ = $0.00497$
$w_{42} \, x_2$ = $0.222 \times 0.08$ = $0.01776$

$w_{43} x_3$ = $0.269 \times 0.09$ = $0.02421$
$w_{44} x_4$ = $0.3 \times 0.10$ = $0.03$
$w_{45} x_5$ = $0.220 \times 0.11$ = $0.0242$
$\therefore$ Max $(0.00497, 0.01776, 0.02421, 0.03, 0.0242) = 0.03$

$$f\left(\max_{j \in N_m} W_{ij} \ X_j\right) = f \ (0.03) = 0.03 \ (\text{Since } 0.03 \in [0, 1])$$

$$\therefore y_4 = 0.03$$

(v)
$w_{51} x_1$ = $0.0355 \times 0.07$ = $0.002485$
$w_{52} x_2$ = $0.111 \times 0.08$ = $0.00888$
$w_{53} x_3$ = $0.1345 \times 0.09$ = $0.012105.$

Now, we partition Q into five elements each by leaving out the first and the last element from the table as $Q_1$ $Q_2$ and $Q_3$ and calculate $P_1, P_2$ and $P_3$ as in the earlier case:

for        $Q_1$              $R_1$
           0.07            0.0071
           0.08            0.0222
           0.09            0.0269
           0.10            0.0300
           0.11            0.0220.

The fuzzy relation equation is

$$P_1 \circ \begin{bmatrix} 0.06 \\ 0.08 \\ 0.09 \\ 0.10 \\ 0.11 \end{bmatrix} = \begin{bmatrix} 0.0071 \\ 0.0222 \\ 0.0269 \\ 0.0300 \\ 0.0220 \end{bmatrix}.$$

Calculate max $w_{ij} x_j$ as follows
        $j \in N_m$
(i)
$w_{11} x_1$ = $0.0142 \times 0.07$ = $0.000994$
$w_{12} x_2$ = $0.0444 \times 0.08$ = $0.003552$
$w_{13} x_3$ = $0.0538 \times 0.09$ = $0.004842$
$w_{14} x_4$ = $0.06 \times 0.10$ = $0.006$

$w_{15} x_5$ $=$ $0.044 \times 0.11$ $=$ $0.00484$

$\therefore$ Max $(0.000994, 0.003552, 0.004842, 0.006, 0.00484) = 0.006$

$$f\left(\max_{j \in N_m} W_{ij} \; X_j\right) = f(0.006) = 0.006 \text{ (Since } 0.006 \in [0, 1])$$

$$\therefore y_1 = 0.006$$

(ii)

$w_{21} x_1$ $=$ $0.01775 \times 0.07$ $=$ $0.0012425$
$w_{22} x_2$ $=$ $0.0555 \times 0.08$ $=$ $0.00444$
$w_{23} x_3$ $=$ $0.06725 \times 0.09$ $=$ $0.0060525$

$$P_3 = \begin{bmatrix} 0.0288 & 0.036 & 0.048 & 0.149 & 0.072 \\ 0.0356 & 0.0445 & 0.05933 & 0.178 & 0.089 \\ 0.0189 & 0.023625 & 0.0315 & 0.0945 & 0.04725 \\ 0.0376 & 0.047 & 0.06266 & 0.0188 & 0.094 \\ 0.0182 & 0.02275 & 0.03033 & 0.091 & 0.0455 \end{bmatrix}.$$

On observing from the table, we see the preference $P_1$, $P_2$ and $P_3$ correspond to the peak hours of the day, $h_5$ that is 10 a.m. with 300 passengers, $h_{11}$ that is 4 p.m. with 381 passengers and $h_{15}$ that is 8 p.m. with 376 number of passengers. Thus this partition gives us three preferences, which coincides with the real data as proved by the working.

We now partition Q arbitrarily, that is the number of elements in each partition is not the same and by a adopting the above method we obtain the following results:

for

| $Q_1$ | $R_1$ |
|---|---|
| 0.06 | 0.0096 |
| 0.07 | 0.0071 |
| 0.08 | 0.0222 |

we have $P_1 = \begin{bmatrix} 0.03 & 0.06 & 0.12 \\ 0.0221875 & 0.044375 & 0.08875 \\ 0.069375 & 0.13875 & 0.2775 \end{bmatrix}.$

For

| $Q_2$ | $R_2$ |
|---|---|
| 0.09 | 0.269 |
| 0.10 | 0.300 |
| 0.11 | 0.220 |
| 0.12 | 0.241 |
| 0.13 | 0.265 |

we have

$$P_2 = \begin{bmatrix} 0.1345 & 0.17933 & 0.1076 & 0.08966 & -.06725 \\ 0.15 & 0.2 & 0.12 & 0.1 & 0.075 \\ 0.11 & 0.14666 & 0.088 & 0.0733 & 0.055 \\ 0.1205 & 0.16066 & 0.0964 & 0.08033 & 0.06025 \\ 0.1325 & 0.17666 & 0.106 & 0.08833 & 0.06625 \end{bmatrix}.$$

For
| $Q_3$ | $R_3$ |
|-------|-------|
| 0.14 | 0.0249 |
| 0.15 | 0.0114 |

we have $P_3 = \begin{bmatrix} 0.083 & 0.166 \\ 0.038 & 0.076 \end{bmatrix}$

and for
| $Q_4$ | $R_4$ |
|-------|-------|
| 0.16 | 0.0381 |
| 0.17 | 0.0288 |
| 0.18 | 0.0356 |
| 0.19 | 0.0189 |
| 0.20 | 0.0376 |
| 0.21 | 0.0182 |

we have

$$\begin{bmatrix} 0.09525 & 0.0635 & 0.047625 & 0.0381 & 0.1905 & 0.03175 \\ 0.072 & 0.048 & 0.036 & 0.0288 & 0.144 & 0.024 \\ 0.08233 & 0.05322 & 0.04111 & 0.03193 & 0.178 & 0.02411 \\ 0.4725 & 0.0312 & 0.023625 & 0.0189 & 0.0945 & 0.001575 \\ 0.0762 & 0.0508 & 0.0381 & 0.03048 & 0.188 & 0.0254 \\ 0.04275 & 0.0285 & 0.0213756 & 0.0171 & 0.091 & 0.001425 \end{bmatrix}.$$

We obtain in the preferences $P_1$, $P_2$, $P_3$ and $P_4$ by partitioning the given data into a set of three elements, a set of five elements, a set of two elements and a set of six elements. On observing from the table, we see that these preferences correspond to the peak hours of the day, $h_3$ that is 8 a.m. with 222 passengers, $h_5$ that is 10 a.m. with 300 passengers, $h_9$ that is 2 p.m. with 249 number of passengers and $h_{11}$ that is 4 p.m. with 381 number of passengers. Thus this partition gives us four preferences which coincides with the real data as proved by the working.

Thus the Government sector can run more buses at the peak hours given and also at the same time restrain the number of buses in the non peak hours we derived the following conclusions:

1. The fuzzy relation equation described given by 1 can give only one preference function P ∘ Q = R but the partition method described by us in this paper can give many number of preferences or desired number of preferences.
2. Since lot of research is needed we feel some other modified techniques can be adopted in the FRE P ∘ Q = R.
3. We have tested our method described in the real data taken from Pallavan Transport Corporation and our results coincides with the given data.
4. We see the number of preference is equal to the number of the partition of Q.
5. Instead of partitioning Q, if we arbitrarily take overlapping subsets of Q certainly we may get the same preference for two or more arbitrary sets.

We see that our method of the fuzzy relation equation can be applied to the peak hour problem in a very successful way. Thus only partitioning of Q can yield non-overlapping unique solution.

Finally, in our method we do not force the error function $E_p$ to become zero, by using many stages or intermittent steps. We accept a value very close to zero for $E_p$ as a preference solution.

For more please refer [98].

## 2.3 Study of the proper proportion of Raw material mix in cement plants using FRE

By the use of fuzzy relational equations and fuzzy neural network for fuzzy relation equation method we study the proper proportion of raw material mix to find the best quality of clinker.

As in other cases we use FRE and when the solution is unavailable by the method of FRE we adopt the neural networks we give the definition and control algorithm and use to solve the problem. We show by our method the cement industries can produce a desired quality of clinker. For more refer [102]

## 2.4 The effect of globalization on Silk weavers who are Bonded labourers using FRE

The strategies of globalization and the subsequent restructuring of economies, including the increased mechanization of labor has had stifling effects on the lives of the silk weavers in the famous Kancheepuram District in Tamil Nadu, India. Here, we study the effects of globalization, privatization and the mechanization of labor, and how this has directly affected (and ruined) the lives of thousands of silk weavers, who belong to a particular community whose tradition occupation is weaving. This research work is based on surveys carried out in the Ayyampettai village near Kancheepuram. The population of this village is around 200 families, and almost all of them are involved in the weaving of silk saris. They are skilled weavers who don't have knowledge of any other trade. Most of them are bond to labor without wages, predominantly because they had reeled into debt for sums ranging from Rs. 1000 upwards. They barely manage to have a square meal a day, and their work patterns are strenuous - they work from 6 a.m. to 7 p.m. on all days, expect the new moon day when they are forbidden from weaving.

Interestingly, their children are not sent to school, they are forced into joining the parental occupation, or into taking petty jobs in order to secure the livelihood. The villagers point to the advent of electric looms and reckon that their lives were much more bearable before this mechanization, at least they used to get better incomes. The wide scale introduction to electric looms / textile machines / power looms, has taken away a lot of their job opportunities. For instance, the machine can weave three silk saris which manually takes fifteen days to weave in just three hours. Also, machine woven silk saris are preferred to hand woven silk saris as a result of which their life is further shattered. Interviews with the weavers revealed the careless and negligent approach of the government to their problem. Here, we study their problem and the effect of globalization on their lives using Fuzzy Relational Equations. We have arrived at interesting conclusions to understand and assay this grave problem.

We have made a sample survey of around 50 families out of the 200 families; who are bonded labourers living in Ayyampettai near Kancheepuram District in Tamil Nadu; have become bonded for Rs.1000 to Rs.2000. They all belong to Hindu community viz. weavers or they traditionally call themselves as Mudaliar caste. Most of the owners are also Mudaliars. They were interviewed using a linguistic questionnaire. Some of the notable facts about their lives are as follows:

1. They do not know any other trade or work but most of them like to learn some other work.

2. They are living now below the poverty line because of the advent of electrical or power looms which has drastically affected their income.

3. The whole family works for over 10 hours with only one day i.e. new moon day in a month being a holiday. On new moon day they don't weave and they are paid by their owners on that day.

4. Only one had completed his school finals. All others have no education for they have to learn the trade while very young.

5. They don't have even a square meal a day.

6. Becoming member of Government Society cannot be even dreamt for they have to pay Rs.3000/- to Rs.5000 to Government and 3 persons should give them surety. So out of the 200 families there was only one was a Government Society member. After the globalization government do not give them any work because marketers prefer machine woven saris to hand woven ones.

7. Owners of the bonded labourers are not able to give work to these labourers.

8. Observations shows that female infanticide must be prevalent in these families as over 80% of the children are only males.

9. The maximum salary a family of 3 to 4 members is around Rs. 2000/- 5% of them alone get this 90% of the families get below Rs.2000 p.m.

10. Paying as rent, electricity, water, etc makes them live below poverty line.

The following attributes are taken as the main point for study:

$B_1$ – No knowledge of any other work has made them not only bonded but live in penury.

$B_2$ – Advent of power looms and globalization (modern textile machinery) has made them still poorer.

$B_3$ – Salary they earn in a month.

| $B_4$ | – | No savings so they become more and more bonded by borrowing from the owners, they live in debts. |
| $B_5$ | – | Government interferes and frees them they don't have any work and Government does not give them any alternative job. |
| $B_6$ | – | Hours / days of work. |

We have taken these six heads $B_1$, $B_2$, … , $B_6$ related to the bonded labourers as the rows of the fuzzy relational matrix.

The main attributes / heads $O_1$, $O_2$, $O_3$, $O_4$ related to the owners of the bonded labourers are :

| $O_1$ | – | Globalization / introduction of modern textile machines. |
| $O_2$ | – | Profit or no loss. |
| $O_3$ | – | Availability of raw goods. |
| $O_4$ | – | Demand for finished goods. |

Using these heads related to owners along columns the fuzzy relational equations are formed using the experts opinions.

The following are the limit sets using the questionnaire :

| $B_1 \geq 0.5$ | Means no knowledge of other work hence live in poverty. |
| $B_2 \geq 0.5$ | Power loom / other modern textile machinery had made their condition from bad to worse. |
| $B_3 \geq 0.5$ | Earning is mediocre. ($B_3 < 0.5$ implies the earning does not help them to meet both ends). |
| $B_4 \geq 0.4$ | No saving no debt. ($B_4 < 0.4$ implies they are in debt). |
| $B_5 \geq 0.5$ | Government interference has not helped. ($B_5 < 0.5$ implies Government Interference have helped). |
| $B_6 \geq 0.4$ | 10 hours of work with no holidays. ($B_6 < 0.4$ implies less than 10 hours of work). |

| $O_1 \geq 0.5$ | The globalizations / government has affected the owners of the bonded labourers drastically ($O_1 < 0.5$ implies has no impact on owners). |
| $O_2 \geq 0.5$ | Profit or no loss ($O_2 < 0.5$ implies total loss). |
| $O_3 \geq 0.6$ | Availability of raw materials. ($O_3 < 0.6$ implies |

shortage of raw material).

$O_4 \geq 0.5$       Just they can meet both ends i.e. demand for finished goods and produced goods balance. ($O_4 < 0.5$ implies no demand for the finished product i.e. demand and supply do not balance).

The opinion of the first expert who happens to be a bonded labor for the two generations aged in seventies is given vital importance and his opinion is transformed into the Fuzzy Relational Equation

$$P = \begin{array}{c} \\ B_1 \\ B_2 \\ B_3 \\ B_4 \\ B_5 \\ B_6 \end{array} \begin{array}{cccc} O_1 & O_2 & O_3 & O_4 \\ \left[ \begin{array}{cccc} .8 & 0 & 0 & 0 \\ .8 & .3 & .3 & 0 \\ .1 & .2 & .3 & .4 \\ 0 & .1 & .1 & .1 \\ .8 & .1 & .2 & .4 \\ .2 & .4 & .4 & .9 \end{array} \right] \end{array} .$$

By considering the profit suppose the owner gives values for Q where $Q^T = [.6, .5, .7, .5]$. Now P and Q are known in the fuzzy relational equation P o Q = R .

Using the max-min principle in the equation P o Q = R.

We get $R^T = [.6, .6, .4, .1, .6, .5]$ In the fuzzy relational equation P o Q = R, P corresponds to the weightages of the expert, Q is the profit the owner expects and R is the calculated or the resultant giving the status of the bonded labourers. When we assume the owners are badly affected by globalizations, but wants to carry out his business with no profit or loss, with moderate or good availability of the raw material and they have enough demand or demand and supply balance we obtain the following attributes related with the bonded labourers. The bonded labourers live in acute poverty as they have no other knowledge of any other work. The power loom has made their life from bad to worst, but the earning is medium with no savings and debts. They do not receive any help from the government, but they have to labor more than ten hours which is given by $[.6, .6, .4, .1, .6, .5]^T$.

Using the same matrix P and taking the expected views of the bonded labourers R to be as $[ .6, .4, .5, .4, .2, .6]^T$ .

Using the equation $P^T$ o R = Q. We obtain Q = $[.6, .4, .4, .6]^T$.

The value of Q states the owners are affected by globalization. They have no profit but loss. They do not get enough raw materials to give to the bonded labor as the market prefers machine woven saris to hand made ones so the demand for the finished goods declines. Thus according to this expert, the main reason for their poverty is due to globalization i.e. the advent of power looms has not only affected them drastically as they do not have the knowledge of any other trade but it has also affected the lives of their owners.

A small owner who owns around ten bonded labor families opinion is taken as the second experts opinion. The weighted matrix P as given by the second expert is:

$$P = \begin{bmatrix} .7 & .1 & 0 & 0 \\ .9 & .2 & .3 & 0 \\ .0 & .1 & .2 & .3 \\ 0 & 0 & .1 & .1 \\ .9 & 0 & .1 & .4 \\ .1 & .2 & .4 & .7 \end{bmatrix}.$$

By considering the profit the owner expects i.e. taking Q = [.6, .5, .7, .5]$^T$
We calculate R using P o Q$^T$ = R
   i.e. R   = [.6, .6, .3, .1, .6, .5]$^T$.

We obtain the following attributes from R related with the bonded labourers.

They live in below poverty, as they have no other trade but the earning is medium with no savings and new debts. They do not get any help from the government, but they have to work more than 10 hours a day which is given by [.6, .6, .3, .1, .6, .5]$^T$ .

Using the same P i.e. the weightages we now find Q giving some satisfactory norms for the bonded labourers.

By taking R = [ .6, .4, .5, .4, .2, .6]$^T$ and using the equation P$^T$ o R = Q,
   i.e. Q = [.6, .2, .4, .6]$^T$

which states the owners are badly affected by globalization. They have no profit but loss they do not get enough raw materials and the demand for the finished goods declines.

The third expert is a very poor bonded labor. The fuzzy relational matrix P given by him is

$$P = \begin{bmatrix} .9 & 0 & 0 & 0 \\ .5 & .3 & .4 & .1 \\ .2 & .2 & .2. & 3 \\ 0 & 0 & .1 & .2 \\ .7 & .2 & .2 & .4 \\ .2 & .3 & .3 & .8 \end{bmatrix}.$$

By considering the profit the owner expects i.e. taking Q = [.6, .5, .7, .5]$^T$ and using the relational equation P o Q = R, we calculate R;

R = [.6, .5, .3, .2, .6, .5]$^T$ .

We obtain the following attributes from R related with the bonded labourers.

This reveals that the bonded labourers standard of living is in a very pathetic condition. They do not have any other source of income or job. Their earning is bare minimum with no savings. Neither the government comes forward to help them nor redeem them from their bondage. In their work spot, they have to slog for 10 hours per day.

Using the same P i.e. the weightages we now find Q by giving some satisfactory norms for the bonded labourers.

By taking R = [ .6, .4, .5, .4, .2, .6]$^T$ and using the equation P$^T$ o R = Q,

i.e. Q = [.6, .3, .4, .6]$^T$.

The value of Q states due to the impact of globalization (modern textile machinery), the owners are badly affected. They are not able to purchase enough raw materials and thus the out put from the industry declines. The owners do not get any profit but eventually end up in a great loss.

The following conclusions are not only derived from the three experts described here but all the fifty bonded labourers opinions are used and some of the owners whom we have interviewed are also ingrained in this analysis.

1.    Bonded labourers are doubly affected people for the advent of globalization (modern textile machinery) has

denied them small or paltry amount, which they are earning in peace as none of them have knowledge of any other trade.

2.  The government has not taken any steps to give or train them on any trade or work or to be more precise they are least affected about their living conditions of them. Some of them expressed that government is functioning to protect the rich and see the rich do not loose anything but they do not even have any foresight about the bonded labourers or their petty owners by which they are making the poor more poorer.

3.  Bonded labourers felt government has taken no steps to eradicate the unimaginable barrier to become members of the government society. They have to pay Rs.3000/- and also they should spell out and get the surety of 3 persons and the three persons demand more than Rs.3000/- each so only they are more comfortable in the hands of their owners i.e. as bonded labourers were at least they exist with some food, though not a square meal a day.

4.  It is high time government takes steps to revive the life of the weavers who work as bonded labourers by training and giving them some job opportunities.

5.  They felt government was killing the very talent of trained weavers by modernization as they have no knowledge of any other trade.

6.  Child labor is at its best in these places as they cannot weave without the help of children. Also none of the children go to school and they strongly practice female infanticide.

7.  Government before introducing these modern textile machineries, should have analyzed the problem and should have taken steps to rehabilitate these weavers. Government has implemented textile machineries without foresight.

This research work will be published with the coauthor T.Narayanamoorthy.

## 2.5  Study of Bonded Labor Problem Using FRE

As we analyze a sample of over 1000 families of bonded labourers living and working in rice mills as bonded labourers in and around Red Hills Area, Chennai, India. It is a shocking information to see that roughly 1630 children work in the rice mill industry of which 40% are in the age group 0-5 years and 50% are in the are group 6-14 years and even the basic primary education is denied to them.

They have become bonded labourers for just less than Rs.5000/-. The children have become slaves and a possession of them in true sense.

We study the role played by politicians, educationalists, social workers, human right workers and above all the legal role and other factors to study and eradicate such cruelty. So in this analysis we have taken up the problem of Bonded labourers using the tool of FREs. We solve the fuzzy relation equation P o Q = R, where P, Q and R are matrices with entries from [0, 1]. Several interesting conclusions are derived from our study.

## 2.6  Data Compression with FREs

Kaoru Hirota and Witold Pedrycz [35] have studied data compression with fuzzy relational equations. They have introduced a concept of fuzzy relation calculus to the problems of image processing. They have discussed fuzzy relation based on data compression. The methodology of data compression hinges on the theory of fuzzy relational equations were the solutions to the specific class of equations give rise to a reconstructed fuzzy relation (image). The main properties of fuzzy relational equations were introduced and analyzed with respect to the resulting reconstruction and compression capabilities.
For the entire work please refer [35].

## 2.7  Applying FRE to Threat Analysis

The threat level posed by various targets in real life depends on various factors some are situation dependent and some related to characteristics of the target being analyzed, such are its formations, its firing status etc. A great deal of work has been done by Institute for simulation and Training (IST) to identify there factors Breached et all have in an ingenuous way applied fuzzy relation equations to threat analysis. They have presented a model to illustrate how to apply a fuzzy relational equation algorithm to perform threat analysis in the context of computer Generated Forces Systems such Mod SAF (Modular semi automated forces) Using fuzzy relational equation the proposed algorithm generates the data from the historic information and its earlier runs. FREs have been successfully applied to threat analysis.

## 2.8 FREs-application to medical diagnosis

Elie Sanchez [84] has studied truth qualification and fuzzy relations in natural languages and its applications to medical diagnosis. A biomedical application in which medical knowledge is expressed in a rule form with AND ed fuzzy propositions in the antecedent illustrates the aggregation of these measures for medical diagnosis assistance Elie Sanchez has illustrated the applications in the field of inflammatory protein variations [84].

A pattern of Medical knowledge consists here of a tableau with linguistic entries that will be interpreted as fuzzy sets, having in mind that different experts might provide somehow different characterizations for a same pattern. This medical knowledge is so translated into fuzzy propositions. A complete description of the problem and its application is found in [84].

## 2.9 A fuzzy relational identification algorithm and its application to predict the behavior to a motor drive system

Fuzzy systems are usually named model free estimators. They estimate input-output relations without the need of an analytical model of how outputs depend on inputs and encode the sampled information in a parallel distributed frame work called fuzzy structure.

Three main types of fuzzy structure are

      (1) Rule Based Systems,
      (2) Fuzzy relational systems and

(3) Fuzzy functional systems.

Let 'o' denote the max-min composition operator. $\overline{X}_1, \overline{X}_2, ..., \overline{X}_n$ denote the input fuzzy sets, $\overline{Y}$ stands for the output fuzzy set and R is the fuzzy relational matrix expressing the system's input-output relationship.

$$\overline{Y} = \overline{X}_1 \bullet \overline{X}_2 \bullet \cdots \bullet \overline{X}_n \bullet R \; . \qquad \text{(A)}$$

Fuzzy relational equations A describes multiple - input – single output fuzzy systems. From a system theory point of view, the following simplified version of (A) can be considered as a single input single output fuzzy system

$$\overline{Y} = \overline{X} \circ R \qquad \text{(B)}$$

and (C) is discretised for each instant

$$\overline{Y}_k^k = \overline{X}_k \circ R_k \; . \qquad \text{(C)}$$

Equation (C) can also be rewritten as

$$\overline{Y}_k(y_k) = \sup_{x_k \in X} \left[ \min \left( \overline{X}_k(x_k), R_k(x_k, y_k) \right) . \qquad \text{(D)}$$

A fuzzy relation R is written as a set of fuzzy rules with fuzzy sets defined on each universe of discourse. For a single-input-single output system (B) defined with n fuzzy sets for $\overline{X}$ and $\overline{Y}$, R is an n × n matrix of possibility measures with each element being denoted as in

$$\text{R (i, j)} = p_{ij}. \qquad \text{(E)}$$

Each matrix element can be translated as a linguistic simple rule like

$$[\text{If } \overline{X}_i, \; then \; \overline{Y}_j] \text{ with possibility b}_{ij} \qquad \text{(F)}$$

and for each condition $\overline{X}_i$ there are n simple rules that form a compound rule. Having given the main fuzzy set operations one can work for the relational identification of a motor drive system and the new fuzzy relational identification algorithm and obtain the method of speed signal of a motor drive system.

## 2.10  Application of genetic algorithms in problems in chemicals industry

Chemical Industries and Automobiles are extensively contributing to the pollution of environment, Carbon monoxide, nitric oxide, ozone, etc., are understood as the some of the factors of pollution from chemical industries. The maintenance of clean and healthy atmosphere makes it necessary to keep the pollution under control which is caused by combustion waste gas. The authors have suggested theory to control waste gas pollution in environment by oil refinery using fuzzy linear programming. To the best of our knowledge the authors [97]are the first one to apply fuzzy linear programming to control or minimize waste gas in oil refinery.

An oil refinery consists of several inter linked units. These units act as production units, refinery units and compressors parts. These refinery units consume high-purity gas production units. But the gas production units produce high-purity gas along with a low purity gas. This low purity gas goes as a waste gas flow and this waste gas released in the atmosphere causes pollution in the environment. But in the oil refinery the quantity of this waste gas flow is an uncertainty varying with time and quality of chemicals used in the oil refinery. Since a complete eradication of waste gas in atmosphere cannot be made; here one aims to minimize the waste gas flow so that pollution in environment can be reduced to some extent. Generally waste gas flow is determined by linear programming method. In the study of minimizing the waste gas flow, some times the current state of the refinery may already be sufficiently close to the optimum. To over come this situation we adopt fuzzy linear programming method.

The fuzzy linear programming is defined by

$$\text{Maximize} \quad z = cx$$
$$\text{Such that} \quad Ax \leq b$$
$$x \leq 0$$

where the coefficients $A$, $b$ and $c$ are fuzzy numbers, the constraints may be considered as fuzzy inequalities with variables $x$ and $z$. We use fuzzy linear programming to determine uncertainty of waste gas flow in oil refinery which pollutes the environment.

Oil that comes from the ground is called "Crude oil". Before one can use it, oil has to be purified at a factory called a "refinery", so as to convert into a fuel or a product for use. The refineries are high-tech factories, they turn crude oil into useful energy products.

During the process of purification of crude oil in an oil refinery a large amount of waste gas is emitted to atmosphere which is dangerous to human life, wildlife and plant life. The pollutants can affect the health in various ways, by causing diseases such as bronchitis or asthma, contributing to cancer or birth defects or perhaps by damaging the body's immune system which makes people more susceptible to a variety of other health risks. Mainly, this waste gas affects Ozone Layer. Ozone (or Ozone Layer) is 10-50 km above the surface of earth. Ozone provides a critical barrier to solar ultraviolet radiation, and protection from skin cancers, cataracts, and serious ecological disruption. Further sulfur dioxide and nitrogen oxide combine with water in the air to form sulfuric acid and nitric acid respectively, causing acid rain. It has been estimated that emission of 70 percentage of sulfur dioxide and nitrogen oxide are from chemical industries.

We cannot stop this process of oil refinery, since oil and natural gas are the main sources of energy. We cannot close down all oil refineries, but we only can try to control the amount of pollution to a possible degree. In this paper, the authors use fuzzy linear programming to reduce the waste gas from oil refinery.

The authors describe the knowledge based system (KBS) that is designed and incorporate it in this paper to generate an on-line advice for operators regarding the proper distribution of gas resources in an oil refinery. In this system, there are many different sources of uncertainty including modeling errors, operating cost, and different opinions of experts on operating strategy. The KBS consists of sub-functions, like first sub-functions, second sub-functions, etc. Each and every sub-functions are discussed relative to certain specific problems.

For example: The first sub-function is mainly adopted to the compressor parts in the oil refineries. Till date they were using stochastic programming, flexibility analysis and process design problems for linear or non-linear problem to compressor parts in oil refinery. Here we adopt the sub function to study the proper distribution of gas resources in an oil refinery and also use fuzzy linear programming (FLP) to minimize the waste gas flow. By the term proper distribution of gas we include the study of both the production of high-purity gas as well as the amount of waste gas flow which causes pollution in environment.

In 1965, Lofti Zadeh [115, 116] wrote his famous paper formally defining multi-valued, or "fuzzy" set theory. He extended traditional set theory by changing the two-values indicator functions i.e., 0, 1 or the crisp function into a multi-valued membership function. The membership function assigns a "grade of membership" ranging from 0 to 1 to each object in the fuzzy set. Zadeh formally defined fuzzy sets, their properties, and various properties of algebraic fuzzy sets. He introduced the concept of linguistic variables which have values that are linguistic in nature (i.e. pollution by waste gas = {small pollution, high pollution, very high pollution}).

Fuzzy Linear Programming (FLP): FLP problems with fuzzy coefficients and fuzzy inequality relations as a multiple fuzzy reasoning scheme, where the past happening of the scheme correspond to the constraints of thee FLP problem. We assign facts (real data from industries) of the scheme, as the objective of the FLP problem. Then the solution process consists of two steps. In the fist step, for every decision variable, we compute the (fuzzy) value of the objective function via constraints and facts/objectives. At the second step an optimal solution to FLP problem is obtained at any point, which produces a maximal element to the set of objective functions (in the sense of the given inequality relation).

The Fuzzy Linear Programming (FLP) problem application is designed to offer advice to operating personnel regarding the distribution of Gas within an oil refinery (Described in Figure 2.10.1) in a way which would minimize the waste gas in environment there by reduce the atmospheric pollution .

GPUI, GPU2 and GPU3 are the gas production units and GGG consumes high purity gas and vents low purity gas. Gas from these production units are sent to some oil refinery units, like sulfur, methanol, etc. Any additional gas needs in the oil refinery must be met by the gas production unit GPU3.

The pressure swing adsorption unit (PSA) separates the GPU2 gas into a high purity product stream and a low purity tail stream (described in the Figure 2.10.1). $C_1$, $C_2$, $C_3$, $C_4$, $C_5$, are compressors. The flow lines that dead –end is an arrow represent vent to flare or fuel gas. This is the wasted gas that is to be minimized. Also we want to minimize the letdown flow from the high purity to the low purity header

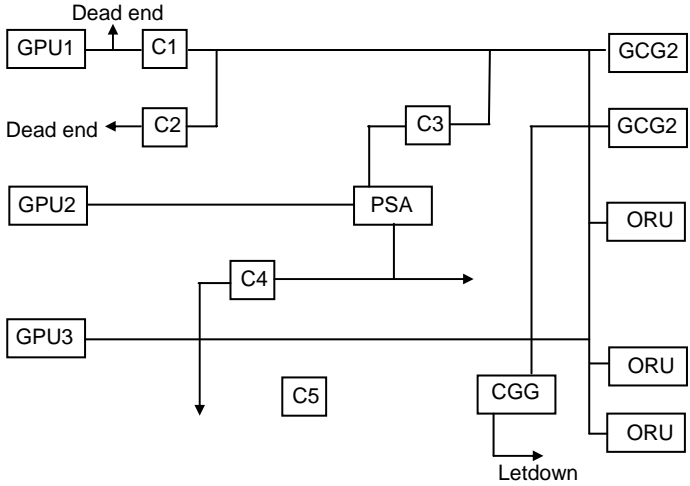

Figure: 2.10.1

FLP is a method of accounting for uncertainty is used by the authors for proper distribution of gas resources, so as to minimize the waste gas flow in atmosphere. FLP allows varying degrees of set membership based on a membership function defined over a range of values. The membership function usually varies from 0 to 1. FLP allow the representation of many different sources of uncertainty in the oil refinery. These sources may (or) may not be probabilistic in nature. The uncertainty is represented by membership functions describing the parameters in the optimization model. A solution is found that either maximizes a given feasibility measure and maximizes the wastage of gas flow. FLP is used here to characterize the neighborhood of solutions that defines the boundaries of acceptable operating states.

Fuzzy Linear Programming (FLP) can be stated as;

$$\left.\begin{array}{c} \max \text{imize } z = cx \\ \text{s.t } Ax \le b \\ x \ge 0 \end{array}\right\} \quad \dots (*)$$

The coefficients A, b and c are fuzzy numbers, the constraints may be considered as fuzzy inequalities. The decision space is defined by the constraints with c, $x \in N$, $b \in R^m$ and $A \in R^m$, where N, $R^m$, and $R^{mxn}$ are reals.

The optimization model chosen by the knowledge based system (KBS) is determined online and is dependent on the refinery units. This optimization method is to reduce the amount of waste gas in pollution.

We aim to

1. The gas ($GCG_2$) vent should be minimized.
2. The let down flow should be minimized and
3. The make up gas produced by the as production unit (GPU3) should be minimized.

Generally the waste gas emitted by the above three ways pollute the environment. The objective function can be expressed as the sum of the individual gas waste flows. The constrains are given by some physical limitations as well as operator entries that describe minimum and maximum desired flows.

The obtained or calculated resultant values of the decision variables are interpreted as changes in the pressure swing adsorption feed, and the rate that gas is imported to CGG and gas production unit (GPU3). But in the optimization model there is uncertainty associated with amount of waste gas from oil refinery, and also some times the current state of the refinery may already be sufficiently close to the optimum.

For example to illustrate the problem, if the fuzzy constraints $x_1$, the objects are taken along the x-axis are shown in the figures 2 and 3, which represent the expression.

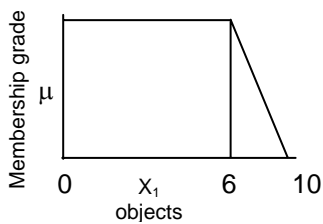

Figure: 2.10.2

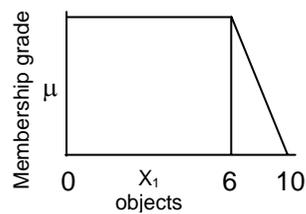

Figure: 2.10.3

$x_1 \leq 8$ (with tolerance p = 2)          (1)

The membership function μ are taken along the y-axis i.e. $\mu(x_1)$ lies in [0, 1] this can be interpreted as the confidence with which this constraint is satisfied (0 for low and 1 for high). The fuzzy inequality constraints can be redefined in terms of their α-cuts.

$$\{S_\alpha \,/\, \alpha \,\varepsilon\, [0, 1]\}, \text{ where } S_\alpha = \{\gamma \,/\, (\mu\,(\gamma) \geq \alpha)\}.$$

The parameter α is used to turn fuzzy inequalities into crisp inequalities. So we can rewrite equation (1)

$$x_1 \leq 6 + 2\,(2)\,(1-\alpha)$$
$$x_1 \leq 6 + 4\,(1-\alpha)$$

where $\alpha \,\varepsilon\, [0, 1]$ expressed in terms of α in this way the fuzzy linear programming problem can be solved parametrically. The solution is a function on α

$$x^* = f(\alpha) \qquad\qquad (2)$$

with the optimal value of the objective function determined by substitution in equation (1).

$$z^* = cx^* = g(\alpha). \qquad\qquad (3)$$

This is used to characterize the objective function. The result covers all possible solutions to the optimization problem for any point in the uncertain interval of the constraints.

The α-cuts of the fuzzy set describes the region of feasible solutions in figures 2 and 3. The extremes ($\alpha = 0$ and $\alpha = 1$) are associated with the minimum and maximum values of $x^*$ respectively. The given equation (2) can also be found this, is used to characterize the objective function. The result covers all possible solutions to the optimization problem for any point in the uncertain interval of the constraints.

Fuzzy Membership Function to Describe Uncertainty: The feasibility of any decision ($\mu_D$) is given by the intersection of the fuzzy set describing the objective and the constraints.

$$\mu_D\,(x) = \mu_z(x)\,{}^{\wedge}\mu_N\,(x)$$

where ^ represents the minimum operator, that is the usual operation for fuzzy set intersection. The value of $\mu_N$ can be easily found by intersecting the membership values for each of the constraints.

$$\mu_N (x) = \mu_1(x) \wedge \mu_2 (x) \wedge \ldots \wedge \mu_m (x).$$

The membership functions for the objective ($\mu_z$) however is not obvious z is defined in (2). Often, predetermined aspiration target values are used to define this function. Since reasonable values of this kind may not be available, the solution to the FLP equation (3) is used to characterize this function.

$$\mu_z(x) = \begin{cases} 1 & \text{if } z(x) \geq b(0) \\ \dfrac{z(x) - b(1)}{b(0) - b(1)} & \text{if } b(1) \leq z(x) \leq b(0) \\ 0 & \text{if } z(x) \leq b(0). \end{cases} \qquad (5)$$

The result is that the confidence value increases as the value of the objective value increases. This is reasonable because the goal is to maximize this function the limits on the function defined by reasonable value is obtained by extremes of the objective value.

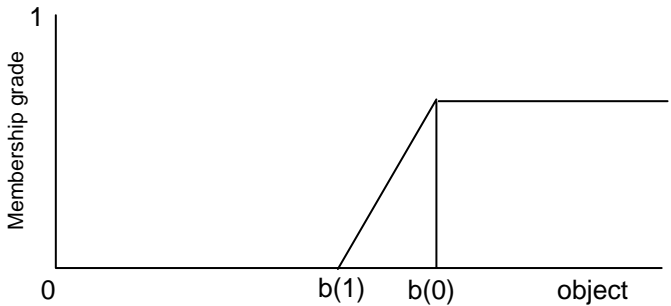

Figure: 2.10.4

These are the results generated by the fuzzy linear programming. Since both $\mu_N$ and $\mu_z$ have been characterized, now our goal is to describe the appropriateness of any operation state. Given any operating x, the feasibility can be specified based on the objective value, the constraints and the estimated uncertainty is got using equation (4). The value of $\mu_D$ are shown as the intersection of the two membership functions.

Defining the decision region based on the intersection we describe the variables and constraints of our problem. The variable $x_1$ represents the amount of gas fed to pressure swing adsorption from the gas production unit. The variable $x_2$ represents the amount of gas production that is sent to CGG. This problem can be represented according to equation (*). The constraints on the problem are subjected to some degree of uncertainty often some violation of the constraints within this range of uncertainty is tolerable. This problem can be represented according to equation (*). Using the given refinery data from the chemical plant.

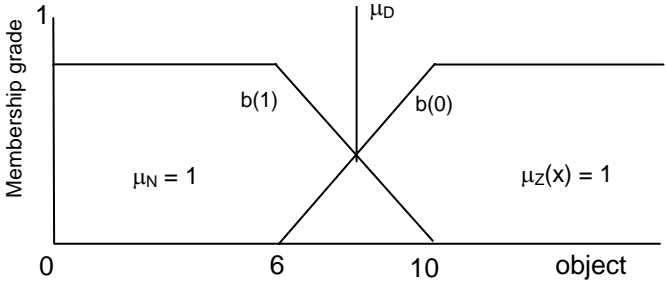

Figure: 2.10.5

$$c = [-0.544 \; 3]$$

$$A = \begin{bmatrix} 1 & 0 \\ 0 & 1 \\ -0.544 & 1 \end{bmatrix},$$

$$b = \begin{bmatrix} 33.652 \\ 23.050 \\ 4.743 \end{bmatrix}$$

Using equation (*) we get

$Zc = -0.544 \; x_1 + 3x_2$ it represents gas waste flow. The gas waste flow is represented by the following three equations:

i. $x_1 + 0x_2 \leq 33.652$ is the total dead – end waste flow gas.

ii. $0x_1 + x_2 = 23.050$ is the total (GCG2) gas consuming gas – treaters waste flow gas.

iii. $-0544 \; x_1 + x_2 \leq 4.743$ is the total let-down waste flow gas.

All flow rates are in million standard cubic feet per day. (i.e. 1 MMSCFD = 0.3277 m$^3$/s at STP). The value used for may be considered to be desired from operator experts opinion. The third constraint represents the minimum let-down flow receiving to keep valve from sticking. The value to this limit cannot be given an exact value, therefore a certain degree of violation may be tolerable. The other constraints may be subject to some uncertainty as well as they represent the maximum allowable values for $x_1$ and $x_2$. In this problem we are going to express all constrains in terms of $\alpha$, $\alpha$, $\varepsilon$ [0, 1]. We have to chose a value of tolerance on the third constraint as $p_3 = 0.1$, then this constraint is represented parametrically as

$$a_3 \, x \leq (b_3 - p_3) + 2p_3 \, (1 - \alpha).$$

For example, if we use crisp optimization problem with the tolerance value $p = 0.1$ we obtain the following result:

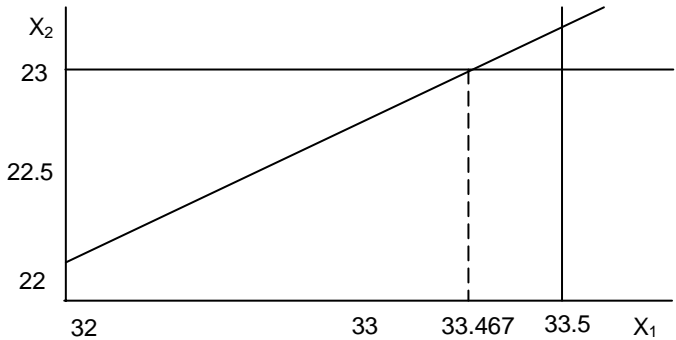

Figure: 2.10.6

where $x_1$ represents the amount of gas fed to PSA from gas production unit which is taken along the x axis, and $x_2$ amount of gas sent to CGG which is taken along the y axis,
we get $x_1 = 33.469$, when $x_2 = 23.050$

$$x^* = \begin{bmatrix} 33.469 \\ 23.050 \end{bmatrix}$$

$z = 50.941$. Finally we compare this result with our fuzzy linear programming method.

We replace two valued indicator function method by fuzzy linear programming.

Fuzzy Linear Programming is used now to maximize the objective function as well as minimize the uncertainty (waste flow gas). For that all of the constraints are expressed in terms of α, α, ∈ [0, 1].

$$a_3 \, x \leq (b_3 - p_3) + 2p_3 \, (1 - \alpha). \quad \alpha \in [0, 1]$$

where $a_3$ is the third row in the matrix A. i.e. $= 0.544x_1 + x_2 \leq 4.843 - 0.2 \, \alpha$, when the tolerance $p_3 = 0.3$, we fix the value of α ε [0.9,1], when the tolerance $p_3 = 0.1$, we see α ε [0.300, 0.600].

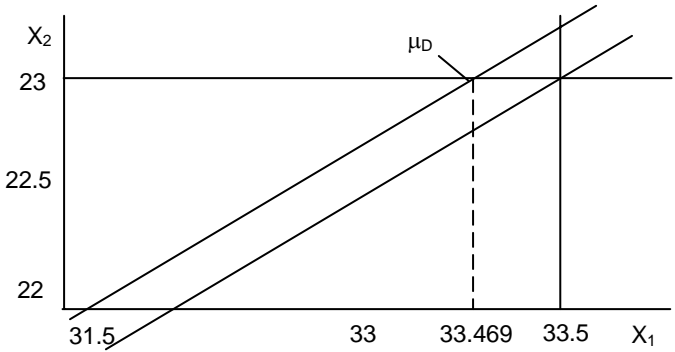

Figure: 2.10.7

where $x_1$ represents the amount of gas fed to PSA from gas production unit which is taken along the x axis, and $x_2$ amount of gas that is sent to CGG which is taken along the y axis,

When $x_2 = 23.050$ and     α = 0.0, we get   $x_1 = 33.469$.
When $x_2 = 23.050$ and     α = 0.4, we get   $x_1 = 33.616$

The set $(\mu_z)$ is defined in equation 5. Fuzzy Linear Programming solution is

$$x^* = f(\alpha) = \begin{bmatrix} 33.469 \\ 23.050 \end{bmatrix}$$

this value is recommended as there is no changes in the operating policy.

So we have to chose the value for α as 0.6 for the tolerance $p_3 = 0.1$, we get the following graph where $x_1$ represents the amount of gas fed to PSA from gas production unit which is taken along the x axis, and $x_2$ amount of gas that sent to CGG which is taken along the y axis,

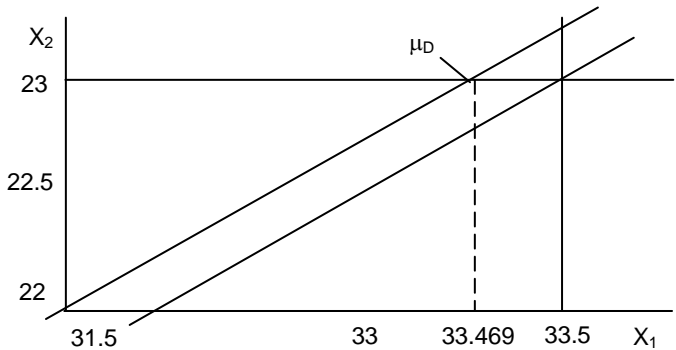

Figure: 2.10.8

when $x_2 = 23.050$ and $\alpha = 0.0$ we get $x_1 = 33.469$
when $x_2 = 23.050$ and $\alpha = 0.6$ we get $x_1 = 33.689$.

The operating region

$$x^* = f(0.6) = \begin{bmatrix} 33.689 \\ 23.050 \end{bmatrix}.$$

Now if the tolerance on the third constraint is increased to $p_3 = 0.2$. This results is the region shown in the following graph. As expected the region has increased to allow a larger range of operating states.

when $x_2 = 23.050$ and $\alpha = 0.0$ we get     $x_1 = 33.285$
when $x_2 = 23.050$ and $\alpha = 0.9$ we get     $x_1 = 33.947$.

The operating region is

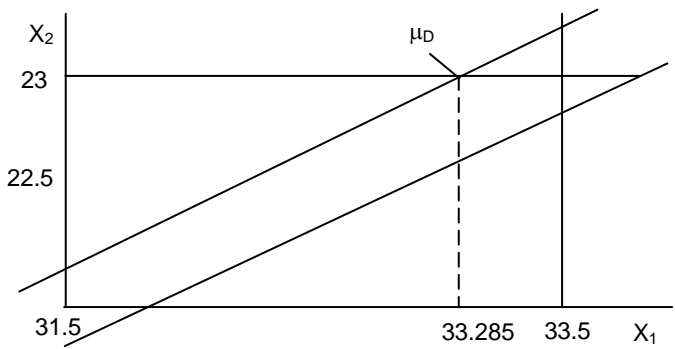

Figure: 2.10.9

where $x_1$ represents the amount of gas fed to PSA from gas production unit which is taken along the x axis, and $x_2$ amount of gas that is sent to CGG which is taken along the y axis.

$$x^* = f(0.9) = \begin{bmatrix} 33.947 \\ 23.050 \end{bmatrix}.$$

The fuzzy linear programming solution is

$$x^* = f(\alpha) = \begin{bmatrix} 33.285 \\ 23.050 \end{bmatrix}$$

$$z^* = 51.043.$$

Finally we have to take $\alpha \varepsilon [0.9, 1.00]$.

Choose $\alpha = 0$ and when the tolerance $p_3 = 0.3$ we get the following graph when $x_2 = 23.050$ we get $x_1 = 33.101$.

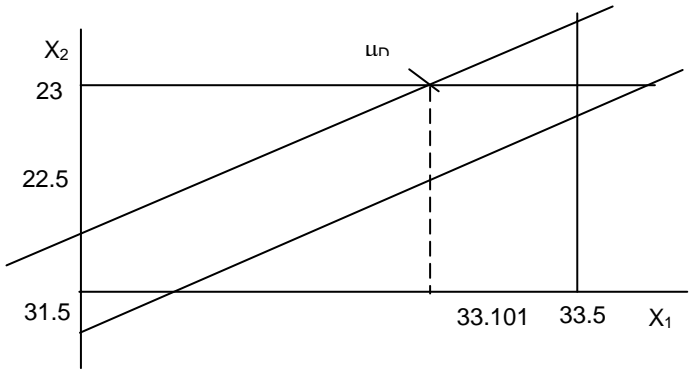

Figure: 2.10.10

where $x_1$ represents the amount of gas fed to PSA from gas production unit; and $x_2$ amount of gas that is sent to CGG.

When $\alpha = 1$ and $x_2 = 23.050$ we get $x_1 = 34.204$. The operating region is

$$x^* = f(1.0) = \begin{bmatrix} 33.204 \\ 23.050 \end{bmatrix}.$$

The fuzzy linear programming solutions are

$$x^* = f(\alpha) = \begin{bmatrix} 33.101 \\ 23.050 \end{bmatrix}.$$

The fuzzy linear programming solutions are

$$z^* = g(\alpha) = 51.143.$$

We chose maximum value from the Fuzzy Linear Programming method i.e. $z^* = 51.143$.

Thus when we work by giving varying membership functions and use fuzzy linear programming we see that we get the minimized waste gas flow value as 33.101 in contrast to 33.464 measured in million standard cubic feet per day and the maximum gas waste flow of system is determined to be 51.143 in contrast to their result of 50.941 measured in million standard cubic feet per day. Since the difference we have obtained is certainly significant, this study when applied to any oil refinery will minimize the waste gas flow to atmosphere considerably and reduce the pollution.

## 2.11 Semantics of implication operators and fuzzy relational products

Here, we analyze the data obtained from the HIV/AIDS patients in a meaningful and natural way by using fuzzy relational operators.

The result, of the assessment by expert a is a relation $R(a)$ given by a matrix of which the $(ij)^{th}$ component $R_{ij}(a)$ is the degree to which the symptom $C_i$ is attributed to the patient $P_j$. The inverse $R^{-1}(a)$ of this relation is a relation from patients to symptoms where $R_{jm}^{-1}(a)$ is the degree to which patient $P_j$ was considered to exemplify symptom $C_m$. Given two such relations $R^{-1}(a)$ and $R(a')$ where expert a' may or may not be the same as a, we can make two extremely interesting comparisons by forming two triangular products of the relations.

The first of these gives us a relation from patients to patients defined as follows: the relation $U(a, a')=R^{-1}(a) \theta R(a')$ (Note: we use the symbol $\theta$ to represent the triangular product of the relation) has for its $(jm)^{th}$ component $U_{jm}(a.a')$ the degree to which the attribution (by expert a) of symptoms to $P_j$ implies their attribution (by expert a') to $P_m$.

To understand this product, suppose assume that the original data is crisp, that is to say binary. Then the attribution by a to $P_j$ of symptom $C_k$ implies its attribution by a' to $P_m$, to the degree 0 or 1 given by the classical table (table 1) for the material implication $R_{jk}^{-1}(a) \rightarrow R_{km}(a')$ namely

Table 2.11.1

| $R_{jk}^{-1}(a)/R_{km}(a')$ | 0 | 1 |
|---|---|---|
| 0 | 1 | 1 |
| 1 | 0 | 1 |

The mean value of this over the number of symptoms is plausibly to be taken as the degree we seek that is, the degree to which attributions by expert a to $P_j$ imply attributions by expert a' to $P_m$. This idea is embodied in the formula

$$U_{jm}(a.a') = (R^{-1}(a) \; \theta \; R(a'))_{jm}$$
$$= (1/N_k) \Sigma_k (R_{jk}^{-1}(a) \rightarrow R_{km}(a')). \qquad (1)$$

After finding the Relational Product matrix U, we find the various α-cuts of this fuzzy relation U. (Note: For a given fuzzy set A defined on a Universal set X and any number $\alpha \in [0,1]$, the α-cut of A is defined as $^{\alpha}A = \{x / A(x) \geq \alpha\}$).

In our model, we denote $P_i$ as the $i^{th}$ patient i=1,2,…,10 and $C_j$ as the $j^{th}$ symptom j=1,2,…,8. we give the raw data in a tabular form (table 2). From the raw data we construct the Relational matrix R(a) where the row corresponds to the symptoms namely disabled ($C_1$), difficult to cope with ($C_2$), dependent ($C_3$), apathetic and unconcerned ($C_4$), blaming oneself ($C_5$), very ill ($C_6$), depressed ($C_7$), and anxious and worried ($C_8$) and the column corresponds to the patients. Also for our study the expert a is same as expert a'.

Table 2.11.2

| $C_j / P_i$ | $P_1$ | $P_2$ | $P_3$ | $P_4$ | $P_5$ | $P_6$ | $P_7$ | $P_8$ | $P_9$ | $P_{10}$ |
|---|---|---|---|---|---|---|---|---|---|---|
| Disabled | 0 | 0 | 1 | 0 | 0 | 0 | 0 | 1 | 0 | 0 |
| Difficult to cope | 0 | 0 | 0 | 0 | 0 | 1 | 0 | 0 | 0 | 0 |

| with | | | | | | | | | | |
|---|---|---|---|---|---|---|---|---|---|---|
| Dependent | 0 | 0 | 0 | 1 | 0 | 0 | 0 | 0 | 0 | 0 |
| Apathetic and Unconcerned | 1 | 1 | 1 | 0 | 1 | 0 | 1 | 0 | 1 | 0 |
| Blaming oneself | 0 | 0 | 1 | 0 | 0 | 0 | 1 | 0 | 0 | 0 |
| Very ill | 0 | 0 | 0 | 1 | 1 | 0 | 0 | 0 | 0 | 1 |
| Depressed | 1 | 0 | 1 | 0 | 1 | 0 | 0 | 1 | 1 | 0 |
| Anxious and worried | 1 | 1 | 0 | 1 | 0 | 0 | 1 | 0 | 0 | 0 |

The Relational matrix R(a) is as follows:

$$
\begin{pmatrix}
0 & 0 & 1 & 0 & 0 & 0 & 0 & 1 & 0 & 0 \\
0 & 0 & 0 & 0 & 0 & 1 & 0 & 0 & 0 & 0 \\
0 & 0 & 0 & 1 & 0 & 0 & 0 & 0 & 0 & 0 \\
1 & 1 & 1 & 0 & 1 & 0 & 1 & 0 & 1 & 0 \\
0 & 0 & 1 & 0 & 0 & 0 & 1 & 0 & 0 & 0 \\
0 & 0 & 0 & 1 & 1 & 0 & 0 & 0 & 0 & 1 \\
1 & 0 & 1 & 0 & 1 & 0 & 0 & 1 & 1 & 0 \\
1 & 1 & 0 & 1 & 0 & 0 & 1 & 0 & 0 & 0
\end{pmatrix}
$$

Now we apply the rule in equation (1) and we get the following product matrix $U = R^{-1}(a) \; \theta \; R(a')$ :

$$
\begin{pmatrix}
1 & .88 & .88 & .75 & .88 & .63 & .88 & .75 & .88 & .63 \\
1 & 1 & .88 & .88 & .88 & .75 & 1 & .75 & .88 & .75 \\
.75 & .63 & 1 & .5 & .75 & .5 & .75 & .75 & .75 & .63 \\
.75 & .75 & .63 & 1 & .75 & .63 & .75 & .63 & .63 & .75 \\
.88 & .75 & .88 & .75 & 1 & .63 & .75 & .75 & .88 & .75 \\
.88 & .88 & .88 & .88 & .88 & 1 & .88 & .88 & .88 & .88 \\
.88 & .88 & .88 & .75 & .75 & .63 & 1 & .63 & .75 & .63 \\
.88 & .75 & 1 & .75 & .88 & .75 & .75 & 1 & .88 & .75 \\
1 & .88 & 1 & .75 & 1 & .75 & .88 & .88 & 1 & .75 \\
.88 & .88 & .88 & 1 & 1 & .88 & .88 & .88 & .88 & 1
\end{pmatrix}
$$

α-cut of U for α = 1.

$$\begin{pmatrix} 1 & 0 & 0 & 0 & 0 & 0 & 0 & 0 & 0 & 0 \\ 1 & 1 & 0 & 0 & 0 & 0 & 1 & 0 & 0 & 0 \\ 0 & 0 & 1 & 0 & 0 & 0 & 0 & 0 & 0 & 0 \\ 0 & 0 & 0 & 1 & 0 & 0 & 0 & 0 & 0 & 0 \\ 0 & 0 & 0 & 0 & 1 & 0 & 0 & 0 & 0 & 0 \\ 0 & 0 & 0 & 0 & 0 & 1 & 0 & 0 & 0 & 0 \\ 0 & 0 & 0 & 0 & 0 & 0 & 1 & 0 & 0 & 0 \\ 0 & 0 & 1 & 0 & 0 & 0 & 0 & 1 & 0 & 0 \\ 1 & 0 & 1 & 0 & 1 & 0 & 0 & 0 & 1 & 0 \\ 0 & 0 & 0 & 1 & 1 & 0 & 0 & 0 & 0 & 1 \end{pmatrix}$$

α -cut of U for α =0.88

$$\begin{pmatrix} 1 & 1 & 1 & 0 & 1 & 0 & 1 & 0 & 1 & 0 \\ 1 & 1 & 1 & 1 & 1 & 0 & 1 & 0 & 1 & 0 \\ 0 & 0 & 1 & 0 & 0 & 0 & 0 & 0 & 0 & 0 \\ 0 & 0 & 0 & 1 & 0 & 0 & 0 & 0 & 0 & 0 \\ 1 & 0 & 1 & 0 & 1 & 0 & 0 & 0 & 1 & 0 \\ 1 & 1 & 1 & 1 & 1 & 1 & 1 & 1 & 1 & 1 \\ 1 & 1 & 1 & 0 & 0 & 0 & 1 & 0 & 0 & 0 \\ 1 & 0 & 1 & 0 & 1 & 0 & 0 & 1 & 1 & 0 \\ 1 & 1 & 1 & 0 & 1 & 0 & 1 & 1 & 1 & 0 \\ 1 & 1 & 1 & 1 & 1 & 1 & 1 & 1 & 1 & 1 \end{pmatrix}$$

From the matrix U for α =1, we observe that the patient $P_2$ is related to the patients $P_1$ and $P_7$ i.e., all the symptoms which are attributed to $P_2$ ($C_4$ and $C_8$) is found in both the patients $P_1$ and $P_7$. So if we want to take steps to improve the health condition (both physical and mental) of the patients $P_1$ and $P_2$ (or $P_2$ and $P_7$), we may find the solutions for the problems of the patient $P_1$ (or $P_7$) then the problems of the patient $P_2$ will be solved because the symptoms which are found in patient $P_2$ is included in the symptoms which are found in the patients $P_1$ and $P_7$.

Similarly the patient $P_9$ is related to the patients $P_1$, $P_3$, and $P_5$. So all the symptoms which are found in patient $P_9$ ($C_4$ and $C_7$) is also found in the patients $P_1$, $P_3$, and $P_5$. So we can say that the patient $P_9$ is come under the category of the patients $P_1$, $P_3$, and $P_5$. In a similar way we can divide all the patients in to different smaller categories so that we can study the problems of the patients clearly and properly.

When we look the matrix U for $\alpha$ =0.88, we observe that the patient $P_1$ is related to the patients $P_2$, $P_3$, $P_5$, $P_7$, and $P_9$. So the symptoms which are found in patient $P_1$ is also found in the patients $P_2$, $P_3$, $P_5$, $P_7$, and $P_9$ with the degree of possibility 0.88. We also observe that the patients $P_6$ and $P_{10}$ are related to every other patient. But we have to note that here the degree of the possibility is only 0.88 so we cannot say that the patients are perfectly related but in the previous case i.e., for $\alpha$ =1 the relation between the patients are perfect. Similar analysis is carried out for all other patients.

Here, we give a relation from symptoms to symptoms which is a triangle product defined as:

$$(R(a)\ \theta\ R^{-1}(a'))_{ik} = (1/N_j)\ \Sigma_j\ (R_{ij}(a) \rightarrow R_{jk}^{-1}(a')).$$

The only difference in this relation is that, when we compare it with the previous one the implication is in the direction $R_{jk}^{-1}(a) \rightarrow R_{km}(a')$ but now the implication is in the reverse order so we can find the triangular product $R(a) \rightarrow R^{-1}(a')$ in a similar way as before by just replacing the position of the given matrix and its transpose.

The Relational product matrix $V = R(a)\ \theta\ R^{-1}(a')$ for this case is follows:

$$\begin{pmatrix}
1 & .8 & .8 & .9 & .9 & .8 & 1 & .8 \\
.9 & 1 & .9 & .9 & .9 & .9 & .9 & .9 \\
.9 & .9 & 1 & .9 & .9 & 1 & .9 & 1 \\
.5 & .4 & .4 & 1 & .6 & .5 & .8 & .7 \\
.9 & .8 & .8 & 1 & 1 & .8 & .9 & .9 \\
.7 & .7 & .8 & .8 & .7 & 1 & .8 & .8 \\
.7 & .5 & .5 & .9 & .6 & .6 & 1 & .6 \\
.6 & .6 & .7 & .9 & .7 & .7 & .7 & 1
\end{pmatrix}$$

$\alpha$ -cut of V for $\alpha$ =1

$$\begin{pmatrix} 1 & 0 & 0 & 0 & 0 & 0 & 1 & 0 \\ 0 & 1 & 0 & 0 & 0 & 0 & 0 & 0 \\ 0 & 0 & 1 & 0 & 0 & 1 & 0 & 1 \\ 0 & 0 & 0 & 1 & 0 & 0 & 0 & 0 \\ 0 & 0 & 0 & 1 & 1 & 0 & 0 & 0 \\ 0 & 0 & 0 & 0 & 0 & 1 & 0 & 0 \\ 0 & 0 & 0 & 0 & 0 & 0 & 1 & 0 \\ 0 & 0 & 0 & 0 & 0 & 0 & 0 & 1 \end{pmatrix}$$

$\alpha$ -cut of V for $\alpha = 0.9$

$$\begin{pmatrix} 1 & 0 & 0 & 1 & 1 & 0 & 1 & 0 \\ 1 & 1 & 1 & 1 & 1 & 1 & 1 & 1 \\ 1 & 1 & 1 & 1 & 1 & 1 & 1 & 1 \\ 0 & 0 & 0 & 1 & 0 & 0 & 0 & 0 \\ 1 & 0 & 0 & 1 & 1 & 0 & 1 & 1 \\ 0 & 0 & 0 & 0 & 0 & 1 & 0 & 0 \\ 0 & 0 & 0 & 1 & 0 & 0 & 1 & 0 \\ 0 & 0 & 0 & 1 & 0 & 0 & 0 & 1 \end{pmatrix}$$

From the matrix V for $\alpha = 1$, it is clear that the symptom very ill ($C_1$) is related to the symptom depressed ($C_7$), so we can say that "the patients who are very ill are depressed" i.e., if a patient having the symptom $C_1$ then the patient must have the symptom $C_7$. Similarly the symptom dependent ($C_3$) is related to the symptoms very ill ($C_6$) and anxious and worried ($C_8$). So we can say that a patient having the symptom $C_3$ must have the symptoms $C_6$ and $C_8$. Similarly the symptom blaming oneself ($C_5$) is related to the symptom apathetic and unconcerned ($C_4$), so we can say that "the patients who are blaming oneself are apathetic and unconcerned".

Also from the matrix V for $\alpha = 0.9$, we observe that both the symptoms difficult to cope with ($C_2$) and dependent ($C_3$) are related to all other symptoms. So a patient having any one of the symptoms $C_2$ and $C_3$ then the patient must have all the other symptoms, but here we should note that the degree of the possibility is only 0.9. Also the symptom disabled ($C_1$) is related to the symptoms apathetic and unconcerned ($C_4$), blaming oneself ($C_5$) and depressed ($C_7$). So we can say that "the patients who are disabled are apathetic and unconcerned, blaming oneself and depressed" with the degree of possibility 0.9. In a similar way we can interpret the other symptoms.

Here, we analyze the data in a very different way. First, we explain the method and then we illustrate it.

Two observers a and a' (who may or may not be in fact the same individual), use checklist (i.e., a list where we enter the data of a particular person) for two persons (patients) $P_j$ and $P_m$ (j = m or j ≠ m) as follows:

<div align="center">
a uses the list on $P_j$<br>
a' uses the list on $P_m$
</div>

$\alpha_{vw}$ = the number of items (symptoms), which a marks v for $P_j$ and a' marks w for $P_m$, where v,w ∈ {0,1}.

The "contingency table" is shown in Table 2.11.3:

<div align="center">Table 2.11.3</div>

| $P_j$ / $P_m$ | 0 | 1 | |
|---|---|---|---|
| 0 | $\alpha_{00}$ | $\alpha_{01}$ | $\alpha_{00} + \alpha_{01} = r_0$ |
| 1 | $\alpha_{10}$ | $\alpha_{11}$ | $\alpha_{10} + \alpha_{11} = r_1$ |
| | $\alpha_{00} + \alpha_{10} = c_0$ | $\alpha_{01} + \alpha_{11} = c_1$ | $r_0 + r_1 = c_0 + c_1 = n$ |

Now we define $x = r_1/n$ and $y = c_1/n$ and we also define $x_{ij}$ = value of x corresponding to the contingency table for $(ij)^{th}$ patient. Similarly we define $y_{ij}$ = value of y corresponding to the contingency table for $(ij)^{th}$ patient. Now we define a new implication operator namely Kleene-Dienes operator in such a way that

$$x \rightarrow y = \max (1 - x, y).$$

Now we define a relation between the patient i and patient j in such a way that $P_i \rightarrow P_j = x_{ij} \rightarrow y_{ij} = \max (1 - x_{ij}, y_{ij})$. After finding the value of $P_i \rightarrow P_j$ for all i,j we can construct the relational matrix W(a, a') as we have formed in section 2.

In our problem, we take five patients and the expert's opinion for these patients is shown below in the form of a table (table 4) where the row corresponds to the symptoms namely very ill ($C_1$), apathetic and unconcerned ($C_2$), depressed ($C_3$), anxious and worried ($C_4$), and disabled ($C_5$) and column corresponds to the patients:

Table 2.11.4

|  | $P_1$ | $P_2$ | $P_3$ | $P_4$ | $P_5$ |
|---|---|---|---|---|---|
| $C_1$ | 0 | 0 | 0 | 0 | 1 |
| $C_2$ | 1 | 1 | 1 | 1 | 1 |
| $C_3$ | 1 | 1 | 0 | 0 | 1 |
| $C_4$ | 1 | 0 | 1 | 1 | 0 |
| $C_5$ | 0 | 1 | 0 | 0 | 0 |

The number '0' in the first column of the first row represents the symptom $C_1$ (very ill) is not found in patient 1 but the first element '1' in the second row represents the symptom $C_2$ (apathetic and unconcerned) is found in patient 1. In a similar way we can interpret the other values. Now using the above table, we construct twenty five different contingency tables for these five patients. For example the contingency table for the pair ($P_1$,$P_3$) is:

| $P_1 / P_3$ | 0 | 1 |
|---|---|---|
| 0 | 2 | 0 |
| 1 | 1 | 2 |

when we compare the above table with table 3, we observe that the corresponding value for $\alpha_{00}$ in this table is 2. Since $\alpha_{00}=2$, the 0number of symptoms which expert marks 0 for $P_1$ and for $P_3$ is 2 i.e., two of the five symptoms are not found in both the patients namely $C_1$ and $C_5$.

Now we have to find the values of $x_{ij}$ and $y_{ij}$ for every pair ($P_i$,$P_j$) i,j=1,2,…5. Using these values and using the Kleene-Dienes operator we can find $P_i \to P_j$ for all i,j in such a way that

$$P_i \to P_j = \max(1 - x_{ij}, y_{ij}).$$

For example, form the contingency table for ($P_1$,$P_3$), we can find that $r_1=3, c_1=2, x_{ij}=r_1/5=0.6$, and $y_{ij}=c_1/5=0.4$. Then by the Kleene-Dienes operator $P_1 \to P_3=0.4$. In a similar way we can find the values for all other pairs ($P_i$, $P_j$). So the Relational product matrix W for our case is as follows:

$$\begin{pmatrix} .6 & .6 & .4 & .4 & .6 \\ .6 & .6 & .4 & .4 & .6 \\ .6 & .6 & .6 & .6 & .6 \\ .6 & .6 & .6 & .6 & .6 \\ .6 & .6 & .4 & .4 & .6 \end{pmatrix}$$

$\alpha$ -cut of W for $\alpha = 0.6$

$$\begin{pmatrix} 1 & 1 & 0 & 0 & 1 \\ 1 & 1 & 0 & 0 & 1 \\ 1 & 1 & 1 & 1 & 1 \\ 1 & 1 & 1 & 1 & 1 \\ 1 & 1 & 0 & 0 & 1 \end{pmatrix}$$

From the above matrix, we observe that the patient $P_3$ is related to everyone i.e., the symptoms found in $P_3$ are also found in all other persons. Also each patient is related to the patient $P_2$ but the patient $P_2$ is related only to the patients $P_1$ and $P_5$ i.e., all the symptoms found in each patient also found in $P_2$ but the symptoms found in patient $P_2$ is found only in the patients $P_1$ and $P_5$ with the degree of possibility 0.6. Also the patient $P_4$ is related to everyone but none of the patients is related to the patient $P_4$ (except $P_3$) with the degree of possibility 0.6. In a similar way we can interpret the other patients.

Here we construct the contingency tables for the symptoms by using the table 4 in the section 4, for example the contingency table for the pair $(C_4, C_3)$ is as follows:

| $C_4 / C_3$ | 0 | 1 |
|---|---|---|
| 0 | 1 | 1 |
| 1 | 2 | 1 |

when we compare the above table with table 3, we observe that the corresponding value for $\alpha_{10}$ in this table is 2. Since $\alpha_{10}=2$, out of five patients two of them having the symptom anxious and worried ($C_4$) but not the symptom depressed ($C_3$) is 2.

Now by the same method as we have seen in section 4 we can find the values of $x_{ij}$ and $y_{ij}$ for every pair $(C_i, C_j)$ i,j=1,2,…5. Using these values and using the Kleene-Dienes operator we can find $C_i \rightarrow C_j$ for all i,j in such a way that

$$C_i \rightarrow C_j = \max(1-x_{ij}, y_{ij}).$$

So the Relational product matrix W in this case is as follows:

$$\begin{pmatrix} .8 & 1 & .8 & .8 & .8 \\ .2 & 1 & .6 & .6 & .2 \\ .4 & 1 & .6 & .6 & .4 \\ .4 & 1 & .6 & .6 & .4 \\ .8 & 1 & .8 & .8 & .8 \end{pmatrix}$$

$\alpha$-cut of W for $\alpha = 1$

$$\begin{pmatrix} 0 & 1 & 0 & 0 & 0 \\ 0 & 1 & 0 & 0 & 0 \\ 0 & 1 & 0 & 0 & 0 \\ 0 & 1 & 0 & 0 & 0 \\ 0 & 1 & 0 & 0 & 0 \end{pmatrix}$$

$\alpha$-cut of W for $\alpha = 0.8$

$$\begin{pmatrix} 1 & 1 & 1 & 1 & 1 \\ 0 & 1 & 0 & 0 & 0 \\ 0 & 1 & 0 & 0 & 0 \\ 0 & 1 & 0 & 0 & 0 \\ 1 & 1 & 1 & 1 & 1 \end{pmatrix}$$

From $\alpha$-cut of W for $\alpha = 1$, we observe that every symptom is related to the symptom apathetic and unconcerned ($C_2$) i.e., a patient having any one of the symptom $C_i$ i=1,2,…5 should have the symptom $C_2$. Since all other entries in the matrix W for $\alpha = 1$ are zero so we cannot find a perfect relation between other symptoms. But from the matrix W for $\alpha = 0.8$, we observe that both the symptoms very ill ($C_1$) and disabled ($C_5$) are related to all other symptoms with the degree of possibility 0.8 i.e., if a patient having any one of the symptoms $C_1$ and $C_5$ then the patient should have all other symptoms with the degree of possibility 0.8. So if we decrease the value of α then we get more relations but the degree of perfection also decreases.

To assign a measure to the degree to which expert say "yes" to items in the checklist for $P_j$ implies expert's saying "yes" to these same items for $P_m$: in briefer words, a measure of the support these fine data give to the statement "if yes-j then yes-m".

In classical logic, "if yes-j then yes-m" is satisfied "by performance" whenever yes-j and yes-m occur together, and "by default" whenever no-j occurs, regardless of the m-answer. Thus all entries in the contingency table (table 2) support the statement except for $\alpha_{10}$. Thus if we weight all items equally, the appropriate classical measure of support for the assertion is m = 1 − ($\alpha_{10}$ / n).

For example, the contingency table corresponding to the pair $(P_2, P_3)$ is:

| $P_2 / P_3$ | 0 | 1 |
|---|---|---|
| 0 | 1 | 1 |
| 1 | 2 | 1 |

So the classical measure of support for this pair is m=1-(2/5) = 0.6. So we can interpret that the symptoms which are present (not present) in patient $P_2$ is present (not present) in patient $P_3$ (except present in $P_2$ and not in $P_3$) with the degree of possibility 0.6.

This is the only method, which can help the expert compare the symptoms in patients and inter-relate the patients. Further some of these methods give the interdependence between the symptoms in patients. (This research work would be appearing with T.Johnson).

Chapter Three

# BASIC NOTIONS AND NEW CONCEPTS ON NEUTROSOPHY

The analysis of most of the real world problems involves the concept of indeterminacy. Here one cannot establish or cannot rule out the possibility of some relation but says that cannot determine the relation or link; this happens in legal field, medical diagnosis even in the construction of chemical flow in industries and more chiefly in socio economic problems prevailing in various countries. So this chapter defines new concepts which paves way for the building of the Neutrosophic relational equation (NRE). This chapter has 6 sections each dealing neutrosophic concepts. First section gives the basic concepts defined by [87-90, 57].

## 3.1 Neutrosophic set

In this section we just recall the concepts about neutrosophic sets given by [87-90].

**DEFINITION 3.1.1:** *Let T, I, F be real standard or non-standard subsets of [ ⁻0, 1⁺ ],*

*with    sup T = t_sup, inf T = t_inf,*
*           sup I = i_sup, inf I = i_inf,*
*           sup F = f_sup, inf F = f_inf,*

*and    n_sup = t_sup+i_sup+f_sup,*
*           n_inf = t_inf+i_inf+f_inf.*

Let U be a universe of discourse, and M a set included in U. An element x from U is noted with respect to the set M as x(T, I, F) and belongs to M in the following way:
    It is t% true in the set, i% indeterminate (unknown if it is) in the set, and f% false, where t varies in T, i varies in I, f varies in F.

Statically T, I, F are subsets, but dynamically T, I, F are functions / operators depending on many known or unknown parameters.

*Example 3.1.1:* Let A and B be two neutrosophic sets. One can say, by language abuse, that any element neutrosophically belongs to any set, due to the percentages of truth/indeterminacy/falsity involved, which varies between 0 and 1 or even less than 0 or greater than 1.

Thus: x (50, 20, 30) belongs to A (which means, with a probability of 50% x is in A, with a probability of 30% x is not in A, and the rest is undecidable); or y (0, 0, 100) belongs to A (which normally means y is not for sure in A); or z (0, 100, 0) belongs to A (which means one does know absolutely nothing about z's affiliation with A).

More general, x ((20-30), (40-45) 4[50-51], {20, 24, 28}) belongs to the set A, which means:

1. with a probability in between 20-30% x is in A (one cannot find an exact approximate because of various sources used);
2. with a probability of 20% or 24% or 28% x is not in A;
3. the indeterminacy related to the appurtenance of x to A is in between 40-45% or between 50-51% (limits included).

The subsets representing the appurtenance, indeterminacy, and falsity may overlap, and n_sup = 30+51+28 > 100 in this case.

**Physics Examples**

a) For example the Schrodinger's Theory says that the quantum state of a photon can basically be in more than one place in the same time, which translated to the neutrosophic set means that an element (quantum state) belongs and does not belong to a set (one place) in the same time; or an element (quantum state) belongs to two different sets (two different places) in the same time. It is a question of "alternative worlds" theory very well represented by the neutrosophic set theory.

In Schroedinger's Equation on the behavior of electromagnetic waves and "matter waves" in quantum theory, the wave function Psi which describes the superposition of possible states may be simulated by a neutrosophic function, i.e. a function whose values are not unique for each argument from the domain of definition (the vertical line test fails, intersecting the graph in more points). Don't we better describe, using the attribute "neutrosophic" than "fuzzy" or any others, a quantum particle that neither exists nor non-exists?

b) How to describe a particle ∽ in the infinite micro-universe that belongs to two distinct places $P_1$ and $P_2$ in the same time? ∽ χ $P_1$ and ∽ ϖ $P_1$ as a true contradiction, or ∽ χ $P_1$ and ∽ χ ·$P_1$.

**Philosophical Examples**
Or, how to calculate the truth-value of Zen (in Japanese) / Chan (in Chinese) doctrine philosophical proposition: the present is eternal and comprises in itself the past and the future?

In Eastern Philosophy the contradictory utterances form the core of the Taoism and Zen/Chan (which emerged from Buddhism and Taoism) doctrines.

How to judge the truth-value of a metaphor, or of an ambiguous statement, or of a social phenomenon which is positive from a standpoint and negative from another standpoint?

There are many ways to construct them, in terms of the practical problem we need to simulate or approach. Below there are mentioned the easiest ones:

**Application**
A cloud is a neutrosophic set, because its borders are ambiguous, and each element (water drop) belongs with a neutrosophic probability to the set. (e.g. there are a kind of separated water drops, around a compact mass of water drops, that we don't know how to consider them: in or out of the cloud).

Also, we are not sure where the cloud ends nor where it begins, neither if some elements are or are not in the set. That's why the percent of indeterminacy is required and the neutrosophic probability (using subsets - not numbers - as components) should be used for better modeling: it is a more organic, smooth, and especially accurate estimation. Indeterminacy is the zone of ignorance of a proposition's value, between truth and falsehood.

**Neutrosophic Set Operations**
One notes, with respect to the sets A and B over the universe U,
 x = x($T_1$, $I_1$, $F_1$) χ A and x = x($T_2$, $I_2$, $F_2$) χ B, by mentioning x's *neutrosophic probability appurtenance*.
And, similarly, y = y(T', I', F') χ B.

**1. Complement of A**
*If x( $T_1$, $I_1$, $F_1$ ) χ A,*
*then x( {1}0$T_1$, {1}0$I_1$, {1}0$F_1$ ) χ C(A).*

**2. Intersection**

*If x( $T_1$, $I_1$, $F_1$ ) $\chi$ A, x( $T_2$, $I_2$, $F_2$ ) $\chi$ B,*
*then x( $T_1 ?T_2$, $I_1 ?I_2$, $F_1 ?F_2$ ) $\chi$ A $\Im$ B.*

**3. Union**

*If x( $T_1$, $I_1$, $F_1$ ) $\chi$ A, x( $T_2$, $I_2$, $F_2$ ) $\chi$ B,*
*then x( $T_1/T_2 0 T_1 ?T_2$, $I_1/I_2 0 I_1 ?I_2$, $F_1/F_2 0 F_1 ?F_2$ ) $\chi$ A $\cup$ B.*

**4. Difference**

*If x( $T_1$, $I_1$, $F_1$ ) $\chi$ A, x( $T_2$, $I_2$, $F_2$ ) $\chi$ B,*
*then x( $T_1 0 T_1 ?T_2$, $I_1 0 I_1 ?I_2$, $F_1 0 F_1 ?F_2$ ) $\chi$ A \ B,*
*because A \ B = A $\Im$ C(B).*

**5. Cartesian Product**

*If x( $T_1$, $I_1$, $F_1$ ) $\chi$ A, y( T', I', F' ) $\chi$ B,*
*then ( x( $T_1$, $I_1$, $F_1$ ), y( T', I', F' ) ) $\chi$ A $\times$ B.*

**6. M is a subset of N**

*If x( $T_1$, $I_1$, $F_1$ ) $\chi$ M $\cup$ x( $T_2$, $I_2$, $F_2$ ) $\chi$ N,*
*where inf $T_1$ [ inf $T_2$, sup $T_1$ [ sup $T_2$, and inf $F_1$ $\mu$ inf $F_2$, sup $F_1$ $\mu$ sup $F_2$.*

We just recall the neutrosophic n-ary Relation:

**DEFINITION 3.1.2:** *Let $A_1$, $A_2$, …, $A_n$ be arbitrary non-empty sets. A Neutrosophic n-ary Relation R on $A_1 \times A_2 \times … \times A_n$ is defined as a subset of the Cartesian product $A_1 \times A_2 \times … \times A_n$, such that for each ordered n-tuple $(x_1, x_2, …, x_n)$(T, I, F), T represents the degree of validity, I the degree of indeterminacy, and F the degree of non-validity respectively of the relation R.*

From the intuitionistic logic, paraconsistent logic, dialetheism, fallibilism, paradoxes, pseudoparadoxes, and tautologies we transfer the "adjectives" to the sets, i.e. to intuitionistic set (set incompletely known), paraconsistent set, dialetheist set, faillibilist set (each element has a percentage of indeterminacy), paradoxist set (an element may belong and may not belong in the same time to the set), pseudoparadoxist set, and tautological set respectively.

Hence, the neutrosophic set generalizes:

- the *intuitionistic set*, which supports incomplete set theories (for 0 < n < 1, 0 [ t, i, f [ 1) and incomplete known elements belonging to a set;
- the *fuzzy set* (for n = 1 and i = 0, and 0 [ t, i, f [ 1);
- the *classical set* (for n = 1 and i = 0, with t, f either 0 or 1);
- the *paraconsistent set* (for n > 1, with all t, i , f < 1+);
- the *faillibilist set* (i > 0);
- the *dialetheist set*, a set M whose at least one of its elements also belongs to its complement C(M); thus, the intersection of some disjoint sets is not empty;
- the *paradoxist set* (t = f = 1);
- the *pseudoparadoxist set* (0 < i < 1, t = 1 and f > 0 or t > 0 and f = 1);
- the *tautological set* (i , f < 0).

Compared with all other types of sets, in the neutrosophic set each element has three components which are subsets (not numbers as in fuzzy set) and considers a subset, similarly to intuitionistic fuzzy set, of "indeterminacy" - due to unexpected parameters hidden in some sets, and let the superior limits of the components to even boil *over 1* (over flooded) and the inferior limits of the components to even freeze *under 0* (under dried).

For example: an element in some tautological sets may have t > 1, called "over included". Similarly, an element in a set may be "over indeterminate" (for i > 1, in some paradoxist sets), "over excluded" (for f > 1, in some unconditionally false appurtenances); or "under true" (for t < 0, in some unconditionally false appurtenances), "under indeterminate" (for i < 0, in some unconditionally true or false appurtenances), "under false" (for f < 0, in some unconditionally true appurtenances).

This is because we should make a distinction between unconditionally true (t > 1, and f < 0 or i < 0) and conditionally true appurtenances (t [ 1, and f [ 1 or i [ 1).

In a *rough set* RS, an element on its boundary-line cannot be classified neither as a member of RS nor of its complement with certainty. In the neutrosophic set a such element may be characterized by x(T, I, F), with corresponding set-values for T, I, F [⁻0, 1⁺].

One first presents the evolution of sets from fuzzy set to neutrosophic set. Then one introduces the neutrosophic components T, I, F which represent the membership, indeterminacy, and non-membership values respectively, where ]⁻0, 1⁺[ is the non-standard unit interval, and thus one defines the neutrosophic set. One gives examples from mathematics, physics, philosophy, and applications of the neutrosophic set. Afterwards, one introduces the neutrosophic set operations (complement, intersection, union, difference, Cartesian product, inclusion, and n-ary relationship), some generalizations and comments on them, and finally the distinctions between the neutrosophic set and the intuitionistic fuzzy set.

The *fuzzy set* (FS) was introduced by L. Zadeh in 1965, where each element had a degree of membership.

The *intuitionistic fuzzy set* (IFS) on a universe X was introduced by K. Atanassov in 1983 as a generalization of FS, where besides the degree of membership $\mu_A(x) \chi[0,1]$ of each element $x \chi X$ to a set A there was considered a degree of non-membership $\nu_A(x) \chi[0,1]$, but such that

$$\ldots x \chi X \ \mu_A(x) + \nu_A(x) \leq 1. \tag{1}$$

According to Deschrijver & Kerre the *vague set* defined by Gau and Buehrer was proven by Bustine & Burillo (1996) to be the same as IFS.

Atanassov defined the *interval-valued intuitionistic fuzzy set* (IVIFS) on a universe X as an object A such that:

$$A = \{(x, M_A(X), N_A(x)), x \chi X\}, \tag{2}$$

with

$$M_A : X \tau Int([0,1]) \text{ and } N_A : X \tau Int([0,1]) \tag{3}$$

and

$$\ldots x \chi X \ \sup M_A(x) + \sup N_A(x) \leq 1. \tag{4}$$

Belnap defined a four-valued logic, with truth (T), false (F), unknown (U), and contradiction (C). He used a bilattice where the four components were inter-related.

In 1995, starting from philosophy (when I fretted to distinguish between *absolute truth* and *relative truth* or between *absolute falsehood* and *relative falsehood* in logics, and respectively between *absolute membership* and *relative membership* or *absolute non-membership* and *relative non-membership* in set theory) I began to use the non-standard analysis. Also, inspired from the sport games (winning, defeating, or tight scores), from votes (pro, contra, null/black votes), from positive/negative/zero numbers, from yes/no/NA, from decision making and control theory (making a decision, not making, or hesitating), from accepted / rejected / pending, etc. and guided by the fact that the law of excluded middle did not work any longer in the modern logics, I combined the non-standard analysis with a tri-component logic/set/probability theory and with philosophy (I was excited by paradoxism in science and arts and letters, as well as by paraconsistency and incompleteness in knowledge). How to deal with all of them at once, is it possible to unity them?

[87-90] proposed the term "neutrosophic" because "neutrosophic" etymologically comes from "neutro-sophy" [French *neutre* < Latin *neuter*, neutral, and Greek *sophia*, skill/wisdom] which means knowledge of neutral thought, and this third/neutral represents the main distinction between "fuzzy" and "intuitionistic fuzzy" logic/set, i.e. the *included middle* component (Lupasco-Nicolescu's logic in philosophy), i.e. the neutral / indeterminate / unknown part (besides the "truth" / "membership" and "falsehood" / "non-membership" components that both appear in fuzzy logic/set). See the Proceedings of the First International Conference on Neutrosophic Logic, The University of New Mexico, Gallup Campus, 1-3 December 2001, at http://www.gallup.unm.edu/~smarandache/FirstNeutConf.htm.

We need to present these set operations in order to be able to introduce the neutrosophic connectors. Let $S_1$ and $S_2$ be two (unidimensional) real standard or non-standard subsets included in the non-standard interval $]^-0, \infty)$ then one defines:

1. *Addition of classical Sets:*
   $S_1/S_2 = \{x \xi x = s_1 + s_2, \text{ where } s_1 \chi S_1 \text{ and } s_2 \chi S_2\}$,
   with inf $S_1/S_2 =$ inf $S_1 +$ inf $S_2$, sup $S_1/S_2 =$ sup $S_1 +$ sup $S_2$;
   and, as some particular cases, we have
   $\{a\}/S_2 = \{x \xi x = a + s_2, \text{ where } s_2 \chi S_2\}$
   with inf $\{a\}/S_2 = a +$ inf $S_2$, sup $\{a\}/S_2 = a +$ sup $S_2$.

2. *Subtraction of classical Sets:*
   $S_1 0 S_2 = \{x \xi x = s_1 - s_2, \text{ where } s_1 \chi S_1 \text{ and } s_2 \chi S_2\}$.
   with inf $S_1 0 S_2 =$ inf $S_1 -$ sup $S_2$, sup $S_1 0 S_2 =$ sup $S_1 -$ inf $S_2$;
   and, as some particular cases, we have
   $\{a\} 0 S_2 = \{x \xi x = a - s_2, \text{ where } s_2 \chi S_2\}$,
   with inf $\{a\} 0 S_2 = a -$ sup $S_2$, sup $\{a\} 0 S_2 = a -$ inf $S_2$;
   also $\{1^+\} 0 S_2 = \{x \xi x = 1^+ - s_2, \text{ where } s_2 \chi S_2\}$,
   with inf $\{1^+\} 0 S_2 = 1^+ -$ sup $S_2$, sup $\{1^+\} 0 S_2 = 100 -$ inf $S_2$.

3. *Multiplication of classical Sets:*
   $S_1 ? S_2 = \{x \xi x = s_1 \exists s_2, \text{ where } s_1 \chi S_1 \text{ and } s_2 \chi S_2\}$.
   with inf $S_1 ? S_2 =$ inf $S_1 \exists$ inf $S_2$, sup $S_1 ? S_2 =$ sup $S_1 \exists$ sup $S_2$;

and, as some particular cases, we have

$\{a\}?S_2 = \{x\xi x=a\exists s_2$, where $s_2\chi S_2\}$,

with inf $\{a\}?S_2 = a * \inf S_2$, sup $\{a\}?S_2 = a \exists \sup S_2$;

also $\{1^+\}?S_2 = \{x\xi x=1\exists s_2$, where $s_2\chi S_2\}$,

with inf $\{1^+\}?S_2 = 1^+ \exists \inf S_2$, sup $\{1^+\}?S_2 = 1^+ \exists \sup S_2$.

4. *Division of a classical Set by a Number:*
   Let k $\chi^*$, then $S_1 2k = \{x\xi x=s_1/k$, where $s_1\chi S_1\}$.

Compared to *Belnap's quadruplet logic*, NS and NL do not use restrictions among the components – and that's why the NS/NL have a more general form, while the middle component in NS and NL (the indeterminacy) can be split in more subcomponents if necessarily in various applications.

**Differences between Neutrosophic Set (NS) and Intuitionistic Fuzzy Set (IFS)**

a) Neutrosophic Set can distinguish between *absolute membership* (i.e. membership in all possible worlds; we have extended Leibniz's absolute truth to absolute membership) and *relative membership* (membership in at least one world but not in all), because NS(absolute membership element)=$1^+$ while NS(relative membership element)=1. This has application in philosophy (see the neutrosophy). That's why the unitary standard interval [0, 1] used in IFS has been extended to the unitary non-standard interval $]^-0, 1^+[$ in NS.

Similar distinctions for *absolute or relative non-membership*, and *absolute or relative indeterminant appurtenance* are allowed in NS.

b) In NS there is no restriction on T, I, F other than they are subsets of $]^-0, 1^+[$, thus: $^-0 [ \inf T + \inf I + \inf F [ \sup T + \sup I + \sup F [ 3^+$.

This non-restriction allows paraconsistent, dialetheist, and incomplete information to be characterized in NS {i.e. the sum of all three components if they are defined as points, or sum of superior limits of all three components if they are defined as subsets can be >1 (for paraconsistent information coming from different sources), or < 1 for incomplete information}, while that information can not be described in IFS because in IFS the components T (membership), I (indeterminacy), F (non-membership) are restricted either to t+i+f=1 or to $t^2 + f^2 [ 1$, if T,

I, F are all reduced to the points t, i, f respectively, or to sup T + sup I + sup F = 1 if T, I, F are subsets of [0, 1].

Of course, there are cases when paraconsistent and incomplete information can be normalized to 1, but this procedure is not always suitable.

c) Relation (3) from interval-valued intuitionistic fuzzy set is relaxed in NS, i.e. the intervals do not necessarily belong to Int[0,1] but to [0,1], even more general to ]-0, 1+[.

d) In NS the components T, I, F can also be *non-standard* subsets included in the unitary non-standard interval ]⁻0, 1⁺[, not only *standard* subsets included in the unitary standard interval [0, 1] as in IFS.

e) NS, like dialetheism, can describe paradoxist elements, NS(paradoxist element) = (1, I, 1), while IFL can not describe a paradox because the sum of components should be 1 in IFS.

f) The connectors in IFS are defined with respect to T and F, i.e. membership and non-membership only (hence the Indeterminacy is what's left from 1), while in NS they can be defined with respect to any of them (no restriction).

g) Component "I", indeterminacy, can be split into more subcomponents in order to better catch the vague information we work with, and such, for example, one can get more accurate answers to the *Question-Answering Systems* initiated by Zadeh. {In Belnap's four-valued logic (1977) indeterminacy is split into Uncertainty (U) and Contradiction (C), but they were inter-related.}

## 3.2 Fuzzy neutrosophic sets

The notion of fuzzy neutrosophic sets are introduced in this section. The two types of fuzzy neutrosophic sets are crisp fuzzy neutrosophic sets and fuzzy neutrosophic sets both are introduced and some of their properties are discussed in this section.

Throughout this book *I* denotes the indeterminacy.

**DEFINITION 3.2.1:** *Let FN = [0, n I/ n ∈ [0, 1]] denotes the fuzzy interval of indeterminacy fuzzy interval of neutrosophy.*

**DEFINITION 3.2.2:** *Let X be a universal set; the map μ : X → {0, 1, I} is called the crisp fuzzy neutrosophic set, that is μ maps elements of X to the set 0, 1 and I.*

*μ(x) = 0; for x ∈ X implies x is a non member.*

$\mu(x) = 1$; for $x \in X$ implies x is a member.
$\mu(x) = I$; for $x \in X$ implies the membership of x is an indeterminacy.

We illustrate this by the following example:

**Example 3.2.1:** Let X be the set of all people in the age group from 1 to 25, which includes whites, blacks and browns. If $\mu$: X $\rightarrow \{0, 1, I\}$ such that only white have membership and black have no membership. Then the browns remain as indeterminate.

$\mu(x) = 0$; implies x is a black.
$\mu(x) = 1$; implies x is a white.
$\mu(x) = I$; implies of x is brown.

So the function or the map $\mu$ is a crisp neutrosophic set.

**DEFINITION 3.2.3:** *Let X be any set; FN = {[0, nI] | n $\in$ [0, 1]}. The map $\mu$: X $\rightarrow$ FN $\cup$ [0, 1] is said to be the fuzzy neutrosophic set of X; clearly $\mu(x)$ can belong to [0, 1] for some x $\in$ X and $\mu(x) = nI$ for n $\in$ [0, 1] for some other x. $\mu$ need not in general be a crisp fuzzy neutrosophic set of X.*

**Example 3.2.2:** X = {set of people suspected for doing a crime}
$\mu$: X $\rightarrow$ [0, 1] $\cup$ FN.

For every x $\in$ X the judge cannot say for certain that he had an hand in the crime for some x $\in$ X, the judge may say his part in performing or in participation in the crime is indeterminate he can also give degrees of indeterminacy of each x $\in$ X.

Thus fuzzy neutrosophic set gives a method by which the degrees of uncertainty is also measured. For instance if we say $\mu(x) = 0.2$ it implies that x has involved in the crime with 0.2 very less participation in that crime is accepted but if $\mu(x) = 0.2I$, it implies that his very involvement is doubtful or indeterminate and we suspect him; but we can not claim that he can be a criminal with the evidences produced. Thus we see the notion of fuzzy neutrosophic sets will certainty find its place in fields like giving judgment on some criminal cases, in medicine, and in any field where the concept of indeterminacy is involved.

Now we proceed on to define further algebraic structures on these crisp fuzzy neutrosophic sets and fuzzy neutrosophic sets.

**DEFINITION 3.2.4:** *Let* $\in$ *[F N* $\cup$ *[0, 1]] denote the set of all closed intervals of neutrosophic members in F N* $\cup$ *[0, 1]. Clearly* $\in$ *[F N* $\cup$ *[0, 1]]* $\subset$ *P [F N* $\cup$ *[0, 1]] where P [F N* $\cup$ *[0, 1]] denotes the set of all neutrosophic subsets and other subsets of FN* $\cup$ *[0, 1].*

The neutrosophic subsets of this type are called neutrosophic in travel valued neutrosophic set.

*Example 3.2.3:* Let $\zeta$[FN $\cup$ [0, 1]] be the set of neutrosophic interval valued neutrosophic subsets. [.03$I$, .7$I$] $\cup$ [.5, .7] is neutrosophic closed interval in $\zeta$[FN $\cup$ [0, 1]].

We define as in case of fuzzy sets the notion of the most important concepts viz. $\alpha$ - cut and its variant a strong $\alpha$ - cut in case of neutrosophic sets. We know given a fuzzy set A defined on X and any number $\alpha \in$ [0 1] the $\alpha$ cut, $^{\alpha}$A and the strong $\alpha$ - cut $^{\alpha+}$A are crisp set.

$$^{\alpha}A = \{x \mid A(x) \geq \alpha\}$$
$$^{\alpha+}A = \{x \mid A(x) > \alpha\}.$$

That is, the $\alpha$-cut (or the strong $\alpha$ - cut) of a fuzzy set A, is the crisp set $^{\alpha}$A (or the crisp set $^{\alpha+}$A) that contains all the elements of the universal set X whose membership grades in A are greater than or equal to (or only greater than) the specified value of $\alpha$.

Now we define for neutrosophic set A defined on X and for any value $\alpha_N \in$ [0, 1] $\cup$ [F N] the $\alpha_N$ – cut $^{\alpha_N} A$ and the strong $\alpha_N$ cut $^{\alpha_N^+} A$ as

$$^{\alpha_N} A = \{x \mid A(x) \geq \alpha_N\}$$
$$^{\alpha_N^+} A = \{x \mid A(x) > \alpha_N\}.$$

Clearly when $\alpha_N \in$ [0 1] $^{\alpha_N} A = {}^{\alpha}A$ and $^{\alpha_N^+} A = {}^{\alpha^+} A$. We also grade the indeterminary as if x, y $\in$ F N. say if

x > y (i.e. x = .7$I$ y = .3$I$ )
x < y (if x = .2$I$ and y = .6$I$ )
if x $\ngtr$ y if (x. 3$I$ and y = .5$I$ )

$\ngtr$ not comparable with y. The incomparability occurs only when x ∈ F N and y ∈ [ 0 1] not with in the set F N or with in the set [0, 1].

Another important property of both $\alpha_N$ – cuts and strong $\alpha_N$ – cuts are that if the pair $\alpha_N$ and $\alpha'_N$ ∈ FN or $\alpha'_N$ and $\alpha_N$ ∈ [ 0 1] 'or' in the mutually exclusive sense we have
if $\alpha_N < \alpha'_N$ then

$$^{\alpha_N} A \supseteq {}^{\alpha'_N} A$$

if $\alpha_N$ ∈ [ 0 1] and $\alpha'_N$ ∈ F N or vice versa we cannot say any thing about it. Also in case $\alpha_N$ and $\alpha'_N$ ∈ F N or $\alpha_N$ and $\alpha'_N$ ∈ [ 0 1] we have

$$^{\alpha_N} A \cap {}^{\alpha'_N} A = {}^{\alpha'_N} A$$

$$^{\alpha_N} A \cup {}^{\alpha'_N} A = {}^{\alpha_N} A$$

Similar results hold good in case of $^{\alpha_N^+} A$ and $^{\alpha_N^{'+}} A$. An obvious consequence of this property is that all $\alpha_N$-cuts and all strong $\alpha_N$-cuts of any neutrosophic set form two distinct families of nested crisp sets.

The 1- cuts $^1 A$ is often called the core of $A$. The $I$-cut $^I A$ is often called the Neutrosophic core (N-core) of $A$. The N-height $h_N (A) \cup h_N (A')$ of a neutrosophic set $A$ is the largest membership grade obtained by any element in that set $\cup h_N (A')$.

Formally $h_N (A) \cup h_N (A') = \sup_{x \in X} A (x) \cup \sup_{x \in X} A' (x)$

where $A$: X → [0, 1] and

$$A': X \to F N$$

if $h_N (A') = \phi$ then the N-height and height coincide otherwise we get the N-height as the union of two height one from height and other from the neutrosophic set FN. A fuzzy set $A$ is called N – normal when $h_N (A) = 1$ or $I$; it is N-subnormal when $h_N(A) < 1$ (or $< I$). The N-height of $A$ may also be viewed as the supremum of $\alpha_N$ for which $^{\alpha_N} A \neq \phi$. Let R be the set of reals R$I$ = {r$I$ / r ∈ R and $I$ be the indeterminacy}. We call this as the set of real neutrosophy contrasts with FN the set of fuzzy neutrosophic set. We venture to define R$I$ the real neutrosophy to define convexity. (R$I$)$^n$ for some n ∈ N is defined similar to R$^n$. This property is viewed as a generalization of the classical concept of convexity of crisp sets.

In order to make the neutrosophic convexity consistent with the classical definition of convexity it is required that $\alpha_N$ – cuts of a convex neutrosophic set be convex for all $\alpha_N \in (0\ 1] \cup FN \setminus \{0\}$, in the classical sense (0 – cut is excluded here since it is always equal to $(RI)^n$ in this case and this includes $-\infty I$ to $\infty I$).

Any property generalized from fuzzy set theory that is preserved in all $\alpha_N$ -cuts for $\alpha_N \in (0\ 1] \cup FN \setminus \{0\}$ in the classical sense is called a N-cutworthy property, if it is preserved in all strong $\alpha_N$ – cuts for $\alpha \in [0,\ 1]$ it is called a N-strong cut worthy property.

The reader is expected to give an example of $\alpha$ - cut worthy property and $\alpha$ - strong cutworthy. The N – standard complement $\overline{A}$ of neutrosophic set $A$ with respect to the universal set X is defined for all $x \in X$ by the equation

$$\overline{A}_N(x) = 1 - A(x) \cup I - A(x);$$

elements of X for which $A(x) = \overline{A}_N(x)$ are called N-equilibrium points of A.

N – standard union and
N – standard intersection
are defined for all $x \in X$ by equations

$$(A \cap B)_N(x) \quad = \quad \min \{A(x), B(x) \mid A, B : x \to [0\ 1]$$
$$\cup \quad \min \{A_N(x), B_N(x) \mid A_N, B_N : x \to FN\}.$$

Similarly

$$(A \cup B)_N(x) \quad = \quad \max \{A(x), B(x) \mid A, B : X \to [0\ 1]\} \cup$$
$$\max \{A_N(x)\ B_N(x) \mid A_N, B_N : X \to FN\}.$$

$$\min \{A(x),\ 1\text{-}A(x)] \quad = \quad 0 \mid A : X \to [01]\} \text{ or}$$
$$\{\min \{A_N(x),\ 1 - A_N(x)\} \quad = \quad 0 / A_N\ X \to FN\}$$

Now as in case of fuzzy subsets we define for neutrosophic subset the notion of the degree of n-subset hood $S_N(AB)$ where A and B are neutrosophic subsets of some universal set X.

$$S_N(A,\ B) = \frac{|A| - \sum \max\{0,\ A(x) - B(x)\}}{|A|}$$

$S_N$ (A, B) may be defined provided A and B are defined on the same interval [0, 1] or FN ('or' in the mutually exclusive sense)

$$0 \le S_N \ (A, B) \le 1 \ or \ 0 \le S_N \ (A, B) \le I \ .$$

$S_N$ (A, B) will remain undefined if A and B are in different sets i.e., one subset defined over [ 0 1 ] and the other over F N. It is left for the reader as an exercise to prove

$$S_N \ (A, B) = \frac{\left| A \cap B \right|}{\left| A \right|}$$

where $\cap$ denotes the standard neutrosophic intersection.

### 3.3 On Neutrosophic lattices

In this section we introduce the notion of neutrosophic lattices and give some of its properties. Three types of neutrosophic lattices are dealt in this section; viz. integral neutrosophic lattice, neutrosophic chain lattice and mixed neutrosophic lattices. Here we mainly define them and illustrate with examples. *I* denotes the concept of indeterminacy.

**DEFINITION 3.3.1:** *Let $N = L \cup \{I\}$ where L is any lattice and I an indeterminate.*

Define the max, min operation on N as follows

Max {x, *I*} = *I* for all x ∈ L\ {1}
Max {1, *I*} = *I*
Min {x, *I*} = *I* for all x ∈ L \ {0}
Min {0, *I*} = 0

We know if x, y ∈ L then max and min are well defined in L.
    N is called the integral neutrosophic lattice.

*Example 3.3.1:* Let N = L ∪ {*I*} given by the following diagram:

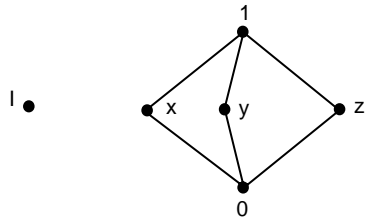

clearly    Min {x, $I$ } = $I$ for all x ∈ L\ {0}
           Min {0, $I$ } = 0
           Max {x, $I$ } = $I$ for all x ∈ L \ {1}
           Max {1, $I$ } = 0.

We see N is an integral neutrosophic lattice and clearly the order of N is 6.

**Remark 3.3.1:** If L is a lattice of order n and N = L ∪ {$I$ } be an integral neutrosophic lattice then order of N is n + 1.
2. For a integral neutrosophic lattice N also {0} is the minimal element and {1} is the maximal element of N.

**DEFINITION 3.3.2:** *Let $C_I$ = {nI | n ∈ [0, 1)} ∪ {1}*
*$C_I$  can be made into a lattice by defining max and min as follows*

           *Min {0, nI} = 0 for all n ∈ [ 0, 1)*
           *Max {1, nI} = 1 for all n ∈ [ 0, 1)*
           *Min {$n_1I$, $n_2I$} = $n_1I$ if $n_1 ≤ n_2$ for all $n_1$, $n_2$ ∈ ( 0, 1)*
           *Max {$n_1I$, $n_2I$} = $n_2I$ if $n_1 ≤ n_2$ for all $n_1$, $n_2$ ∈ (0, 1)*

*Clearly $C_I$ is a lattice called the neutrosophic chain lattice.*
     It is however important to note that we do not have any relation between integral neutrosophic lattice and neutrosophic chain lattice.
     Next we proceed onto define mixed neutrosophic lattice and pure neutrosophic lattice.

**DEFINITION 3.3.3:** *Let N = $L_I$, $I$ ∪ {0, 1} where $L_I$ is any lattice*
*$L_II$ = {xI / x ∈ $L_I$ \ {0, 1}}.*
*N = $L_II$ ∪ {0, 1} is a lattice under the following min. max operations*
           *Min {$x_1I$, $x_2 I$} = $x_1 I$ if and only if*
           *Min {$x_1$, $x_2$} = $x_1$*
           *Max {$x_1I$, $x_2 I$} = $x_2 I$ if and only if max {$x_1$, $x_2$}= $x_2$*
           *Min {0, xI} = 0 and*
           *Max {1, xI} = 1.*
*Clearly N is a lattice called the pure neutrosophic lattice.*

**Remark 3.3.2:** All neutrosophic chain lattices are pure neutrosophic lattices. Clearly integral neutrosophic lattices are not pure neutrosophic lattices. All neutrosophic chain lattices in general need not be pure neutrosophic lattices.

***Example 3.3.2:*** Consider

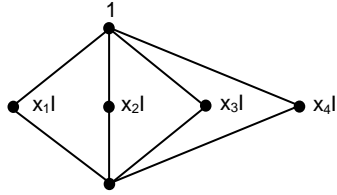

Clearly N = L *I* ∪ {0, 1} where L

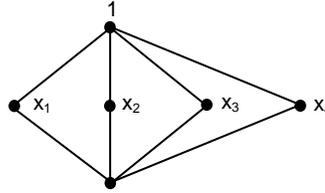

is a pure neutrosophic lattice but is not a neutrosophic chain lattice.

**DEFINITION 3.3.4:** *Let L be a lattice N = L ∪ L I is said to be special neutrosophic lattice where LI = {xI / x ∈ L \ {0 1}, if on N is defined min and max as follows:*

$$Min \{x, xI\} = x \ I \ for \ all \ x \in L \setminus \{0 \ 1\}.$$
$$Max \{x, xI\} = xI \ for \ all \ x \in L \setminus \{0 \ 1\}$$
$$Min \{x_1, x_2I\} = x_1 \ if \ x_1 < x_2$$
$$Min \{x_1, x_2I\} = x_2I \ if \ x_2 < x_1$$
$$Min \{x_1I, x_2I\} = x_1I \ iff \ x_1 < x_2$$
$$Min \{0, xI\} = 0$$
$$Max \{x_1, x_2I\} = x_2 \ I \ iff \ x_1 < x_2$$
$$Max \{x_1, x_2I\} = x_1 \ iff \ x_1 > x_2$$
$$Max \{x_1I, x_2I\} = x_2I \ iff \ x_1 < x_2$$
$$Max \{1, xI\} = 1.$$

*Min (x, y} and max {x, y} for all x, y ∈ L is defined in the usual way. N with these max-min function is called the special neutrosophic lattice.*

Next we proceed on to define mixed special neutrosophic lattices.

**DEFINITION 3.3.5:** *Let N = $L_1$ ∪ $L_2$ I where $L_1$ and $L_2$ are two distinct lattices. Define max and min on N as follow:.*

$$Min \{x_1, x_2I\} = x_2I, \ x_1 \in L_1 \setminus \{0\} \ x_2 \in L_2 \setminus \{0\}$$

$Min \{x_1, 0\} = 0$
$Min \{0, x_2I\} = 0$
$Max \{x_1, x_2I\} = x_2I$
$Max \{1, x_2I\} = 1$
$Max \{x_1, 1\} = 1$

*Min and max on elements of $L_1$ and $L_2$ are done as in case of lattices. Then we call N the mixed special neutrosophic lattice.*

***Example 3.3.3:*** Let $N = L \cup L \, I$ where L is

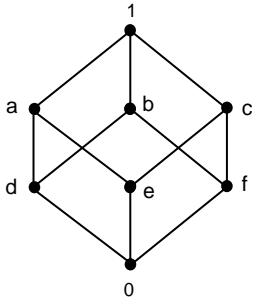

N is a special neutrosophic lattice with 14 elements in it.

***Example 3.3.4:*** Let $N = L \cup L \, I$, where L is

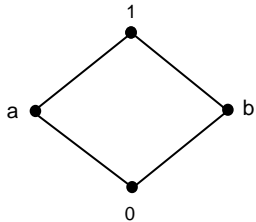

N is a special neutrosophic lattice with $N = \{0, 1, aI, bI, a, b\}$ with six elements.

***Example 3.3.5:*** Let $N = L_1 \cup L_2 \, I$ where $L_1$ is

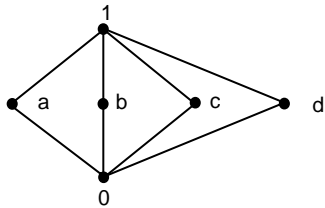

And $L_2$ is

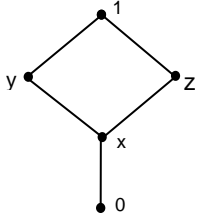

N = {0, 1, a, b, c, d, x$I$, y$I$, z$I$} is a mixed special neutrosophic lattice with 9 elements in it.

Now we proceed on to define the notion of types of neutrosophic Boolean algebra.

**DEFINITION 3.3.6:** *Let B be a Boolean algebra of order greater than 2. N = B $\cup$ I is defined as the integral neutrosophic Boolean algebra. Clearly we do not have the notion of neutrosophic chain Boolean algebra. The order of an integral neutrosophic Boolean algebra is $2^n + 1$ for all integers n, n $\geq$ 2.*

**DEFINITION 3.3.7:** *Let N = L$_1$ I $\cup$ {0, 1} be a pure neutrosophic lattice, N is called as pure neutrosophic Boolean algebra provided L$_1$ is a Boolean algebra of order greater than or equal to four. On similar lines if we replace in the definition of the special neutrosophic lattice N = L $\cup$ L I L by a Boolean algebra then we call N the special neutrosophic Boolean algebra.*

The notion of mixed special neutrosophic Boolean algebra is defined in an analogous way.

We can define direct product of neutrosophic lattices as in case of direct product of lattices. We know a lattice L is called complete if $\wedge$H and $\vee$H exist for any subset H $\subset$ L. We define complete integral neutrosophic lattice, complete special neutrosophic lattice and complete mixed special neutrosophic lattice as lattices for which $\wedge$H and $\vee$H exist for any subsets of any type of neutrosophic lattices.

### 3.4 Neutrosophic notions : Basic concepts

To define the concept of neutrosophic binary relations, neutrosophic relational equations and other related concepts we need the basic notion called the neutrosophic function and other related neutrosophic concepts.

**DEFINITION 3.4.1:** *A set is defined by a function $N_A$ called the neutrosophic function, that declares which element of X are members of the set, which are not members and which are indeterminate. Set A is defined by the neutrosophic function $N_A$ as follows:*

$$N_A(x) = \begin{cases} 1 & if & x \in A \\ 0 & if & x \notin A \\ I & if \ we \ cannot \ decide \end{cases}$$

*where I denotes the indeterminate symbol. That is the neutrosophic function $N_A$ maps elements of X to elements of the set $[0, 1] \cup I$ which is formally expressed by*

$$N_A : X \rightarrow [0, 1] \cup I.$$

*For each $x \in X$ when $N_A(x) = 1$, x is declared to be a member of A when When $N_A(x) = 0$, x is declared to be a non member of A when $N_A(x) = I$, x cannot be determined whether it is a member of A or a non member of A.*

If every member of a set A is also a member of set B (i.e. if $x \in A$ implies $x \in B$) then A is called a subset of B, and this is written as $A \subseteq B$.

Every set is a subset of itself, and every set is a subset of the universal set. *If* $A \subseteq B$ and $B \subseteq A$ then A and B contain the same members. They are then called equal sets that is denoted by $A = B$. We have the following relation between fuzzy characteristic function and neutrosophic function.

**THEOREM 3.4.2:** *Every characteristic function is a neutrosophic function. But all neutrosophic functions in general need not be a characteristic function.*

*Proof:* Since every characteristic function is obviously neutrosophic function with no indeterminate part associated with it. On the other hand if $N_A$ is a neutrosophic function with at least a x $\in$ X such that $N_A$ (x) = $I$ then clearly $N_A$ is not a characteristic function.

Several fuzzy sets representing linguistic concepts such as low, medium high and so on with indeterminates are used to define states of a variable. Such a variable is usually called a neutrosophic variable.

**THEOREM 3.4.3:** *Every fuzzy variable is always a neutrosophic variable, but all neutrosophic variables in general are not fuzzy variables.*

*Proof:* By the very definition every fuzzy variable $N_A$ (x) $\neq$ $I$ for any x in X so is a neutrosophic variable. But in case of a neutrosophic variable if we have a x $\in$ X such that $N_A(X) = I$ then we see it is not a fuzzy variable.

Let N {[0, 1] $\cup$ $I$} denote the set of all neutrosophic sets that can be defined within the universal set [0, 1] $\cup$ $I$, N{[0, 1] $\cup$ $I$} is called the neutrosophic power set of [0, 1] $\cup$ $I$.

Level 2 neutrosophic sets are those, which have their membership functions A : N {X} $\rightarrow$ [0 , 1] $\cup$ $I$. where N {x} denotes the neutrosophic power set of X. (the set of all neutrosophic sets of X).

Level 2 neutrosophic sets are generalized to level 3 fuzzy sets by using a universal whose elements are level 2 neutrosophic sets. Higher level neutrosophic sets can be obtained recursively in the same way.

$$A : N (X) \rightarrow N \{[0, 1] \cup I \}$$

other combinations are also possible.

Given a neutrosophic set A defined on X and any value $\alpha \in$ [0, 1] $\cup$ $I$, the $\alpha$ - cut $^\alpha$A and the strong $\alpha$ - cut $^{\alpha+}$A are the neutrosophic crisp sets .

$$^\alpha A = \{ x \mid A (x) \geq \alpha \}$$
$$^{\alpha+}A = \{ x \mid A (x) > \alpha \}.$$

The neutrosophic crisp sets are crisp set if $\alpha \in$ [0, 1] if $\alpha = I$ then the set is a neutrosophic crisp set.

The set of all levels $\alpha \in$ [0, 1] $\cup$ $I$ that represent distinct $\alpha$ - cuts of a given set A is called a neutrosophic level set of A. Formally $\Lambda$ (A) = { $\alpha$ | A (x) = $\alpha$ for some x $\in$ X; $\alpha \in$ [ 0, 1] $\cup$ $I$}. where $\Lambda$ denotes the level set of neutrosophic set A defined on X.

Unlike in fuzzy sets we cannot in case of neutrosophic sets always say for any $\alpha_1, \alpha_2 \in [0, 1] \cup I$ $\alpha_1 < \alpha_2$ (or $\alpha_2 < \alpha_1$); for if $\alpha_1 = I$ and $\alpha_2 \in [0, 1]$ we cannot compare them. This is a marked difference between fuzzy crisp sets with $\alpha$ - cuts and neutrosophic sets with $\alpha$ - cuts.

Now we have to define as usual $1 -$ cut $^1A$ to be core of A; $^IA$ to be indeterminate core of A. A neutrosophic set a is called normal when h (A) = 1 is called subnormal when h (A) < 1 and neutrosophic if h (A) = $I$ where

$$h(A) = \sup_{x \in X} A(x);$$

A (x) are neutrosophic sets.

Several results true in case of fuzzy sets can be easily extended in case of neutrosophic sets. For more in this direction refer [43].

The representations of an arbitrary neutrosophic set A in terms of the special neutrosophic set is carried out in the following way. Let $\alpha \in [0, 1] \cup I$ we convert each of the $\alpha$ - cuts to a special neutrosophic set $_\alpha A$ defined for each x $\in$ X as

$$_\alpha A(x) = \alpha. \, ^\alpha A(x)$$
$$\text{if } \alpha = I, \, IA(x) = I^\alpha. \, A(x).$$

**DEFINITION 3.4.4:** *We say a crisp function f: X $\rightarrow$ Y is neutrosified when it is extended to act on neutrosophic sets defined on X and Y. That is the neutrosified function for which the same symbol is usually used has the form*

$$f: N(X) \rightarrow N(Y)$$

*and its inverse function $f^{-1}$, has the form*

$$f^{-1}: N(Y) \rightarrow N(X)$$

*[Here N(X) denotes the neutrosophic power set of X]. A principle for neutrosifying crisp functions is called an extension principle.*

**Conventions about neutrosophic sets**

Let A, B $\in$ $N$(X) i.e.,

$$A : X \rightarrow [0, 1] \cup I$$
$$B : X \rightarrow [0, 1] \cup I$$

$(A \cap B) (x) = \min \{A (x), B (x)$ if $A (x) = I$ or $B (x) = I$

then $(A \cap B) (x)$ is defined to be $I$ i.e., $\min \{A (x), B (x)) = I$ $(A \cup B) (x) = \max \{A (x), B (x)$ if one of $A (x) = I$ or $B (x) = I$ then $(A \cup B) (x) = I$ i.e., $\max \{A (x), B (x)\} = I$.

Thus it is pertinent to mention here that if one of $A (x) = I$ or $B(x) = I$ then $(A \cup B) (x) = (A \cap B) x$. i.e., is the existence of indeterminacy $\max \{A (x), B (x)\} = \min \{A (x), B(x)\} = I$.

$\overline{A} (x) = 1- A (x)$; if $A (x) = I$ then $\overline{A} (x) = A (x) = I$.

### 3.5 Neutrosophic matrices and fuzzy neutrosophic matrices

In this section we define the concept of neutrosophic matrices and fuzzy neutrosophic matrices and the operations on these matrices are also given.

**DEFINITION 3.5.1:** *Let $M_{nxm} = \{(a_{ij}) \ / \ a_{ij} \in K(I)\}$, where $K (I)$, is a neutrosophic field. We call $M_{nxm}$ to be the neutrosophic matrix.*

***Example 3.5.1:*** Let $Q(I) = \langle Q \cup I \rangle$ be the neutrosophic field.

$$M_{4\times 3} = \begin{pmatrix} 0 & 1 & I \\ -2 & 4I & 0 \\ 1 & -I & 2 \\ 3I & 1 & 0 \end{pmatrix}$$

is the neutrosophic matrix, with entries from rationals and the indeterminacy $I$. We define product of two neutrosophic matrices whenever the product is defined as follows:

Let

$$A = \begin{pmatrix} -1 & 2 & -I \\ 3 & I & 0 \end{pmatrix}_{2\times 3} \text{ and } B = \begin{pmatrix} -I & 1 & 2 & 4 \\ 1 & I & 0 & 2 \\ 5 & -2 & 3I & -I \end{pmatrix}_{3x4}$$

$$AB = \begin{bmatrix} -6I+2 & -1+4I & -2-3I & I \\ -2I & 3+I & 6 & 12+2I \end{bmatrix}_{2 \times 4}$$

(we use the fact $I^2 = I$).

Now we proceed onto define the notion of fuzzy integral neutrosophic matrices and operations on them, for more about these refer [103].

**DEFINITION 3.5.2:** *Let $N = [0, 1] \cup I$ where $I$ is the indeterminacy. The $m \times n$ matrices $M_{m \times n} = \{(a_{ij}) / a_{ij} \in [0, 1] \cup I\}$ is called the fuzzy integral neutrosophic matrices. Clearly the class of $m \times n$ matrices is contained in the class of fuzzy integral neutrosophic matrices.*

**Example 3.5.2:** Let $\quad A = \begin{pmatrix} I & 0.1 & 0 \\ 0.9 & 1 & I \end{pmatrix}$.

A is a $2 \times 3$ integral fuzzy neutrosophic matrix.

We define operation on these matrices. An integral fuzzy neutrosophic row vector is a $1 \times n$ integral fuzzy neutrosophic matrix. Similarly an integral fuzzy neutrosophic column vector is a $m \times 1$ integral fuzzy neutrosophic matrix.

**Example 3.5.3:** $A = (0.1, 0.3, 1, 0, 0, 0.7, I, 0.002, 0.01, I, 1, 0.12)$ is a integral row vector or a $1 \times 11$, integral fuzzy neutrosophic matrix.

**Example 3.5.4:** $B = (1, 0.2, 0.111, I, 0.32, 0.001, I, 0, 1)^T$ is an integral neutrosophic column vector or B is a $9 \times 1$ integral fuzzy neutrosophic matrix.
    We would be using the concept of fuzzy neutrosophic column or row vector in our study.

**DEFINITION 3.5.3:** *Let $P = (p_{ij})$ be a $m \times n$ integral fuzzy neutrosophic matrix and $Q = (q_{ij})$ be a $n \times p$ integral fuzzy neutrosophic matrix. The composition map $P \bullet Q$ is defined by $R = (r_{ij})$ which is a $m \times p$ matrix where $r_{ij} = \max_{k} \min(p_{ik} q_{kj})$ with the assumption $\max(p_{ij}, I) = I$ and $\min(p_{ij}, I) = I$ where $p_{ij} \in [0, 1]$. $\min(0, I) = 0$ and $\max(1, I) = 1$.*

*Example 3.5.5:* Let

$$P = \begin{bmatrix} 0.3 & I & 1 \\ 0 & 0.9 & 0.2 \\ 0.7 & 0 & 0.4 \end{bmatrix}, \ Q = (0.1, I, 0)^T$$

be two integral fuzzy neutrosophic matrices.

$$P \bullet Q = \begin{bmatrix} 0.3 & I & 1 \\ 0 & 0.9 & 0.2 \\ 0.7 & 0 & 0.4 \end{bmatrix} \bullet \begin{bmatrix} 0.1 \\ I \\ 0 \end{bmatrix} = (I, I, 0.1).$$

*Example 3.5.6:* Let

$$P = \begin{bmatrix} 0 & I \\ 0.3 & 1 \\ 0.8 & 0.4 \end{bmatrix}$$

and

$$Q = \begin{bmatrix} 0.1 & 0.2 & 1 & 0 & I \\ 0 & 0.9 & 0.2 & 1 & 0 \end{bmatrix}.$$

One can define the max-min operation for any pair of integral fuzzy neutrosophic matrices with compatible operation.

Now we proceed onto define the notion of fuzzy neutrosophic matrices. Let $N_s = [0, 1] \cup nI \mathbin{/} n \in (0, 1]\}$; we call the set $N_s$ to be the fuzzy neutrosophic set.

**DEFINITION 3.5.3:** *Let $N_s$ be the fuzzy neutrosophic set. $M_{n \times n} = \{(a_{ij}) \mathbin{/} a_{ij} \in N_s\}$ we call the matrices with entries from $N_s$ to be the fuzzy neutrosophic matrices.*

*Example 3.5.6:* Let $N_s = [0,1] \cup \{nI \mathbin{/} n \in (0,1]\}$ be the set

$$P = \begin{bmatrix} 0 & 0.2I & 0.31 & I \\ I & 0.01 & 0.7I & 0 \\ 0.31I & 0.53I & 1 & 0.1 \end{bmatrix}$$

P is a $3 \times 4$ fuzzy neutrosophic matrix.

***Example 3.5.7:*** Let $N_s = [0, 1] \cup \{nI / n \in (0, 1]\}$ be the fuzzy neutrosophic matrix. $A = [0, 0.12I, I, 1, 0.31]$ is the fuzzy neutrosophic row vector:

$$B = \begin{bmatrix} 0.5I \\ 0.11 \\ I \\ 0 \\ -1 \end{bmatrix}$$

is the fuzzy neutrosophic column vector.

Now we proceed on to define operations on these fuzzy neutrosophic matrices.

Let $M = (m_{ij})$ and $N = (n_{ij})$ be two $m \times n$ and $n \times p$ fuzzy neutrosophic matrices. $M \bullet N = R = (r_{ij})$ where the entries in the fuzzy neutrosophic matrices are fuzzy indeterminates i.e. the indeterminates have degrees from 0 to 1 i.e. even if some factor is an indeterminate we try to give it a degree to which it is indeterminate for instance $0.9I$ denotes the indeterminacy rate is high where as $0.01I$ denotes the low indeterminacy rate. Thus neutrosophic matrices have only the notion of degrees of indeterminacy. Any other type of operations can be defined on the neutrosophic matrices and fuzzy neutrosophic matrices. The notion of these matrices will be used to define neutrosophic relational equations and fuzzy neutrosophic relational equations.

## 3.6 Characteristics and significance of the newer Paradigm Shift using indeterminacy

The concept of a scientific paradigm was introduced by Thomas Kuhn ion his important highly influential book. The Structure of Scientific Revolutions [45] it is defined as a art of theories, standards principles and methods that are taken for granted by the scientific community in a given field. Using this concept Kuhn characterizes scientific paradigm as a process of normal science, based upon a particular paradigm, are interwoven with periods of paradigm shifts which are referred to by Kuhn as scientific revolutions.

In his book Kuhn illustrates the notion of a paradigm shift by many well – documented examples from the history of Science.

Some of the most visible paradigm shifts are associated with the names of Copernicus (astronomy), Newton (mechanics), Lavoisier (chemistry) Darwin (biology), Maxwell (electromagnetism), Einstein (mechanics) and Godel (mathematics).

Although paradigm shifts vary from one another in their scope pace and other features they share a few general characteristics. Each paradigm shift is initiated by emerging problems that are difficult or impossible to deal with the current paradigm (paradoxes, anomalies etc.). Each paradigm when proposed is initially rejected in various forms (it is ignored ridiculed, attacked etc) by most scientists in the given field. Those who usually support the new paradigm are either very young or very new to the field and consequently not very influential. Since the paradigm is initially not well developed the position of its proponents is weak. The paradigm eventually gains its status on pragmatic grounds by demonstrating that it is more successful than the existing paradigm in dealing with problems that are generally recognized as acute. As a rule, the greater the scope of a paradigm shift the longer it takes for the new paradigm to be generally accepted.

The same need was expressed by Zadeh [1962] three years before he actually proposed the new paradigm of mathematics based upon the concept of a fuzzy set. When the new paradigm was proposed by Zadeh the usual process of a paradigm shift began. The concept of a fuzzy set which underlies this new paradigm was initially ignored ridiculed or attacked by many, while it was supported only by a few, mostly young and not influential. In spite of the initial lack of interest skepticism or even open hostility. The new paradigm persevered with virtually no support in the 1960's matured significantly and gained some support in 1970s and began to demonstrate its superior pragmatic utility in the 1980s.

The paradigm shift is still ongoing and it will likely take much longer than usual to complete it. This is not surprising since the scope of the paradigm shift is enormous. The new paradigm does not affect any particular field of science but the very foundation of science. In fact it challenges the most sacred element of the foundations. The Aristotelian two valued logic, which for millennia has been taken for granted and viewed as inviolable. The acceptance of such a radical challenge is surely difficult for most scientists, it requires an open mind, enough time, and considerable effort to properly comprehend the meaning and significance of the paradigm shift involved.

At this time we can recognize at least four features that make the new paradigm superior to the classical paradigm.

The new paradigm allows us to express irreducible observation and measurement uncertainties in their various manifestation and make these uncertainties intrinsic to empirical data. Such data which are based on graded distinctions among states of relevant variables are usually called fuzzy data. When fuzzy data is processed their intrinsic uncertainties their by the results are more meaningful.

This new paradigm offers for greater resources for managing complexity and controlling computational cost. This new paradigm has considerably greater expressive power consequently it can effectively deal with a broader class of problems. This has greater capability to capture human common sense reasoning decision making and other aspects of human cognition.

Now still new concept is the neutrosophy, the neutrosophic set and the related concept. This still newer paradigm allows us to express the indeterminacies involved in the analysis of empirical data. Indeterminacy is more powerful in a way than uncertainties. When neutrosophic data are processed their indeterminacies are processed as well and the consequent results are more meaningful.

Further this newer paradigm offers far greater resources for managing indeterminacy for in any study be it scientific or otherwise the role played by the indeterminacy factor is significant but till date we have been ignoring this factor. The newer paradigm has greater power in analyzing the indeterminacy present it several human problems be it legal, medical, chemical or industrial. The newer paradigm has a greater capability to capture the reality for in reality a lot of indeterminacy is involved.

Chapter Four

# NEUTROSOPHIC RELATIONAL EQUATIONS AND ITS APPLICATIONS

We have introduced several properties about neutrosophic concepts in the earlier chapter. Study of neutrosophic relational equations (NREs) will find its applications whenever indeterminacy plays a vital role. Thus in problems which involves indeterminacy certainly NREs would be more appropriate than FREs.

This chapter has eleven sections which introduces and applies NRE to several real world problems as well analyses and studies for the first time; optimization of NRE using max product composition; solution using lattices and so on.

## 4.1 Binary neutrosophic Relation and their properties

In this section we introduce the notion of neutrosophic relational equations and fuzzy neutrosophic relational equations and analyze and apply them to real-world problems, which are abundant with the concept of indeterminacy. We also mention that most of the unsupervised data also involve at least to certain degrees the notion of indeterminacy.

Throughout this section by a neutrosophic matrix we mean a matrix whose entries are from the set $N = [0, 1] \cup I$ and by a fuzzy neutrosophic matrix we mean a matrix whose entries are from $N' = [0, 1] \cup \{nI / n \in (0,1]\}$.

Now we proceed on to define binary neutrosophic relations and binary neutrosophic fuzzy relation.

A binary neutrosophic relation $R_N(x, y)$ may assign to each element of X two or more elements of Y or the indeterminate $I$. Some basic operations on functions such as the inverse and composition are applicable to binary relations as well. Given a neutrosophic relation $R_N(X, Y)$ its domain is a neutrosophic set on $X \cup I$ domain R whose membership function is defined by $\text{domR}(x) = \max_{y \in X \cup I} R_N(x, y)$ for each $x \in X \cup I$.

That is each element of set $X \cup I$ belongs to the domain of R to the degree equal to the strength of its strongest relation to any member of set $Y \cup I$. The degree may be an indeterminate $I$ also. Thus this is one of the marked difference between the binary fuzzy relation and the binary neutrosophic relation. The range of $R_N(X,Y)$ is a neutrosophic relation on Y, ran R whose membership is defined by $\text{ran } R(y) = \max_{x \in X} R_N(x, y)$ for each $y \in$ Y, that is the strength of the strongest relation that each element of Y has to an element of X is equal to the degree of that element's membership in the range of R or it can be an indeterminate $I$.

The height of a neutrosophic relation $R_N(x, y)$ is a number h(R) or an indeterminate $I$ defined by $h_N(R) = \max_{y \in Y \cup I} \max_{x \in X \cup I} R_N(x, y)$. That is $h_N(R)$ is the largest membership grade attained by any pair $(x, y)$ in R or the indeterminate $I$.

A convenient representation of the neutrosophic binary relation $R_N(X, Y)$ are membership matrices $R = [\gamma_{xy}]$ where $\gamma_{xy} \in R_N(x, y)$. Another useful representation of a binary neutrosophic relation is a neutrosophic sagittal diagram. Each of the sets X, Y represented by a set of nodes in the diagram, nodes corresponding to one set are clearly distinguished from nodes representing the other set. Elements of X' × Y' with non-zero membership grades in $R_N(X, Y)$ are represented in the diagram by lines connecting the respective nodes. These lines are labeled with the values of the membership grades.

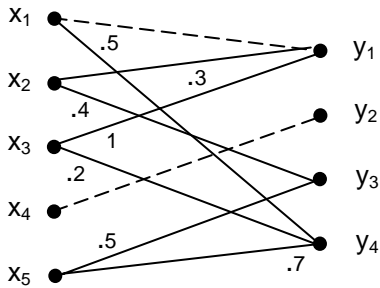

**FIGURE: 4.1.1**

An example of the neutrosophic sagittal diagram is a binary neutrosophic relation $R_N(X, Y)$ together with the membership neutrosophic matrix is given below.

$$\begin{array}{c} & \begin{array}{cccc} y_1 & y_2 & y_3 & y_4 \end{array} \\ \begin{array}{c} x_1 \\ x_2 \\ x_3 \\ x_4 \\ x_5 \end{array} & \left[ \begin{array}{cccc} I & 0 & 0 & 0.5 \\ 0.3 & 0 & 0.4 & 0 \\ 1 & 0 & 0 & 0.2 \\ 0 & I & 0 & 0 \\ 0 & 0 & 0.5 & 0.7 \end{array} \right] \end{array}$$

The inverse of a neutrosophic relation $R_N(X, Y) = R(x, y)$ for all $x \in X$ and all $y \in Y$. A neutrosophic membership matrix $R^{-1} = [\,r_{yx}^{-1}\,]$ representing $R_N^{-1}(Y, X)$ is the transpose of the matrix R for $R_N(X, Y)$ which means that the rows of $R^{-1}$ equal the columns of R and the columns of $R^{-1}$ equal rows of R. Clearly $(R^{-1})^{-1} = R$ for any binary neutrosophic relation.

Consider any two binary neutrosophic relation $P_N(X, Y)$ and $Q_N(Y, Z)$ with a common set Y. The standard composition of these relations which is denoted by $P_N(X, Y) \bullet Q_N(Y, Z)$ produces a binary neutrosophic relation $R_N(X, Z)$ on $X \times Z$ defined by $R_N(x, z) = [P \bullet Q]_N(x, z) = \max_{y \in Y} \min[P_N(x, y), Q_N(x, y)]$ for all $x \in$ X and all $z \in Z$.

This composition which is based on the standard $t_N$-norm and $t_N$-co-norm, is often referred to as the max-min composition. It can be easily verified that even in the case of binary neutrosophic relations $[P_N(X, Y) \bullet Q_N(Y, Z)]^{-1} = Q_N^{-1}(Z, Y) \bullet P_N^{-1}(Y, X)$. $[P_N(X, Y) \bullet Q_N(Y, Z)] \bullet R_N(Z, W) = P_N(X, Y) \bullet [Q_N(Y, Z) \bullet R_N(Z, W)]$, that is, the standard (or max-min) composition is associative and its inverse is equal to the reverse composition of the inverse relation. However, the standard composition is not commutative, because $Q_N(Y, Z) \bullet P_N(X, Y)$ is not well defined when $X \neq Z$. Even if $X = Z$ and $Q_N(Y, Z) \circ P_N(X, Y)$ are well defined still we can have $P_N(X, Y) \circ Q(Y, Z) \neq Q(Y, Z) \circ P(X, Y)$.

Compositions of binary neutrosophic relation can the performed conveniently in terms of membership matrices of the relations. Let $P = [p_{ik}]$, $Q = [q_{kj}]$ and $R = [r_{ij}]$ be membership matrices of binary relations such that R = P ° Q. We write this using matrix notation

$$[r_{ij}] = [p_{ik}] \circ [q_{kj}]$$

where $r_{ij} = \max_{k} \min (p_{ik}, q_{kj})$.

A similar operation on two binary relations, which differs from the composition in that it yields triples instead of pairs, is known as the relational join. For neutrosophic relation $P_N (X, Y)$ and $Q_N (Y, Z)$ the relational join $P * Q$ corresponding to the neutrosophic standard max-min composition is a ternary relation $R_N (X, Y, Z)$ defined by $R_N (x, y, z) = [P * Q]_N (x, y, z) = \min [P_N (x, y), Q_N (y, z)]$ for each $x \in X$, $y \in Y$ and $z \in Z$.

This is illustrated by the following Figure 4.1.2.

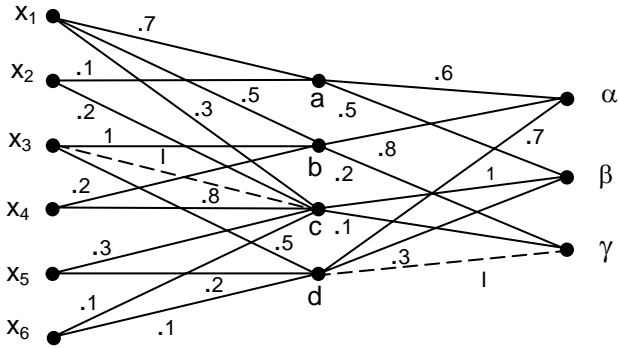

**FIGURE: 4.1.2**

In addition to defining a neutrosophic binary relation there exists between two different sets, it is also possible to define neutrosophic binary relation among the elements of a single set X.

A neutrosophic binary relation of this type is denoted by $R_N(X, X)$ or $R_N (X^2)$ and is a subset of $X \times X = X^2$.

These relations are often referred to as neutrosophic directed graphs or neutrosophic digraphs. [103]

Neutrosophic binary relations $R_N (X, X)$ can be expressed by the same forms as general neutrosophic binary relations. However they can be conveniently expressed in terms of simple diagrams with the following properties.

I.      Each element of the set X is represented by a single node in the diagram.

II.     Directed connections between nodes indicate pairs of elements of X for which the grade of membership in R is non zero or indeterminate.

III.    Each connection in the diagram is labeled by the actual membership grade of the corresponding pair in R or in indeterminacy of the relationship between those pairs.

The neutrosophic membership matrix and the neutrosophic sagittal diagram is as follows for any set X = {a, b, c, d, e}.

$$\begin{bmatrix} 0 & I & .3 & .2 & 0 \\ 1 & 0 & I & 0 & .3 \\ I & .2 & 0 & 0 & 0 \\ 0 & .6 & 0 & .3 & I \\ 0 & 0 & 0 & I & .2 \end{bmatrix}$$

Neutrosophic membership matrix for x is given above and the neutrosophic sagittal diagram is given below.

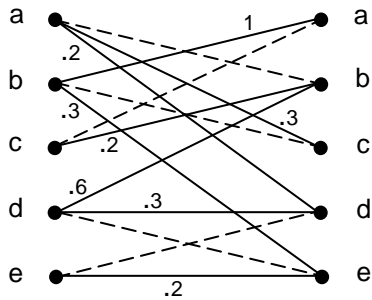

Figure: 4.1.3

Neutrosophic diagram or graph is left for the reader as an exercise.

The notion of reflexivity, symmetry and transitivity can be extended for neutrosophic relations $R_N$ (X, Y) by defining them in terms of the membership functions or indeterminacy relation.

*Thus $R_N$ (X, X) is reflexive if and only if $R_N$ (x, x) = 1 for all x ∈ X. If* *this is not the case for some x ∈ X the relation is irreflexive.*

*A weaker form of reflexivity, if for no x in X; $R_N(x, x) = 1$ then we call the relation to be anti-reflexive referred to as $\in$-reflexivity, is sometimes defined by requiring that*

$$R_N(x, x) \geq \in \text{ where } 0 < \in < 1.$$

*A fuzzy relation is symmetric if and only if*

$$R_N(x, y) = R_N(y, x) \text{ for all } x, y, \in X.$$

*Whenever this relation is not true for some x, y $\in$ X the relation is called asymmetric. Furthermore when $R_N(x, y) > 0$ and $R_N(y, x) > 0$ implies that $x = y$ for all $x, y \in X$ the relation R is called anti-symmetric.*

*A fuzzy relation $R_N(X, X)$ is transitive (or more specifically max-min transitive) if*

$$R_N(x, z) \geq \max_{y \in Y} \min [R_N(x, y), R_N(y, z)]$$

*is satisfied for each pair $(x, z) \in X^2$. A relation failing to satisfy the above inequality for some members of X is called non-transitive and if $R_N(x, x) < \max_{y \in Y} \min [RN(x, y), RN(y, z)]$ for all $(x, x) \in X^2$, then the relation is called anti-transitive*

*Given a relation $R_N(X, X)$ its transitive closure $\overline{R}_{NT}(x, X)$ can be analyzed in the following way.*

The transitive closure on a crisp relation $R_N(X, X)$ is defined as the relation that is transitive, contains

$$R_N(X, X) < \max_{y \in Y} \min [R_N(x, y) R_N(y, z)]$$

for all $(x, x) \in X^2$, then the relation is called anti-transitive. Given a relation $R_N(x, x)$ its transitive closure $\overline{R}_{NT}(X, X)$ can be analyzed in the following way.

The transitive closure on a crisp relation $R_N(X, X)$ is defined as the relation that is transitive, contains $R_N$ and has the fewest possible members. For neutrosophic relations the last requirement is generalized such that the elements of transitive closure have the smallest possible membership grades, that still allow the first two requirements to be met.

Given a relation $R_N(X, X)$ its transitive closure $\overline{R}_{NT}(X, X)$ can be determined by a simple algorithm.

Now we proceed on to define the notion of neutrosophic equivalence relation.

**DEFINITION 4.1.1:** *A crisp neutrosophic relation $R_N(X, X)$ that is reflexive, symmetric and transitive is called an neutrosophic equivalence relation. For each element x in X, we can define a crisp neutrosophic set $A_x$ which contains all the elements of X that are related to x by the neutrosophic equivalence relation.*

*Formally $A_x = [ y / (x, y) \in R_N (X, X)]$. $A_x$ is clearly a subset of X. The element x is itself contained in $A_x$, due to the reflexivity of R because R is transitive and symmetric each member of $A_x$ is related to all other members of $A_x$. Further no member of $A_x$ is related to any element of X not included in $A_x$. This set $A_x$ is clearly referred to as an neutrosophic equivalence class of $R_N (X, x)$ with respect to x. The members of each neutrosophic equivalence class can be considered neutrosophic equivalent to each other and only to each other under the relation R.*
But here it is pertinent to mention that in some X; (a, b) may not be related at all to be more precise there may be an element a $\in$ X which is such that its relation with several or some elements in X \ {a} is an indeterminate. The elements which cannot determine its relation with other elements will be put in as separate set.

A neutrosophic binary relation that is reflexive, symmetric and transitive is known as a neutrosophic equivalence relation.

Now we proceed on to define Neutrosophic intersections neutrosophic t – norms ($t_N$ – norms)

Let A and B be any two neutrosophic sets, the intersection of A and B is specified in general by a neutrosophic binary operation on the set N = [0, 1] $\cup$ I, that is a function of the form

$$i_N: N \times N \to N.$$

For each element x of the universal set, this function takes as its argument the pair consisting of the elements membership grades in set A and in set B, and yield the membership grade of the element in the set constituting the intersection of A and B. Thus,

$$(A \cap B) (x) = i_N [A(x), B(x)] \text{ for all } x \in X.$$

In order for the function $i_N$ of this form to qualify as a fuzzy intersection, it must possess appropriate properties, which ensure that neutrosophic sets produced by $i_N$ are intuitively acceptable as meaningful fuzzy intersections of any given pair of neutrosophic sets. It turns out that functions known as $t_N$- norms, will be introduced and analyzed in this section. In fact the class of $t_N$-norms is now accepted as equivalent to the class of neutrosophic fuzzy intersection. We will use the terms $t_N$ – norms and neutrosophic intersections inter changeably.

Given a $t_N$ – norm, $i_N$ and neutrosophic sets A and B we have to apply:

$$(A \cap B) (x) = i_N [A (x) , B (x)]$$

for each $x \in X$, to determine the intersection of A and B based upon $i_N$.

However the function $i_N$ is totally independent of x, it depends only on the values A (x) and B(x). Thus we may ignore x and assume that the arguments of $i_N$ are arbitrary numbers a, b $\in$ [0, 1] $\cup$ $I$ = N in the following examination of formal properties of $t_N$-norm.

A neutrosophic intersection/ $t_N$-norm $i_N$ is a binary operation on the unit interval that satisfies at least the following axioms for all a, b, c, d $\in$ N = [0, 1] $\cup$ $I$.

$$
\begin{array}{ll}
1_N & i_N (a, 1) = a \\
2_N & i_N (a, I) = I \\
3_N & b \leq d \text{ implies} \\
& i_N (a, b) \leq i_N (a, d) \\
4_N & i_N (a, b) = i_N (b, a) \\
5_N & i_N (a, i_N(b, d)) = i_N (a, b), d).
\end{array}
$$

We call the conditions $1_N$ to $5_N$ as the axiomatic skeleton for neutrosophic intersections / $t_N$ – norms. Clearly $i_N$ is a continuous function on N \ $I$ and $i_N$ (a, a) < a $\forall a \in$ N \ $I$

$$i_N (I\, I) = I.$$

If $a_1 < a_2$ and $b_1 < b_2$ implies $i_N (a_1, b_1) < i_N (a_2, b_2)$. Several properties in this direction can be derived as in case of t-norms.

The following are some examples of $t_N$ –norms

1.      $i_N (a, b) = \min (a, b)$
        $i_N (a, I) = \min (a, I) = I$ will be called as standard neutrosophic intersection.

2. $i_N$ (a, b) = ab for a, b ∈ N \ *I* and $i_N$ (a, b) = *I* for a, b ∈ N
   where a = *I* or b = *I* will be called as the neutrosophic
   algebraic product.
3. Bounded neutrosophic difference.
   $i_N$ (a, b) = max (0, a + b – 1) for a, b ∈ *I*
   $i_N$ (a, *I*) = *I* is yet another example of $t_N$ – norm.
   4. Drastic neutrosophic intersection
   5.

$$i_N \text{ (a, b) } = \begin{cases} a \text{ when } b = 1 \\ b \text{ when } a = 1 \\ I \text{ when } a = I \\ \quad \text{ or } b = I \\ \quad \text{ or } a = b = I \\ 0 \quad \text{ otherwise} \end{cases}$$

As *I* is an indeterminate adjoined in $t_N$ – norms. It is not easy to
give then the graphs of neutrosophic intersections. Here also we
leave the analysis and study of these $t_N$ – norms (i.e. neutrosophic
intersections) to the reader.

The notion of neutrosophic unions closely parallels that of
neutrosophic intersections. Like neutrosophic intersection the
general neutrosophic union of two neutrosophic sets A and B is
specified by a function

$\mu_N$: N × N → N where N = [0, 1] ∪ *I*.

The argument of this function is the pair consisting of the
membership grade of some element x in the neutrosophic set A
and the membership grade of that some element in the
neutrosophic set B, (here by membership grade we mean not only
the membership grade in the unit interval [0, 1] but also the
indeterminacy of the membership). The function returns the
membership grade of the element in the set A ∪ B.

Thus (A ∪ B) (x) = $\mu_N$ [A (x), B(x)] for all x ∈ X. Properties
that a function $\mu_N$ must satisfy to be initiatively acceptable as
neutrosophic union are exactly the same as properties of functions
that are known. Thus neutrosophic union will be called as
neutrosophic t-co-norm; denoted by $t_N$ – co-norm.

A neutrosophic union / $t_N$ – co-norm $\mu_N$ is a binary operation
on the unit interval that satisfies at least the following conditions
for all a, b, c, d ∈ N = [0, 1] ∪ *I*

$$C_1 \qquad \mu_N(a, 0) = a$$
$$C_2 \qquad \mu_N(a, I) = I$$
$$C_3 \qquad b \leq d \text{ implies}$$
$$\mu_N(a, b) \leq \mu_N(a, d)$$
$$C_4 \qquad \mu_N(a, b) = \mu_N(b, a)$$
$$C_5 \qquad \mu_N(a, \mu_N(b, d))$$
$$= \quad \mu_N(\mu_N(a, b), d).$$

Since the above set of conditions are essentially neutrosophic unions we call it the axiomatic skeleton for neutrosophic unions / $t_N$-co-norms. The addition requirements for neutrosophic unions are

i. $\mu_N$ is a continuous functions on $N \setminus \{I\}$
ii. $\mu_N(a, a) > a$.
iii. $a_1 < a_2$ and $b_1 < b_2$ implies $\mu_N(a_1, b_1) < \mu_N(a_2, b_2)$; $a_1, a_2, b_1, b_2 \in N \setminus \{I\}$

We give some basic neutrosophic unions.
Let $\mu_N : [0, 1] \times [0, 1] \to [0, 1]$

$$\mu_N(a, b) = \max(a, b)$$
$$\mu_N(a, I) = I \text{ is called as the standard}$$
$$\text{neutrosophic union.}$$
$$\mu_N(a, b) = a + b - ab \text{ and}$$
$$\mu_N(a, I) = I.$$

This function will be called as the neutrosophic algebraic sum.

$$\mu_N(a, b) = \min(1, a + b) \text{ and } \mu_N(a, I) = I$$

will be called as the neutrosophic bounded sum. We define the notion of neutrosophic drastic unions

$$\mu_N(a, b) = \begin{cases} a \text{ when } b = 0 \\ b \text{ when } a = 0 \\ I \text{ when } a = I \\ \quad \text{ or } b = I \\ 1 \text{ otherwise.} \end{cases}$$

Now we proceed on to define the notion of neutrosophic Aggregation operators. Neutrosophic aggregation operators on neutrosophic sets are operations by which several neutrosophic sets are combined in a desirable way to produce a single neutrosophic set.

Any neutrosophic aggregation operation on n neutrosophic sets ($n \geq 2$) is defined by a function $h_N: N^n \to N$ where $N = [0, 1] \cup I$ and $N^n = \underbrace{N \times ... \times N}_{n-times}$ when applied to neutrosophic sets $A_1$, $A_2,..., A_n$ defined on X the function $h_N$ produces an aggregate neutrosophic set A by operating on the membership grades of these sets for each $x \in X$ (Here also by the term membership grades we mean not only the membership grades from the unit interval [0, 1] but also the indeterminacy $I$ for some $x \in X$ are included). Thus

$$A_N (x) = h_N (A_1 (x), A_2 (x),..., A_n(x))$$

for each $x \in X$.

## 4.2  Optimization of NRE with max-product composition

The study of the neutrosophic relational equations

$$x \text{ o } A = b$$

where $A = (a_{ij})_{m \times n}$ neutrosophic matrix with entries from [0 1] $\cup$ FN, $b = (b_1, ..., b_n)$ $b_i \in$ [0 1] $\cup$ F N and 'o' is the max-min composition.

Even in case of fuzzy relation equation the resolution of the equation x o A = b is an interesting on going research topic. [61] instead of finding all solutions of x o A = b, they assume f (x) to be the user criterion function they solve the non linear programming model with fuzzy relation constraints

$$\begin{align} &\text{Min f(x)} \\ \text{such that} \quad &x \text{ o } A = b \end{align} \qquad (1)$$

A minimizer will be the best solution. At this juncture when we seek to transform this model in to a neutrosophic model we are not fully aware of how the solution looks like will we get one best solution or several such solutions. All these are suggested as problems in chapter 5. Thus contrary to fuzzy relation constrains our problem

"min f(x)
such that          x o A = b"

subjects to neutrosophic relation constrains also.

It is left once again for the reader to establish that when the solution set of the neutrosophic relation equation x o A = b is not empty the existence of a minimal solution. (In fact we have a finite number of minimal solution it can be established in an analogous way as in FRE under certain constraints).

The optimal solution can be obtained as in case of fuzzy relation equations. The study and analysis of finding solution to the problem (1) with a non linear objective function will be termed as a non linear optimization problem with neutrosophic relation constrains.

The process of obtaining the solution is similar to that done in the case of fuzzy relation equation. This study itself can be written as a book we leave the interested reader to develop. However we would be bringing forth a separate book on this topic.

### 4.3 Method of solution to NRE in a complete Brouwerian lattices

Xue-ping Wang [110] has given a method of solution to fuzzy relation equations in a complete Brouwerian lattice. Here we give a method of solution to neutrosophic relation equations in a complete Brouwerian neutrosophic lattice.

Unfortunately, how to solve a fuzzy relation equation in a complete Brouwerian lattice is still an open problem so it is an open problem to study neutrosophic relation equation in a complete Brouwerian neutrosophic lattice.

We only propose a few problems as even in case of fuzzy relation equations the method of solution is open.

### 4.4 Multi objective optimization problem with NRE constraints

We define analogous to multi objective optimization problems with fuzzy relation equation the procedure of solution to neutrosophic relational equation.

Let $X_n = \{[0\ 1] \cup FN\}^m$. Let A = m × n neutrosophic fuzzy matrix $I = \{1, 2,…, m\}$ and $J = \{1, 2,…, n\}$. b = $[b_j]_{1 \times n}$ such that $a_{ij} \in [0, 1] \cup FN$ for all $i \in I$ and $j \in J$. Given A and b, a system of neutrosophic relation equation defined

$$x \text{ o } A = b \qquad\qquad (1)$$

where 'o' represents the max-min composition defined in chapter 3 of this book for neutrosophic min max and relational equations.

A solution to (1) is a neutrosophic vector $x = (x_1, \ldots, x_m)$, $x_i \in [0\ 1] \cup FN$ such that

$$\max_{i \in I} [\min (x_i, a_{ij})] = b_j \ \forall \ j \in J \qquad (2)$$

In other words the optimization problem we are interested in, has the following form

$$\begin{aligned} &\text{Min } \{f_1(x), f_2(x), \ldots, f_k(x)\} \\ &\text{Such that} \quad x \ o \ A = b \end{aligned} \qquad (3)$$

$x_i \in [0\ 1] \cup [F\ N]$, $i \in I$ where $f_k(x)$ is an objective function, $k \in K = \{1, 2, \ldots p\}$. In case of FRE, the properties of efficient points were investigated by [110].

Here $f_k$'s are the neutrosophic criterion vectors. The problem (3) is restated as follows. Let $X_N$ be the feasible domain

$$X_N = \{x \in R^n \cup \{IR\}^n \ / \ x \ o \ A = b\}, \ x_i \in R^n \cup \{IR\}^n\}.$$

For each $x \in X_N$ we say x is a solution vector and define

$$z = \{f_1(x), f_2(x), \ldots, f_p(x)\}$$

to be its neutrosophic criterion vector. Moreover define

$$Z_N = \{\{z \in R^p \cup \{IR\}^p \ / \ (f_1(x), f_2(x), \ldots, f_p(x)) \text{ for some } x \in X_N\},$$

z may be indeterminate or may be in $R^P_N$.

**DEFINITION 4.4.1:** *A point $\bar{x} \in X_N$ is an efficient or a pare to optimal solution to problem (3) if and only if there does not exist any $x \in X_N$ such that $f_k(x) \leq f_k(\bar{x})$ for all $k \in K$, and $f_k(x) < f_k(\bar{x})$ for at least one k; otherwise $\bar{x}$ is an inefficient solution.*

**DEFINITION 4.4.2:** *Let $Z^1, Z^2 \in Z_N$ be two neutrosophic criterion vectors. Then $Z^1$ dominates $Z^2$ if and only if $Z^1 \leq Z^2$ and $Z^1 \neq Z^2$. That is $Z^1_k \leq Z^2_k$ for all $k \in K$ and $Z^1_k < Z^2_k$ for at least one k.*

*Note:* Here comparison is possible only within criterion vectors or neutrosophic criterion vectors if one is just a criterion vector and other a neutrosophic criterion vector comparison is not possible.

**DEFINITION 4.4.3:** *Let $z \in Z_N$. Then $\bar{z}$ is non-dominated if and only if there does not exists any $z \in Z_N$ that dominates. $\bar{z}$ is a dominated neutrosophic criterion vector.*

The idea of dominance is applied to the neutrosophic criterion vectors whereas the idea of efficiency is applied to the solution vectors. A point $\bar{x} \in X_N$ is efficient if its criterion vector is non-dominated in $Z_N$. The set of all effective points is called the neutrosophic efficient set or neutrosophic p are to optimal set. The set of all non-dominated neutrosophic criterion vectors is called the non-dominated neutrosophic set.

As in case of FREs a system of NREs may be manipulated in a way such that the required computational effort of the proposed genetic algorithm is reduced. Due to the requirement of x o A = b some components of every solution neutrosophic vector or vector may have to assume a specific value of an indeterminate value. As in case of FREs the genetic operators are applied to the reduced problem.

We at the first stage divide the problem into two components i.e. the components of vector b is divided into b' and b", b' having fuzzy values and b" having neutrosophic values so that b' and b" can be ordered and the matrix A is correspondingly rearranged. Notice the corresponding $\hat{x}$, and $\hat{x}$" which are the maximum solution can be obtained by the following formula

$$\hat{x}' = A \text{ @ b'} = \left[ \overset{n'}{\underset{j=1}{\wedge}} (a'_{ij} \text{ @ } b'_j) \right]$$

where '$\wedge$' is the min operator

$$\hat{x}'' = A \text{ @ b''} = \left[ \overset{n''}{\underset{j=1}{\wedge}} (a''_{ij} \text{ @ } b''_j) \right] \quad (4)$$

$$a'_{ij} \text{ @ } b'_j = \begin{cases} 1 \; if \; a'_{ij} \le b'_j \\ b'_j \; if \; a'_j > b'_j \end{cases}$$

$$\text{a"}_{ij} @ \text{b"}_j = \begin{cases} 1 \; if \; a"_{ij} \le b"_j \\ b"_j \; if \; a"_j > b"_j \end{cases} \tag{5}$$

These $\hat{x}$', and $\hat{x}$" can be used to check whether the feasible domain is empty.

If max [min ($\hat{x}'_i$, $a'_{ij}$)] = b'$_j$ $\forall_j \in$ J

(6')

and max [min ($\hat{x}"_i$, $a"_{ij}$)] = b"$_j$ $\forall_j \in$ J (6")

then $\hat{x}$' and $\hat{x}$" is the maximum solution (1) $\hat{x} = \hat{x}' \cup \hat{x}$" is the solution; other wise problem (3) is infeasible.

**DEFINITION 4.4.4:** If *a value change in some elements of a given neutrosophic relation matrix, A has no effect on the solution of the corresponding neutrosophic relation equations this value change is called an N-equivalence operation.*

**DEFINITION 4.4.5:** *Given a system of neutrosophic relation equation.*

*(1) a N-pseudo – characteristic matrix P = [P$_{ij}$] is defined as*

$$P_{ij} = \begin{cases} 1 & if \; a'_{ij} > b'_j \\ 0 & if \; a'_{ij} = b'_j \\ -1 & if \; a'_{ij} < b'_j \\ I & if \; a"_{ij} > b"_j \\ 0 & if \; a"_{ij} = b"_j \\ -I & if \; a"_{ij} < b"_j \end{cases}$$

*The N-pseudo characteristic matrix will be referred to as a NP – matrix.*

**DEFINITION 4.4.6:** *Given the maximum solution $\hat{x} = \hat{x}' \cup \hat{x}$" if there exists some i $\in$ I and some j $\in$ J such that $\hat{x}'_i \wedge a'_{ij} = b'_j$ and $\hat{x}"_i \wedge a"_{ij} = b"_j$ then the corresponding a'$_{ij}$ and a"$_{ij}$ of matrix A is called a N-critical element for $\hat{x} = \hat{x}' \cup \hat{x}$".*

The results and definitions related to this can be carried out in an analogous way as in case of multi objective optimization problems with fuzzy relation equation constraints.

### 4.5 Neural Neutrosophic relational system with a new learning algorithm

A neutrosophic t-norm (Triangular norm denoted by Nt – norm or $t_N$ – norm) is a function.

$t_N$: $[0\ 1] \cup I \times [0\ 1] \cup I \rightarrow [0\ 1] \cup I$ satisfying for any a, b, c, d $\in$ $[0, 1] \cup I$ the following conditions

   i.     $t_N$ (a, 1) = a and $t_N$ (a, 0) = 0
   ii.    If a $\leq$ d then it implies $t_N$ (a, b) $\leq$ $t_N$ (a, d)
   iii.   $t_N$ (a, $I$) = $I$
   iv.    $t_N$ (a, b) = $t_N$ (b, a)
   v.     $t_N$ (a, $t_N$ (b, d)) = $t_N$ ($t_N$ (a, b), d).

Moreover $t_N$ is called Neutrosophic Archimedean if and only if

   i.     $t_N$ is a continuous function.
   ii.    $t_N$ (a $t_N$ (a, a)) < a $\forall$ a $\in$ (0, 1).

It is left as an exercise for the reader to study and analyze the properties of $t_N$-norm. The general structure of a conventional neutrosophic neuron (N-neuron) can be shown in the figure.

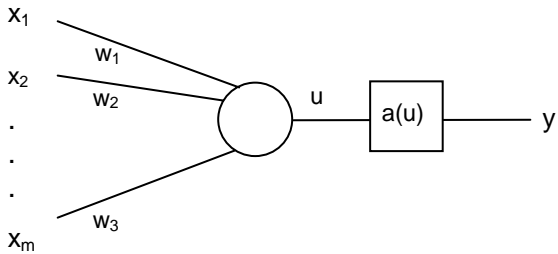

Figure: 1.9.1

   The equation that describes this kind of N-neuron is as follows

$$y = a\left[\sum_{i=1}^{m} w_i x_i + \partial\right] \qquad (1)$$

where a is a non linearity $\partial$ is the threshold and $w_i$ are weights which can also be indeterminates i.e., $I$ that can change on line with the aid of a learning process.

The compositional neuron has the same structure with the neuron of equation (1) but it can still be described by the equation

$$y = a\left[\underset{J \in N}{S}(x_j w_j)\right] \tag{2}$$

where S is a fuzzy union operator as $S_N$-norm, $t_N$ is a neutrosophic intersection operator and a is the activation function

$$a(x) = \begin{cases} 0 & x \in (-\infty, 0) \\ x & x \in [0, 1] \\ 1 & x \in (0, \infty) \\ 0 & x \in (-\infty I, 0) \\ xI & x \in [0, I] \\ I & x \in (0, \infty I) \end{cases}$$

The analogous properties work in this situation as in case of neural fuzzy relational system with a new learning algorithm.

## 4.6 Unattainable solution of NRE

Hideyuki Imai et al [39] have obtained a necessary and sufficient conditions for the existence of a partially attainable and an unattainable solution.

Here we study the problem in the context of NRE. Let U and V be nonempty sets and let N(U), N(V), N(U × V) be the collection of neutrosophic sets of U, V and U × V respectively. The equation

$$X \circ A = B \tag{*}$$

is a neutrosophic relational equation where A $\in$ N (U × V) and B $\in$ N(V) are given and X $\in$ N(U) is unknown and 'o' denotes the usual $\wedge$-$\vee$ composition. A neutrosophic set X satisfying the equation above is called a solution of the equation.

If $\mu_X : U \to FN \cup [0\ 1]$

$\mu_A : U \times V \to FN \cup [0\ 1]$ and $\mu_B : U \times V \to FN \cup [0\ 1]$ are their neutrosophic membership functions. Now equation (*) is as follows.

$V(\mu_X(u) \wedge \mu_A(u, v) = \mu_B(v))$ ($\forall v \in V$). The solution set in case of FRE has been investigated by several researchers and several important conclusions have been derived by the analysis of the solution set in case of NRE is at its dormant state. It will be an interesting and innovative research to study and analyze the solution set in case of NRE.

We propose several problems about the solution set in case of NRE. When both the sets U and V are both finite sets, show that the solution set is completely determined by the greatest solution and the set of minimal solutions. Now we proceed on to define the neutrosophic membership function (N-membership function) of the neutrosophic sets X, Y, $\in N(U)$.

**DEFINITION 4.6.1:** *Let $\mu_x$ and $\mu_y$ be the neutrosophic membership functions of the neutrosophic set X, Y $\in N(U)$ respectively. Then the partial order $\leq$ the join $\vee$ the meet $\wedge$ are defined as follows.*

$$\mu_X \leq \mu_y \Leftrightarrow \text{for all } u \in U \ (\mu_X(u) \leq \mu_y(u))$$
$$\mu_X \vee \mu_Y : \ni u \rightarrow \mu_X(u) \vee \mu_Y(u) \in FN \cup [0\ 1]$$
$$\mu_X \wedge \mu_Y : \ni u \rightarrow \mu_X(u) \wedge \mu_Y(u) \in FN \cup [0\ 1].$$

*Note that $\mu_X < \mu_Y$ is equivalent to $X \subset Y$ for X, Y $\in N(U)$. The greatest solution and minimal solution for NRE minimal solution are defined as in case of FRE. We can in case of NRE also define the attainable and unattainable solution.*

As they can be derived as a matter of routine it is left as an exercise for the reader.

## 4.7 Specificity Shift in solving NRE

The specificity shift method can be classified as an approach situated between analytical and numerical method of solving FRE. When we study or impose specificity shift for solving NRE we call them neutrosophic specificity shift (N-specificity shift).

$$X \ \square \ R = y \qquad (1)$$

which denotes a FRE with max t-composition where t is assumed to be a continuous t-norm while X R and y are viewed as fuzzy set and a fuzzy relation defined in finite universe of discourse.

We take X, y and R to be neutrosophic sets and the norm t is assumed to be a continuous $t_N$-norm. [14] have studied in case of FRE. Here we just sketch the method analogous to that given by [14] for FRE in case of NRE. Given a collection of neutrosophic data $((x(1), y(1)), (x(2), y(2),…, (x(N), y(N))$ where $x(k) \in \{[0\ 1] \cup FN\}^n$ and $y(k) \in \{[0\ 1] \cup FN\}^m$.

Determine a neutrosophic relation R satisfying the collection of the relational constrains (Neutrosophic relational equations).

$$x(k) \ \Box \ R = y\ (k)$$

Expressing the above equation in terms of corresponding membership functions of $x(k)$, $y(k)$ and R we drive

$$y_j\ (k) = \bigvee_{i=1}^{n} (x_i\ (k)\ t\ r_{ij})$$

where $k = 1, 2, …, N; j = 1, 2, …, m$.

The problem has to be talked using interpolation under the special assumptions and conditions. Further one cannot always be guaranteed of a solution. The solution may or may not exist.

It is appropriate to study whether the maximal neutrosophic (fuzzy) relation with the membership function equal to

$$R = \bigcap_{k=1}^{N} (x(k) \rightarrow y(k))$$

exist. This task is also left as an exercise to the reader.

It is pertinent to mention here that the computations involve an intersection of the individual neutrosophic relations determined via a pseudo complement associated with the $t_N$ –norm. The study and analysis of the problem under the neutrosophic setup happens to be very different in only certain situations. When the notion of indeterminacy plays a vital role in the problem we see solutions so obtained may be entirely different from the usual fuzzy relation equations. As use of the optimization methods leads to better approximate solutions yet the entire procedure with indeterminacy could be quite time consuming. Further more as there could be a multiplicity of solutions such approaches usually identify only one of them leaving the rest of them unknowns or as a indeterminacy. The entire investigated setting can be changed from the fuzzy model to the neutrosophic model with appropriate modification to ingrain the element of indeterminacy.

### 4.8 NRE with deneutrofication algorithm for the largest solution

Kagei [42] has provided algorithm for solving a new fuzzy relational equation including defuzzificaiton. Here we extend it in case of a new neutrosophic relational equation and adopt in our study the notion of deneutrofication. An input neutrosophic set is first transformed into an internal neutrosophic set by a neutrosophic relation. The internal neutrosophic relation is obtained from neutrosophic input and deneutrosified output.

The problem is classified under two types in case of fuzzy setup. In type *I* problem there exists nontrivial largest solution. For type *II* the largest solution is trivial. The neutrosophic relational equations which are used analogous to fuzzy relational equation is

$$q_\lambda^N = p_\lambda^N \circ R_N \text{ where } p_\lambda^N \text{ and } q_\lambda^N$$

are neutrosophic sets and X and Y respectively ($\lambda$ is an index of the equation), $R_N$ is a neutrosophic relation from X to Y to be solved and '$\circ_N$' is a neutrosophic composition operator defined analogous to a fuzzy composition operator.

Study the existence or non-existence of the largest solution of $R_N$ when the neutrosophic set is defined as a mapping from a nonempty set to a complete Brouwerian neutrosophic lattice. Some researchers used fuzzy sets $q_\lambda$ as output data. We can use the neutrosophic sets $q_\lambda^N$ as output data. Study and analyze the output data. In certain cases the output data can itself be a deneutrosified data.

Now we shall indicate how to solve the neutrosophic relation $R_N$ when the system includes deneutrosophic processes.

In these systems the neutrosophic membership values have to be compared with each other for the deneutrosofication of internal neutrosophic sets. Therefore a complete totally ordered set must be used this is done only as two steps one with fuzzy values and another only with neutrosophic values as both these values cannot be compared with each other. Thus the procedure used for fuzzy sets can be adopted, by replacing U by $U_N$ a union the complete subsets of the unit interval [0 1] and FN and all analogous operations are performed to attain a solution.

Thus we assume we may arrive at a unique solution which is the unique solution given by fuzzy sets and the other given by the neutrosophic sets by giving a unique output.

## 4.9 Solving neutrosophic relation equation with a linear objective function

[24] has given an optimization model with a linear objective function subject to a system of FRE. Solving such equations is not an easy job due to the non-convexity of its feasible domain. Now we define analogous methods for NRE.

Let $A_N = (a_{ij})$; $a_{ij} \in FN \cup [0\ 1]$ be an $m \times n$ dimensional fuzzy matrix and $b_N = (b_1,\ldots, b_n)^T$, $b_{ij} \in [0\ 1] \cup FN$ be an n-dimensional neutrosophic vector, then the following system of FRE is defined by $A_N$ and $b_N$

$$x \ o \ A_N = b_N \qquad (1)$$

where 'o' denotes the usual min max composition defined in Chapter 3. In other words we try to find a neutrosophic vector $x = (x_1,\ldots, x_m)^T$ with $x_i \in [0\ 1] \cup FN$ such that

$$\max \min (x_i, a_{ij}) = b_j \quad (j = 1, 2, \ldots, n) \qquad (2)$$

The resolution of the NRE (1) is an interesting and a undeveloped research topic at it is very dormant state as even in case of FRE the study of the problem and the research is only on going.

Let $C_N = (c_1,\ldots, c_m)^T \in R^m \cup \{IR\}^m$ be a m-dimensional vector where $c_i$ represents the weight (or cost) associated with variable $x_i$ for $i = 1, 2,\ldots,n$. We consider the following optimization problem.

Minimize $$Z_n = \sum_{i=1}^{n} c_i x_i \qquad (3)$$

such that $x \ o \ A_N = b_N$, $x_i \in [0\ 1] \cup FN$ compare it with regular linear programming and fuzzy linear programming this linear optimization problem subject to NRE is very different from FRE and are of different nature.

Note the feasible domain of problem (3) is the solution set of system (1). We denote it by $X (A_N \ b_N) = \{x = (x_1, \ldots, x_m)^T / (x_1, \ldots, x_m)^T \in R^m \cup (IR)^m$ such that $x \ o \ A_N = b_N x_i \in [0\ 1] \cup FN\}$.

To characterize $X (A_N, b_N)$ we define $I = (1, 2, \ldots, m)$, $J = \{1, 2, \ldots, n\}$ and $X_N = \{x \in R^m \cup (IR)^m / x_i \in FN \cup [0\ 1]$ for all

$i \in I$}. For $x^1$, $x^2 \in X_N$ we say $x^1 \leq x^2$ if and only if $x_i^1 \leq x_i^2$, $\forall i \in I$. In this way '$\leq$' forms a partial order relation on $X_N$ and $(X_N, \leq)$ becomes a neutrosophic lattice. We call $\hat{x} \in X_N$ $(A_N, b_N)$ a maximum solution, if $x \leq \hat{x}$ for all $x \in X_N$ $(A_N$ $b_N)$. Similarly $\hat{x} \in X_N$ $(A_N, b_N)$ is called a minimum solution if $x \leq \hat{x}$ for all $x \in X_N$ $(A_N, b_N)$. When $X_N$ $(A_N, b_N) \neq \phi$ it can be completely determined by one maximum solution and a finite number of minimum solutions.

The maximum solution can be obtained by assigning

$$\hat{x} \quad = A_N \ @ \ b_N = \left[ \overset{n}{\underset{j=1}{\wedge}} (a_{ij} \ @ \ b_j) \right]_{i \in I}$$

where $a_{ij} \ @ \ b_j = \begin{cases} 1 \ if \ a_{ij} \leq b_j \\ b_j \ if \ a_{ij} > b_j \\ I \ if \ a_{ij} \ is \ incomparable \ with \ b_j \end{cases}$

Suppose we denote the set of all minimum solution by $\widetilde{X}_N$ $(A_N, b_N)$ then

$$\widetilde{X}_N (A_N, b_N) = \bigcup_{\bar{x} \in \bar{X}_N (A_N, b_N)} \left\{ x \in X_N \mid \breve{x} \leq x \leq \hat{x} \right\}.$$

Now we take a close look at $X_N$ $(A_N, b_N)$. All analogous results in the Neutrosophic setup is left as exercise for the reader.

### 4.10 Some properties of minimal solution for a NRE.

NRE occurs in practical problems for instance in neutrosophic reasoning. Therefore it is necessary to investigate properties of the set of solutions [39] have given a necessary and sufficient conditions for existence of a minimal solution of a FRE defined on an infinite index sets.

Here we suggest analogous results in case of NRE. Let I and J be the index set and let $A_N = (a_{ij})$ be a coefficient matrix $b_N = (b_j)$ be a constant vector where $i \in I$ and $j \in J$ the equation

$$x \ o \ A_N = b_N \qquad (1)$$

or $\qquad \vee (x_i \wedge a_{ij}) = b_j$

for all j ∈ J is called the Neutrosophic relation equation where o denotes the sup-min composition and all $x_i$, $b_j$, $a_{ij}$ are in [0 1] ∪ FN. An x which satisfies (1) is called a solution of equation (1)

**DEFINITION 4.10.1:** *Let $(P_N, \leq)$ be a partially ordered set and $X_N \subset P_N$. A minimal element of $X_N$ is an element $p \in X_N$ such that $x < p$ for $x \in X_N$. The greatest element of $X_N$ is an element $g \in X_N$ such that $x \leq g$ for all $x \in X_N$.*

**DEFINITION 4.10.2:** *Let $a_N = (a_{ij})$ and $b_N = (b_{ij})$ be neutrosophic vectors. Then the partial order $\leq$ the join $\vee$ and the meet $\wedge$ are defined as follows.*

$$a_N \leq b_N \Leftrightarrow a_i \leq b_i \text{ for all } i \in I \ a_N \vee b_N \ \underset{-}{\Delta},$$

$$(a_i \vee b_i) \ a_N \wedge b_N \ \underset{-}{\Delta} \ a_i \wedge b_i).$$

**DEFINITION 4.10.3:** *Let $\{[0\ 1]^l \cup FN^l, \leq\}$ be a poset with the partial order given in definition 4.10.2. and let $x_N \subset [0\ 1]^l \cup [FN]^l$ be the solution set of equation (1). The greatest element of $x_N$ is minimal element of $x_N$ and $x^0_N$ is called the greatest solution; which denote a minimal solution and a set of minimal solutions of Equation (1) respectively.*

**DEFINITION 4.10.4:** *For a, b, ∈ [0 1] ∪ FN*

$$a_N \ \alpha \ b_N \ \underset{-}{\Delta} \ \begin{cases} 1 \text{ if } a_N \leq b_N \\ b_N \text{ otherwise} \\ I \text{ if } a_N \text{ and } b_N \text{ are incomparable} \end{cases}$$

*Moreover*

$$A_N @ b^{-1}_N \ \underset{-}{\Delta} \ \left[ \underset{j \in J}{\wedge} a_{ij} o b_j \right]$$

*where $b^{-1}_N$ denotes the transposition of vector $b_N$.*

**DEFINITION 4.10.5:** *The solution $x = (x_{ij}) \in x_N$ is N-attainable for $j \in J$ if there exists $i_j \in I$ such that $x_i \wedge a_{ij} = b_j$ and the solution $x = (x_i) \in \aleph_N$ is unattainable for $j \in J$ if $x_i \wedge a_{ij} < b_j$ for all $i \in I$.*

**DEFINITION 4.10.6:** *The solution $x \in \aleph_N$ is called a N-attainable solution if x is N-attainable for $j \in J$, the solution is called an*

*N-unattainable solution if x is $(\aleph_N)_j^{(-1)}$ for $j \in J$ and the solution*
*$x \in X$ is called a N-partially attainable solution if $x \in \aleph_N$ is*
*neither an N-attainable solution nor an N-unattainable solution.*
*In other words, $x \in \aleph_N$ is an N-attainable solution if and only if*
*$x \in (\aleph_N)_j^{(+)}$ $x \in \aleph_N$ is an N-unattainable solution if and only if*
*$x \in (\aleph_N)_j^{(-1)}$ $x \in \aleph_N$ is a N-partially attainable solution if and*
*only if $x \in \aleph_N - (\aleph_N)_j^{(+)} - (\aleph_N)_j^{(-1)}$.*

The set of all N-partially attainable solution is denoted by $(X_N)_j^{(*)}$. All properties related to attainable, unattainable partially attainable solutions can be also be derived with appropriate modifications in case of N-attainable, N-unattainable and N-partially attainable solution. Some of these results are proposed as problems in chapter V.

## 4.11  Applications of NRE to Real World Problems

In this section we hint the applications of NRE to various real world problems like flow rate in chemical plants, transportation problem, study of bonded labor problem study of interrelations among HIV/ AIDS affected patients and use of genetic algorithms in chemical problems.

### 4.11.1 Use of NRE in chemical engineering

The use of FRE for the first line has been used in the study of flow rates in chemical plants. In this study we are only guaranteed of a solution but when we use NRE in study of flow rates in the chemical plants we are also made to understand that certain flow rates are indeterminates depending on the leakage, chemical reactions and the new effect due to chemical reactions which may change due to change in the density/ viscosity of the fluid under study their by changing the flow rates while analyzing as a mathematical model. So in the study of flow rates in chemical plants some indeterminacy are also related with it. FRE has its own limitation for it cannot involve in its analysis the indeterminacy factor.

We have given analysis in chapter 2 using FRE. Now we suggest the use of NRE and bring out its importance in the determination of flow rates in chemical plants.

Consider the binary neutrosophic relations $P_N$ (X, Y) $Q_N$ (Y, Z) and R (X, Z) which are defined on the sets X, Y and Z. Let the membership matrices of P, Q and R be denoted by P = $[p_{ij}]$, Q = $[q_{jk}]$ and R = $[r_{ij}]$ respectively where $p_{ij}$ = $P(x_i, y_j)$, $q_{jk}$ = Q $(y_j, r_k)$ and $r_{ik}$ = R $(x_i, z_k)$ for i∈$I$ = $N_n$, j∈J = $N_m$ and k ∈ K = $N_k$ entries of P, Q and R are taken for the interval [0 1] × FN. The three neutrosophic matrices constrain each other by the equation

$$P \text{ o } Q = R \qquad (1)$$

where 'o' denotes the max-min composition (1) known as the Neutrosophic Relational Equation (NRE) which represents the set of equation

$$\max p_{ij} \, q_{jk} = r_{ik}. \qquad (2)$$

For all i ∈ $N_n$ and k ∈ $N_s$. After partitioning the matrix and solving the equation (1) yields maximum of $q_{jk} < r_{ik}$ for some $q_{jk}$, then this set of equation has no solution so to solve equation (2) we invent and redefine a feed – forward neural networks of one layer with n-neurons with m inputs. The inputs are associated with $w_{ij}$ called weights, which may be real, or indeterminates from R$I$. The neutrosophic activation function $f_N$ is defined by

$$f_N(a) = \begin{cases} 0 \; if \; a < 0 \\ a \; if \; a \in [01] \\ 1 \; if \; a > 1 \\ aI \; if \; a \in FN \\ I \; if \; aI > I \\ 0 \; if \; in \; aI, a < 0. \end{cases}$$

The out put $y_i = f_N$ (max $w_{ij} x_j$). Now the NRE is used to estimate the flow rates in a chemical plant. In places where the indeterminacy is involved the expert can be very careful and use methods to overcome indeterminacy by adopting more and more constraints which have not been given proper representation and their by finding means to eliminate the indeterminacy involved in the weights.

In case of impossibility to eliminate these indeterminacy one can use the maximum caution in dealing with these values which are indeterminates so that all types of economic and time loss can be met with great care. In the flow rate problem the use of NRE mainly predicts the presence of the indeterminacy which can be minimized using $f_N$; where by all other in-descrepancies are given due representation.

We suggest the use of NRE for when flow rates are concerned in any chemical plant the due weightage must be given the quality of chemicals or raw materials which in many cases are not up to expectations, leakage of pipe, the viscosity or density after chemical reaction time factor, which is related with time temperature and pressure for which usually due representations, is not given only ideal conditions are assumed. Thus use of NRE may prevent accident, economic loss and other conditions and so on.

### 4.11.2 Use of NRE to determine the peak hours of the day for transport system.

In an analogous way to modified fuzzy relation equations we given a sketch of the modified form of Neutrosophic relation equations and analyze the passenger preference for a particular hour of a day. Since the very term preference is a fuzzy term and sometimes even an indeterminate one we at the out set are justified in using these neutrosophic relational equations. Let $P_N o Q_N = R_N$ be any neutrosophic relational equation where P, Q and R are neutrosophic matrices. We reduce the NRE, $P_N o Q_N = R_N$ into NREs

$$P_N^1 \; o \; Q_N^1 = R_N^1, \, P_N^2 \; o \; Q_N^2 = R_N^2, ..., P_N^S \; o \; Q_N^S = R_N^S$$

where $Q_N = Q_N^1 \bigcup Q_N^2 \bigcup ... \bigcup Q_N^S$ such that $Q_N^i \bigcap Q_N^j = \phi$ if $i \neq j$. Hence by this method we get S-preferences.

We briefly describe the modified or new NRE used here. We know the NRE can be represented by

$$P_N o Q_N = R_N \tag{1}$$

where 'o' is the max product composition with $P_N = [p_{ij}]$, $Q_N = [q_{jk}]$, $R_N = [r_{ik}]$ with $i \in N_n$, $j \in N_m$ and $k \in K_s$. We want to determine $P_N$ (All $P_N$, $Q_N$ and $R_N$ are neutrosophic fuzzy matrices i.e. the entries of $P_N$, $Q_N$ and $R_N$ are from [0 1] $\cup$ F N).

Equation 1 gives

$$\max_{j \in N_m} p_{ij} \ q_{jk} = r_{ik} \qquad (2)$$

for all i $\in$ $N_n$ and k $\in$ $N_s$ to solve equation (2) for $p_{ij}$ we use the linear activation function $f_N$ for all a $\in$ $R_N$ by

$$f_N (a) = \begin{cases} 0 & if \ a < 0 \\ a & if \ a \in [01] \\ aI & if \ a \in [0I] \\ 1 & if \ a > 1 \\ I & if \ a > I. \end{cases}$$

We as in case of FRE work with equation 1 so that the error function becomes very close to zero. The solution is then expressed by the weight $w_{ij}$ as $p_{ij} = w_{ij}$ for all i $\in$ $N_n$ and j $\in$ $N_m$; $P_N = (w_{ij})$ is the neutrosophic n $\times$ n matrix.

Thus by transforming the single neutrosophic relation equation into a collection of NREs we analyze each

$$P_N^i \ oQ_N^i = R_N^i$$

for $1 \leq i \leq S$ and obtain the preferences. The preferences so calculated can also be indeterminates when the corresponding weights are indeterminates i.e. from FN. For the data discussed in chapter two of this book we have taken a nice preferences. We may have data were the preferences of the day may be very fluctuating say for one day the number of passengers in the same route for a particular hour may be 144 and for the same hour in the same route for the same hour on some other day the number of passenger may be in a single digit in such cases we have to use indeterminacy and work with NRE for average will not serve the purpose for when they ply on that hour it will be a very heavy loss so that if that hour for a particular route is kept as an indeterminate the government may choose to run a bus or not to run to save money or avoid loss. Thus when the indeterminacy is present we are in a position to use the indeterminacy factor and accordingly work with caution so that unnecessary economic loss is averted. Thus NREs prove itself to be useful when the real data under analysis is a fluctuating one for a particular period.

### 4.11.3 Use of NRE to the study the problems of bonded labor

Now we apply it to the problems faced by the silk weavers who are bonded labourers using NRE. We take the attributes associated with this problem only as given in chapter 2. Now we use the NRE $P_N$ o $Q_N = R_N$ where $P_N$, $Q_N$ and $R_N$ are neutrosophic matrices.

Now we determine the neutrosophic matrix associated with the attributes relating the bonded labourers and the owners using NRE.

$$P_N = \begin{array}{c} \\ B_1 \\ B_2 \\ B_3 \\ B_4 \\ B_5 \\ B_6 \end{array} \begin{array}{cccc} O_1 & O_2 & O_3 & O_4 \\ \left[\begin{array}{cccc} .6 & 0 & .3I & O \\ .7 & .4 & .3 & .I \\ .3 & .4 & .3 & .3 \\ .3I & 0 & .3 & .4I \\ .8 & .4I & .2 & .4 \\ 0 & .4 & .5 & .9 \end{array}\right] \end{array}$$

Suppose the $Q^T_N = [.6, .5, .7, .9]$. Now $P_N$ and $Q_N$ are known in the neutrosophic relational equation $P_N$ o $Q_N = R_N$. Using the max-min principle in the equation $P_N$ o $Q_N = R_N$ we get $R^T_N = \{.6, .8I, .4, .4I, .6, .9\}$.

In the neutrosophic relational equation $P_N$ o $Q_N = R_N$, $P_N$ corresponds to the weightages of the expert, $Q_N$ is the profit the owner expects and $R_N$ is the calculated or the resultant giving the status of the bonded labourers. Now we have taken a neutrosophic vector in which the demand for the finished goods is in the indeterminate state for we see people at large do not seek now a days hand woven materials to the machine woven ones for the reasons very well known to them so we see we cannot say how far the demand is for such goods so we give it the highest indeterminate value. Now we consider the next highest value for globalization followed by availability of raw goods and finally profit or no loss.

Using these we obtain the neutrosophic resultant vector $R^T_N = \{.6, .8I, .4, .4I, .6, .9\}$ where hours of days work is the highest and the advent of power looms and globalization has made them still poorer is an higher indeterminacy followed by no knowledge of any other work has made them only bonded but live in penury and government interferes and frees them they don't have any work and government does not give them any alternative job remains next maximum.

Thus we have given only one illustration one can work with several of the possible neutrosophic vectors and derive the resultants. Several experts opinion can be taken and their neutrosophic resultant vectors can be determined. We have just given illustration of one case.

Now we give yet another relation between 8 symptoms and 10 patients; $P_1$, $P_2$,…, $P_{10}$. The 8 symptoms which are taken are $S_1$,…, $S_8$ are given as follows.

|  |  |  |
|---|---|---|
| $S_1$ | - | Disabled |
| $S_2$ | - | Difficult to cope with |
| $S_3$ | - | Dependent on others |
| $S_4$ | - | Apathetic and unconcerned |
| $S_5$ | - | Blaming oneself |
| $S_6$ | - | Very ill |
| $S_7$ | - | Depressed |
| $S_8$ | - | Anxious and worried |

We study the 10 patients $P_1$,…, $P_{10}$ related to 8 symptoms which they suffer. It is pertinent to mention that all the patients may not suffer all types of diseases / symptoms some of the diseases / symptoms they suffer may be an indeterminate.

Now using the experts opinion who is the ward doctor we give the related neutrosophic matrix $P_N$ with weights $w_{ij}$. $P_N$ is a $8 \times 10$ matrix.

|  | $P_1$ | $P_2$ | $P_3$ | $P_4$ | $P_5$ | $P_6$ | $P_7$ | $P_8$ | $P_9$ | $P_{10}$ |
|---|---|---|---|---|---|---|---|---|---|---|
| $S_1$ | 0 | 0 | .2I | .5 | 0 | 0 | .6 | .7 | 0 | .5I |
| $S_2$ | 0 | 0 | 0 | 0 | 0 | .2 | 1 | 0 | .9 | .6 |
| $S_3$ | .5I | 0 | 0 | 0 | .9 | 0 | 0 | 0 | 0 | 0 |
| $S_4$ | .7 | 0 | 0 | .8I | 0 | .3 | 0 | .8 | 0 | 0 |
| $S_5$ | 0 | .8I | .3 | 0 | .7 | 1 | 0 | .3 | .7I | 0.7 |
| $S_6$ | .3 | .7 | .0 | .3 | 0 | 0 | 0 | 1 | 1 | 0 |
| $S_7$ | .9 | .4 | 0 | 0 | .8I | .9 | 0 | 0 | 0 | .4 |
| $S_8$ | .2I | 0 | 0 | 0 | 0 | 0 | .7I | 0 | .2 | .3 |

Let $Q_N$ denote the set of 8 symptoms / diseases i.e., $Q_N$ the neutrosophic vector is a $1 \times 8$ neutrosophic matrix.

Consider the neutrosophic equation $Q_N \circ P_N = R_N$; clearly when $Q_N$ and $P_N$ are known we can easily solve the neutrosophic relational equation and obtain the neutrosophic vector $R_N$.

Let $Q_N = (.3, .7, .5I .3, 0, .3, .2, .3I)$ be the neutrosophic vector given by the expert he feels the dependence is an indeterminate concept to some extent and difficult to cope with is present in most patients and in fact all patients suffer from depression. Using the neutrosophic equation $Q_N \circ P_N = R_N$ we calculate $R_N$ as follows:

$$R_N = (.5I, .3, .3, .3, .5I, 0, .7, .3, .7, .6)$$

which shows for the given input the relations with the patients. For in case of the first patient a combination of symptoms given by $Q_N$ results in an indeterminate the same happens to be true for the $5^{th}$ patients so for the given set of symptoms / disease given by $Q_N$ has nil influence on the patient $P_6$. But the same set of combination of symptoms has the maximum influence on the patients $P_7$ and $P_9$ with their membership grade equal to 0.7. The patient $P_{10}$ has 0.6 membership grade for the same neutrosophic vector $Q_N$. Thus the doctor can feed in any combination of the neutrosophic vectors and get the relative influence on the patients. We have illustrated this for a particular $Q_N$, an interested reader can work with any desired $Q_N$.

Thus we see how NRE can be used in the medical field to compare the relative effect on the patients apply neutrosophic linear programming defined by

$$\text{Maximize } z = cx, \, Ax \leq b \text{ such that } x \leq 0$$

where the coefficients A, b, c, $\in$ FN $\cup$ [0 1] the constraints may be considered as neutrosophic inequalities with variables x and z. Construct a neutrosophic linear programming to determine the uncertainty for any real world problem.

Chapter Five

# SUGGESTED PROBLEMS

Here we suggest a few problems for the reader to solve. Some of the problems will help one to build neutrosophic models.

1. Define for the neutrosophic relational equation x o A = b, where A = $(a_{ij})_{mxn}$ is a neutrosophic matrix, (b = $(b_1 \ldots b_n)$, $b_i \in [0\ 1] \cup I$. The users criterion function f(x) so as to form.

    i. Non-linear programming model with neutrosophic relation constrains.

    Min     f(x)

    Such that

    x o A = b                               (1)

    ii. Does minimizer of equation (1) prove a best solution to the user based on the objective function f(x).

2. Apply this model in medical diagnosis like symptom/ disease model or death wish of terminally ill patients.

3. Find a method of solution to neutrosophic relation equations in a complete Brouwerian neutrosophic lattice. Does a solution exist?

    (Hint: A solution to this problem will also give a solution to the fuzzy relation equations in a complete Brouwerian lattice).

4. Whether there exists a minimal elements in the solution set? Can one determine the minimal elements in the neutrosophic relation equation?

5. Define and describe the multi objective optimization problems with neutrosophic relation equation constraints.

6. Obtain properties of $t_N$-norm, compare a t-norm and a $t_N$ norm.

7. Develop an efficient learning algorithm for neural neutrosophic relational system with a new learning algorithm.

8. Define equality index for NRE mentioned in problem 7.

9. Study the solution set of NRE defined in chapter 4.

10. Investigate the solution set of the NRE described in problem 9, when

   i. Both U and V are finite
   ii. One of U or V is infinite.

11. If both the sets U and V are finite; does it imply all solutions are attainable for V for any NRE. (In case of FRE when both U and V are finite all solutions are attainable for V).

12. Find condition for X to be an unattainable solution in case of NRE.

13. Obtain any other interesting notions about attainable and unattainable solutions in case of NRE.

14. Define neutrosophic interpolation problem given in the form of the maximal neutrosophic relation with the membership function equal to

$$R = \bigcap_{K=1}^{N} (x_N(k) \to y_N(k)) .$$

15. Compare the neutrosophic model with a fuzzy model given in Problems 9 and 10 for any specific problem.

16. Form an algorithm to tackle largest solution of type I problem described in Kagei in case new neutrosophic relational equations.

17. Give the quasi-largest solution to type II problem in case of NRE analogous to those given by Kagei.

18. Prove If $X_N (A_N, b_N) \neq \phi$ then $I_j \neq \phi$ for all $j \in J$.

19. If $x \in X_N (A_N, b_N)$ then for each $j \in J$ does there exists $i \in I$ such that $x_{i_o} \wedge a_{i_o j} = b_j$ and $x_i \wedge a_{ij} \leq b_j$ for all $i \in I$.

20. Find an algorithm for finding an optimal solution of problem,

$$\text{minimize } Z_N = \sum_{i=1}^{m} C_i x_i$$

such that

$$x_N \text{ o } A_N = b_N,$$
$$x_i \in [0\ 1] \cup FN.$$

21. $\aleph_N \neq \phi \Leftrightarrow A_N @ b_N^{-1} \in \aleph_N$ then is $A_N @ b_N^{-1}$ the greatest solution ?

22. When both index sets I and J are finite $\aleph_N \neq \phi$ implies $\aleph_N^o \neq \phi$ then does $x \in \aleph_N \Leftrightarrow (\bar{x}_N \in \aleph_N^o, \bar{x}_N \leq x < \hat{x}_N)$ ?

23. If $\hat{x}_N$ is the greatest solution of the NRE and J is a finite set will $\hat{x} \in (X_N)_j^{(+)} \Leftrightarrow X^0 \neq \phi$ ?

24. Construct a neutrosophic linear programming analogous to fuzzy linear programming defined by

$$\text{Maximize } z = cx$$

such that

$$Ax \leq b$$
$$x \leq 0$$

where the coefficients A, b, c, $\in$ FN $\cup$ [0 1] the constraints may be considered as neutrosophic inequalities with variables x and z.

25.   Use the neutrosophic linear programming to solve the problem of control of waste gas pollution in environment by oil refinery. (The analogous problem done in chapter two using fuzzy linear programming).

26.   Define neutrosophic relational products analogous to fuzzy relational products given in (chapter 2)
       Study using the new definition given in Problem () the relation between systems/diseases and its relation with the patients.

# BIBLIOGRAPHY

In this book we have collected and utilized a comprehensive list of research papers that have dealt with the field of FRE. The bibliography list given below is exhaustive and lists all these research papers that we could get access to and which were referred to during the course of our work on the current book. The reference quality of this book is enhanced by its very extensive bibliography, annotated by various notes and comments sections making the book broadly and easily accessible.

1.      Adamopoulos, G.I., and Pappis, C.P., Some Results on the Resolution of Fuzzy Relation Equations, *Fuzzy Sets and Systems*, 60 (1993) 83-88. *[9, 33, 117]*

2.      Adlassnig, K.P., Fuzzy Set Theory in Medical Diagnosis, *IEEE Trans. Systems, Man, Cybernetics*, 16 (1986) 260-265. *[117, 122]*

3.      Akiyama, Y., Abe, T., Mitsunaga, T., and Koga, H., A Conceptual Study of Max-composition on the Correspondence of Base Spaces and its Applications in Determining Fuzzy Relations, *Japanese J. of Fuzzy Theory Systems*, 3 (1991) 113-132. *[64]*

4.      Bezdek, J.C., *Pattern Recognition with Fuzzy Objective Function Algorithm*, Plenum Press, New York, 1981. *[9,54]*

5.      Birkhoff, G., *Lattice Theory*, American Mathematical Society, 1979. *[61-63, 126, 162-63]*

6.      Blanco, A., Delgado, M., and Requena, I., Solving Fuzzy Relational Equations by Max-min Neural Network, *Proc. 3$^{rd}$ IEEE Internet Conf. On Fuzzy Systems*, Orlando (1994) 1737-1742. *[88-90, 102, 134]*

7.      Buckley, J.J., and Hayashi, Y., Fuzzy Neural Networks: A Survey, *Fuzzy Sets and Systems*, 66 (1994) 1-13. *[101]*

8.  Cechiarova, K., Unique Solvability of Max-Min Fuzzy Equations and Strong Regularity of Matrices over Fuzzy Algebra, *Fuzzy Sets and Systems*, 75 (1995) 165-177. *[83]*

9.  Cheng, L., and Peng, B., The Fuzzy Relation Equation with Union or Intersection Preserving Operator, *Fuzzy Sets and Systems*, 25 (1988) 191-204. *[40, 44]*

10. Chung, F., and Lee, T., A New Look at Solving a System of Fuzzy Relational Equations, *Fuzzy Sets and Systems*, 99 (1997) 343-353. *[9, 54]*

11. Czogala, E., Drewniak, J., and Pedrycz, W., Fuzzy relation Applications on Finite Set, *Fuzzy Sets and Systems*, 7 (1982) 89-101. *[9, 33-34, 117-18, 141]*

12. De Jong, K.A., An Analysis of the Behavior of a Class of Genetic Adaptive Systems, *Dissertation Abstracts Internet*, 86 (1975) 5140B. *[54]*

13. Di Nola, A., and Sessa, S., On the Set of Composite Fuzzy Relation Equations, *Fuzzy Sets and Systems*, 9 (1983) 275-285. *[71]*

14. Di Nola, A., Pedrycz, W., and Sessa, S., Some Theoretical Aspects of Fuzzy Relation Equations Describing Fuzzy System, *Inform Sci.*, 34 (1984) 261-264. *[87, 266]*

15. Di Nola, A., Pedrycz, W., Sessa, S., and Wang, P.Z., Fuzzy Relation Equation under a Class of Triangular Norms: A Survey and New Results, *Stochastica*, 8 (1984) 99-145. *[87, 141]*

16. Di Nola, A., Relational Equations in Totally Ordered Lattices and their Complete Resolution, *J. Math. Appl.*, 107 (1985) 148-155. *[33]*

17. Di Nola, A., Sessa, S., Pedrycz, W., and Sanchez, E., *Fuzzy Relational Equations and their Application in Knowledge Engineering*, Kluwer Academic Publishers, Dordrecht, 1989. *[44, 55-56, 64, 71-72, 75, 102, 117]*

18. Di Nola, A., Pedrycz, W., Sessa, S., and Sanchez, E., Fuzzy Relation Equations Theory as a Basis of Fuzzy Modeling: An Overview, *Fuzzy Sets and Systems*, 40 (1991) 415-429. *[9, 54]*

19.    Di Nola, A., On Solving Relational Equations in Brouwerian Lattices, *Fuzzy Sets and Systems*, 34 (1994) 365-376. [55]

20.    Drewniak, J., *Fuzzy Relation Calculus*, Univ. Slaski, Katowice, 1989. *[140-41]*

21.    Drewniak, J., Equations in Classes of Fuzzy Relations, *Fuzzy Sets and Systems*, 75 (1995) 215-228. *[140, 146]*

22.    Dubois, D., and Prade, H., Fuzzy Relation Equations and Causal Reasoning, *Fuzzy Sets and Systems*, 75 (1995) 119-134. *[126]*

23.    Dumford, N., and Schwartz, J.T., *Linear Operators Part I*, Interscience Publishers, New York, 1958. *[116]*

24.    Fang, S.C., and Li, G., Solving Fuzzy Relation Equations with a Linear Objective Function, *Fuzzy Sets and Systems*, 103 (1999) 107-113. *[33, 49, 117-18, 122, 268]*

25.    Fang, S.C., and Puthenpurn, S., Linear Optimization and Extensions: Theory and Algorithm, Prentice-Hall, New Jersey, 1993. *[117]*

26.    Galichet, S., and Foulloy, L., Fuzzy Controllers: Synthesis and Equivalences, IEEE Trans. Fuzzy Sets, (1995) 140-145. *[110]*

27.    Gavalec, M., Solvability and Unique Solvability of Max-min Fuzzy Equations. *Fuzzy Sets and Systems*, 124 (2001) 385-393. *[81-83]*

28.    Gottwald, S., Approximately Solving Fuzzy Relation Equations: Some Mathematical Results and Some Heuristic Proposals, *Fuzzy Sets and Systems*, 66 (1994) 175-193. *[148-49, 153]*

29.    Gottwald, S., Approximate Solutions of Fuzzy Relational Equations and a Characterization of t-norms that Define Matrices for Fuzzy Sets, *Fuzzy Sets and Systems*, 75 (1995) 189-201. *[87, 146, 148-49, 153, 156]*

30.    Guo, S.Z., Wang, P.Z., Di Nola, A., and Sessa, S., Further Contributions to the Study of Finite Fuzzy Relation

Equations, *Fuzzy Sets and Systems*, 26 (1988) 93-104. *[33, 117]*

31.  Gupta, M.M., and Qi, J., Design of Fuzzy Logic Controllers based on Generalized T-operators, Fuzzy Sets and Systems, 40 (1991) 473-489. *[64]*

32.  Gupta, M.M., and Qi, J., Theory of T-norms and Fuzzy Inference, *Fuzzy Sets and Systems*, 40 (1991) 431-450. *[64]*

33.  Gupta, M.M., and Rao, D.H., On the Principles of Fuzzy Neural Networks, *Fuzzy Sets and Systems*, 61 (1994) 1-18. *[64, 87-88, 101-02]*

34.  Higashi, M., and Klir, G.J., Resolution of Finite Fuzzy Relation Equations, *Fuzzy Sets and Systems*, 13 (1984) 65-82. *[9, 34, 49, 59, 68, 85-86, 117-18, 124, 141]*

35.  Hirota, K., and Pedrycz, W., Specificity Shift in Solving Fuzzy Relational Equations, *Fuzzy Sets and Systems*, 106 (1999) 211-220. *[196]*

36.  Holland, J., *Adaptation in Natural and Artificial Systems*, The University of Michigan Press, Ann Arbor, 1975. *[68]*

37.  Hong, D.H., and Hwang, S.Y., On the Compositional Rule of Inference under Triangular Norms, *Fuzzy Sets and Systems*, 66 (1994) 25-38. *[64]*

38.  Hormaifar, A., Lai, S., and Qi, X., Constrained Optimization via Genetic Algorithm, *Simulation*, 62 (1994) 242-254. *[50, 52, 54]*

39.  Imai, H., Kikuchi, K., and Miyakoshi, M., Unattainable Solutions of a Fuzzy Relation Equation, *Fuzzy Sets and Systems*, 99 (1998) 193-196. *[67-68, 70, 122, 126. 265, 270]*

40.  Jenei, S., On Archimedean Triangular Norms, *Fuzzy Sets and Systems*, 99 (1998) 179-186. *[64]*

41.  Joines, J.A., and Houck, C., On the Use of Non-stationary Penalty Function to Solve Nonlinear Constrained Optimization Problems with Gas, *Proc. 1$^{st}$ IEEE Internal Conf. Evolutionary Computation*, 2 (1994) 579-584. *[50, 52, 54]*

42.     Kagei, S., Fuzzy Relational Equation with Defuzzification Algorithm for the Largest Solution, *Fuzzy Sets and Systems*, 123 (2001) 119-127. *[75, 267]*

43.     Klir, G.J., and Yuan, B., *Fuzzy Sets and Fuzzy Logic: Theory and Applications*, Prentice-Hall, Englewood Cliffs NJ, 1995. *[16, 21-22, 25, 30, 41, 43-44, 64, 86, 97, 242]*

44.     Kosko, B., *Neural Networks and Fuzzy Systems: Dynamical Approach to Machine Intelligence*, Prentice-Hall, Englewood Cliffs NJ, 1992. *[68]*

45.     Kuhn, T., *The Structure of Scientific Revolutions*, Univ. of Chicago Press, 1962. *[247]*

46.     Kurano, M., Yasuda, M., Nakatami, J., and Yoshida, Y., A Limit Theorem in Some Dynamic Fuzzy Systems, *Fuzzy Sets and Systems*, 51(1992) 83- 88. *[68]*

47.     Kurano, M., Yasuda, M., Nakatami, J., and Yoshida, Y., A Fuzzy Relational Equation in Dynamic Fuzzy Systems, *Fuzzy Sets and Systems*, 101 (1999) 439-443. *[113-14]*

48.     Kuratowski, K., *Topology I*, Academic Press New York 1966. *[114]*

49.     Lee, C.C., Theoretical and Linguistic Aspects of Fuzzy Logic Controller, *Automatica*, 15 (1979) 553-577. *[108]*

50.     Lee, C.C., Fuzzy Logic in Control Systems: Fuzzy Logic Controller, Part I and II, *IEEE Trans. Systems, Man and Cybernetics*, 20 (1990) 404 – 405. *[108]*

51.     Lettieri, A., and Liguori, F., Characterization of Some Fuzzy Relation Equations Provided with one Solution on a Finite Set, *Fuzzy Sets and Systems*, 13 (1984) 83-94. *[60]*

52.     Li, G. and Fang, S.G., *On the Resolution of Finite Fuzzy Relation Equations,* OR Report No.322, North Carolina State University, Raleigh, North Carolina, 1986. *[9, 117, 122]*

53.     Li, X., and Ruan, D., Novel Neural Algorithm Based on Fuzzy S-rules for Solving Fuzzy Relation Equations Part I, *Fuzzy Sets and Systems*, 90 (1997) 11-23. *[87-88, 92, 105]*

54.     Li, X., and Ruan, D., Novel Neural Algorithms Based on Fuzzy S-rules for Solving Fuzzy Relation Equations Part II, *Fuzzy Sets and Systems*, 103 (1999) 473-486. *[87-88]*

55.     Li, X., and Ruan, D., Novel Neural Algorithm Based on Fuzzy S-rules for Solving Fuzzy Relation Equations Part III, *Fuzzy Sets and Systems*, 109 (2000) 355-362. *[87-88, 101, 105]*

56.     Li, X., Max-min Operator Network and Fuzzy S-rule, *Proc. 2nd National Meeting on Neural Networks*, Nanjing, 1991. *[68]*

57.     Liu, F., and Smarandache, F., *Intentionally and Unintentionally. On Both, A and Non-A, in Neutrosophy*. http://lanl.arxiv.org/ftp/math/papers/0201/0201009.pdf *[222]*

58.     Loetamonphing, J., and Fang, S.C., Optimization of Fuzzy Relation Equations with Max-product Composition, *Fuzzy Sets and Systems*, 118 (2001) 509-517. *[32, 34, 38, 40]*

59.     Loetamonphing, J., Fang, S.C., and Young, R.E., Multi Objective Optimization Problems with Fuzzy Relation Constraints, *Fuzzy Sets and Systems*, 127 (2002) 147-164. *[57, 60]*

60.     Lu, J., An Expert System Based on Fuzzy Relation Equations for PCS-1900 Cellular System Models, *Proc. South-eastern INFORMS Conference*, Myrtle Beach SC, Oct 1998. *[49]*

61.     Lu, J., and Fang, S.C., Solving Nonlinear Optimization Problems with Fuzzy Relation Equation Constraints, *Fuzzy Sets and Systems*, 119 (2001) 1-20. *[48-49, 52, 54-55, 60, 259]*

62.     Luo, C.Z., Reachable Solution Set of a Fuzzy Relation Equation, *J. of Math. Anal. Appl.*, 103 (1984) 524-532. *[49, 68, 70]*

63.     Luoh, L., Wang, W.J., Liaw, Y.K., New Algorithms for Solving Fuzzy Relation Equations, *Mathematics and Computers in Simulation*, 59 (2002) 329-333. *[85-86]*

64.     Michalewicz, Z. and Janikow, Z., Genetic Algorithms for Numerical Optimization, *Stats. Comput.*, 1 (1991) 75-91. *[50]*

65.     Michalewicz, Z., and Janikow, Z., *Handling Constraints in Genetic Algorithms*, Kaufmann Publishers, Los Altos CA, 1991. *[50, 52, 54]*

66.     Michalewicz, Z., *Genetic Algorithms + Data Structures = Evolution Programs*, Springer, New York, 1994. *[50]*

67.     Miyakoshi, M., and Shimbo, M., Solutions of Fuzzy Relational Equations with Triangular Norms, *Fuzzy Sets and Systems*, 16 (1985) 53-63. *[101]*

68.     Miyakoshi, M., and Shimbo, M., Sets of Solution Set Invariant Coefficient Matrices of Simple Fuzzy Relation Equations, *Fuzzy Sets and Systems*, 21 (1987) 59-83. *[68]*

69.     Miyakoshi, M., and Shimbo, M., Sets of Solution Set Equivalent Coefficient Matrices of Fuzzy Relation Equation, *Fuzzy Sets and Systems*, 35 (1990) 357-387. *[125]*

70.     Mizumoto, M., and Zimmermann, H.J., Comparisons of Fuzzy Reasoning Methods, *Fuzzy Sets and Systems*, 8 (1982) 253-283. *[110, 123]*

71.     Mizumoto, M., Are Max and Min Suitable Operators for Fuzzy Control Methods?, *Proc. 6$^{th}$ IFSA World Congress I*, Sao Paulo, Brazil, (1995) 497-500. *[68]*

72.     Neundorf, D., and Bohm, R., Solvability Criteria for Systems of Fuzzy Relation Equations, *Fuzzy Sets and Systems*, 80 (1996) 345-352. *[156]*

73.     Ngayen, H.T., A Note on the Extension Principle for Fuzzy Sets, *J. Math Anal. Appl.*, 64 (1978) 369-380. *[70]*

74.     Pavlica, V., and Petrovacki, D., About Simple Fuzzy Control and Fuzzy Control Based on Fuzzy Relational Equations, *Fuzzy Sets and Systems*, 101 (1999) 41-47. *[108]*

75.     Pedrycz, W., Fuzzy Relational Equations with Generalized Connectives and their Applications, *Fuzzy Sets and Systems*, 10 (1983) 185-201. *[9, 66, 148]*

76.     Pedrycz, W., *Fuzzy Control and Fuzzy Systems*, Wiley, New York, 1989. *[9, 66, 110-13]*

77.     Pedrycz, W., Inverse Problem in Fuzzy Relational Equations, *Fuzzy Sets and Systems*, 36 (1990) 277-291. *[9, 66]*

78.     Pedrycz, W., Processing in Relational Structures: Fuzzy Relational Equations, *Fuzzy Sets and Systems*, 25 (1991) 77-106. *[9, 66, 89]*

79.     Pedrycz, W., s-t Fuzzy Relational Equations, *Fuzzy Sets and Systems*, 59 (1993) 189-195. *[9, 64, 66]*

80.     Pedrycz, W., Genetic Algorithms for Learning in Fuzzy Relational Structures, *Fuzzy Sets and Systems*, 69 (1995) 37-52. *[9, 66]*

81.     Praseetha, V.R., *A New Class of Fuzzy Relation Equation and its Application to a Transportation Problem*, Masters Dissertation, Guide: Dr. W. B. Vasantha Kandasamy, Department of Mathematics, Indian Institute of Technology, April 2000. *[167]*

82.     Prevot, M., Algorithm for the Solution of Fuzzy Relation, *Fuzzy Sets and Systems*, 5 (1976) 38-48. *[9, 33, 117]*

83.     Ramathilagam, S., *Mathematical Approach to the Cement Industry problems using Fuzzy Theory*, Ph.D. Dissertation, Guide: Dr. W. B. Vasantha Kandasamy, Department of Mathematics, Indian Institute of Technology, Madras, November 2002. *[9, 54]*

84.     Sanchez, E., Resolution of Composite Fuzzy Relation Equation, *Inform. and Control*, 30 (1976) 38-48. *[9, 31-33, 44, 55-56, 68, 75, 79, 87, 101, 117, 124, 140-41, 152, 156-57, 163, 197]*

85.     Schitkowski, K., *More Test Examples for Non-linear Programming Codes*, Lecture Notes in Economics and Mathematical Systems: 282, Springer, New York, 1987. *[167]*

86.     Sessa, S., Some Results in the Setting of Fuzzy Relation Equation Theory, *Fuzzy Sets and Systems*, 14 (1984) 217-248. *[87, 101]*

87.     Smarandache, F., *Collected Papers III,* Editura Abaddaba, Oradea, 2000.
        http://www.gallup.unm.edu/~smarandache/CP3.pdf *[222]*

88.     Smarandache, F., *A Unifying Field in Logics: Neutrosophic Logic. Neutrosophy, Neutrosophic Set, Neutrosophic Probability and Statistics*, third edition, Xiquan, Phoenix, 2003. *[222]*

89.     Smarandache, F., Definitions Derived from Neutrosophics, In Proceedings of the *First International Conference on Neutrosophy, Neutrosophic Logic, Neutrosophic Set, Neutrosophic Probability and Statistics*, University of New Mexico, Gallup, 1-3 December 2001. *[222]*

90.     Smarandache, F., Dezert, J., Buller, A., Khoshnevisan, M., Bhattacharya, S., Singh, S., Liu, F., Dinulescu-Campina, Lucas, C., and Gershenson, C., *Proceedings of the First International Conference on Neutrosophy, Neutrosophic Logic, Neutrosophic Set, Neutrosophic Probability and Statistics*, The University of New Mexico, Gallup Campus, 1-3 December 2001. *[222]*

91.     Stamou, G.B., and Tzafestas, S.G., Neural Fuzzy Relational Systems with New Learning Algorithm, *Mathematics and Computers in Simulation*, 51 (2000) 301-314. *[64]*

92.     Stamou, G.B., and Tzafestas, S.G., Resolution of Composite Fuzzy Relation Equations based on Archimedean Triangular Norms, *Fuzzy Sets and Systems*, 120 (2001) 395-407. *[40, 63-64]*

93.     Steuer, R.E., *Multiple Criteria Optimization Theory: Computation and Applications*, Wiley, New York, 1986. *[59, 63]*

94.     Sugeno, M., *Industrial Applications of Fuzzy Control*, Elsevier, New York, 1985. *[108, 140]*

95.     Ukai, S., and Kaguei, S., Automatic Accompaniment Performance System using Fuzzy Inference, *Proc. Sino Japan Joint Meeting: Fuzzy sets and Systems*, Beijing, E1-5 (1990) 1-4. *[76]*

96.     Ukai, S., and Kaguei, S., Automatic Generation of Accompaniment Chords using Fuzzy Inference, *J. Japan Soc. Fuzzy Theory Systems*, 3 (1991) 377-381. *[76]*

97.     Vasantha Kandasamy, W.B., Neelakantan, N.R., and Kannan, S.R., Operability Study on Decision Tables in a Chemical Plant using Hierarchical Genetic Fuzzy Control Algorithms, *Vikram Mathematical Journal*, 19 (1999) 48-59. *[199]*

98.     Vasantha Kandasamy, W.B., and Praseetha, R., New Fuzzy Relation Equations to Estimate the Peak Hours of the Day for Transport Systems, *J. of Bihar Math. Soc.,* 20 (2000) 1-14. *[188]*

99.     Vasantha Kandasamy, W.B., Neelakantan, N.R., and Kannan, S.R., Replacement of Algebraic Linear Equations by Fuzzy Relation Equations in Chemical Engineering, In *Recent Trends in Mathematical Sciences*, Proc. of Int. Conf. on Recent Advances in Mathematical Sciences held at IIT Kharagpur on Dec. 20-22, 2001, published by Narosa Publishing House, (2001) 161-168. *[167]*

100.    Vasantha Kandasamy, W.B., and Balu, M. S., Use of Weighted Multi-Expert Neural Network System to Study the Indian Politics, *Varahimir J. of Math. Sci.*, 2 (2002) 44-53. *[167]*

101.    Vasantha Kandasamy, W.B., and Mary John, M., Fuzzy Analysis to Study the Pollution and the Disease Caused by Hazardous Waste From Textile Industries, *Ultra Sci.*, 14 (2002) 248-251. *[167]*

102.    Vasantha Kandasamy, W.B., Neelakantan, N.R., and Ramathilagam, S., Use of Fuzzy Neural Networks to Study the Proper Proportions of Raw Material Mix in Cement Plants, *Varahmihir J. Math. Sci.*, 2 (2002) 231-246. *[188]*

103.    Vasantha Kandasamy, W.B., and Smarandache, F., *Fuzzy Cognitive Maps and Neutrosophic Cognitive Maps*, Xiquan, Phoenix, 2003. *[167]*

104.    Wagenknecht, M., and Hatmasann, K., On the Existence of Minimal Solutions for Fuzzy Equations with Tolerances, *Fuzzy Sets and Systems*, 34 (1990) 237-244. *[124]*

105.    Wagenknecht, M., On Pseudo Transitive Approximations of fuzzy Relations, *Fuzzy Sets and Systems*, 44 (1991) 45-55. *[68, 145]*

106.    Wang, H.F., An Algorithm for Solving Iterated Complete Relation Equations, *Proc. NAFIPS*, (1988) 242-249. *[64, 117]*

107.    Wang, H.F., Wu, C.W., Ho, C.H., and Hsieh, M.J., Diagnosis of Gastric Cancer with Fuzzy Pattern Recognition, *J. Systems Engg.*, 2 (1992) 151-163. *[49]*

108.    Wang, W.F., A Multi Objective Mathematical Programming Problem with Fuzzy Relation Constraints, *J. Math. Criteria Dec. Anal.*, 4 (1995) 23-35. *[9, 58-59, 63]*

109.    Wang, X., Method of Solution to Fuzzy Relation Equations in a Complete Brouwerian Lattice, *Fuzzy Sets and Systems*, 120 (2001) 409-414. *[55, 162]*

110.    Wang, X., Infinite Fuzzy Relational Equations on a Complete Brouwerian Lattice, *Fuzzy Sets and Systems*, 138 (2003) 657-666. *[260-61]*

111.    Winston, W.L., *Introduction to Mathematical Programming: Applications and Algorithms*, Daxbury Press, Belmont CA, 1995. *[9, 54, 167]*

112.    Yager, R.R., Fuzzy Decision Making including Unequal Objective, *Fuzzy Sets and Systems*, 1 (1978) 87-95. *[33, 58]*

113.    Yager, R.R., On Ordered Weighted Averaging Operators in Multi Criteria Decision Making, *IEEE Trans. Systems, Man and Cybernetics*, 18 (1988) 183-190. *[33, 133]*

114.    Yen, J., Langari, R., and Zadeh, L.A., *Industrial Applications of Fuzzy Logic and Intelligent Systems*, IEEE Press, New York 1995. *[167]*

115.    Zadeh, L.A., Similarity Relations and Fuzzy Orderings, *Inform. Sci.*, 3 (1971) 177-200. *[110, 140, 166, 201]*

116.    Zadeh, L.A., A Theory of Approximate Reasoning, *Machine Intelligence*, 9 (1979) 149- 194. *[110, 166, 201]*

117.    Zimmermann, H.J., *Fuzzy Set Theory and its Applications*, Kluwer, Boston, 1988. *[9, 33, 49, 57, 62, 108-12, 117]*

# INDEX

Fuzzy partial ordering, 26
Fuzzy relation, 10 to 12
Fuzzy relational equations FRE, 9
Fuzzy S-rules, 87
Fuzzy union, 29

**G**

Genetic algorithm, (GA), 49-51

**H**

Hasse's diagram, 26

**I**

Implication, 110
Inefficient solution, 261
Inefficient solution, 59
Infimum, 25
Integral fuzzy neutrosophic column vector, 245
Integral fuzzy neutrosophic matrix, 244
Integral fuzzy neutrosophic row vector, 245
Integral neutrosophic Boolean algebra, 239-240
Integral neutrosophic lattice, 236
Interpolation nodes, 73
Intuitionistic set, 226
Irreflexive, 13, 16

**L**

Linear contraction, 54
Linear extraction, 54
Lower bound, 25
Lukeasiewicz implication, 146-149

**M**

Maximal solution matrix (max SM), 45
Maximal solution operation (max SO), 44
Maximum relation, 21-22
Maximum solution, 41, 85
Max-min composition, 12, 56

Max-product composition, 32
Max-t-composition, 71
Mean solution matrix  (mean SM), 45-46
Mean solution, 41
Membership matrix, 11
Membership neutrosophic matrix, 251-52
Minimal solution matrix (min SM), 45
Minimal solution operator (min SO), 44
Minimal solution, 21-22, 41
Minimum solution, 85
Mixed neutrosophic lattice, 236-237
Modus ponens tautology, 110

**N**

Nat topology, 90
N-attainable, 271
N-critical, 263
N-cut worthy, 234
N-equilibrium points, 234-235
N-equivalence operation, 263-264
Neutrosifying crisp function, 243
Neutrosophic algebraic product, 256-257
Neutrosophic algebraic sum, 258-259
Neutrosophic Archimedean, 264
Neutrosophic chain Boolean algebra, 239-240
Neutrosophic chain lattice, 236
Neutrosophic convexity, 234
Neutrosophic core (N-core), 233
Neutrosophic criterion vectors, 261
Neutrosophic digraph, 252
Neutrosophic directed graphs, 252
Neutrosophic equivalence class, 255
Neutrosophic equivalence relation, 255
Neutrosophic equivalent, 255-256
Neutrosophic fields, 244
Neutrosophic function, 240, 242-243
Neutrosophic lattices, 235-236
Neutrosophic matrix, 244
Neutrosophic members,, 232
Neutrosophic n-ary relation, 225
Neutrosophic probability appurtenance, 224
Neutrosophic relational equations (NRE), 249

**R**

Reflexive, 13, 16
Relative pseudo-complement, 163
Rough set, 226

**S**

Sagittal diagram, 14
Similarity relation, 17
Smooth derivative, 88
Solution matrix, 40
Solvability degree, 154
Special neutrosophic Boolean algebra, 239-240
Special neutrosophic lattice, 236-237
Specificity shift, 71
Standard intersection, 97
Standard neutrosophic intersection, 256-257
Standard union, 98
Sup i, 26
Super point, 54
Supremum, 25
Sup-t-FRE, 40
Symmetric, 13

**T**

Tautological set, 226
t-co-norm, 97-98, 101
Ternary relation, 12, 13
$t_N$-co-norm, 257
t-norms, 30, 43, 157
Tolerance relation, 22
Transitive, 15
Truth degree, 153-154
Two-layered max-min operation, 90

**U**

Unattainable solution, 67, 69, 70, 124
Upper bound, 25

# About the Authors

**Dr.W.B.Vasantha Kandasamy** is an Associate Professor in the Department of Mathematics, Indian Institute of Technology Madras, Chennai, where she lives with her husband Dr.K.Kandasamy and daughters Meena and Kama. Her current interests include Smarandache algebraic structures, fuzzy theory, coding/ communication theory. In the past decade she has guided eight Ph.D. scholars in the different fields of non-associative algebras, algebraic coding theory, transportation theory, fuzzy groups, and applications of fuzzy theory of the problems faced in chemical industries and cement industries. Currently, six Ph.D. scholars are working under her guidance. She has to her credit 255 research papers of which 203 are individually authored. Apart from this, she and her students have presented around 294 papers in national and international conferences. She teaches both undergraduate and post-graduate students and has guided over 41 M.Sc. and M.Tech. projects. She has worked in collaboration projects with the Indian Space Research Organization and with the Tamil Nadu State AIDS Control Society. She has authored a Book Series, consisting of ten research books on the topic of Smarandache Algebraic Structures which were published by the American Research Press.

She can be contacted at vasantha@itm.ac.in
You can visit her work on the web at: http://mat.iitm.ac.in/~wbv

---

**Dr.Florentin Smarandache** is an Associate Professor of Mathematics at the University of New Mexico, Gallup Campus, USA. He published over 60 books and 80 papers and notes in mathematics, philosophy, literature, rebus. In mathematics his research papers are in number theory, non-Euclidean geometry, synthetic geometry, algebraic structures, statistics, and multiple valued logic (fuzzy logic and fuzzy set, neutrosophic logic and neutrosophic set, neutrosophic probability). He contributed with proposed problems and solutions to the Student Mathematical Competitions. His latest interest is in information fusion were he works with Dr.Jean Dezert from ONERA (French National Establishment for Aerospace Research in Paris) in creasing a new theory of plausible and paradoxical reasoning (DSmT).

He can be contacted at smarand@unm.edu